\documentclass[11pt,a4paper]{amsart}
\usepackage[T1]{fontenc}
\usepackage{lmodern,microtype}
\usepackage[margin=25mm]{geometry}
\usepackage{amsmath,amssymb,amsthm,mathtools,booktabs,enumitem,tabularx}
\usepackage{xcolor}
\usepackage[hidelinks,hypertexnames=false]{hyperref}
\allowdisplaybreaks[1]
\newtheorem{theorem}{Theorem}[section]
\newtheorem{lemma}[theorem]{Lemma}
\newtheorem{proposition}[theorem]{Proposition}
\newtheorem{corollary}[theorem]{Corollary}

\theoremstyle{remark}
\newtheorem{remark}[theorem]{Remark}
\theoremstyle{plain}
\newtheorem*{lemma*}{Lemma}
\newcommand{\R}{\mathbb R}
\newcommand{\C}{\mathbb C}\newcommand{\NN}{\mathbb N}
\newcommand{\cC}{\mathcal C}
\newcommand{\avg}[2]{\overline{#1}_{#2}}
\newcommand{\ipA}[2]{\left(#1,#2\right)_A}
\newcommand{\eps}{\varepsilon}\newcommand{\grad}{\nabla}
\newcommand{\dif}{\,\mathrm d}\newcommand{\dV}{\,\mathrm dV_g}
\newcommand{\osc}{\operatorname{osc}}
\newcommand{\LS}{\operatorname{deg}_{\mathrm{LS}}}
\numberwithin{equation}{section}
\hypersetup{pdftitle={The topological degree of the regular SU(3) Toda system},pdfauthor={}}
\title[The degree of the $SU(3)$ Toda system]{The Topological Degree of the Regular $SU(3)$ Toda System}
\author{
    Youngae Lee
    \and
	Haoyu Li
	\and
	Juncheng Wei
	\and
	Wen Yang
    \and
    Lei Zhang
}
\date{}
\newcommand{\qform}{\mathfrak b}
\newcommand{\typeA}{\mathsf A}
\newcommand{\typeB}{\mathsf B}
\begin{document}
\begin{abstract}
We resolve the long-standing open problem of computing the
Leray--Schauder degree of the regular $SU(3)$ Toda system on a closed
oriented surface $M$. If $d_{m,k}$ denotes the degree in the block
$(4\pi m,4\pi(m+1))\times(4\pi k,4\pi(k+1))$, then
\[
 \sum_{m,k\geq0}d_{m,k}x^my^k
 =(1-x)^{\chi(M)-1}(1-y)^{\chi(M)-1}(1-xy)^{\chi(M)},
\]
where $\chi(M)$ is the Euler characteristic of $M$. Nonvanishing
of $d_{m,k}$ thus gives an existence criterion in terms of the
topology of $M$ and the parameter block.

To solve this problem, we use generic perturbations and local
second-order information on the weights to exclude collisions of
shadow centers and splitting at mixed blow-up points. The shadow
degrees and the Morse index then give the degree jump at each
critical value.
\end{abstract}
\maketitle

\tableofcontents

\section{Introduction}

Let $(M,g)$ be a connected and closed Riemannian surface. Let $h_1,h_2$ be smooth positive functions on $M$, let $\rho_1,\rho_2>0$ and write $\rho=(\rho_1,\rho_2)$. Write $\dV$ for the area element of $g$ and let
\(A=(a_{ij})=\begin{pmatrix}
2 & -1\\
-1 & 2
\end{pmatrix}\) be the Cartan matrix of $SU(3)$. We normalize the metric so that \(|M|=1\). The regular \(\mathrm{SU}(3)\) Toda system is then given by
\begin{equation}\tag{1.1}\left\{\label{1.1}
\begin{aligned}
&-\Delta_g u_i
=\displaystyle\sum_{j=1}^{2}a_{ij}\rho_j
\left(
\frac{h_j e^{u_j}}
{\int_M h_j e^{u_j}\dV}-1
\right)\mbox{ in }M, \\ &
\displaystyle\int_M u_i\dV=0,
\end{aligned}
\qquad i=1,2,\right.
\end{equation}
Let
\[
 \mathring H^1(M)
 =\left\{v\in H^1(M):\int_Mv\dV=0\right\},
 \qquad \mathcal H=[\mathring H^1(M)]^2.
\]
For $u=(u_1,u_2)\in\mathcal H$, we define the operator
\[
 K_{\rho,i}(u):=(-\Delta_g)^{-1}
 \left[\sum_{j=1}^2a_{ij}\rho_j
 \left(\frac{h_je^{u_j}}{\int_Mh_je^{u_j}\dV}-1\right)\right],
 \qquad i=1,2,
\]
where $(-\Delta_g)^{-1}$ is the mean-zero inverse of $-\Delta_g$.
With $K_\rho=(K_{\rho,1},K_{\rho,2})$ and $I$ the identity on
$\mathcal H$, we write $F_\rho:=I-K_\rho$.

For $m,k\in\mathbb Z_{\geq0}$, set
\[
  \mathcal B_{m,k}
  :=(4\pi m,4\pi(m+1))\times(4\pi k,4\pi(k+1)).
\]
The a priori bounds of Lin--Wei--Yang--Zhang
\cite[Theorem~1.1]{LWYZ18}, together with elliptic estimates, imply that
for each $\rho\in\mathcal B_{m,k}$ the Leray--Schauder degree of
$F_\rho$ is well defined on every sufficiently large ball centered
at the origin in $\mathcal H$.

By the uniform a priori bounds along homotopies (see Section~\ref{sec:noncritical-degree})
and the homotopy invariance of the Leray--Schauder degree, the degree
is constant on $\mathcal B_{m,k}$ and independent of the choice of
positive weights. We denote this common value by $d_{m,k}$.

Our main result determines the degrees $d_{m,k}$ by giving an explicit
formula for their generating function.

\setcounter{equation}{1}
\begin{theorem}\label{th1.1} Let $h_1,h_2\in C^\infty(M)$ be strictly positive, and let $\chi(M)$
denote the Euler characteristic of $M$. Then the degrees $d_{m,k}$
satisfy
\begin{equation}\label{1.2}
  D(x,y):=\sum_{m,k\geq0}d_{m,k}x^my^k
  =(1-x)^{\chi(M)-1}(1-y)^{\chi(M)-1}(1-xy)^{\chi(M)}.
\end{equation}
\end{theorem}

\medskip

Write $\chi=\chi(M)$ and set $b_j=\binom{j-\chi}{j}$ for $j\geq0$
and $b_j=0$ for $j<0$. Comparing coefficients of $x^my^k$ in
\eqref{1.2} gives
\begin{equation}\label{1.3}
  d_{m,k}
  =\sum_{j=0}^{\min(m,k)}
    (-1)^j\binom{\chi}{j}b_{m-j}b_{k-j}.
\end{equation}

Theorem~\ref{th1.1} gives a degree formula for the Toda system in the
spirit of the degree calculations for the scalar mean field equation
\[
  -\Delta_g v
  =\lambda\left(\frac{he^v}{\int_M he^v\dV}-1\right),
  \qquad \int_M v\dV=0,
\]
where $h$ is smooth and positive. The study of the degree for this
equation goes back to the work of Y.~Y.~Li \cite{Li99}.
Chen and Lin first obtained sharp
estimates for multiple bubbles \cite{CL02} and then used these
estimates to compute the topological degree \cite{CL03}. Degree counting has since become an effective method for proving the existence of solutions to mean-field-type equations and systems. The method relies on detailed blow-up analysis, whose origins can be traced to the work of Brezis and Merle \cite{BM91}, the quantization of local masses at
regular blow-up points by Li and Shafrir \cite{LS94}, and the classification of entire finite-mass solutions of the regular Liouville equation by Chen and Li \cite{ChenLi91}. For a scalar mean field equation with collapsing singularities, Lee, Lin and Yang constructed bubbling solutions whose mass does not concentrate completely \cite{LLY19}. The same blow-up theory also applies to other questions, including the local uniqueness of bubbling solutions. Bartolucci, Jevnikar, Lee and Yang proved uniqueness of bubbling solutions with prescribed blow-up points under nondegeneracy hypotheses \cite{BJLY}.

Among mean field equations and systems, the Toda system exhibits particularly complicated blow-up behavior. Early studies of Toda systems on compact surfaces include the work of Jost, Lin and Wang, which studies bubbling for the Toda system and proved existence at some critical parameters \cite{JLW06}. For the $SU(3)$ case, Lin, Wei and Zhao obtained sharp estimates for fully bubbling solutions \cite{LWZhao12} and Lin, Wei and Zhang identified a finite set of possible local mass pairs at blow-up for singular systems \cite{LWZ}.

Variational existence results include the noncoercive regime treated by Malchiodi and Ruiz \cite{MR13}, the treatment of singular sources under topological and noncriticality hypotheses by Battaglia, Jevnikar, Malchiodi and Ruiz \cite{BJMR15} and the topological join construction of Jevnikar, Kallel and Malchiodi for surfaces of arbitrary genus in specified parameter ranges \cite{JKM15}. Malchiodi and Ruiz determined the parity of the Leray--Schauder degree on diagonal parameter blocks for the regular $SU(3)$
system \cite{MR15}. 

A natural question is how to compute the topological degree of the $SU(3)$ Toda system. Lee, Lin, Wei and Yang used a shadow system to study the degree jump at the first critical value of one parameter \cite{LLWY18}. Lee, Lin, Yang and Zhang treated the corresponding one-point blow-up problem for rank-two Toda systems with simple singularities \cite{LLYZ20}. For the nonsymmetric $B_2$ and $G_2$ systems, Lee, Lin and Yan obtained degree formulas in parameter ranges involving multiple bubbles \cite{LLY22}. The above results give partial information on the degree. In this paper, we derive a \emph{complete} formula for the Leray--Schauder degree of the regular $SU(3)$ Toda system in every noncritical parameter block. In a forthcoming paper \cite{LLWYZ}, we also derive the degree-counting formula for the \emph{singular} $SU(3)$ Toda system.

Section~\ref{rev:sec:global} proves Theorem~\ref{th1.1} by combining
the degree below the first critical value with the degree jumps
across subsequent critical values.
Sections~\ref{rev:sec:scalar}--\ref{rev:sec:ordinary} establish the
results used in that proof, including the blow-up classification,
the shadow degree formula, and the local degree computations for
mixed and ordinary configurations.

\section{Setting and proof of Theorem 1.1}
\label{rev:sec:global}\label{global-sec:localization}

\subsection{Functional setting and notation}\label{sec:main}
Recall the Toda system \eqref{1.1} and the map $F_\rho$ on $\mathcal H$ in the introduction. Here, $(-\Delta_g)^{-1}$ denotes the inverse of $-\Delta_g$ on the space of zero-mean functions. It is well-defined based on Poincar\'e's inequality and Lax--Milgram. Thus \eqref{1.1} is equivalent to
\[
 F_\rho(u)=u-K_\rho(u)=0.
\]

We normalize the Green function $G$ by
\[
 -\Delta_g G(\cdot,p)=\delta_p-1\mbox{ in }M,
 \qquad
 \int_M G(\cdot,p)\dV=0.
\]
Let $d_g$ be the geodesic distance on $(M,g)$. The associated
Robin function is defined by
\[
 r_g(p):=\lim_{x\to p}
 \left(G(x,p)+\frac{1}{2\pi}\log d_g(x,p)\right),
\]
A \emph{density} on $M$ is a nonnegative integrable function $f$.
Its associated measure is $f\dV$. For a density $f$, we use the operator notation
$Gf(x)=\int_M G(x,p)f(p)\,\mathrm dV_g(p)$.

In what follows, $\nu$ denotes the sequence index and limits are taken as $\nu\to\infty$ unless otherwise stated.
We suppress this index on solutions and associated parameters,
densities and centers whenever the meaning is clear, and pass
to subsequences without relabeling.

Given a bounded sequence of parameter vectors
$\rho_\nu=(\rho_{1,\nu},\rho_{2,\nu})^T\in(0,\infty)^2$,
let $u_\nu=(u_{1,\nu},u_{2,\nu})^T$ be a solution of
\eqref{1.1} with parameter vector $\rho_\nu$ for each $\nu$.
We use the notation and Green representation
\begin{equation}\label{2.1}
 \begin{aligned}
 C =A^{-1}=\frac13\begin{pmatrix}2&1\\1&2\end{pmatrix},
 \qquad
 f_{i,\nu}=\frac{\rho_{i,\nu}h_i e^{u_{i,\nu}}}
                 {\int_M h_i e^{u_{i,\nu}}\dV},\qquad
 u_\nu =A\begin{pmatrix}Gf_{1,\nu}\\Gf_{2,\nu}\end{pmatrix}.
 \end{aligned}
\end{equation}

\subsection{Local coordinates}\label{sec:notation}
In local coordinates $z=x+iy$, we write $\dif A=dx\,dy$ and use $\dif s$ for Euclidean arclength. We denote coordinate disks centered at the origin by $B_r$ or $D_r$. The plane $\R^2$ is also taken as $\C$. We choose all isothermal charts of $M$. Let $\boldsymbol n$ denote the outward unit normal to the boundary of the domain under consideration. In these coordinates $v_z=(v_x-iv_y)/2$ and $\Delta_0$ is the Euclidean Laplacian.

In an isothermal coordinate with $g=e^\phi|dz|^2$, let $K_g$ denote the Gaussian curvature and write the regular part of the Green kernel as $R(z,\zeta)$, so that $G(z,\zeta)=-(2\pi)^{-1}\log|z-\zeta|+R(z,\zeta)$.
For a sequence as in \eqref{2.1}, we set
\begin{align}\label{def:Fnu}
F_{i,\nu}:=e^\phi f_{i,\nu}\,\,\mbox{ and }\,\,U_{i,\nu}:=\log F_{i,\nu}.
\end{align}
We write $U_\nu=(U_{1,\nu},U_{2,\nu})^T$. For a parameter vector $\rho$, write $f^\rho=(f_1^\rho,f_2^\rho)^T$
for the forcing term in local coordinates. In terms of $U_\nu$, the
Toda system \eqref{1.1} takes the form
\begin{equation}\label{2.2}
  -\Delta_0 U_\nu
  = A(F_{1,\nu},F_{2,\nu})^T + f^{\rho_\nu},
  \qquad
  f_i^\rho
  = -e^\phi(A\rho)_i - \Delta_0(\log h_i+\phi),
  \quad i=1,2.
\end{equation}

\subsection{The degree on noncritical blocks}
\label{sec:noncritical-degree}

A family of Toda solutions is called \emph{unbounded} if the
supremum of its $\mathcal H$-norms is infinite. For a sequence
with weights $h_{i,\nu}$, set
\[
 I_{i,\nu}=\int_M h_{i,\nu}e^{u_{i,\nu}}\dV,\qquad i=1,2.
\]
In \eqref{2.1} and \eqref{2.2}, we use $h_{i,\nu}$ for the
weight of the $\nu$-th equation.

\begin{lemma}\label{lem:toda-concentration}
Let $(u_\nu)$ be a sequence of Toda solutions with parameters
$\rho_\nu$ and smooth positive weights $h_{i,\nu}$. Assume
\[
 \|u_\nu\|_{\mathcal H}\to\infty,\qquad
 \rho_\nu\to\rho_*\in(0,\infty)^2,\qquad h_{i,*}>0,
\]
and
\[
 \|h_{i,\nu}-h_{i,*}\|_{C^k(M)}\to0
 \qquad\text{for every integer }k\geq0,\quad i=1,2.
\]
Along a subsequence, there are nonnegative measures $\mu_i$,
a finite nonempty set $S\subset M$ and a limit $u_*$ such that
\[
 f_{i,\nu}\dV\rightharpoonup\mu_i,\qquad
 \mu_i(M)=\rho_{i,*},\qquad
 u_\nu\to u_*\quad\text{in }C^k_{\mathrm{loc}}(M\setminus S)
 \quad\text{for every integer }k\geq0.
\]
For every $q\in S$,
\begin{equation}\label{2.3}
 \frac{1}{4\pi}\bigl(\mu_1(\{q\}),\mu_2(\{q\})\bigr)
 \in\{(1,0),(0,1),(1,2),(2,1),(2,2)\}.
\end{equation}
Moreover, for $i=1,2$,
\[
 \gamma_i(q):=\frac{1}{2\pi}\sum_{j=1}^2a_{ij}\mu_j(\{q\})>2
 \text{ for some }q\in S
 \quad\Longrightarrow\quad
 I_{i,\nu}\to\infty,\qquad
 \mu_i=\sum_{p\in S}\mu_i(\{p\})\delta_p.
\]
\end{lemma}

\begin{proof}
\textit{Step 1. Local convergence.}
Since $|M|=1$ and $u_{i,\nu}$ has zero mean, Jensen's inequality
gives \begin{align}\label{I_lower}I_{i,\nu}\geq\min_M h_{i,\nu}\geq c>0.\end{align} The mean-zero
elliptic estimate and \eqref{2.1} give
\begin{align}\label{ineq:unbdd}
 \|u_\nu\|_{\mathcal H}
 \leq C\max_i\|f_{i,\nu}\|_{L^\infty(M)}
 \leq C\exp\bigl(\max_i\sup_M u_{i,\nu}\bigr).
\end{align}
In particular, $\max_i\sup_M u_{i,\nu}\to\infty$.
The mass identity $\int_M f_{i,\nu}\dV=\rho_{i,\nu}$ gives,
after passage to a subsequence,
\[
 f_{i,\nu}\dV\rightharpoonup\mu_i,\qquad
 \mu_i(M)=\rho_{i,*}.
\]

Fix $p>1$ and $\varepsilon_0>0$ with $p\varepsilon_0<2\pi$.
Set
\begin{align}\label{def:S}
 \lambda_\nu:=(f_{1,\nu}+f_{2,\nu})\dV,\qquad
 \lambda:=\mu_1+\mu_2,\qquad
 S:=\{q\in M:\lambda(\{q\})\geq\varepsilon_0\}.
\end{align}
The total mass bound gives $\#S<\infty$. For $x\notin S$,
choose a geodesic ball $B_g(x,2r)$ such that
\[
 \limsup_{\nu\to\infty}\lambda_\nu(B_g(x,2r))
 \leq\lambda(\overline{B_g(x,2r)})<\varepsilon_0.
\]
For $y\in B_g(x,r)$, \eqref{2.1} and the lower bound for $G$
give
\[
 u_{i,\nu}(y)\leq C+
 \frac1\pi\int_{B_g(x,2r)}
       \log^+\frac1{d_g(y,z)}\,\dif\lambda_\nu(z),
 \qquad \log^+t=\max\{\log t,0\}.
\]
Here the Green kernel off $B_g(x,2r)$ is uniformly bounded for
$y\in B_g(x,r)$. Jensen's inequality for the normalized
restriction of $\lambda_\nu$ to $B_g(x,2r)$ gives
\[
 \int_{B_g(x,r)}e^{p u_{i,\nu}}\dV
 \leq C\sup_{z\in B_g(x,2r)}
 \int_{B_g(x,r)}
 \exp\!\left(\frac{p\varepsilon_0}{\pi}
               \log^+\frac1{d_g(y,z)}\right)\dV(y)
 \leq C.
\]
Together with \eqref{I_lower}, we see that  $f_{i,\nu}$ is bounded in $L^p(B_g(x,r))$.
H\"older's inequality in the Green representation gives a
uniform bound for $|u_{i,\nu}|$ on $B_g(x,r/2)$. Interior
elliptic estimates and a diagonal subsequence give the asserted
local convergence on $M\setminus S$. If $S=\varnothing$, a
finite covering bounds $\max_i\sup_M u_{i,\nu}$, a contradiction.

\medskip
\noindent\textit{Step 2. Local masses.}
Fix $q\in S$. Choose an isothermal coordinate disk $D$ with
$D\cap S=\{q\}$, and fix $r>0$ with $\overline B_r(q)\subset D$.
Let
\[
 -\Delta_0\eta_{i,\nu}=f_i^{\rho_\nu}\quad\text{in }D,\qquad
 \eta_{i,\nu}=0\quad\text{on }\partial D.
\]
Let $\eta_{i,*}$ solve the same problem with forcing
$f_i^{\rho_*}$ evaluated at $h_{i,*}$. For $z=q+ry$, we define
\begin{equation}\label{eq:local-quantization-rescaling}
\begin{aligned}
 v_{i,\nu}(y)
 &=U_{i,\nu}(q+ry)-\eta_{i,\nu}(q+ry)
     +2\log r+\eta_{i,\nu}(q),\\
 H_{i,\nu}(y)
 &=\exp\bigl(\eta_{i,\nu}(q+ry)-\eta_{i,\nu}(q)\bigr).
\end{aligned}
\end{equation}
By \eqref{def:Fnu} and \eqref{2.2}, on $B_1$ we have
\[
 -\Delta_y v_{i,\nu}
 =\sum_{j=1}^2 a_{ij}H_{j,\nu}e^{v_{j,\nu}},\qquad
 H_{i,\nu}e^{v_{i,\nu}}=r^2F_{i,\nu}(q+ry).
\]
Fix $0<\alpha<1$. The classical Schauder estimate
\cite[p.~101, Eq.~(6.42)]{GT01} gives
\[
 \|\eta_{i,\nu}-\eta_{i,*}\|_{C^{2,\alpha}(\overline D)}
 \leq C\|f_i^{\rho_\nu}-f_i^{\rho_*}\|_{C^{0,\alpha}(\overline D)}
 \longrightarrow 0.
\]
The convergence of the forcing follows from \eqref{2.2}.
Hence $H_{i,\nu}\to H_{i,*}>0$ in $C^1(\overline B_1)$, where
$H_{i,*}$ is defined by the same formula with $\eta_{i,*}$.

By \eqref{2.1} and \eqref{def:Fnu},
\[
 U_{i,\nu}
 =u_{i,\nu}+\log h_{i,\nu}+\phi+
   \log(\rho_{i,\nu}/I_{i,\nu}).
\]
The last term is constant in space and bounded from above.
The local bounds off $S$ and the bounds for $\eta_{i,\nu}$ give
\[
 \sup_K v_{i,\nu}\leq C_K
 \quad(K\Subset B_1\setminus\{0\}),\qquad
 \operatorname{osc}_{\partial B_1}v_{i,\nu}\leq C.
\]
The change of variables gives, for $0<s\leq1$,
\[
 \int_{B_s}H_{i,\nu}e^{v_{i,\nu}}\dif A
 =\int_{B_{rs}(q)}F_{i,\nu}\dif A
 =\int_{B_{rs}(q)}f_{i,\nu}\dV
 \leq\rho_{i,\nu}.
\]
If $\max_i\sup_{B_1}v_{i,\nu}$ were bounded along a subsequence,
both $F_{i,\nu}$ would be bounded on $B_r(q)$. This contradicts
$\lambda(\{q\})>0$. Thus $\max_i\sup_{B_1}v_{i,\nu}\to+\infty$.
The hypotheses of \cite[Lemma~3.1(i)]{LLWY18} are satisfied.
The change of variables identifies the local masses with
$\mu_i(\{q\})$ and gives \eqref{2.3}.

\medskip
\noindent\textit{Step 3. Complete concentration of a component.}
For a finite measure $\mu$, we write
$G\mu(x)=\int_M G(x,y)\,\dif\mu(y)$. By \eqref{2.1} and weak
convergence,
\[
 u_{i,*}=\sum_{j=1}^2a_{ij}G\mu_j
 \qquad\text{on }M\setminus S.
\]
Indeed, the identity holds first in distributions and then agrees
with the smooth local limit outside of $S$. Fix $q\in S$ with
$\gamma_i(q)>2$. Fix an isothermal coordinate disk $D$ centered at $q=0$. Then
\[
\begin{aligned}
 \frac1{2\pi}\int_0^{2\pi}u_{i,*}(se^{\mathrm i\theta})\,\dif\theta
 &=\frac1{2\pi}\sum_{j=1}^2\int_M
   \left(\int_0^{2\pi}a_{ij}G(se^{\mathrm i\theta},y)\,\dif\theta\right)
   \dif\mu_j(y)\\
 &=\gamma_i(q)\log(1/s)
   -\frac1{(2\pi)^2}\sum_{j=1}^2\int_{D\setminus\{q\}}
   \left(\int_0^{2\pi}a_{ij}\log|se^{\mathrm i\theta}-y|
         \,\dif\theta\right)\dif\mu_j(y)+O(1)\\
 &=\gamma_i(q)\log(1/s)
   -\frac1{2\pi}\sum_{j=1}^2\int_{D\setminus\{q\}}
     a_{ij}\log\max\{s,|y|\}\,\dif\mu_j(y)+O(1)\\
 &=\gamma_i(q)\log(1/s)+o(\log(1/s)).
\end{aligned}
\]
Here, for $0<s<\delta<1$ with $B_\delta(q)\subset D$,
\[
 \left|\sum_{j=1}^2\int_{D\setminus\{q\}}
 a_{ij}\log\max\{s,|y|\}\,\dif\mu_j(y)\right|
 \leq C_\delta+\log(1/s)\sum_{j=1}^2
 |a_{ij}|\,\mu_j(B_\delta(q)\setminus\{q\}).
\]
Divide by $\log(1/s)$ and let $s\downarrow0$. Then let
$\delta\downarrow0$, using
$\mu_j(B_\delta(q)\setminus\{q\})\to0$.

If $I_{i,\nu}$ were bounded along a subsequence, Fatou's lemma
and $\min_M h_{i,\nu}\geq c>0$ would give $e^{u_{i,*}}\in L^1(M)$.
For $2<\gamma'<\gamma_i(q)$, Jensen's inequality on circles gives
\[
 \frac1{2\pi}\int_0^{2\pi}e^{u_{i,*}(se^{\mathrm i\theta})}
       \,\dif\theta\geq s^{-\gamma'}
 \qquad\text{for all sufficiently small }s>0.
\]
This contradicts integrability. Hence $I_{i,\nu}\to\infty$.
Since $u_{i,\nu}$ is locally bounded off $S$, $f_{i,\nu}$ tends
locally uniformly to zero there. Thus $\mu_i$ is supported on $S$.
\end{proof}

\begin{remark}
By (\ref{ineq:unbdd}), the unboundedness in $\mathcal H$-norm implies the unboundedness of $\|f_{i,\nu}\|_{L^\infty(M)}$ and of $\max_{i}\sup_M u_{i,\nu}$.
\end{remark}

\begin{remark}\label{r:CompleteComcentration}
At least one component concentrates completely. Indeed, fix $q\in S$.
By \eqref{2.3}, some $i\in\{1,2\}$ satisfies $\gamma_i(q)>2$.
Lemma~\ref{lem:toda-concentration} gives $I_{i,\nu}\to\infty$.
For every $K\Subset M\setminus S$, the local bound for $u_{i,\nu}$ gives
\[
 \sup_K f_{i,\nu}
 \leq
 \frac{\rho_{i,\nu}\|h_{i,\nu}\|_{L^\infty(M)}
       \exp(\sup_K u_{i,\nu})}{I_{i,\nu}}
 \leq \frac{C_K}{I_{i,\nu}}
 \longrightarrow0.
\]
Thus $\mu_i$ is supported on $S$ and $\mu_i=\sum_{p\in S}\mu_i(\{p\})\delta_p$.
\end{remark}

\begin{lemma}
Let $\rho:[0,1]\to\mathcal B_{m,k}=(4\pi m,4\pi(m+1))\times(4\pi k,4\pi(k+1))$
be a continuous path and let $h_i(t,x)>0$ be smooth on $[0,1]\times M$.
The solutions are uniformly bounded in $\mathcal H$ along the path.
\end{lemma}

\begin{proof}
Suppose the solutions along the path are unbounded. Choose $t_\nu$
and corresponding solutions with $\|u_\nu\|_{\mathcal H}\to\infty$.
After passage to a subsequence, $t_\nu\to t_*$. Continuity of $\rho$
and smoothness of the weights give the convergence hypotheses of
Lemma~\ref{lem:toda-concentration} with limits $\rho(t_*)$ and
$h_i(t_*,\cdot)$.

Choose $q\in S$. By \eqref{2.3}, there exists $i\in\{1,2\}$
such that $\gamma_i(q)>2$. Lemma~\ref{lem:toda-concentration}
then implies that $\mu_i$ is supported on $S$, so $\rho_{i,*}=\mu_i(M)=\sum_{p\in S}\mu_i(\{p\})\in4\pi\mathbb N$.
This contradicts $\rho_*\in\mathcal B_{m,k}$.
\end{proof}

We verify that $(t,u)\mapsto K_{\rho(t)}(u)$, with weights
$h_i(t,\cdot)$, is continuous and maps bounded sets into
relatively compact subsets of $\mathcal H$. On bounded subsets of $\mathcal H$, the Moser--Trudinger inequality
\cite{Moser71} bounds $e^{u_i}$ in every finite $L^p(M)$.
Together with the positive lower bound for the denominators,
elliptic estimates and the compact embedding
$W^{2,p}(M)\hookrightarrow H^1(M)$ for $p>2$ give compactness
of $K_\rho$ uniformly along the path. Strong convergence in
$\mathcal H$ gives convergence in measure, and the uniform
$L^{2p}$ bound gives convergence of the exponentials in $L^p$.
This also gives continuity of the homotopy.

A common large ball therefore defines the Leray--Schauder degree
along the path. Join parameters by a line segment in
$\mathcal B_{m,k}$ and join positive weights by linearly
interpolating their logarithms. Homotopy invariance
\cite[Theorems~8.1 and~8.2]{Deimling85} shows that the degree
is constant in each block and independent of the positive weights.

\subsection{Blow-up limits and the shadow system}\label{sec:shadow-setting}
We call a concentration point \emph{ordinary} or \emph{mixed}
if its local mass pair is $(4\pi,0)$ or $(8\pi,4\pi)$,
respectively. We denote these types by $\typeA$ and $\typeB$.

\begin{lemma}\label{le2.1}
Suppose that $(u_\nu)$ is unbounded, that
$\rho_{2,\nu}\to\rho_2>0$ with $\rho_2\notin4\pi\mathbb N$ and that $\rho_{1,\nu}\to\rho_1^*>0$.
Along a subsequence, the following hold.

\begin{enumerate}[label=\textup{(\roman*)},leftmargin=*]
\item There exist nonnegative integers $a,b$ and a concentration set
consisting of $a$ type $A$ points $P=(p_1,\ldots,p_a)$ and
$b$ type $B$ points $Q=(q_1,\ldots,q_b)$. Let $n_*\dV$ be the
absolutely continuous part in the sense of Lebesgue decomposition of the weak limit of $f_{2,\nu}\dV$.
Then
\begin{equation}\label{2.4}
 \left\{
 \begin{aligned}
 \rho_1^*&=4\pi(a+2b),
 \qquad \sigma:=\sigma_b=\rho_2-4\pi b>0,\\
 f_{1,\nu}\dV&\rightharpoonup
 4\pi\sum_{i=1}^a\delta_{p_i}+8\pi\sum_{j=1}^b\delta_{q_j},\\
 f_{2,\nu}\dV&\rightharpoonup
 4\pi\sum_{j=1}^b\delta_{q_j}+n_*\dV.
 \end{aligned}
 \right.
\end{equation}
Moreover, there exists $w_*\in C^\infty(M)$ such that
\begin{equation}\label{2.5}
 n_*=\frac{\sigma H_Pe^{w_*}}{\int_M H_Pe^{w_*}\dV}
 \quad\text{where}\quad
 H_P=h_2e^{-4\pi\sum_{i=1}^aG(\cdot,p_i)}.
\end{equation}

\item The following limits hold in
$C^\infty_{\mathrm{loc}}(M\setminus(P\cup Q))$:
\begin{equation}\label{2.6}\left\{
 \begin{aligned}
 u_{1,\nu}&\longrightarrow
 8\pi\sum_{i=1}^aG(\cdot,p_i)
 +12\pi\sum_{j=1}^bG(\cdot,q_j)-\tfrac12w_*,\\
 u_{2,\nu}&\longrightarrow
 w_*-4\pi\sum_{i=1}^aG(\cdot,p_i).
 \end{aligned}\right.
\end{equation}

\item The triple $(w_*,P,Q)$ satisfies the shadow system consisting of
\begin{equation}\label{2.7}
 -\Delta_gw_*
 =2\sigma
 \left(\frac{H_Pe^{w_*}}{\int_MH_Pe^{w_*}\dV}-1\right),
 \qquad \int_Mw_*\dV=0,
\end{equation}
and the center equations. To write them, put
\[
 g_1=\log h_1+4\pi r_g,
 \qquad
 V=2\log h_1+\log h_2+12\pi r_g.
\]
The center equations are
\begin{equation}\label{2.8}\left\{
 \begin{aligned}
 F_i^A&:=
 \nabla g_1(p_i)-\frac12\nabla w_*(p_i)
 +8\pi\sum_{l\ne i}\nabla_xG(p_i,p_l)
 +12\pi\sum_{j=1}^b\nabla_xG(p_i,q_j)=0,
 \quad 1\leq i\leq a,\\
 F_j^B&:=
 \nabla V(q_j)
 +12\pi\sum_{i=1}^a\nabla_xG(q_j,p_i)
 +24\pi\sum_{l\ne j}\nabla_xG(q_j,q_l)=0,
 \quad 1\leq j\leq b,
 \end{aligned}\right.
\end{equation}
\end{enumerate}
\end{lemma}

We postpone the proof to Subsection \ref{sec:limiting-equations}.

Fix $\rho_2\in(4\pi k,4\pi(k+1))$ and  $\rho_1=4\pi m$
with $m\geq1$. Lemma~\ref{le2.1} identifies the two possible types
of concentration points for a sequence with such a limit.
If there are $a$ ordinary and $b$ mixed points, then $a+2b=m$ and
$b\leq k$. By Remark \ref{r:CompleteComcentration}, the first component concentrates completely. The residual
scalar equation for the second component has parameter
$\sigma_b=\rho_2-4\pi b$ and a weight with a quadratic zero at each
ordinary point. This equation and the equations for the centers form
the shadow system.

In follows, for $N\geq1$, we define the ordered configuration space as
\[
 \operatorname{Conf}_N(M)
 =\{(x_1,\ldots,x_N)\in M^N:x_i\ne x_j\text{ for }i\ne j\}.
\]
We let $\operatorname{Conf}_0(M)$ consist of the empty configuration.
We regard centers of each type as unordered. Their configuration space is
\[
 \mathcal C_{a,b}(M)
 =\operatorname{Conf}_{a+b}(M)/(S_a\times S_b).
\]
Here $S_a$ and $S_b$ permute the type $A$ and type $B$ centers,
respectively.

\subsection{The shadow degree and its orientation}\label{sec:degree-conventions}
We use the spaces
\[
 H^r_0(M)=\{v\in H^r(M):\textstyle\int_Mv\dV=0\},
 \qquad X_s=H^s_0(M),\quad s\ge6\ \text{an integer}.
\]
For the shadow degree, the scalar variable lies in $X_s$ and we
write its equation using the mean-zero Green inverse. We place the
scalar equation before the center equations and orient each center
using the orientation of $M$. Notice that any Riemannian surface is orientable (see \cite[I.2.1, p.~14]{FK92}). For a nonnegative weight $V$ of positive mass, put
$$
 \avg vV=\frac{\int_MVv\dV}{\int_MV\dV},\qquad
 \ipA\eta\xi=\int_M\nabla\eta^TA^{-1}\nabla\xi\dV.
$$
Here $A$ is the Cartan matrix in \eqref{1.1}.
The norm induced by $(\cdot,\cdot)_A$ is equivalent to the
$H^1$ norm on $\mathcal H$.
The Morse index of a quadratic form $\qform$ is the dimension of
its negative spectral space. The local index of a nondegenerate
finite-dimensional field $F$ is $\operatorname{sgn}\det DF$.

We compute the Toda degree in $\mathcal H$. The centers are unordered within each type. Their permutation factor $a!b!$ is included in \eqref{1.5}.

We write $Z=(w,P,Q)$ for an unordered shadow solution. If the
scalar linearization at fixed centers is invertible, we use $w(P)$
for the local branch of solutions of the scalar equation
\eqref{2.7} given by the implicit function theorem. The scalar
equation depends on $P$ and does not depend on $Q$. We write
$$I(P):=\int_M H_Pe^{w(P)}\dV\,\,\mbox{ and } \,\,n_*^g:=\sigma_b H_Pe^{w(P)}/I(P)$$ for its residual density on
$M$. Its Euclidean density in a chart is $e^\phi n_*^g$. Write $m_*$ for the Morse index of the scalar linearization and
$F$ for the center field obtained by
substituting $w(P)$ into the center equations \eqref{2.8}. The symbol
$d_S(Z)$ denotes the local shadow degree, and $d^S_{a,b,k}$ the
degree for the entire type. We write $d_S$ when the solution is fixed.

With this orientation, we prove the local index formula
\begin{equation}\label{eq:shadow-local-index}
 d_S(Z)=(-1)^{m_*}\operatorname{sgn}\det DF(Z)
\end{equation}
below at a nondegenerate coupled shadow solution whose fixed-center scalar block is invertible. Here $DF(Z)$ is the derivative of the reduced field $F$ in the center variables $(P,Q)$ at the centers of $Z$. This formula is proved at the end of Section~\ref{sec:limiting-equations}.

\subsection{Choice of weights and the shadow degree}\label{sec:wall-setting}
Write $\rho=\rho_2$ for the fixed second parameter and set
$\sigma_b=\rho-4\pi b$. With $w$ in place of $w_*$, the shadow
equations \eqref{2.7}--\eqref{2.8} take the form
\begin{equation}\label{3.40}
\left\{
\begin{aligned}
 F_0&=-\Delta_gw-2\sigma_b
 \left(\frac{H_Pe^w}{\int_MH_Pe^w\dV}-1\right)=0,
 \qquad \int_Mw\dV=0,\\
 F_i^A&=\nabla g_1(p_i)-\frac12\nabla w(p_i)
 +8\pi\sum_{l\ne i}\nabla_xG(p_i,p_l)
 +12\pi\sum_j\nabla_xG(p_i,q_j)=0,\\
 F_j^B&=\nabla V(q_j)+12\pi\sum_i\nabla_xG(q_j,p_i)
 +24\pi\sum_{l\ne j}\nabla_xG(q_j,q_l)=0.
\end{aligned}
\right.
\end{equation}
The coordinate conventions are those of Section~\ref{sec:notation}. Notice that $F_0$ depends only on $w$ and $P$.s

At a mixed center $q_j$, choose an isothermal coordinate $z$
centered at $q_j$, with $g=e^\phi|dz|^2$ and set
\[
\begin{aligned}
 \psi_{1,j}(z)
 &:=\log h_1(z)+\phi(z)+12\pi R(z,q_j)
+8\pi\sum_iG(z,p_i)+12\pi\sum_{\ell\ne j}G(z,q_\ell)-\tfrac12w(z),\\
 \mathcal T_j
 &:=\psi_{1,j,zz}(q_j)+\psi_{1,j,z}(q_j)^2,\\
 \kappa_m(q)
 &:=(\Delta_g\log h_1)(q)+8\pi m-\rho-2K_g(q).
\end{aligned}
\]
We use the following condition to exclude a collision of two
ordinary centers at $q_j$:
\begin{equation}\label{3.41}
 \kappa_m(q_j)>4e^{-\phi(q_j)}|\mathcal T_j|.
\end{equation}

\begin{theorem}\label{th3.14}
Fix integers $m\geq1$ and $k\geq0$ and let
$\rho\in(4\pi k,4\pi(k+1))$. There exists a pair of smooth positive
weights $(h_1,h_2)$ with the following properties for every type
\[
 a+2b=m,\qquad 0\leq b\leq k.
\]
\begin{enumerate}[label=\textup{(\roman*)},leftmargin=*]
\item The solution set of the shadow system \eqref{3.40} is compact,
with a priori bounds uniform over these types.
\item The solution set consists of finitely many nondegenerate
solutions. At each solution, the scalar linearization with the
centers fixed is invertible.
\item There is a common constant $\delta>0$ such that \eqref{3.41}
holds in the form
\[
 \kappa_m(q_j)-4e^{-\phi(q_j)}|\mathcal T_j|\geq\delta
\]
at every mixed center of every solution.
\end{enumerate}
\end{theorem}

The proof of Theorem~\ref{th3.14} gives a single choice of weights for
all admissible types. For these weights, the compact shadow homotopy
gives the following degree formula. We will prove Theorem \ref{th3.14} in Subsection \ref{subsection:ProofofTheorem2.6}.

With the help of Theorem \ref{th3.14}, we show that

\begin{theorem}\label{th3.15}
For the weights of Theorem~\ref{th3.14},
\begin{equation}\label{3.47}
 d^S_{a,b,k}
 =\binom{\chi}{a+b}\binom{a+b}{a}
 \sum_{j=0}^a\binom aj b_{k-b-j},
 \qquad \chi=\chi(M).
\end{equation}
\end{theorem}

The proof is postponed to Subsection \ref{subsection:ProofofTheorem2.6}.

We use $\chi$, $b_j$ and $D(x,y)$ as in Theorem~\ref{th1.1}.
All generating functions are formal power series. The generalized
binomial theorem gives
\begin{equation}\tag{1.5}\label{1.4}
 B_\chi(y)=\sum_{k\geq0}b_ky^k=(1-y)^{\chi-1}.
\end{equation}
Here $[y^n]F(y)$ denotes the coefficient of $y^n$
in the formal power series $F(y)$, with $[y^n]F(y)=0$
for $n<0$. The centers are unordered within each type.

For the weights chosen in Theorem~\ref{th3.14}, the shadow degree
for $a+2b=m$ is given by Theorem~\ref{th3.15} as
\begin{equation}\label{1.5}
 d^S_{a,b,k}
 =\binom{\chi}{a+b}\binom{a+b}{a}
 [y^{k-b}]\bigl(B_\chi(y)(1+y)^a\bigr).
\end{equation}

\subsection{Scales and normalization at a wall}\label{sec:scale-conventions}

Fix integers $m\geq1$ and $k\geq0$ and fix
$\rho_2\in(4\pi k,4\pi(k+1))$. Recall that $P=(p_i)_{i=1}^a$ and
$Q=(q_j)_{j=1}^b$ are the ordinary and mixed centers, with
$a+2b=m$ and $b\leq k$. Let $Z^0=(w^0,P^0,Q^0)$ be a solution of \eqref{3.40}. Set
\[
 \sigma=\sigma_b=\rho_2-4\pi b,\qquad
 f^0=\frac{\sigma H_{P^0}e^{w^0}}
           {\int_M H_{P^0}e^{w^0}\dV}.
\]

With the centers fixed, the partial derivative of the scalar
component $F_0$ is
\begin{align}
 \partial_w F_0(w^0,P^0) & \colon X_s\longrightarrow H^{s-2}_0(M),\nonumber\\
 \partial_w F_0(w^0,P^0)[\xi]
 &=-\Delta_g\xi
  -2f^0\left(\xi-\frac1{\sigma}\int_M f^0\xi\dV\right).\nonumber
\end{align}
Assume that $\partial_w F_0(w^0,P^0)$ is invertible.
This derivative is a compact perturbation of $-\Delta_g$
between the indicated mean-zero spaces. It therefore has
Fredholm index zero and is invertible under this assumption.

By the implicit function theorem, there are neighborhoods
$\mathcal U\subset \operatorname{Conf}_a(M)$ of $P^0$ and $\mathcal V\subset X_s$
of $w^0$, and a $C^1$ map $P\mapsto w(P)$ from $\mathcal U$
to $\mathcal V$ with $w(P^0)=w^0$. The local solution set of
the scalar equation \eqref{2.7} is the graph
\[
 \{(v,P)\in\mathcal V\times\mathcal U:F_0(v,P)=0\}
 =\{(w(P),P):P\in\mathcal U\}.
\]
Since $F_0$ is independent of the positions $Q$ for fixed $b$, the same map $w(P)$ applies for all nearby $Q$. For $P$ near $P^0$, we have
\begin{equation}\label{7.1}
 \begin{gathered}
 I(P)=\int_M H_Pe^{w(P)}\dV,\qquad
 n_*^g(P)=\frac{\sigma H_Pe^{w(P)}}{I(P)},\\
 -\Delta_g w(P)=2(n_*^g(P)-\sigma),\qquad
 \int_M w(P)\dV=0.
 \end{gathered}
\end{equation}

Set $\rho^0=(4\pi m,\rho_2)$, $\mu=\rho_1-4\pi m$ and
$F_\mu=F_{(4\pi m+\mu,\rho_2)}$. Replacing $w_*$ by $w(P)$ in \eqref{2.6}, we define
$(u_1^0,u_2^0)$ by
\begin{equation}\label{7.2}
 \begin{aligned}
 u_1^0&=8\pi\sum_iG(\cdot,p_i)+12\pi\sum_jG(\cdot,q_j)
         -\tfrac12w(P),\\
 u_2^0&=w(P)-4\pi\sum_iG(\cdot,p_i).
 \end{aligned}
\end{equation}
For a common parameter $t>0$, set
\begin{equation}\label{7.4}
 c_1=2\log t,\qquad c_2=\log(\sigma/I(P)),\qquad
 \delta=\sqrt t,\qquad \ell=\log(1/t).
\end{equation}
Here $\delta$ is the reference outer scale. The outer densities
in an isothermal chart are $e^{c_i}h_i e^{u_i^0+\phi}$. In particular,
\[
 e^{c_2}h_2e^{u_2^0+\phi}=e^\phi n_*^g=:n_*^{\rm loc}.
\]
For $b=0$, the same normalization applies with $\sigma=\rho_2$.

At a mixed center $q_j$, use a centered chart $z$ and write $r=|z|$.
Let $B_j\in\C$ be the shape parameter in $z^2-\epsilon_jB_j$.
For $\epsilon_j,\beta_j>0$, the local model densities are
\begin{equation}\label{7.8}
 \begin{aligned}
 F_j&=\epsilon_j^2+\beta_j\epsilon_j^2r^2
                          +|z^2-\epsilon_jB_j|^2,\\
 Q_j&=\beta_j\epsilon_j^4+4\epsilon_j^2r^2
                          +\beta_j\epsilon_j^2|z^2+\epsilon_jB_j|^2,\\
 D_j^{\rm G}&=4\beta_j\epsilon_j^4,\qquad
 n_{1,j}=\frac{4Q_j}{F_j^2},\qquad
 n_{2,j}=\frac{4D_j^{\rm G}F_j}{Q_j^2},\\
 U_j&=(\log n_{1,j},\log n_{2,j}).
 \end{aligned}
\end{equation}

\begin{lemma*}[Matching of the local scales]
At each ordinary center $p_i$, take an isothermal chart with $\phi(p_i)=0$ and $\nabla\phi(p_i)=0$. The matched width
$s_i$ and its density are given by
\begin{equation}\label{7.5}
 \begin{aligned}
 \psi_i^A&=\log h_1+\phi+8\pi R(\cdot,p_i)
 +8\pi\sum_{l\ne i}G(\cdot,p_l)
 +12\pi\sum_jG(\cdot,q_j)-\tfrac12w(P),\\
 s_i&=\frac t2 e^{\psi_i^A(p_i)/2},\qquad
 m_{A,i}(z)=\frac{4s_i^2}{(s_i^2+|z|^2)^2}.
 \end{aligned}
\end{equation}

At each mixed center $q_j$, set
\begin{equation}\label{7.6}
 \psi_j^B=\log h_1+\phi+12\pi R(\cdot,q_j)
 +8\pi\sum_iG(\cdot,p_i)
 +12\pi\sum_{l\ne j}G(\cdot,q_l)-\tfrac12w(P).
\end{equation}
Take the centered isothermal chart with $\phi(q_j)=0$ and
$\partial_z\psi_j^B(q_j)=0$. The matched parameters are
\begin{equation}\label{7.7}
 \begin{gathered}
 \epsilon_j=\frac t4 e^{\psi_j^B(q_j)/2},\qquad
 \delta_j=\sqrt{\epsilon_j},\qquad
 \beta_j=n_*^{\rm loc}(q_j)>0,\\
 a_j^2=\frac{\beta_j\epsilon_j^2(1+|B_j|^2)}4.
 \end{gathered}
\end{equation}
Here $\delta_j$ and $a_j$ are the outer and inner widths,
respectively. The charts can be chosen smoothly in the centers.
For centers near $(P^0,Q^0)$ and $B_j$ in bounded sets, we have
$s_i\asymp t$, $\delta_j\asymp\sqrt t$, and $a_j\asymp t$
uniformly as $t\to0$.
\end{lemma*}

\begin{proof}
Rescale each centered chart so that $\phi(0)=0$. Under
$z=f(\zeta)=\zeta+c\zeta^2$, the metric factor becomes
$\widetilde\phi=\phi\circ f+2\log|f'|$. At $p_i$, take
$c=-\partial_z\phi(0)/2$ to obtain $\nabla\widetilde\phi(0)=0$.

At $q_j$, the first outer density is
$t^2e^{\psi_j^B(z)}|z|^{-6}$. In the new chart,
\[
 \widetilde\psi_j^B=\psi_j^B\circ f+2\log|f'|
                         -6\log|f(\zeta)/\zeta|,
 \qquad
 \partial_\zeta\widetilde\psi_j^B(0)=\partial_z\psi_j^B(0)-c.
\]
Take $c=\partial_z\psi_j^B(0)$. These changes depend smoothly
on the centers.

At $p_i$, the expansion of $G$ in Section~\ref{sec:notation}
gives the first outer density $t^2e^{\psi_i^A(z)}r^{-4}$.
As $r\to0$ with $s_i/r\to0$,
$m_{A,i}(z)\sim4s_i^2r^{-4}$. Matching the coefficients gives
$4s_i^2=t^2e^{\psi_i^A(p_i)}$.

At $q_j$, as $r\to0$ with $\sqrt{\epsilon_j}/r\to0$,
\eqref{7.8} gives
\[
 F_j\sim r^4,\qquad Q_j\sim4\epsilon_j^2r^2,\qquad
 n_{1,j}\sim16\epsilon_j^2r^{-6}.
\]
Thus $16\epsilon_j^2=t^2e^{\psi_j^B(q_j)}$. For fixed $z\ne0$
and fixed $\beta_j,B_j$, as $\epsilon_j\to0$ we have
\[
 n_{2,j}(z)\longrightarrow
       \frac{\beta_j}{(1+\beta_jr^2/4)^2}.
\]
Match the value at zero of this limiting profile to
$n_*^{\rm loc}(q_j)$, which is positive since $q_j\notin P$.
Normalize the inner width by $n_{2,j}(0)=4/a_j^2$. Evaluating
\eqref{7.8} at zero gives
\[
 \frac4{a_j^2}
 =\frac{16}{\beta_j\epsilon_j^2(1+|B_j|^2)}.
\]
The uniform scale comparisons follow from \eqref{7.5} and \eqref{7.7}.
\end{proof}

For $b=0$, write $\psi_i=\psi_i^A$ and $d_i(P)=s_i/t$.
At a mixed center, write $d_j(P,Q)=\epsilon_j/t$.
In single-point formulas, we omit the center subscript on
the mixed parameters and widths. The spectral calculation
uses a reference scale $s$ comparable to $t$.

In each centered chart, write $r=|z|$ and set
\begin{equation}\label{7.9}
 \begin{aligned}
 D_i(z)&=-2\log\left(1+\frac{s_i^2}{r^2}\right),\\
 L_{1,j}(z)&=2\log\epsilon_j+\log16-6\log r
             +\log\left(1+\frac{\beta_jr^2}{4}\right),\\
 L_{2,j}(z)&=\log\beta_j
             -2\log\left(1+\frac{\beta_jr^2}{4}\right),\\
 L_j&=(L_{1,j},L_{2,j})^T,\qquad
 n_{{\rm diff},j}=e^{L_{2,j}}.
 \end{aligned}
\end{equation}
Choose radial cutoffs $\chi_i^A$ and $\chi_j^B$ in the centered
charts about $p_i$ and $q_j$, respectively. Their supports are
pairwise disjoint, and each cutoff equals one on a smaller disk.
Define the ansatz and its mean-zero version by
\begin{equation}\label{eq:ansatz}
 \begin{aligned}
 v^{\rm app}
  &=\binom{u_1^0}{u_2^0}+\sum_i\chi_i^A\binom{D_i}{-D_i/2}
       +\sum_j\chi_j^B(U_j-L_j),\\
 u^{\rm app}
 &=v^{\rm app}-\int_Mv^{\rm app}\dV.
 \end{aligned}
\end{equation}
We seek a solution of the form $u=u^{\rm app}+\eta$.

For an all-ordinary shadow solution $Z^0=(w^0,P^0,\varnothing)$,
use $f^{\rho^0}$ from \eqref{2.2} and the chart of \eqref{7.5} to set
\begin{equation}\label{8.19}
 \kappa_i(P)=-f_1^{\rho^0}(p_i)
 =\Delta_g\log h_1(p_i)+8\pi m-\rho_2-2K_g(p_i).
\end{equation}
Define the logarithmic mass coefficient by
\begin{equation}\label{8.23}
 \begin{aligned}
 l(P,w(P))&=\sum_i e^{\psi_i(p_i)}\kappa_i(P),\\
 L_0&=\frac\pi2 l(P^0,w^0)
     =2\pi\sum_i d_i(P^0)^2\kappa_i(P^0).
 \end{aligned}
\end{equation}

\begin{lemma}\label{le8.4}
For weights as in Theorem~\ref{th3.14}, an arbitrarily small
perturbation in any fixed smooth norm gives $l\ne0$ at every
shadow solution of type $(m,0)$. The shadow solutions, their
nondegeneracy, and \eqref{3.41} at the mixed centers can be preserved.
\end{lemma}

The proof is given in Section~\ref{sec:ordinary-mass}.
We choose the weights of Theorem~\ref{th3.14} with the additional
nonvanishing condition of Lemma~\ref{le8.4}. Homotopy invariance
at noncritical parameters gives the same degree for the original
weights.

Lemma~\ref{le2.1} gives the weak limits of the Toda densities
and the shadow equations for every concentrating sequence.
With the chosen weights, each limiting shadow solution
satisfies \eqref{3.41} at its mixed centers.

\subsection{Actual sequences and the two local degrees}\label{sec:two-local-degrees}
\subsubsection{Configurations with mixed points}
Use the wall and scale conventions of Section~\ref{sec:scale-conventions}.
For $b\geq1$, the strict inequality from the choice of weights
gives \eqref{7.35} in the normalized charts \eqref{7.6}--\eqref{7.7}.
With fixed noncritical $\rho_2$ and invertible fixed-center scalar
linearization, Theorem~\ref{th7.7} places each actual sequence with
\eqref{2.4}--\eqref{2.6} in the family of Proposition~\ref{pr7.2}
after modulation. The localization estimates \eqref{7.50}--\eqref{7.53}
and \eqref{7.62} also hold for $\mu=0$.

At each mixed center $q_j$, the local degree calculation uses the
following quantities in the normalized chart
\begin{equation}\label{7.21}
 \begin{split}
 \kappa_j&=\Delta_g\log h_1(q_j)+8\pi m-\rho_2-2K_g(q_j),\\
 T_j&=\partial_z^2\psi^B_j(q_j),\qquad
 \mathbf t_j=(\Re T_j,-\Im T_j),\\
 \mathcal M(B)&=\int_{\C}\frac{|v+B|}{(1+|v|^2)^2}\dif A_v,\\
 \mathcal H(B)&=\mathcal M(B)-B\cdot\nabla\mathcal M(B).
 \end{split}
\end{equation}
The regular field $\psi^B_j$ in \eqref{7.6} is $\psi_{1,j}$ with
$w=w(P)$. In the normalized chart of \eqref{7.7}, we have
$\kappa_j=\kappa_m(q_j)$ and $T_j=\mathcal T_j$ in the notation of
\eqref{3.41}. The proof of Theorem~\ref{th7.4} gives the positive
integral representation of $\mathcal H(B)$ and the local degree formula.

\begin{corollary}
\label{co7.8}
Under the assumptions of Theorem~\ref{th7.7}, suppose the limiting
solution of \eqref{3.40} is nondegenerate. Let $B_{j,*}$ solve
\eqref{4.33} in the charts of \eqref{7.7} and put
$$
 S_* =\sum_{j=1}^b d_j(P^0,Q^0)\kappa_j\mathcal H(B_{j,*})>0,
 \qquad \epsilon_j=td_j(P,Q).
$$
Every Toda sequence concentrating at this shadow solution satisfies
\begin{equation}\label{7.63}
 \frac\mu t\longrightarrow S_*,\qquad
 \|\eta\|_{H^1}=O(t^{2/3}|\log t|).
\end{equation}
Concentration requires $\mu>0$. For small $\mu>0$, all sufficiently
concentrated solutions lie in a bounded open set on whose boundary the
Toda map has no zero. Its degree is
\begin{equation}\label{7.64}
 (-1)^ad_S.
\end{equation}
The estimates are uniform over a finite set of shadow solutions.
\end{corollary}

\subsubsection{All-ordinary configurations}
Fix a separated nondegenerate ordinary shadow solution $(P^0,w^0)$
with invertible fixed-center scalar linearization. Use the chart and
scale conventions of Section~\ref{sec:scale-conventions} and the
projected family of Theorem~\ref{th8.2}. The next theorem gives
the local degree of this family.

\begin{theorem}
\label{th8.5}
Assume $l(P^0,w^0)\ne0$. For small $\mu\ne0$, the bubbling
neighborhood contains solutions only for $\mu L_0>0$. On that side,
\begin{equation}\label{8.24}
 \deg_{\mathrm{LS},\mathrm{loc}}=(-1)^m\operatorname{sgn}(\mu)d_S.
\end{equation}
The positive-side degree minus the negative-side degree is
\begin{equation}\label{8.25}
 (-1)^m d_S.
\end{equation}
The localization result below applies these conclusions to every
ordinary blow-up sequence with this shadow limit.
\end{theorem}

For an actual sequence with \eqref{8.28}, Theorem~\ref{th8.9} gives
\eqref{8.39}--\eqref{8.40}. After modulation, the sequence belongs
to the projected family, including at $\mu=0$.
Proposition~\ref{pr8.10} gives the zero-multiplier mass expansion
\eqref{8.43} and the scale comparison \eqref{8.44} on the allowed
side. Together with the local degree computation, these results
give the following statement about the degree neighborhood.

\begin{corollary}
\label{co8.11}
Fix a nondegenerate ordinary shadow solution with invertible
scalar linearization at fixed centers and $l\ne0$. Concentration
at this solution requires $\mu L_0>0$. After modulation, every
such sequence eventually lies in the neighborhood of
Theorem~\ref{th8.5}.

With the wall map $F_\mu$ of Section~\ref{sec:scale-conventions},
write $Z^0=(w^0,P^0,\varnothing)$ and let $\mathcal U_{Z^0,\mu}$
be this bounded open neighborhood for $\mu L_0>0$ and the
empty set otherwise. Its boundary is zero-free, and
\begin{equation}\label{8.45}
 \LS(F_\epsilon,\mathcal U_{Z^0,\epsilon},0)
 -\LS(F_{-\epsilon},\mathcal U_{Z^0,-\epsilon},0)
 =(-1)^m d_S,\qquad 0<\epsilon<\epsilon_0.
\end{equation}
The estimates are uniform over finitely many shadow solutions.
\end{corollary}

\subsubsection{The local jump and exclusion of wall sequences}
The two local results can now be compared at the same wall.

\begin{theorem}
\label{th9.1}
Use the chosen weights and the wall map $F_\mu$ of
Section~\ref{sec:scale-conventions}.
Then the following hold.
\begin{enumerate}[label=\textup{(\roman*)},leftmargin=*]
\item Every unbounded sequence of solutions with $\mu\to0$ has
a subsequence concentrating at a shadow solution $Z$ of admissible
type $(a,b)$.
\item For small $\mu\ne0$, there are bounded open sets
$\mathcal U_{Z,\mu}$, possibly empty, such that
$$
 0\notin F_\mu(\partial\mathcal U_{Z,\mu}),\qquad
 u_\nu\text{ concentrates at }Z
 \ \Longrightarrow\ u_\nu\in\mathcal U_{Z,\mu_\nu}\text{ eventually},
$$
and
\begin{equation}\label{9.1}
 \deg(F_\epsilon,\mathcal U_{Z,\epsilon},0)
 -\deg(F_{-\epsilon},\mathcal U_{Z,-\epsilon},0)=(-1)^ad_S(Z).
\end{equation}
\item The solution set $F_0^{-1}(0)$ is bounded.
\end{enumerate}
These assertions hold uniformly over the finitely many shadow solutions.
\end{theorem}

\begin{proof}
\textup{(i)} follows from Lemma~\ref{le2.1}. The translation
Pohozaev identity gives the center equations \eqref{2.8}.

\textup{(ii)} Corollaries~\ref{co7.8} and \ref{co8.11} give the
bounded open sets with zero-free boundaries. If a solution sequence
with $\mu_\nu\to0$ and $\mu_\nu\ne0$ concentrates at $Z$, it eventually
lies in $\mathcal U_{Z,\mu_\nu}$. Set $\mathcal U_{Z,\mu}=\varnothing$
on the side where concentration at $Z$ is excluded.

Fix $Z$ of type $(a,b)$ and take $\epsilon>0$ sufficiently small. Set
\[
 d_T^s(Z)=\deg(F_{s\epsilon},\mathcal U_{Z,s\epsilon},0)
 \quad(s=\pm1),\qquad C_Z=(-1)^a d_S(Z).
\]
Write $d_T^\pm=d_T^{\pm1}$ and let $d_T$ denote the degree on the
contributing side. If $b\geq1$, Corollary~\ref{co7.8} gives
\[
 d_T^+(Z)=C_Z,\qquad d_T^-(Z)=0.
\]
If $b=0$, then $a=m$ and $L_0\ne0$. Set $s=\operatorname{sgn}L_0$.
Theorem~\ref{th8.5} and Corollary~\ref{co8.11} give
\[
 d_T^s(Z)=sC_Z,\qquad d_T^{-s}(Z)=0.
\]
Thus $d_T^+(Z)-d_T^-(Z)=s(sC_Z)=C_Z$ in the ordinary case. This
proves \eqref{9.1} in both cases.

\textup{(iii)} Suppose $F_0^{-1}(0)$ is unbounded. By \textup{(i)},
a subsequence of solutions at $\mu=0$ concentrates at a shadow
solution $Z$ of type $(a,b)$, with scale $t>0$ and $t\to0$.
Here $\rho_2$ is fixed and noncritical. For the chosen weights,
$Z$ is separated and nondegenerate. Its scalar linearization at
fixed centers is invertible.

If $b\geq1$, the sequence satisfies \eqref{2.4}--\eqref{2.6}, and
the limiting mixed centers satisfy \eqref{7.35}.
Theorem~\ref{th7.7} puts the sequence in the projected family after
modulation. The estimates apply since $|\mu|=0\leq C_0t$.
The centers and the parameters $B_j$ converge to their limiting
values. All multipliers vanish. Equation~\eqref{7.25} gives
\[
 \mu=tS_*+o(t),\qquad S_*>0,
\]
where $S_*$ is as in Corollary~\ref{co7.8}. At $\mu=0$, division
by $t$ gives $0=S_*+o(1)$, a contradiction.

If $b=0$, the sequence satisfies \eqref{8.28} and $L_0\ne0$.
Theorem~\ref{th8.9} puts it in the projected family also at $\mu=0$.
Apply \eqref{8.18} to the actual solutions with zero sources, as in
the proof of Proposition~\ref{pr8.10}. This gives \eqref{8.43},
\[
 \mu=L_0t^2\log(1/t)+O(t^2).
\]
At $\mu=0$, division by $t^2\log(1/t)>0$ gives
$0=L_0+O(1/\log(1/t))$. This contradicts $L_0\ne0$ and proves
\textup{(iii)}. The neighborhoods and estimates can be chosen
uniformly over the finitely many shadow solutions.
\end{proof}

After the additional perturbation of Lemma~\ref{le8.4},
Theorem~\ref{th9.1} gives, at each shadow solution $Z$,
\begin{equation}\label{1.6}
 d_T^+(Z)-d_T^-(Z)=(-1)^a d_S(Z).
\end{equation}
\subsection{Separation, excision and the degree jump}
\begin{samepage}
\begin{proposition}
\label{pr9.3}
For the chosen weights there is a bounded open set $\mathcal B$
such that $0\notin F_\mu(\partial\mathcal B)$ for small $|\mu|$.
For each sign and every sufficiently small $\epsilon>0$, there are
pairwise disjoint bounded open sets $\mathcal U_{Z,\pm\epsilon}$
outside $\mathcal B$ such that
\begin{enumerate}[label=\textnormal{(L\arabic*)}]
\item $F_{\pm\epsilon}^{-1}(0)\subset
       \mathcal B\cup\bigcup_Z\mathcal U_{Z,\pm\epsilon}$;
\item $\deg(F_\epsilon,\mathcal U_{Z,\epsilon},0)
       -\deg(F_{-\epsilon},\mathcal U_{Z,-\epsilon},0)=(-1)^ad_S(Z)$.
\end{enumerate}
\end{proposition}
\end{samepage}
\begin{proof}
\textit{Step 1.} We first choose $\mathcal B$.
By Theorem~\ref{th9.1}, $F_0^{-1}(0)$ is bounded. Since $K$ is compact
on bounded sets, $F_0^{-1}(0)$ is compact. Then, we can set
\[
 \mathcal B:=\{u\in\mathcal H:\|u\|_A<R\}\mbox{ with }
 \|u\|_A^2
 =\int_M\sum_{i,j=1}^2 (A^{-1})_{ij}
   \langle\nabla_g u_i,\nabla_g u_j\rangle_g\,dV_g
\]
and $R>0$ sufficiently large. No we prove $0\notin F_\mu(\partial\mathcal B)$. Otherwise, $\partial\mathcal B$ contained zeros of $F_\mu$ along a sequence
$\mu\to0$, a subsequence would converge to a zero of $F_0$ on
$\partial\mathcal B$. This contradicts the choice of $\mathcal B$.

\textit{Step 2.} We separate $\mathcal B$ into three parts and compute the degree on them respectively.

\noindent\textit{(i) Convergence to the outer pair.}
Shrink the center domains in Corollaries~\ref{co7.8} and \ref{co8.11}
so that the corresponding scalar branch graphs have compact, pairwise
disjoint closures in the space of unordered triples $(w,P,Q)$. 
Let $v_\nu^{\rm app}$ be $v^{\rm app}$ in \eqref{eq:ansatz}
evaluated at the parameters of $u_\nu$. 
Fix $Z$ and take $u_\nu=u_\nu^{\rm app}+\eta_\nu
\in\mathcal U_{Z,\mu_\nu}$ with $\mu_\nu\to0$. After passing to a
subsequence, we have
\[
 t_\nu\to0,\qquad \|\eta_\nu\|_{H^1}\to0,\qquad
 (w_\nu,P_\nu,Q_\nu)\to(\bar w,\bar P,\bar Q).
\]
The limiting triple belongs to the corresponding graph closure.
Fix $r_0,s_0,\epsilon_0>0$ and write $r=|z|$ in the moving charts.
At an ordinary center, \eqref{7.9} gives
\[
 |D_i|=2\log(1+s_i^2/r^2)\le C(1+|\log r|)
 \qquad (0<r<r_0,\quad 0<s_i\le s_0).
\]
At a mixed center, fix $B_0>0$ and $0<\beta_-\le\beta_+<\infty$.
For
\[
 0<r<r_0,\qquad 0<\epsilon_j\le\epsilon_0,\qquad
 |B_j|\le B_0,\qquad \beta_-\le\beta_j\le\beta_+,
\]
the formulas in \eqref{7.8} give the uniform bounds
\[
 F_j\asymp(\epsilon_j+r^2)^2,\qquad
 Q_j/\epsilon_j^2\asymp r^2+\epsilon_j^2.
\]
By \eqref{7.8} and \eqref{7.9},
\[
 U_j-L_j
 =-\sum_{i=1}^2
 \left(
   a_{i1}\log\frac{F_j}{r^4}
   +a_{i2}\log\frac{Q_j}{\epsilon_j^2r^2(4+\beta_j r^2)}
 \right)e_i.
\]
Both ratios have a uniform positive lower bound. They are bounded
above by fixed powers of $1+r^{-1}$. Hence
$|U_j-L_j|\le C(1+|\log r|)$ under the same conditions.
As $t_\nu\to0$, the local corrections in \eqref{eq:ansatz}
tend to zero uniformly on compact subsets away from the limiting
centers. Dominated convergence in the moving charts gives
\[
 v_\nu^{\rm app}-u_\nu^0\longrightarrow 0
 \qquad\text{in }L^1(M).
\]
The outer pair $u_\nu^0$ has zero mean. Thus
\[
 \int_M v_\nu^{\rm app}\dV
 =\int_M(v_\nu^{\rm app}-u_\nu^0)\dV\longrightarrow 0.
\]
Let $\bar u^0$ be the pair in \eqref{2.6} for
$(\bar w,\bar P,\bar Q)$. It follows that
\[
 u_\nu\longrightarrow\bar u^0
 \qquad\text{in }L^2_{\rm loc}(M\setminus(\bar P\cup\bar Q)).
\]

\noindent\textit{(ii) Pairwise disjointness of the local domains.}
Suppose two distinct local domains intersect along a sequence
$\mu_\nu\to0$. Apply \textup{(i)} to a sequence in their intersection.
The two outer pairs agree away from both sets of limiting centers.
The coefficients $-4$ and $-6$ of the logarithmic terms in the first
component determine the ordinary and mixed centers, respectively.
The second component then determines the scalar limit. Thus
\[
 (\bar w_1,\bar P_1,\bar Q_1)=(\bar w_2,\bar P_2,\bar Q_2).
\]
This contradicts the disjointness of the two graph closures.

\noindent\textit{(iii) Separation from $\mathcal B$.}
Suppose $u_\nu\in\mathcal U_{Z,\mu_\nu}\cap\mathcal B$ along a
sequence $\mu_\nu\to0$. Since $\mathcal B$ is bounded in $\mathcal H$,
we may pass to a subsequence. By \textup{(i)},
\[
 u_\nu\rightharpoonup u_*\quad\text{in }\mathcal H,\qquad
 u_*=\bar u^0\quad\text{a.e.\ on }M\setminus(\bar P\cup\bar Q).
\]
Since $m\geq1$, the first component of $\bar u^0$ has a logarithmic
singularity. Hence $u_*\notin\mathcal H$, a contradiction. Therefore
\[
 \mathcal U_{Z,\mu}\cap\mathcal B=\varnothing
 \qquad\text{for all sufficiently small }|\mu|.
\]

\textit{Step 3.} We cover the zeros of $F_{\pm\epsilon}$.
Suppose zeros at $\mu=\pm\epsilon$ lie outside the union in
\textup{(L1)} along a sequence $\epsilon\downarrow0$. A bounded
subsequence has a further subsequence converging to a zero of $F_0$
in $\mathcal B$.

By Theorem~\ref{th9.1}, an unbounded subsequence has a further
subsequence in one of the sets $\mathcal U_{Z,\pm\epsilon}$.
Both cases contradict the choice of the sequence, proving
\textup{(L1)}. Equation~\eqref{9.1} gives \textup{(L2)}.
\end{proof}

\begin{theorem}\label{th9.4}
For every $m\ge1$ and $k\ge0$,
\begin{equation}\label{9.3}
 d_{m,k}-d_{m-1,k}
 =\sum_{a+2b=m}(-1)^a d^S_{a,b,k}.
\end{equation}
We set $d^S_{a,b,k}=0$ for $b>k$.
\end{theorem}
\begin{proof}
By Proposition~\ref{pr9.3}, additivity and homotopy invariance on
$\mathcal B$ give
$$
 \begin{aligned}
 d_{m,k}&=\deg(F_\epsilon,\mathcal B,0)
       +\sum_Z\deg(F_\epsilon,\mathcal U_{Z,\epsilon},0),\\
 d_{m-1,k}&=\deg(F_{-\epsilon},\mathcal B,0)
       +\sum_Z\deg(F_{-\epsilon},\mathcal U_{Z,-\epsilon},0),\\
 \deg(F_\epsilon,\mathcal B,0)&=\deg(F_{-\epsilon},\mathcal B,0).
 \end{aligned}
$$
Hence \eqref{9.1} gives
\[
 d_{m,k}-d_{m-1,k}
   =\sum_Z(-1)^{a(Z)}d_S(Z)
   =\sum_{a+2b=m}(-1)^a d^S_{a,b,k}.
\]
For each type, the sum of the local indices equals the shadow
degree on $X_s$ by \eqref{eq:shadow-local-index},
\eqref{eq:shadow-tubular-homotopy} and additivity. The value
in \eqref{3.47} is the degree on this same space. The comparison
in \eqref{eq:scalar-degree-space-comparison} uses the same frozen
operator in two Sobolev spaces. The identity $G_{\rm CL}=G$ identifies
the Chen--Lin weight, which has a quadratic zero at each ordinary
mark. This yields the scalar degree \eqref{eq:frozen-scalar-degree}.
Combining it with the configuration factor in \eqref{3.47} identifies
the jump with the explicit sum over types used in the generating function.
The shadow degree is unchanged by the perturbation in
Lemma~\ref{le8.4}. Indeed, uniform compactness of the shadow solution
sets in the weight neighborhood of Theorem~\ref{th3.14} places any
sequence of zeros for convergent weights in the finite union of the
continuation neighborhoods. There the implicit function theorem gives
exactly one continuation of each zero and preserves its local sign.
The same argument applies along the weight segment of the perturbation.
The sum over $Z$ counts each unordered shadow solution once.
The permutation factors have already been included in the
configuration Euler characteristic in \eqref{3.47}.
The admissible weight path within each noncritical block,
established in Section~\ref{sec:noncritical-degree}, identifies both Toda degrees
with those for the original weights. The auxiliary weights may
therefore be chosen separately for each wall.
\end{proof}

\subsection{The degree below the first critical value}
\begin{lemma}\label{le9.2}
For every $k\geq0$, $d_{0,k}=b_k$.
\end{lemma}
\begin{proof}
Fix $r=\rho_2$ and $0\leq\varepsilon=\rho_1\leq s<4\pi$. Put
$$
 f_i=\frac{\rho_i h_ie^{u_i}}{\int_Mh_ie^{u_i}\dV},\quad
 v_i=G*f_i,\qquad (u_1,u_2)=(2v_1-v_2,2v_2-v_1).
$$
The lower bound for $G$ and Jensen's inequality give
$v_2\geq-Cr$ and $f_1\leq C\varepsilon e^{2v_1}$.
Choose $p>2$ and $p\varepsilon_0<2\pi$. For $0<\varepsilon\leq\varepsilon_0$,
$$
 \begin{aligned}
 \int_M e^{2pv_1}\dV
 &\leq\sup_y\int_Me^{2p\varepsilon G(x,y)}\dV(x)\leq C,\\
 \|f_1\|_p&\leq C\varepsilon,
 \qquad \|v_1\|_{W^{2,p}}\leq C\varepsilon,\\
 -\Delta_gw&=2r\left(
 \frac{h_2e^{-v_1}e^w}{\int_Mh_2e^{-v_1}e^w\dV}-1\right),
 \qquad w=2v_2,\\
 h_2e^{-v_1}&\longrightarrow h_2\quad\text{in }C^{1,\alpha},
 \qquad 2r\notin8\pi\mathbb N.
 \end{aligned}
$$
For compactness near $\varepsilon=0$, the scalar weights
$h_2e^{-v_1}$ have a uniform positive lower bound, are uniformly
Lipschitz, and converge in $C^1$ to $h_2$. If $w$ were unbounded
along a sequence with $\varepsilon\downarrow0$, the scalar Green
estimate would give finitely many concentration points and bounded
oscillation on isolating circles. At each point, subtract the
bounded smooth Poisson potential for the scalar background and
apply Proposition~\ref{pr3.13} with $J=0$. In the $2r$ scalar
normalization, each positive point mass is $8\pi$. The circular
mean at such a point has leading coefficient $4$, so the
normalization integral diverges by the same circular Jensen
argument as in Section~\ref{sec:noncritical-degree}. All the mass would concentrate,
contrary to $2r\notin8\pi\mathbb N$. Thus $w$ stays bounded.
Compactness inside the first noncritical block applies when
$\varepsilon$ is bounded away from zero. Together these bounds
give a common degree domain for $\varepsilon\in[0,s]$.

At $\varepsilon=0$, let
$$
 T_r(v)=rG\left(\frac{h_2e^v}{\int_Mh_2e^v\dV}-1\right),
 \qquad S(u_1,u_2)=(u_1+u_2/2,u_2).
$$
Then
$$
 SF_{0,r}S^{-1}(a,b)=(a,b-2T_r(b)),\qquad
 d_{0,k}=\deg(I)\deg(I-2T_r)=b_k
$$
by \cite[Theorem~1.1]{CL15}. The scalar degree uses the same
mean-zero space and Green inverse, with positive weight $h_2$
and parameter $2r\in(8\pi k,8\pi(k+1))$. As in the proof of
Theorem~\ref{th3.15}, a mark of order zero gives the smooth
scalar case. The shear is a conjugacy by an invertible linear
map, so its domain and range orientations cancel in the degree.
\end{proof}

\subsection{Proof of the degree formula}\label{sec:main-proof}

\begin{proof}[Proof of Theorem \ref{th1.1}]
Section~\ref{sec:main} writes the Toda equation as a compact
perturbation of the identity on $\mathcal H$.
Section~\ref{sec:noncritical-degree} gives a priori bounds along
compact paths of parameters and weights within each noncritical
block. The degree is therefore defined on every sufficiently
large ball and is constant in the block and independent of the
smooth positive weights.
Lemma~\ref{le9.2} gives $d_{0,k}=b_k$ by an admissible homotopy to $\rho_1=0$.

Fix $m\geq1$, $k\geq0$, and
$\rho_2\in(4\pi k,4\pi(k+1))$, and consider
$F_\mu=F_{4\pi m+\mu,\rho_2}$ near the wall $\mu=0$.
By Lemma~\ref{le2.1}, an unbounded sequence of zeros with
$\mu\to0$ has a subsequence concentrating at a shadow solution
$Z=(w,P,Q)$ of type $(a,b)$. The first mass identity in
\eqref{2.4} gives $a+2b=m$. The residual mass
$\sigma_b=\rho_2-4\pi b>0$ gives $b\leq k$.
Theorem~\ref{th3.14} provides one choice of weights for all
these types, with finitely many nondegenerate shadow solutions,
invertible scalar blocks at fixed centers, and the strict mixed
inequality \eqref{3.41}. By Lemma~\ref{le8.4} we may also arrange
that $L_0\ne0$ at every all-ordinary shadow solution. The
continuation in the proof of Theorem~\ref{th9.4} preserves
the shadow degrees under this perturbation. The noncritical
weight homotopy allows us to compute the two adjacent Toda
degrees with this choice.

For $b\geq1$, Theorem~\ref{th7.7} and Corollary~\ref{co7.8}
place actual concentrating sequences in the mixed neighborhood
when $\mu\ne0$. Theorem~\ref{th7.4} gives local degree
$(-1)^a d_S(Z)$ on the positive side, and the negative side
contributes zero.
For $b=0$, Theorem~\ref{th8.9}, Proposition~\ref{pr8.10}, and
Corollary~\ref{co8.11} place actual sequences in the ordinary
neighborhood on the side $\mu L_0>0$. On that side,
Theorem~\ref{th8.5} gives local degree
$(-1)^m\operatorname{sgn}(\mu)d_S(Z)$; the opposite side
contributes zero. Thus the positive-side degree minus the
negative-side degree is $(-1)^a d_S(Z)$ in both cases, as in
\eqref{9.1}.

The localization results also apply at $\mu=0$. For actual
solutions the multipliers vanish. In the mixed case \eqref{7.25}
then gives $\mu/t\to S_*>0$, and in the ordinary case
\eqref{8.43} gives $\mu/(t^2\log(1/t))\to L_0\ne0$.
Both limits exclude concentration on the wall. This gives the
bounded wall zero set in Theorem~\ref{th9.1}.

Proposition~\ref{pr9.3} now gives the pairwise disjoint local
domains outside the bounded ball $\mathcal B$. Together with
$\mathcal B$ they contain all zeros on each side for sufficiently
small $|\mu|>0$. The boundary of $\mathcal B$ is zero-free
throughout the crossing. Its two degree contributions are
therefore equal. Additivity and excision leave the sum of the
local jumps over unordered shadow solutions. The indices
\eqref{eq:shadow-local-index} agree with the local degrees of the
shadow map \eqref{eq:shadow-tubular-homotopy}, so grouping the
sum by type gives \eqref{9.3}. Theorem~\ref{th3.15} evaluates
each type degree by \eqref{3.47}.

Set $J_0(y)=B_\chi(y)$ and
$J_m(y)=\sum_{k\geq0}(d_{m,k}-d_{m-1,k})y^k$ for $m\geq1$.
Then
$$
 (1-x)D(x,y)=\sum_{m\geq0}J_m(y)x^m.
$$
Substitute \eqref{3.47} into \eqref{9.3}. For fixed $a,b$, the
sum over $k$ starts at $k=b$ and gives
$y^bB_\chi(y)(1+y)^a$, since $b_{k-b-j}=0$ for negative indices.
Include $(a,b)=(0,0)$ through the initial row $J_0(y)$. Then
\begin{align*}
 (1-x)D(x,y)
 &=B_\chi(y)\sum_{a,b\geq0}
   \binom\chi{a+b}\binom{a+b}a
   (-x(1+y))^{a}(x^2y)^b\\
 &=B_\chi(y)(1-x(1+y)+x^2y)^{\chi}\\
 &=B_\chi(y)(1-x)^\chi(1-xy)^\chi.
\end{align*}
For the second equality, we first sum at fixed $a+b$ and then
apply the binomial theorem. The two substituted terms have positive
order in $x$, so each coefficient receives only finitely many
contributions, including when $\chi<0$.
The factor $1-x$ is invertible in $\mathbb Z[[x,y]]$.
Dividing by it gives \eqref{1.2}. In the product of the three
factors, the coefficient of $(xy)^j$ from $(1-xy)^\chi$ is
$(-1)^j\binom{\chi}{j}$, and the remaining coefficients are
$b_{m-j}$ and $b_{k-j}$ by \eqref{1.4}. Extracting coefficients
gives \eqref{1.3}, including on both coordinate axes.
\end{proof}

\subsubsection*{Technical ingredients of the degree computation}

The computation of the shadow degree in \eqref{3.47} uses the
homotopy of Proposition~\ref{pr3.8} and the local index formula
\eqref{eq:shadow-local-index}. The comparison
\eqref{eq:scalar-degree-space-comparison} identifies the degrees of the
frozen-center scalar map on the two Sobolev spaces. The identification
of the Green kernels in the proof of Theorem~\ref{th3.15} then gives
the scalar factor \eqref{eq:frozen-scalar-degree} from the degree
formula of Chen and Lin \cite[Theorem~1.1]{CL15}.

For the mixed case $b\geq1$, Proposition~\ref{pr4.19} gives the
profiles and annular bounds used in Section~\ref{rev:sec:matching}.
The proof uses the center selection in
\eqref{eq:prescribed-center-selection} and the exhaustion of the local
masses by the two profiles. The inner-scale smallness estimate
\eqref{eq:inner-first-density-small} is used in \eqref{5.24}
and subsequently in \eqref{5.56}. Theorem~\ref{th6.1} gives the model spectral and
inverse bounds used in both local constructions. In the mixed case,
the comparison with the completed weights in the rough parameter
range uses \eqref{eq:completed-weight-capacity}. The weighted
dilation moment \eqref{eq:mixed-dilation-weighted-moment} and the
cross-term estimates \eqref{eq:mixed-dilation-cross-terms} enter
the bound \eqref{7.61}.

Theorems~\ref{th7.7} and~\ref{th8.9} locate actual concentrating
sequences in the mixed and all-ordinary cases, respectively.
These estimates, together with the parameter modulation of
Lemma~\ref{lem:parameter-modulation}, give the interior localization
on the contributing sides stated in Corollaries~\ref{co7.8}
and~\ref{co8.11}. With the multipliers set to zero,
\eqref{7.25} excludes mixed concentration on the wall.
In the all-ordinary case, Proposition~\ref{pr8.10} gives
the corresponding exclusion through \eqref{8.43}.

The local blow-up analysis uses the selection, quantization, and
Pohozaev results of
\cite[Proposition~3A, Theorem~1.3 and \textup{(1-9)}]{LWYZ18}
and the scalar concentration and least-mass results of
\cite[Propositions~2.1, 3.1 and Theorems~A, 1]{Tar23}. The
entire profile classification and the scalar kernel description
come from \cite[Theorems~1.1 and~1.3]{LWY} and
\cite[Theorem~1]{DEM}, respectively.

\subsection{Consequences}
At $\rho_1=8\pi$, the contribution of type $(2,0)$ and, for $k\geq1$,
that of type $(0,1)$ give
$$
 d_{2,k}-d_{1,k}
 =\binom\chi2(b_k+2b_{k-1}+b_{k-2})+\chi b_{k-1}.
$$
The last summand comes from the mixed type $(0,1)$.

\begin{corollary}\label{co9.5}
$d_{m,k}>0$ on every surface of
positive genus. Thus the system has a solution for all positive
noncritical parameters and smooth positive weights.
\end{corollary}
\begin{proof}
If $\chi=-q<0$, each term in
$$
 d_{m,k}=\sum_{j=0}^{\min(m,k)}
 \binom{q+j-1}{j}\binom{m-j+q}{q}\binom{k-j+q}{q}
$$
is positive. If $\chi=0$, every coefficient of $D$ is one.
The nonzero degree yields a solution.
\end{proof}

\begin{corollary}\label{co9.6}
For a torus, $d_{m,k}=1$ in every noncritical block.
\end{corollary}
This follows from \eqref{1.2} with $\chi=0$.

\section{The blowup limit and the shadow equations}\label{rev:sec:scalar}

 \subsection{Derivation of the shadow system}
 \label{sec:limiting-equations}
\begin{proof}[Proof of Lemma~\ref{le2.1}]
\textit{Step 1.}
Pass to a subsequence with $\|u_\nu\|_{\mathcal H}\to\infty$
and apply Lemma~\ref{lem:toda-concentration} with the fixed weights
$h_i$. If $I_{2,\nu}$ diverged along a subsequence, then $\mu_2$
would be supported on $S$. By \eqref{2.3}, this would give
$\rho_2\in4\pi\mathbb N$. Hence $I_{2,\nu}$ is bounded above
and bounded away from zero.

At any $q\in S$, one of the two coefficients $\gamma_i(q)$
is greater than two. The bound for $I_{2,\nu}$ excludes $i=2$,
so the first component concentrates completely. After passing
to a subsequence, set
\[
 I_{2,\nu}\to I_{2,*}\in(0,\infty),\qquad
 L:=\frac{\rho_2}{I_{2,*}}>0,\qquad
 n_*:=Lh_2 e^{u_{2,*}}\quad\text{on }M\setminus S.
\]
By local convergence, $\mu_2=n_*\dV$ on $M\setminus S$. Since $S$ is finite, $n_*\dV$ is the absolutely continuous part
of $\mu_2$. In particular, $n_*\in L^1(M)$.

Integrability now rules out the mass pairs $(0,4\pi)$, $(4\pi,8\pi)$
and $(8\pi,8\pi)$. For $i=1,2$, let $u_{i,*}$ denote the local
limit of $u_{i,\nu}$ away from the concentration set.
Near a blow-up point $q$ with local masses $(m_1,m_2)$,
the Green representation and the lower bound for $G$ give
\[
 u_{2,*}(x)\geq
 -\frac{2m_2-m_1}{2\pi}\log d_g(x,q)-C.
\]
Indeed, the residual potential $2Gn_*$ is bounded below since
$G\geq-C$ and $n_*\geq0$.
Since $n_*=Lh_2e^{u_{2,*}}$ with $L>0$, we obtain
$n_*(x)\geq c\,d_g(x,q)^{-\alpha}$, where $\alpha=4,6,4$
for $(m_1,m_2)=(0,4\pi),(4\pi,8\pi),(8\pi,8\pi)$,
respectively, and $c>0$.
Each case contradicts $n_*\in L^1(M)$.
Thus only types $A$ and $B$ remain.
The local masses at the type $A$ and $B$ points give \eqref{2.4},
with $\sigma=\int_M n_*\dV>0$.

Set $w_*=2Gn_*$.
Green's representation yields \eqref{2.6}. Its second line gives
$n_*=L H_Pe^{w_*}$. Since $\int_M n_*\dV=\sigma$, we obtain
\eqref{2.5}.
The normalization of $G$ also gives \eqref{2.7}.
The weight $H_P$ extends smoothly to $M$, vanishes quadratically
at each type $A$ point, and is positive at each type $B$ point.
Since $n_*\in L^1(M)$, local exponential estimates followed by
elliptic bootstrapping yield $w_*\in C^\infty(M)$.

\medskip
\noindent\textit{Step 2.}
Fix a concentration point $q$ with local masses $(m_1,m_2)$.
Use the local notation of \eqref{2.2}. Choose the coordinate $z$
so that $z(q)=0$. The regular part from Section~\ref{sec:notation}
then gives
\[
 G(z,\zeta)=-\frac{1}{2\pi}\log|z-\zeta|+R(z,\zeta).
\]
By the definition of the Robin function and the symmetry of $R$,
\[
 r_g(z)=R(z,z)+\frac{\phi(z)}{4\pi},
 \qquad
 \nabla r_g(q)
 =2\nabla_xR(q,q)+\frac{\nabla\phi(q)}{4\pi}.
\]
Set $c_1=2m_1-m_2$ and $c_2=2m_2-m_1$, and write
\[
 u_{i,*}(z)=-\frac{c_i}{2\pi}\log|z|+v_i(z),
 \qquad
 \psi_i=\log h_i+\phi+v_i,
 \qquad i=1,2.
\]
Here $v_i$ is smooth near $q$. Since
$U_{i,\nu}-(u_{i,\nu}+\log h_i+\phi)$ is constant in space,
the spatial derivatives of $U_{i,\nu}$ converge away from $q$
to those of $u_{i,*}+\log h_i+\phi$. The regular part of
this last function at $q$ is $\psi_i$.
Apply the local translation Pohozaev identity on a disk about
$q$ with the symmetric multiplier $C=A^{-1}$ in \eqref{2.2}.
The smooth forcing term vanishes in the small-disk limit.
The differentiated Green kernel and the total mass bound give
a bound for $\nabla U_\nu$ in $L^p$ for $1<p<2$, so the forcing
term is $O(r^{2-2/p})$. First let $\nu\to\infty$ at fixed radius.
The boundary density term then vanishes as $r\to0$, since the
limiting regular densities are bounded there.
Write the logarithmic coefficient vector as
$\alpha=A(m_1,m_2)^T/(2\pi)$, so the limiting gradient
of component $i$ is $-\alpha_i\boldsymbol n/r+\nabla\psi_i(q)+O(r)$.
The pure logarithmic stress has zero integral by angular symmetry.
In the mixed terms the two normal contractions cancel by
symmetry of $C$, leaving
$-2\pi\sum_{i,j}C_{ij}\alpha_i\nabla\psi_j(q)$.
Thus, using $CA=I$, the remaining identity is
\[
 m_1\nabla\psi_1(q)+m_2\nabla\psi_2(q)=0.
\]
By \eqref{2.4}, this balance reduces to $\nabla\psi_1(p_i)=0$
at a type $A$ point and to
$2\nabla\psi_1(q_j)+\nabla\psi_2(q_j)=0$ at a type $B$ point.

At a type $A$ point $p_i$, the self-interaction term in
$u_{1,*}$ is $8\pi G(\cdot,p_i)$. The Robin function identity gives
\[
 \nabla\bigl(\phi+8\pi R(\cdot,p_i)\bigr)(p_i)
 =4\pi\nabla r_g(p_i).
\]
Substituting the limit in \eqref{2.6} into
$\nabla\psi_1(p_i)=0$ therefore yields the type $A$ equation
in \eqref{2.8}.

At a type $B$ point $q_j$, the $w_*$ terms cancel in
$2u_{1,*}+u_{2,*}$. The self-interaction term in this combination
is $24\pi G(\cdot,q_j)$. Using
\[
 \nabla\bigl(3\phi+24\pi R(\cdot,q_j)\bigr)(q_j)
 =12\pi\nabla r_g(q_j)
\]
in the balance
$2\nabla\psi_1(q_j)+\nabla\psi_2(q_j)=0$
yields the type $B$ equation in \eqref{2.8}.
The Pohozaev balance uses the total local masses at $q_j$
and thus also applies when the mixed profile splits.
\end{proof}

Use the unordered configuration space $\mathcal C_{a,b}(M)$
of Section~\ref{sec:shadow-setting}, and use the orientation and
weighted averages of Section~\ref{sec:degree-conventions}.
Let $Z=(w,P,Q)$ be a solution of the shadow system
\eqref{2.7}--\eqref{2.8}, with $w$ in place of $w_*$, and set
$f=\sigma_bH_Pe^w/\int_MH_Pe^w\dV$. The derivative of the
scalar fixed-point equation in $X_s$ is $I-\mathcal B_f$, where
\[
 \mathcal B_fv=2G\bigl[f(v-\avg v f)\bigr].
\]
At a shadow solution, the coefficient $f$ is smooth. The operator
is compact and self-adjoint with respect to the Dirichlet inner
product on $H^1_0(M)$. Its bilinear form is
$2\int_M f(v-\avg v f)(\xi-\avg \xi f)\dV$. Thus the Morse index
$m_*$ of the scalar form is the total multiplicity of the
eigenvalues of $\mathcal B_f$ greater than one. The eigenvectors
for a nonzero eigenvalue are smooth, since $\mathcal B_f$ gains
two Sobolev derivatives. The generalized eigenspaces for such an
eigenvalue also agree in $H^1_0(M)$ and $X_s$. In a Jordan chain
$(\mathcal B_f-\lambda)v=v_0$ with $v_0$ smooth and $\lambda\ne0$,
the relation $v=\lambda^{-1}(\mathcal B_fv-v_0)$ allows us to
bootstrap the regularity of $v$.
Consequently $m_*$ gives the scalar Leray--Schauder sign in $X_s$
when $I-\mathcal B_f$ is invertible. Indeed, in the compact linear
sign formula of \cite[Theorem~8.10]{Deimling85}, the exponent
is the sum of algebraic multiplicities of real eigenvalues of the
compact operator greater than one. The common generalized eigenspaces
just established identify this exponent with $m_*$. The hypothesis of
invertibility excludes the eigenvalue one, and the identity-minus-compact
derivative can be used for the local sign by the linearization argument
in \cite[Theorem~8.1, \textup{(D3)}]{Deimling85}: the differentiable remainder is
$o(\|v\|_{X_s})$, whereas the invertible derivative has a positive
linear lower bound.

Write the derivative of the coupled fixed-point and center
equations in $X_s$ and local center coordinates as
$\left(\begin{smallmatrix}A_s&B_s\\ C_s&D_s\end{smallmatrix}\right)$,
with $A_s=I-\mathcal B_f$. If $A_s$ is invertible, the implicit
function theorem gives the scalar branch $w=w(P)$ and the
derivative of the reduced center field is
$D_{(P,Q)}F=D_s-C_sA_s^{-1}B_s$. In fact,
\[
 \begin{pmatrix}A_s&B_s\\C_s&D_s\end{pmatrix}
 =
 \begin{pmatrix}I&0\\C_sA_s^{-1}&I\end{pmatrix}
 \begin{pmatrix}A_s&0\\0&D_s-C_sA_s^{-1}B_s\end{pmatrix}
 \begin{pmatrix}I&A_s^{-1}B_s\\0&I\end{pmatrix}.
\]
At a nondegenerate solution the middle factor is invertible.
Scaling the off-diagonal blocks in the two triangular factors
to zero gives a path of invertible identity-minus-compact
operators. The local product formula gives
\eqref{eq:shadow-local-index}. We use the $H^1_0$ realization for
the scalar form. The point evaluations in the coupled derivative
act on $X_s$.

\subsection{Collisions among ordinary concentration points}
We consider the shadow system \eqref{2.7}--\eqref{2.8}
with $a=N\geq2$, $b=0$, and $\rho:=\sigma_0=\rho_2$,
where $\rho>0$ and $\rho\notin4\pi\mathbb N$.
The equations reduce to
\begin{equation}\label{2.9}
\begin{cases}
 \displaystyle \Delta_gw+2\rho\left(\frac{H_Pe^w}{\int_MH_Pe^w\dV}-1\right)=0,
 \quad \displaystyle\int_Mw\dV=0,\\[6pt]
 \displaystyle \grad g_1(p_i)-\frac12\grad w(p_i)
       +8\pi\sum_{j\ne i}\grad_xG(p_i,p_j)=0,
       \quad 1\leq i\leq N,
\end{cases}
\end{equation}
where the $p_i$ are pairwise distinct.  For a sequence of solutions, we set
$$
 I_\nu=\int_MH_{P_\nu}e^{w_\nu}\dV,
 \qquad f_\nu=\rho H_{P_\nu}e^{w_\nu}/I_\nu.
$$
Thus $f_\nu\dV$ is a positive measure of total mass $\rho$ and
\begin{equation}\label{2.10}
 w_\nu(x)=2\int_MG(x,y)f_\nu(y)\dV(y).
\end{equation}
Pass to a common subsequence along which all centers converge and
$f_\nu\dV$ converges weakly to a finite measure. In the local mass
definition below, we take the inner limit at radii for which the
boundary circle has zero mass under the limiting measure. The
outer limit is the mass of that measure at $q$. A collision
class at $q$ consists of all indices whose centers converge to
$q$. We denote its cardinality by $m$.

\begin{theorem}
\label{th2.2}
For a sequence solving \eqref{2.9}, suppose that $m$ of the centers
$p_{1,\nu},\cdots,p_{N,\nu}$ tend to $q$ and set
$$
 4\pi n=\lim_{r\downarrow0}\lim_{\nu\to\infty}\int_{B_r(q)}f_\nu\dV.
$$
Then $n\in\mathbb Z_{\geq0}$ and
\begin{equation}\label{2.11}
 n^2-(m+1)n+m(m-1)=0,\qquad n<\frac{m+1}2.
\end{equation}
Consequently $(m,n)=(0,0)$, $(1,0)$ or $(2,1)$. 
\end{theorem}

\begin{remark}
Under the assumptions of Theorem \ref{th2.2}, each collision consists
of a pair of centers with total local mass $4\pi$, and there is no
blow-up elsewhere. Several distinct pairs may occur, with no
restriction on their relative scales.
\end{remark}

\begin{proof}[Proof of Theorem~\ref{th2.2}]
\textit{Step 1. Mass quantization and the residual equation.}

Choose an isothermal chart with $g=e^\phi|dz|^2$ near $q$. All formulas
in this section use Euclidean coordinates. Let $\psi$ be a fixed smooth
solution of
$$
 \Delta_0\psi=(2\rho-4\pi N)e^\phi
$$
on the coordinate disk and define
$$
 \widehat u_\nu=w_\nu-\log I_\nu
       -4\pi\sum_{j=1}^NG(\cdot,p_{j,\nu})-\psi.
$$
Then $\widehat u_\nu$ satisfies
\begin{equation}\label{2.12}
 \Delta_0\widehat u_\nu+2\rho\widehat h e^{\widehat u_\nu}
       =4\pi\sum_{p_{j,\nu}\text{ in the chart}}\delta_{p_{j,\nu}},
 \qquad \widehat h=e^\phi h_2e^\psi>0
\end{equation}
in the sense of distributions. The weight $\widehat h$ is independent
of the moving centers, and
$$
 \rho\widehat h e^{\widehat u_\nu}\dif A_z=f_\nu\dV.
$$
For each limiting point $q$, denote its local mass by
$$
M_q:=\lim_{r\downarrow0}\lim_{\nu\to\infty}
\int_{B_r(q)}f_\nu\dV.
$$
Let $D$ be a disk isolating a point $q$ with positive local mass.
Choose an annulus containing $\partial D$ and disjoint from the
concentration set. Cover it by finitely many disks on which the
mass is sufficiently small. The local bounds from the Green
representation then give a uniform oscillation bound on $\partial D$.
On this disk, subtract the logarithms at the poles converging to $q$,
as in the proof of Proposition~\ref{pr3.13}. The resulting smooth
positive factor includes all exterior poles. By separation, it has
a uniform positive lower bound and a uniform Lipschitz bound.
It converges locally uniformly along a subsequence. The normalized
total mass bounds the local mass. Proposition~\ref{pr3.13} gives
$M_q\in4\pi\mathbb N$ in the presence of moving poles and in
the regular case. Its proof uses the factor of two in \eqref{2.12}.
If no mass concentrates at $q$, then $M_q=0$.

As in the proof of Lemma~\ref{le2.1}, applying Jensen's inequality
to the mean-zero function
$w_\nu-4\pi\sum_jG(\cdot,p_{j,\nu})$ gives
$I_\nu\geq\min_M h_2>0$. On a compact set $K$ away from the
point masses, cover $K$ by smaller disks inside disks on which
$\int f_\nu\dV$ is uniformly small. In the representation
$w_\nu=2Gf_\nu$, the potential from outside each larger disk is
bounded on the smaller one. For the near potential, apply Jensen's
inequality to $f_\nu$ divided by its mass on the disk. This bounds the
$L^p$ norm of its exponential for any fixed $p>1$ provided the disk
mass is sufficiently small. The logarithmic singularity then has
exponent less than $2$ in the Jensen integral. Since the weight $H_{P_\nu}$
is uniformly bounded and $I_\nu$ is bounded below, $f_\nu$ is
locally bounded in $L^p$. The Green kernel against this $L^p$
bound and the bounded far field give $\sup_K w_\nu\le C_K$.
If $I_\nu\to+\infty$, it follows that
$\|f_\nu\|_{L^\infty(K)}\le C_K/I_\nu\to0$. The whole mass $\rho$
would then be a sum of multiples of $4\pi$, contrary to
noncriticality. Therefore, after a subsequence, $I_\nu\to
I_0\in(0,\infty)$ and the limiting density away from the atoms
is positive where the limiting weight is positive.

Suppose $m$ centers converge to $q$ and the local mass is $M_q=4\pi n$.
The Green kernel has a global lower bound, so its representation gives
$$
 w_0(z)\geq-\frac{M_q}{\pi}\log|z| -C
$$
near $q$. The limiting weight vanishes to order $2m$, so the limiting
density $f_0$ is bounded below by a positive multiple of
$$
 |z|^{2m-M_q/\pi}=|z|^{2m-4n}.
$$
Its integrability forces
\begin{equation}\label{2.13}
 n<\frac{m+1}{2}.
\end{equation}

Since $n$ is an integer, $m-2n\geq0$. Near $q$, subtracting
$8\pi nG(\cdot,q)$ from $w_0$ leaves a distributional scalar equation
with a smooth weight vanishing to order $2(m-2n)$ and an atom-free
density. On a sufficiently small disk the density has arbitrarily
small mass. Applying the Green representation and Jensen's inequality as above,
we find that the exponential of the remainder belongs to
$L^p_{\mathrm{loc}}$ for some $p>1$. The bounded harmonic and smooth parts contribute bounded
factors, and a logarithmic kernel with coefficient equal to the
small mass has integrable $p$-th exponential. The distributional
equation and local Poisson regularity then give $W^{2,p}_{\rm loc}$
for the remainder. Its smooth weighted equation gives smoothness
by elliptic bootstrapping. In a coordinate disk we can write
\begin{equation}\label{2.14}
 w_0(z)=-4n\log|z|+\eta(z),
 \qquad \eta\in C^2,
 \qquad f_0(z)=O(|z|^{2m-4n}).
\end{equation}

\medskip
\noindent\textit{Step 2. The radial Pohozaev identity.}

Take $q=0$ in the fixed chart, and set $J:=\{j\in\{1,\ldots,N\}:p_{j,\nu}\to q\}$,
so that $|J|=m$. After shrinking the chart if necessary, the other
centers stay a fixed positive distance from it. In this chart, define
\(
 V_\nu(z):=\frac{\rho}{I_\nu}e^{\phi(z)}H_{P_\nu}(z).
\)
Then 
$$
 F_\nu(z)=f_\nu(z)e^{\phi(z)}
       =V_\nu(z)e^{w_\nu(z)},
\qquad
 \grad\log V_\nu(z)=2\sum_{i\in J}\frac{z-p_{i,\nu}}{|z-p_{i,\nu}|^2}+O(1),
$$
where the $O(1)$ term is uniformly bounded with respect to $\nu$ and $z$. Constants involving $I_\nu$ disappear after differentiation. The scalar equation is
$$
 -\Delta_0w_\nu=2F_\nu-2\rho e^\phi.
$$

For any fixed small $r$ such that $\partial B_r(0)$ avoids the limiting
concentration set, the radial Pohozaev identity for the scalar equation reads
\begin{align}
 &\int_{\partial B_r}\left[
 (z\cdot\grad w_\nu)\partial_{\boldsymbol n} w_\nu
 -\frac12(z\cdot{\boldsymbol n})|\grad w_\nu|^2
 +2(z\cdot{\boldsymbol n})F_\nu\right]\dif s \notag\\
 &\quad=\int_{B_r}\left[
 4F_\nu+2F_\nu z\cdot\grad\log V_\nu
 +2\rho e^\phi z\cdot\grad w_\nu\right]\dif A_z.
 \label{2.15}
\end{align}
Every term is integrable before we take the limit because $F_\nu$
vanishes quadratically at each $p_{i,\nu}$.

We first let $\nu\to\infty$ and then $r\downarrow0$.
By \eqref{2.14},
\begin{align}
\lim_{\nu\to\infty}\int_{\partial B_r}
[(z\cdot\grad w_\nu)\partial_{\boldsymbol n}w_\nu
-\frac12(z\cdot\boldsymbol n)|\grad w_\nu|^2]\,\dif s
&=\pi(4n)^2+O(r),\nonumber\\
\lim_{\nu\to\infty}2\int_{\partial B_r}
(z\cdot\boldsymbol n)F_\nu\,\dif s&=O(r^{2m-4n+2}).\nonumber
\end{align}
Hence the limit of the left-hand side of \eqref{2.15}
is $\pi(4n)^2=M_q^2/\pi$.
On the right-hand side, we have
$$
\begin{aligned}
2\int_{B_r}F_\nu z\cdot
\left(\grad\log V_\nu
-2\sum_{i\in J}\frac{z-p_{i,\nu}}{|z-p_{i,\nu}|^2}\right)\dif A_z
&=O(r),\\
2\rho\int_{B_r}e^\phi z\cdot\grad w_\nu\,\dif A_z
&=O(r^{3-2/p}),\qquad 1<p<2,
\end{aligned}
$$
uniformly in $\nu$, by the bounded total mass and the uniform
$W^{1,p}$ bound for $w_\nu$.

In the singular weight term, write $z=(z-p_{i,\nu})+p_{i,\nu}$.
Up to terms that vanish in the stated limit, the right-hand side of
\eqref{2.15} becomes
\begin{equation}\label{2.16}
 4(m+1)M_q+
 4\lim_{r\downarrow0}\lim_{\nu\to\infty}
 \sum_{i\in J}p_{i,\nu}\cdot
 \int_{B_r}F_\nu(z)\frac{z-p_{i,\nu}}{|z-p_{i,\nu}|^2}\dif A_z.
\end{equation}

The Green representation \eqref{2.10} gives, for fixed $r$,
$$
 \int_{B_r}F_\nu(z)\frac{z-p_{i,\nu}}{|z-p_{i,\nu}|^2}\dif A_z
       =\pi\grad w_\nu(p_{i,\nu})+O_r(1).
$$
The contribution from outside $B_r$ is bounded for fixed $r$, and
the regular part of the Green kernel is uniformly smooth. Multiplying
the remainder by $p_{i,\nu}\to0$ makes it vanish in the first limit.

On the other hand, the center equation in \eqref{2.9} gives
\begin{equation}\label{2.17}
 \grad w_\nu(p_{i,\nu})
 =-8\sum_{\substack{j\in J\\j\ne i}}
       \frac{p_{i,\nu}-p_{j,\nu}}{|p_{i,\nu}-p_{j,\nu}|^2}+O(1).
\end{equation}
Using the pairing identity
$$
 \frac{p_i\cdot(p_i-p_j)+p_j\cdot(p_j-p_i)}{|p_i-p_j|^2}=1
$$
and summing \eqref{2.17}, we get
$$
\begin{aligned}
\sum_{i\in J}p_{i,\nu}\cdot\grad w_\nu(p_{i,\nu})
&=-8\sum_{i\in J}\sum_{\substack{j\in J\\j\ne i}}
\frac{p_{i,\nu}\cdot(p_{i,\nu}-p_{j,\nu})}
     {|p_{i,\nu}-p_{j,\nu}|^2}
+O\!\left(\sum_{i\in J}|p_{i,\nu}|\right)\\
&=-8\binom{m}{2}+o(1)
=-4m(m-1)+o(1).
\end{aligned}
$$ Substituting in \eqref{2.16} proves
$$
 \frac{M_q^2}{\pi}=4(m+1)M_q-16\pi m(m-1).
$$
With $M_q=4\pi n$, this is \eqref{2.11}. Its discriminant in $n$ is
$$
 (m+1)^2-4m(m-1)=4-3(m-1)^2.
$$
It is negative for $m\geq3$. If $m=2$, the roots are $n=1,2$,
and \eqref{2.13} excludes $n=2$. If $m=1$, the roots are $n=0,2$,
and integrability leaves only $n=0$. If $m=0$, the roots are $n=0,1$,
and integrability again leaves only $n=0$. This proves Theorem~\ref{th2.2}.
\end{proof}

\subsection{The residual center equation}

Now let $N=2$ and suppose $p_{1,\nu},p_{2,\nu}\to q$. By Theorem~\ref{th2.2},
$$
 f_\nu\dV\rightharpoonup4\pi\delta_q+f_0\dV,
 \qquad \int_Mf_0\dV=\rho-4\pi>0.
$$
The Green representation and cancellation of the limiting zero of order four give
\begin{equation}\label{2.18}
 w_\nu\longrightarrow8\pi G(\cdot,q)+w_*,
 \qquad
 \Delta_gw_*+2(\rho-4\pi)
 \left(\frac{h_2e^{w_*}}{\int_Mh_2e^{w_*}\dV}-1\right)=0,
 \quad \int_Mw_*\dV=0.
\end{equation}
The convergence is $C^2$ locally away from $q$. The function $w_*$ is smooth on all of $M$.

\begin{lemma}\label{le2.3}
The collision point $q$ in \eqref{2.18} satisfies
\begin{equation}\label{2.19}
 \grad_g\bigl(2\log h_1+\log h_2+12\pi r_g\bigr)(q)=0.
\end{equation}
Thus, for every collision limit $(w_*,q)$ of the two-center shadow system,
$w_*$ solves the scalar equation in \eqref{2.18} with parameter $\rho-4\pi$.
The point $q$ satisfies the critical point condition \eqref{2.19}.
\end{lemma}

\begin{remark}
The necessary conditions \eqref{2.18} and \eqref{2.19} on the
collision limit $(w_*,q)$ follow from the two-center shadow system,
independently of the existence of mixed Toda bubbling solutions.
\end{remark}

\begin{proof}[Proof of Lemma \ref{le2.3}]
We use a coordinate calculation to specify the regular-part conventions.
Put $q=0$ in an isothermal chart and write
$$
 G(x,p)=-\frac1{2\pi}\log|x-p|+R(x,p),
 \qquad r_g(p)=R(p,p)+\frac{\phi(p)}{4\pi}.
$$
Here $R$ denotes the regular part of the Green function in the
chosen isothermal chart and satisfies $R(x,p)=R(p,x)$. Thus
$$
\grad(R(p,p))=2\grad_xR(p,p).
$$
Set $R_x=\grad_xR(0,0)$, $W=\grad w_*(0)$ and
$$
S_\nu=\grad w_\nu(p_{1,\nu})+\grad w_\nu(p_{2,\nu})
$$
for short. Adding the two center equations shows that $S_\nu$ has a limit and gives
\begin{equation}\label{2.20}
 \lim_{\nu\to\infty} S_\nu=4\grad g_1(0)+32\pi R_x.
\end{equation}
The singular derivatives of $G(p_1,p_2)$ cancel exactly when we add
the two center equations.

The translation Pohozaev identity for the scalar equation reads
\begin{align*}
 \int_{\partial B_r}\left[
 \partial_{\boldsymbol n} w_\nu\,\grad w_\nu
 -\frac12|\grad w_\nu|^2 {\boldsymbol n}+2F_\nu {\boldsymbol n}\right]\dif s=\int_{B_r}\left[2F_\nu\grad\log V_\nu +2\rho e^\phi\grad w_\nu\right]\dif A_z.
\end{align*}
The limit of the left-hand side is
$$
 -8\pi(8\pi R_x+W),
$$
because $w_0=-4\log|z|+8\pi R(z,0)+w_*(z)$. The terms
$2\int_{\partial B_r}F_\nu\boldsymbol n\,\dif s$ and
$2\rho\int_{B_r}e^\phi\grad w_\nu\,\dif A_z$ vanish in the limit.

Define
$$
I_{i,\nu}(r)=\int_{B_r}F_\nu(z)
\frac{z-p_{i,\nu}}{|z-p_{i,\nu}|^2}\,\dif A_z,
\qquad
Q_{i,\nu}(r)=\grad w_\nu(p_{i,\nu})-\frac1\pi I_{i,\nu}(r).
$$
For fixed $r$, Green's formula gives
$$
\begin{aligned}
\lim_{\nu\to\infty}Q_{i,\nu}(r)
&=8\pi R_x+W-\frac1\pi\int_{B_r}
f_0(z)e^{\phi(z)}\frac{z}{|z|^2}\,\dif A_z\\
&\longrightarrow 8\pi R_x+W
\qquad\text{as }r\downarrow0,
\end{aligned}
$$
since the integral is $O(r)$ by the boundedness of $f_0$ near $0$.

Since
$$
 \grad\log V_\nu
 =\grad(\log h_2+\phi)-4\pi\sum_{i=1}^2\grad_xR(z,p_{i,\nu})
       +2\sum_{i=1}^2\frac{z-p_{i,\nu}}{|z-p_{i,\nu}|^2},
$$
the limit of the right-hand side is
$$
 8\pi\grad(\log h_2+\phi)
       -128\pi^2R_x+4\pi\lim S_\nu-8\pi W.
$$
Equating both sides yields
\begin{equation}\label{2.21}
 \lim_{\nu\to\infty} S_\nu=-2\grad(\log h_2+\phi)+16\pi R_x.
\end{equation}
Combining \eqref{2.20} and \eqref{2.21},
$$
 2\grad g_1+\grad\log h_2+\grad\phi+8\pi R_x=0.
$$
Finally, $g_1=\log h_1+4\pi r_g$ and
$$
 8\pi R_x=4\pi\grad r_g-\grad\phi
$$
give \eqref{2.19}. The terms containing the residual scalar gradient
$W$ cancel exactly.
\end{proof}

\subsection{The midpoint of a colliding pair}

\begin{lemma}\label{lem:pair-midpoint}
Under the assumptions of Theorem~\ref{th2.2}, suppose that $N=2$
and $p_{1,\nu},p_{2,\nu}\to q$. In a fixed isothermal chart centered
at $q$, set $F_\nu=e^\phi f_\nu$ and
$$
 c_\nu=\frac{p_{1,\nu}+p_{2,\nu}}2,
 \qquad \delta_\nu=\frac{|p_{1,\nu}-p_{2,\nu}|}2.
$$
After passing to a subsequence, the following hold.
\begin{enumerate}[label=\textup{(\roman*)}]
\item The rescaled densities satisfy
$$
 \delta_\nu^2F_\nu(c_\nu+\delta_\nu y)\,\dif A_y
 \rightharpoonup4\pi\delta_0
 \qquad\text{locally weakly as measures on }\R^2.
$$
\item There are local maxima $x_\nu$ of $F_\nu$ such that, with
$s_\nu=F_\nu(x_\nu)^{-1/2}$,
$$
 x_\nu=c_\nu+o(\delta_\nu),\qquad s_\nu=o(\delta_\nu).
$$
Moreover,
$$
 \log\!\bigl(s_\nu^2F_\nu(x_\nu+s_\nu y)\bigr)
 \longrightarrow -2\log\left(1+\frac{|y|^2}{4}\right)
 \qquad\text{in }C^2_{\mathrm{loc}}(\R^2).
$$
\end{enumerate}
\end{lemma}

\begin{proof}
\textit{Step 1.} We show that the sequence rescaled by $\delta_\nu$
has a blow-up point in $\R^2$.

Write
$p_{1,\nu}=c_\nu+\delta_\nu e_\nu$ and
$p_{2,\nu}=c_\nu-\delta_\nu e_\nu$, where $|e_\nu|=1$ and,
after passing to a subsequence, $e_\nu\to e$. For
$z=c_\nu+\delta_\nu y$, set
$$
 V_\nu(y):=\widehat u_\nu(z)
 -2\sum_{j=1}^2\log|z-p_{j,\nu}|+6\log\delta_\nu,
 \qquad
 \widetilde F_\nu(y)=\delta_\nu^2F_\nu(z).
$$
By \eqref{2.12}, on the rescaled coordinate disk $D_\nu$ we have
$$
 -\Delta_y V_\nu=2\rho\widehat h(c_\nu+\delta_\nu y)
 |y-e_\nu|^2|y+e_\nu|^2e^{V_\nu}=2\widetilde F_\nu.
$$
Apply \cite[Proposition~3.1, \textup{(2.15)--(2.21)}]{Tar23} to
$\xi_\nu(z)=\widehat u_\nu(z)-2\sum_{j=1}^2\log|z-p_{j,\nu}|$.
Here $\xi_\nu$ is smooth because, after cancellation,
it differs from $w_\nu$ by smooth terms and a spatial constant.
Translate the fixed disk to the midpoint. Choose its radius small
enough that the boundary stays in a fixed isolating annulus.
The two distinct zeros have strengths one and separation
$2\delta_\nu$. The smooth factor $2\rho\widehat h$ is uniformly
positive and Lipschitz and converges in the translated coordinates.
Add a bounded constant to $\xi_\nu$ and compensate in the smooth
factor so that the limiting factor equals one at the origin.
These changes preserve the source measure and the oscillation.
The point mass at $q$ is $8\pi$ in the measure $2F_\nu\,\dif A$.
Since the coefficient is bounded above, a positive point mass
forces $\xi_\nu$ to blow up in every smaller disk. On a fixed
annulus about $q$, \eqref{2.18} and $I_\nu\to I_0\in(0,\infty)$
give boundedness and bounded boundary oscillation. Thus there is no
other limiting blow-up point on the disk. These are the isolated
least-mass hypotheses in that proposition. In its notation, the
pole scale is the maximum distance to the midpoint, namely
$\delta_\nu$, and the logarithmic shift is $6\log\delta_\nu$.

\textit{Step 2.} The proposition gives one atom in the rescaled
measure, of mass $8\pi$ in $2\widetilde F_\nu\,\dif A$, with
convergence to $-\infty$ for $V_\nu$ locally away from that atom.
If the atom were at a limiting zero, at least two zeros would
have to converge to that point. The two limiting zeros here
are the distinct points $e$ and $-e$. Hence for some
$z_0\ne\pm e$ we have
\[
 \widetilde F_\nu\,\dif A\rightharpoonup4\pi\delta_{z_0}
\]
locally as measures. The density tends to zero locally uniformly
away from $z_0$ because the rescaled weight in the equation is
bounded on fixed compact sets.

\textit{Step 3.} We use the first center equation to locate $z_0$.

For $y$ in a fixed compact set away from $z_0$, the part of the
Green integral with $|y-\zeta|<\eta$ is $O(\eta)$ by the local
density bound, and the part with $|\zeta|>L$ is $O(L^{-1})$ by
the total mass bound. Weak convergence against the remaining
truncated kernels gives, locally uniformly away from $z_0$,
$$
 \grad V_\nu(y)\longrightarrow
 -4\frac{y-z_0}{|y-z_0|^2}.
$$
At $y=e_\nu$, the chain rule and the smooth terms in $V_\nu$ give
$\delta_\nu\grad w_\nu(p_{1,\nu})
=\grad V_\nu(e_\nu)+O(\delta_\nu)$. Since
$p_{1,\nu}-p_{2,\nu}=2\delta_\nu e_\nu$, the logarithmic part of $G$
gives $8\pi\delta_\nu\grad_xG(p_{1,\nu},p_{2,\nu})
=-2e_\nu+O(\delta_\nu)$. Thus the first center equation in
\eqref{2.9}, multiplied by $\delta_\nu$, reads
$-\tfrac12\grad V_\nu(e_\nu)-2e_\nu+O(\delta_\nu)=0$.
Since $e_\nu\to e\ne z_0$, the locally uniform convergence of
$\grad V_\nu$ gives
$$
 2\frac{e-z_0}{|e-z_0|^2}-2e=0.
$$
Taking norms gives $|e-z_0|=1$. Substituting this equality back gives
$z_0=0$. This proves \textup{(i)}.

\textit{Step 4.} We choose a local maximum of $F_\nu$ and determine
the scale around it.

Let $y_\nu$ maximize $\widetilde F_\nu$ on
$\overline B_{1/4}$. By \textup{(i)} and the local density bound
away from $0$, we have $y_\nu\to0$ and
$\widetilde F_\nu(y_\nu)\to+\infty$. Thus $y_\nu$ is an interior
maximum for large $\nu$, and $x_\nu=c_\nu+\delta_\nu y_\nu$ is a
local maximum of $F_\nu$ with $x_\nu=c_\nu+o(\delta_\nu)$.
For $s_\nu=F_\nu(x_\nu)^{-1/2}$, we have
$s_\nu/\delta_\nu=\widetilde F_\nu(y_\nu)^{-1/2}\to0$.

\textit{Step 5.} We identify the limiting profile at scale $s_\nu$.

Set $U_\nu(y)=\log(s_\nu^2F_\nu(x_\nu+s_\nu y))$.
In these coordinates, every fixed disk eventually avoids both
sources. Thus \eqref{2.12} gives
$$
 \Delta_y U_\nu+2e^{U_\nu}
 =s_\nu^2(\Delta_0\log\widehat h)(x_\nu+s_\nu y).
$$
On each fixed $B_{2R}$, the choice of $x_\nu$ gives $U_\nu(0)=0$
and $U_\nu\leq0$ for large $\nu$. Let $\chi_\nu$ solve
$\Delta\chi_\nu=\Delta U_\nu$ on $B_{2R}$ with zero boundary
values. The equation gives
$\|\chi_\nu\|_{L^\infty(B_{2R})}\leq C_R$. Harnack's inequality
applied to the nonnegative harmonic function
$C_R-U_\nu+\chi_\nu$ bounds $U_\nu$ from below on $B_R$.
Interior elliptic estimates then give a limit $U$ in
$C^2_{\mathrm{loc}}(\R^2)$ along a subsequence, with
$\Delta U+2e^U=0$ and $U(0)=\max U=0$. For every fixed $R$,
changing variables gives
$$
 \limsup_{\nu\to\infty}\int_{B_R}e^{U_\nu}\,\dif A
 \leq\lim_{\nu\to\infty}\int_{B_{1/4}}\widetilde F_\nu\,\dif A
 =4\pi.
$$
Let $R\to\infty$ and apply the rank-one classification
\eqref{eq:scalar-liouville} with $\gamma=0$. The maximum at
$0$ and the value $e^{U(0)}=1$ give $p=0$ and $s=4$.
Thus $U(y)=-2\log(1+|y|^2/4)$, as in \textup{(ii)}.
\end{proof}

\subsection{The general case: successive pair collisions}

\begin{lemma}\label{lem:successive-pair-collisions}
Let $(w_\nu,P_\nu,Q_\nu)$ solve \eqref{2.7}--\eqref{2.8} with fixed
type $(a,b)$ and scalar parameter $\sigma=\sigma_b>0$,
$\sigma\notin4\pi\NN$. Put
$I_\nu=\int_M H_{P_\nu}e^{w_\nu}\dV$ and
$f_\nu=\sigma H_{P_\nu}e^{w_\nu}/I_\nu$.
After passing to a subsequence, all centers have limits and
$f_\nu\dV\rightharpoonup\mu$.
\begin{enumerate}[label=\textup{(\roman*)}]
\item If $u$ ordinary centers and $v$ mixed centers tend to $q$ and
the local mass at $q$ is $4\pi n=\mu(\{q\})$, then
\[
 (u,v,n)\in\{(0,0,0),(1,0,0),(0,1,0),(2,0,1)\}.
\]
\item Let $\mathcal R=\{r_1,\ldots,r_\ell\}$ be the set of common limits
of the colliding pairs of ordinary centers, i.e., the $(u,v,n)=(2,0,1)$ case. Write
$p_i^0=\lim_{\nu\to\infty}p_{i,\nu}$ and
$q_j^0=\lim_{\nu\to\infty}q_{j,\nu}$, and set
\[
 P'=\{p_i^0:1\leq i\leq a\}\setminus\mathcal R,
 \qquad
 Q'=\{q_j^0:1\leq j\leq b\}\cup\mathcal R.
\]
Then
\[
 w_\nu\longrightarrow 8\pi\sum_{h=1}^{\ell}G(\cdot,r_h)+w'
 \quad\text{in }C^2_{\mathrm{loc}}(M\setminus\{r_1,\ldots,r_\ell\}),
\]
where $(w',P',Q')$ solves \eqref{2.7}--\eqref{2.8} with type
$(a-2\ell,b+\ell)$ and scalar parameter $\sigma-4\pi\ell>0$.
\end{enumerate}
\end{lemma}

\begin{proof}
\textit{Step 1.} We classify the collisions.
Since only ordinary centers enter $H_{P_\nu}$ (see \eqref{2.5}),
the quantization, normalization, and integrability arguments from
\eqref{2.12} through \eqref{2.13} in the proof of
Theorem~\ref{th2.2} apply with $\rho=\sigma$, $N=a$, and $m=u$.
Thus, along a subsequence,
\[
 I_\nu\to I_0\in(0,\infty),\qquad
 n\in\mathbb Z_{\geq0},\qquad n<\frac{u+1}{2}.
\]
In a fixed isothermal chart centered at $q$, write $p_i$ and $q_j$
for the centers tending to $q$. We omit the sequence index $\nu$ in
the following calculation. The remaining centers stay a positive
distance from $q$, so the gradients of the corresponding Green function
terms are uniformly bounded in a fixed neighborhood of $q$. Set
$S_\nu=\sum_{i=1}^u p_i\cdot\grad w_\nu(p_i)$.
In \eqref{2.8}, sum $p_i\cdot F_i^A=0$ over $i$ and
$q_j\cdot F_j^B=0$ over $j$. Using
$\grad_xG(x,y)=-(x-y)/(2\pi|x-y|^2)+O(1)$ and pairing terms of the
same type by $x\cdot(x-y)+y\cdot(y-x)=|x-y|^2$, we obtain
\begin{align*}
 S_\nu&=-8\binom u2-12\sum_{i,j}
       \frac{p_i\cdot(p_i-q_j)}{|p_i-q_j|^2}+o(1),\\
 0&=-6\sum_{i,j}\frac{q_j\cdot(q_j-p_i)}{|p_i-q_j|^2}
       -12\binom v2+o(1).
\end{align*}
The two fractions sum to $1$ for each $(i,j)$, so
\[
 \lim_{\nu\to\infty}S_\nu=-8\binom u2-12uv-24\binom v2.
\]
Applying the radial Pohozaev identity \eqref{2.15} with $m=u$ and
using Green's representation as after \eqref{2.16}, we obtain
\[
 16\pi n^2=16\pi(u+1)n+4\pi\lim_{\nu\to\infty}S_\nu.
\]
Substituting the limit gives
\begin{equation}\label{2.22}
 n^2-(u+1)n+u(u-1)+3uv+3v(v-1)=0.
\end{equation}
The discriminant $4-3(u+2v-1)^2$ gives $u+2v\leq2$. Substituting
these possibilities in \eqref{2.22} and using $n<(u+1)/2$ yields
\textup{(i)}.

\textit{Step 2.} We identify the limiting shadow system.
By \textup{(i)}, with $\sigma'=\sigma-4\pi\ell$, we have
\[
 \mu=4\pi\sum_h\delta_{r_h}+f_0\dV,\qquad
 \int_M f_0\dV=\sigma'>0,\qquad
 w_0=8\pi\sum_hG(\cdot,r_h)+w',\quad w'=2Gf_0.
\]
Green's representation and local elliptic estimates give
$w_\nu\to w_0$ in $C^2_{\mathrm{loc}}$ away from the collision points.
The limiting weight is $H_{P'}e^{-8\pi\sum_hG(\cdot,r_h)}$.
Multiplying it by $e^{w_0}$ and using $\int_M f_0\dV=\sigma'$ gives
\[
 f_0=\frac{\sigma}{I_0}H_{P'}e^{w'}
    =\frac{\sigma'H_{P'}e^{w'}}{\int_M H_{P'}e^{w'}\dV},
 \qquad -\Delta_g w'=2(f_0-\sigma').
\]
The argument for \eqref{2.14} gives smoothness of $w'$ across the
collision points. In the center equations of \eqref{2.8}, the
contribution of the pair merging at $r_h$ to the limiting
potential at a remaining ordinary center and at an original
mixed center, respectively, is
\[
 \left(2\cdot8\pi-\tfrac12\cdot8\pi\right)G(\cdot,r_h)
 =12\pi G(\cdot,r_h),\qquad
 2\cdot12\pi G(\cdot,r_h)=24\pi G(\cdot,r_h).
\]
At a new mixed center $r_h$, apply the local translation Pohozaev
calculation of Lemma~\ref{le2.3} with the external fields
\[
 g_{{\rm ext},\nu}=g_1+8\pi\sum_{p\in P_{{\rm ext},\nu}}G(\cdot,p)
      +12\pi\sum_{q\in Q_\nu}G(\cdot,q),\qquad
 \log h_{2,{\rm ext},\nu}
 =\log h_2-4\pi\sum_{p\in P_{{\rm ext},\nu}}G(\cdot,p),
\]
where $P_{{\rm ext},\nu}$ consists of all ordinary centers except
the two tending to $r_h$. These fields converge smoothly near
$r_h$ because every center included in them stays a positive
distance from a sufficiently small fixed disk. For the local use
of Lemma~\ref{le2.3}, take isothermal coordinates centered at $r_h=0$
and write the local limit as
$w_0=-4\log|z|+8\pi R(z,0)+w_{\rm reg}(z)$.
The contributions of the other collapsing pairs to $w_{\rm reg}$
are smooth on this disk. Let $R_x=\grad_x R(0,0)$,
$W=\grad w_{\rm reg}(0)$ and
$S_\nu=\sum_{i=1}^2\grad w_\nu(p_{i,\nu})$, with the two indices
now denoting this pair. Adding the two Euclidean center equations
cancels the singular gradients and gives
\[
 \lim_\nu S_\nu=4\grad g_{\rm ext}(0)+32\pi R_x.
\]
In the scalar translation identity, the smooth part of the
weight gradient has limit contribution
$8\pi\grad(\log h_{2,\rm ext}+\phi)(0)-64\pi^2R_x$.
For each singular term, subtract the integral over the test disk
from the Green representation of $\grad w_\nu(p_{i,\nu})$, as
in the calculation of $Q_{i,\nu}(r)$ in Lemma~\ref{le2.3}.
For fixed radius, the exterior contribution converges by separation
from $p_{i,\nu}$. As the radius then tends to zero, this difference
tends to $8\pi R_x+W$ because the remaining density is bounded at
$0$. Thus the two singular terms together contribute
$4\pi\lim_\nu S_\nu-64\pi^2R_x-8\pi W$. The boundary stress tends
to $-8\pi(8\pi R_x+W)$. The boundary density term tends to zero
because the order-four zero cancels the singular exponential,
and the bounded background forcing contributes $O(r)$ by the
integrability of the differentiated Green kernel. This gives
\[
 \lim_\nu S_\nu=-2\grad(\log h_{2,\rm ext}+\phi)(0)+16\pi R_x.
\]
Together these two identities give
$\grad(2g_{\rm ext}+\log h_{2,\rm ext}+4\pi r_g)(0)=0$,
since $4\pi\grad r_g=8\pi R_x+\grad\phi$ at $0$.
Each other colliding pair contributes two terms, each with coefficient
$12\pi$, at its common limit. Hence the limiting equation is
\[
 \begin{aligned}
 0
 &=\grad(2g_{\rm ext}+\log h_{2,\rm ext}+4\pi r_g)(r_h)\\
 &=\grad\!\left[V+12\pi\sum_{p\in P'}G(\cdot,p)
     +24\pi\sum_{q\in Q'\setminus\{r_h\}}G(\cdot,q)\right](r_h).
 \end{aligned}
\]
This is the mixed-center equation in \eqref{2.8} at $r_h$, which
completes the proof of \textup{(ii)}.
\end{proof}
\subsection{A necessary condition for collision along the homotopy}

\begin{lemma}
Fix $c<0$ and $\rho>0$ with $\rho\notin4\pi\NN$. Let
$(w_\nu,p_{1,\nu},p_{2,\nu})$ solve the scalar equation in
\eqref{2.9} with $N=2$ and the center equations
\[
 \grad g_1(p_{i,\nu})+c\grad w_\nu(p_{i,\nu})
 +8\pi t_\nu\grad_xG(p_{i,\nu},p_{j,\nu})=0,
 \qquad \{i,j\}=\{1,2\}.
\]
Suppose that $t_\nu\to t_*$, $p_{1,\nu},p_{2,\nu}\to q$ and the
local mass $M_q$ is positive. Then $M_q=4\pi$ and
\begin{equation}\label{2.23}
 t_*=-2c.
\end{equation}
In particular, for $t_\nu\in[0,1]$, a collision with positive local
mass is excluded if $c<-1/2$ and requires $t_*=1$ if $c=-1/2$.
\end{lemma}

\begin{proof}
The local equation \eqref{2.12} has two moving sources, each of
strength one. The two-source case of Proposition~\ref{pr3.13} gives
$M_q=4\pi n$ for some integer $n\geq1$. Since
$\rho\notin4\pi\NN$, the normalization argument following
\eqref{2.12} gives $I_\nu\to I_0\in(0,\infty)$ along a subsequence.
In a coordinate disk centered at $q$, the two quadratic zeros
merge into a zero of order four, while the Green representation
gives $w_0(z)\geq-4n\log|z|-C$. Thus the residual density satisfies
$f_0(z)\geq C^{-1}|z|^{4-4n}$. Its local integrability requires
$4-4n>-2$, so $n<3/2$. Hence $n=1$ and $M_q=4\pi$.

In isothermal coordinates centered at $q$, set
$S_\nu=\sum_{i=1}^2p_{i,\nu}\cdot\grad w_\nu(p_{i,\nu})$.
Taking the scalar product of the $i$th center equation with
$p_{i,\nu}$ and summing over $i$, we use
$\grad_xG(x,y)=-(x-y)/(2\pi|x-y|^2)+O(1)$ and
$x\cdot(x-y)+y\cdot(y-x)=|x-y|^2$ to obtain
$cS_\nu=4t_\nu+o(1)$. The radial Pohozaev identity
\eqref{2.15}, evaluated as in \eqref{2.16} with $m=2$, gives
$16\pi=48\pi+4\pi\lim_{\nu\to\infty}S_\nu
=48\pi+16\pi t_*/c$. This yields \eqref{2.23}.
\end{proof}

\section{Degree of the shadow systems}\label{rev:sec:shadow}

We use the shadow equations and center conventions of
Section~\ref{sec:shadow-setting} and the spaces and orientation
of Section~\ref{sec:degree-conventions}. The mixed-center
inequality \eqref{3.41} is used in the compact homotopy of
Section~\ref{typed-sec:rotation}.

\subsection{An identity for two colliding centers}\label{sec:stress-chapter}

\subsubsection{Two colliding centers}

Fix $\rho>0$ with $\rho\notin4\pi\NN$. We consider \eqref{2.9}
with $N=2$, replacing $-1/2$ by $c<0$ in the center equations.
For $P=(p_+,p_-)$, write $H_{p_+,p_-}=H_P$. The system is
\begin{equation}\label{3.1}
\begin{cases}
 \displaystyle \Delta_gw+2\rho\left(
 \frac{H_{p_+,p_-}e^w}{\int_MH_{p_+,p_-}e^w\dV}-1\right)=0,
 \qquad \int_Mw\dV=0,\\[4pt]
 \displaystyle \nabla g_1(p_i)+c\nabla w(p_i)
 +8\pi\nabla_xG(p_i,p_j)=0,
 \qquad i=\pm,\quad j=\mp,
\end{cases}
\end{equation}
where $p_+\ne p_-$.

Suppose $(w_\nu,p_{+,\nu},p_{-,\nu})$ solves \eqref{3.1} with
$c=c_\nu\to-1/2$ and $p_{\pm,\nu}\to q$. In an isothermal coordinate
$z=x+iy$ centered at $q$, with $g=e^\phi|dz|^2$, write
\[
 p_{\pm,\nu}=s_\nu\pm d_\nu,
 \qquad s_\nu\to0,\quad d_\nu\ne0,\quad d_\nu\to0.
\]
We use Euclidean gradients and the derivative convention
of Section~\ref{sec:notation} in this coordinate. Put
\[
 f_\nu=\rho\frac{H_{p_{+,\nu},p_{-,\nu}}e^{w_\nu}}
                  {\int_MH_{p_{+,\nu},p_{-,\nu}}e^{w_\nu}\dV},
 \qquad F_\nu=f_\nu e^\phi.
\]
We pass to a subsequence. The scalar analysis in
Section~\ref{rev:sec:scalar} then gives the following limits.
\begin{enumerate}
\item We have $f_\nu\dV\rightharpoonup4\pi\delta_q+f_0\dV$ and
$w_\nu\to8\pi G(\cdot,q)+w_*$ locally smoothly on $M\setminus\{q\}$,
where $w_*$ solves the residual equation in \eqref{2.18} and
\[
 f_0=(\rho-4\pi)\frac{h_2e^{w_*}}{\int_Mh_2e^{w_*}\dV}.
\]
In the coordinate disk, $F_\nu\to F_0=f_0e^\phi$ locally smoothly
away from $q$, and $F_0$ is smooth and positive near $q$.

\item Under the rescaling $z=s_\nu+|d_\nu|y$, the densities satisfy
\begin{equation}\label{3.2}
 |d_\nu|^2F_\nu(s_\nu+|d_\nu|y)\dif A_y
 \rightharpoonup4\pi\delta_0
 \qquad\text{locally weakly as measures on }\C.
\end{equation}
In particular, $\int_{B_{\eta|d_\nu|}(s_\nu)}F_\nu\dif A_z\to4\pi$
for every fixed $\eta>0$.
\end{enumerate}

The quantization, normalization, and integrability arguments from
\eqref{2.12} to \eqref{2.13} also apply when $c_\nu$ varies.
They exclude concentration away from $q$ and give $M_q=4\pi n$
with $n\in\{0,1\}$. Set
$S_\nu=\sum_{i=+,-}p_{i,\nu}\cdot\nabla w_\nu(p_{i,\nu})$.
Pairing the center equations in \eqref{3.1} gives
$c_\nu S_\nu=4+o(1)$, so $S_\nu\to-8$. The radial Pohozaev identity
\eqref{2.15}--\eqref{2.16} gives
\[
 16\pi n^2=48\pi n+4\pi\lim_{\nu\to\infty}S_\nu
          =48\pi n-32\pi,
 \qquad n=1.
\]
Green's representation then gives the first limit. The first two
steps of Lemma~\ref{lem:pair-midpoint} use only the scalar equation.
In Step~3, with $e_\nu=d_\nu/|d_\nu|$, the rescaled center equation is
$c_\nu\nabla V_\nu(e_\nu)-2e_\nu+O(|d_\nu|)=0$.
Taking the limit with $c_\nu\to-1/2$ gives the same midpoint condition
and hence \eqref{3.2}.

Write
\[
 G(z,p)=-\frac1{2\pi}\log|z-p|+R(z,p),
 \qquad r_g(z)=R(z,z)+\frac{\phi(z)}{4\pi}.
\]
The subscripts $z$ and $p$ refer to the first and second arguments
of $R$, respectively. Unless otherwise specified, mixed derivatives
such as $R_{zp}$ and $R_{z\bar p}$ are evaluated at $(q,q)$.
We define
\begin{equation}\label{3.3}
 C(z)=2\log h_1(z)+\log h_2(z)+12\pi r_g(z).
\end{equation}
With $q$ fixed, set
\begin{equation}\label{3.4}
 \psi_1(z)=\log h_1(z)+\phi(z)+12\pi R(z,q)-\frac12w_*(z),
 \qquad \psi_2(z)=\log h_2(z)+\phi(z)+w_*(z).
\end{equation}
At $q$ put
\begin{equation}\label{3.5}
 \kappa_1(q)=(\Delta_g\log h_1)(q)+16\pi-\rho-2K_g(q),
 \qquad \mathcal T(q)=\psi_{1,zz}(q)+\psi_{1,z}(q)^2.
\end{equation}

The following theorem gives a criterion that excludes collisions
as $c\to-1/2$ with $c\leq-1/2$.

\begin{theorem}\label{th3.1}
Under the preceding assumptions, $\nabla C(q)=0$. With
$\zeta_\nu=d_\nu^2$, we have
\begin{equation}\label{3.6}
 -2\mathcal T(q)\zeta_\nu
 -\tfrac12e^{\phi(q)}\kappa_1(q)|\zeta_\nu|
 =\frac{-c_\nu-\tfrac12}{c_\nu^2}+o(|\zeta_\nu|).
\end{equation}
A collision with $c_\nu\equiv-1/2$ requires
\begin{equation}\label{3.7}
 e^{\phi(q)}|\kappa_1(q)|=4|\mathcal T(q)|.
\end{equation}
If $c_\nu\leq-1/2$, the condition below excludes a collision at $q$.
\begin{equation}\label{3.8}
 e^{\phi(q)}\kappa_1(q)>4|\mathcal T(q)|.
\end{equation}
\end{theorem}

We postpone the proof of Theorem~\ref{th3.1} until after
Lemma~\ref{le3.3}.

\begin{remark}\label{re3.2}
The condition $e^{\phi(q)}|\kappa_1(q)|\ne4|\mathcal T(q)|$
excludes collisions at $c=-1/2$. The sign condition \eqref{3.8}
also excludes collisions along sequences approaching $-1/2$ from below.
\end{remark}

\subsubsection{A complex differential identity}

Set $U_\nu=\log F_\nu$. In the sense of distributions, the scalar equation in
\eqref{3.1} and the normalization of the Green function give
\begin{equation}\label{3.9}
 \Delta_0U_\nu+2e^{U_\nu}
 =\kappa(z)+4\pi(\delta_{p_{+,\nu}}+\delta_{p_{-,\nu}}),
 \qquad
 \kappa=\Delta_0(\log h_2+\phi)+(2\rho-8\pi)e^\phi.
\end{equation}
Away from the sources, set $T_\nu=U_{\nu,zz}-U_{\nu,z}^2/2$.
By \eqref{3.9} and $\Delta_0=4\partial_z\partial_{\bar z}$,
$U_{\nu,z\bar z}=\kappa/4-e^{U_\nu}/2$. Hence
\begin{equation}\label{3.10}
 \partial_{\bar z}T_\nu
 =\partial_z U_{\nu,z\bar z}-U_{\nu,z}U_{\nu,z\bar z}
 =\tfrac14(\kappa_z-\kappa U_{\nu,z}).
\end{equation}
We construct $Q_\nu$ with
$\partial_{\bar z}Q_\nu=\partial_{\bar z}T_\nu$ away from the sources.

Choose concentric coordinate disks $D_1\Subset D_2\Subset D_3$ about
$q$. On $D_2$, Green's representation gives
\begin{equation}\label{3.11}
 U_{\nu,z}(z)=B_\nu(z)+\sum_{i=+,-}\frac1{z-p_{i,\nu}}
 -\frac1{2\pi}\int_{D_3}\frac{F_\nu(y)}{z-y}\dif A_y.
\end{equation}
The kernels contributing to $B_\nu$ are smooth in $z$ on
$\overline{D_2}$. Weak convergence therefore gives $B_\nu\to B_0$
smoothly there. In the limit, the two source poles cancel the
term $-2/(z-q)$ from $4\pi\delta_q$. We obtain
\begin{equation}\label{3.12}
 U_{0,z}(z)=B_0(z)-\frac1{2\pi}
 \int_{D_3}\frac{F_0(y)}{z-y}\dif A_y,\qquad U_0=\log F_0.
\end{equation}

For a function $a$ on $D_2$, use the Cauchy transform
\[
 \cC a(z)=\frac1\pi\int_{D_2}\frac{a(\xi)}{z-\xi}\dif A_\xi,
 \qquad \partial_{\bar z}\cC a=a\quad\text{in }D_2.
\]
For $y\in D_3$, define
\begin{align}
 E_\kappa(z,y)
 &=\cC\left[\frac{\kappa(\cdot)-\kappa(y)}{\cdot-y}\right](z),
 \label{3.13}\\
 A_\kappa(z,y)
 &=\kappa(y)\frac{\bar z-\bar y}{z-y}+E_\kappa(z,y),
 \qquad z\ne y.
 \label{3.14}
\end{align}
Using $\partial_{\bar z}(z-y)^{-1}=\pi\delta_y$ and
$(\bar z-\bar y)\delta_y=0$, we obtain
\begin{equation}\label{3.15}
 \partial_{\bar z}A_\kappa(z,y)=\frac{\kappa(z)}{z-y}
 \qquad\text{in distributions}.
\end{equation}
The bound
\[
 \left|\frac{\kappa(\xi)-\kappa(y)}{\xi-y}\right|
 \leq \|\nabla\kappa\|_{L^\infty(D_3)}
\]
and dominated convergence show that the quotient depends continuously
on $y$ as an element of $L^p(D_2)$. For $p>2$, the Cauchy transform maps $L^p(D_2)$ to
$C^{0,\alpha}(\overline{D_1})$ for $0<\alpha<1-2/p$.
Thus $E_\kappa$ is bounded and jointly continuous on
$\overline{D_1}\times\overline{D_3}$.

Set
\begin{equation}\label{3.16}
\begin{split}
 Q_\nu(z)=\cC\left[\frac{\kappa_z-\kappa B_\nu}{4}\right](z)
 -\frac14\sum_{i=+,-}A_\kappa(z,p_{i,\nu})+\frac1{8\pi}\int_{D_3}F_\nu(y)A_\kappa(z,y)\dif A_y.
\end{split}
\end{equation}
By \eqref{3.11} and \eqref{3.15},
$\partial_{\bar z}Q_\nu=\partial_{\bar z}T_\nu$ away from the sources.
The integral in \eqref{3.16} is continuous in $z$ by dominated
convergence, since $A_\kappa$ is bounded and continuous off $y=z$.
Define the finite part by
\[
 Q_\nu^{\rm fin}(p_{i,\nu})
 =\lim_{z\to p_{i,\nu}}\left[
 Q_\nu(z)+\frac{\kappa(p_{i,\nu})}{4}
 \frac{\overline{z-p_{i,\nu}}}{z-p_{i,\nu}}
 \right].
\]
Set
\begin{equation}\label{3.17}
 Q_0(z)=\cC\left[\frac{\kappa_z-\kappa B_0}{4}\right](z)
 +\frac1{8\pi}\int_{D_3}F_0(y)A_\kappa(z,y)\dif A_y.
\end{equation}
The two source terms in \eqref{3.16} cancel the contribution of
$4\pi\delta_q$ in the limit. This yields uniform convergence $Q_\nu\to Q_0$ on
compact subsets of $D_1\setminus\{q\}$ by combining local boundedness of
the densities near $y=z$ with weak convergence away from $y=z$.

\begin{lemma}\label{le3.3}
For both $i=+,-$,
\begin{equation}\label{3.18}
 Q_\nu^{\rm fin}(p_{i,\nu})
 =Q_0(q)+\frac{\kappa(q)}4\frac{\bar d_\nu}{d_\nu}+o(1).
\end{equation}
\end{lemma}

\begin{proof}
The term containing $B_\nu$ in \eqref{3.16} converges uniformly.
By joint continuity of $E_\kappa$, the two source terms and the
concentrated mass contribute in the limit
\[
 -\frac14E_\kappa(q,q)-\frac14E_\kappa(q,q)
 +\frac{4\pi}{8\pi}E_\kappa(q,q)=0.
\]
The integral against $F_0$ gives the corresponding part of $Q_0(q)$.

Put $\tau_\nu=\bar d_\nu/d_\nu$. Since
$p_{i,\nu}-p_{j,\nu}=\pm2d_\nu$, the angular term from the other
source is $-\kappa(q)\tau_\nu/4+o(1)$. For the integral we have
\[
\begin{aligned}
\frac1{8\pi}\int_{D_3}F_\nu(y)\kappa(y)
\frac{\overline{p_{i,\nu}-y}}{p_{i,\nu}-y}\dif A_y=\frac12\kappa(q)\tau_\nu +\frac1{8\pi}\int_{D_3}F_0(y)\kappa(y)
       \frac{\bar q-\bar y}{q-y}\dif A_y+o(1).
\end{aligned}
\]
Indeed, on $B_{\nu,\eta}=B_{\eta|d_\nu|}(s_\nu)$ with
$0<\eta<1/2$, the angular quotient is $\tau_\nu+O(\eta)$.
For fixed small $r>0$, \eqref{3.2} and weak convergence give
\[
 \int_{B_{\nu,\eta}}F_\nu\dif A\longrightarrow4\pi,\qquad
 \int_{B_r(q)\setminus B_{\nu,\eta}}F_\nu\dif A
 \longrightarrow\int_{B_r(q)}F_0\dif A=O(r^2).
\]
On $D_3\setminus B_r(q)$, the kernel converges uniformly. Weak convergence
then gives the integral against $F_0$ there.
Take $\nu\to\infty$ first and then $\eta,r\downarrow0$ to obtain
the displayed identity. The angular contribution is therefore
$-\kappa(q)\tau_\nu/4+\kappa(q)\tau_\nu/2=\kappa(q)\tau_\nu/4$,
which proves \eqref{3.18}.
\end{proof}

\begin{proof}[Proof of Theorem~\ref{th3.1}]
\textit{Step 1.} We expand $T_\nu$ at the two source points.

Write $z=p_i+s$, where $p_i=p_{i,\nu}$, and
$U_\nu(z)=\log|s|^2+v_{i,\nu}(z)$ with $v_{i,\nu}$ smooth.
The regular part of \eqref{3.9} gives
$v_{i,\nu,z\bar z}(p_i)=\kappa(p_i)/4$ since $F_\nu(p_i)=0$.
For fixed $\nu$, Taylor expansion as $s\to0$ gives
\begin{equation}\label{3.19}
 T_\nu(z)=-\frac{3}{2s^2}+\frac{a_{i,\nu}}s
 -\frac{\kappa(p_i)}4\frac{\bar s}s
 -\frac12a_{i,\nu}^2+o(1),
 \qquad a_{i,\nu}=-v_{i,\nu,z}(p_i).
\end{equation}
The difference $T_\nu-Q_\nu$ is holomorphic away from the
sources. Its angular terms cancel. After we subtract the poles
in \eqref{3.19}, the remaining function is bounded near each source.
Thus there is a holomorphic function $J_\nu$ on $D_1$ such that
\begin{equation}\label{3.20}
 T_\nu(z)=J_\nu(z)+Q_\nu(z)
 +\sum_{i=+,-}\left(-\frac{3}{2(z-p_{i,\nu})^2}
 +\frac{a_{i,\nu}}{z-p_{i,\nu}}\right).
\end{equation}

We suppress $\nu$ in $p_i,d,c,a_i$. The center equation in
\eqref{3.1} then gives
\begin{equation}\label{3.21}
 a_i=\frac{\alpha(c)}{p_i-p_j}-L_i^c,
 \qquad \alpha(c)=-1-\frac2c,
\end{equation}
where
\begin{equation}\label{3.22}
 L_i^c=(\log h_2+\phi)_z(p_i)-\frac{g_{1,z}(p_i)}c
 -4\pi R_z(p_i,p_i)
 -\left(4\pi+\frac{8\pi}c\right)R_z(p_i,p_j).
\end{equation}
At $c=-1/2$, this becomes
\begin{equation}\label{3.23}
 L_i=2g_{1,z}(p_i)+(\log h_2+\phi)_z(p_i)
 -4\pi R_z(p_i,p_i)+12\pi R_z(p_i,p_j).
\end{equation}
Since $r_{g,z}(q)=2R_z(q,q)+\phi_z(q)/(4\pi)$, we have
$L_i^c\to C_z(q)$ as $p_i,p_j\to q$ and $c\to-1/2$.
Taking the finite part of \eqref{3.20} at $p_i$ and using
\eqref{3.19} gives
\[
 J_\nu(p_i)+Q_\nu^{\rm fin}(p_i)
 =-\frac12a_i^2+\frac{3}{2(p_i-p_j)^2}-\frac{a_j}{p_i-p_j}.
\]
Substituting \eqref{3.21}, we obtain
\begin{align}
 J_\nu(p_+)+Q_\nu^{\rm fin}(p_+)
 &=\frac{3+2\alpha-\alpha^2}{8d^2}
 +\frac{\alpha L_+^c+L_-^c}{2d}-\frac12(L_+^c)^2,
 \label{3.24}\\
 J_\nu(p_-)+Q_\nu^{\rm fin}(p_-)
 &=\frac{3+2\alpha-\alpha^2}{8d^2}
 -\frac{L_+^c+\alpha L_-^c}{2d}-\frac12(L_-^c)^2.
 \label{3.25}
\end{align}

\medskip
\noindent\textit{Step 2.} We prove $C_z(q)=0$ and derive the equation
for $d_\nu^2$.

Let $\Gamma\subset D_1$ be a fixed positively oriented circle
about $q$. On $\Gamma$, $T_\nu-Q_\nu\to T_0-Q_0$ uniformly, where
$T_0=U_{0,zz}-U_{0,z}^2/2$. Since
$\Delta_0U_0+2e^{U_0}=\kappa$, the calculation in
\eqref{3.10}--\eqref{3.17} gives
$\partial_{\bar z}(T_0-Q_0)=0$ inside $\Gamma$. By
\eqref{3.20}--\eqref{3.21},
\[
 \int_\Gamma(T_\nu-Q_\nu)\dif z
 =2\pi i(a_++a_-)
 =-2\pi i(L_+^c+L_-^c)\longrightarrow-4\pi i C_z(q).
\]
The integral tends to zero. Hence $C_z(q)=0$ and $\nabla C(q)=0$.

In \eqref{3.20}, the terms containing $\alpha(c_\nu)$ sum to
$\alpha(c_\nu)/[(z-p_{+,\nu})(z-p_{-,\nu})]$. Since
$\alpha(c_\nu)\to3$, these terms cancel the limiting double poles.
The remaining pole terms tend to zero since $L_\pm^{c_\nu}\to0$.
Thus $J_\nu\to T_0-Q_0$ uniformly on $\Gamma$. The maximum
modulus principle then gives
\begin{equation}\label{3.26}
 J_\nu(p_{+,\nu}),\ J_\nu(p_{-,\nu})
 \longrightarrow J_0(q)=T_0(q)-Q_0(q).
\end{equation}

Set
\begin{equation}\label{3.27}
 L_z=C_{zz}(q)-24\pi R_{zp},\qquad
 L_{\bar z}=C_{z\bar z}(q)-24\pi R_{z\bar p}.
\end{equation}
Expanding \eqref{3.22} about the midpoint gives
\begin{equation}\label{3.28}
 L_{+,\nu}^{c_\nu}-L_{-,\nu}^{c_\nu}
 =2L_zd_\nu+2L_{\bar z}\bar d_\nu+o(|d_\nu|).
\end{equation}
Indeed, writing $F(z,p)$ for the right-hand side of \eqref{3.23},
we have $F(z,z)=C_z(z)$, $F_p(q,q)=12\pi R_{zp}$, and
$F_{\bar p}(q,q)=12\pi R_{z\bar p}$. Hence
\[
\begin{aligned}
 (F_z-F_p)(q,q)&=C_{zz}(q)-2F_p(q,q)=L_z,\\
 (F_{\bar z}-F_{\bar p})(q,q)
 &=C_{z\bar z}(q)-2F_{\bar p}(q,q)=L_{\bar z}.
\end{aligned}
\]
These are the coefficients in the expansion of
$F(s+d,s-d)-F(s-d,s+d)$. The dependence on $c$ in
\eqref{3.22} is smooth. Replacing $-1/2$ by $c_\nu$ in this
difference therefore contributes $O(|c_\nu+1/2|\,|d_\nu|)=o(|d_\nu|)$.

Averaging \eqref{3.24} and \eqref{3.25}, we have
\[
 \frac{(\alpha-1)(L_+^c-L_-^c)}{4d}
 =L_z+L_{\bar z}\frac{\bar d}{d}+o(1),
 \qquad (L_\pm^c)^2=o(1).
\]
Lemma~\ref{le3.3} and \eqref{3.26} give the average of the left-hand
sides as $T_0(q)+\kappa(q)\bar d/(4d)+o(1)$. Thus
\begin{equation}\label{3.29}
 T_0(q)-L_z+
 \left(\frac{\kappa(q)}4-L_{\bar z}\right)\frac{\bar d_\nu}{d_\nu}
 =\frac{3+2\alpha(c_\nu)-\alpha(c_\nu)^2}{8d_\nu^2}+o(1).
\end{equation}
Multiplication by $d_\nu^2$ gives
\begin{equation}\label{3.30}
 \mathcal A(q)\zeta_\nu+\mathcal B(q)|\zeta_\nu|
 =\frac{3+2\alpha(c_\nu)-\alpha(c_\nu)^2}{8}+o(|\zeta_\nu|),
\end{equation}
where
\[
 \mathcal A(q)=T_0(q)-C_{zz}(q)+24\pi R_{zp},\qquad
 \mathcal B(q)=\frac{\kappa(q)}4-C_{z\bar z}(q)+24\pi R_{z\bar p}.
\]
By \eqref{3.21},
\begin{equation}\label{3.31}
 \frac{3+2\alpha(c)-\alpha(c)^2}{8}
 =-\frac{2c+1}{2c^2}=\frac{-c-\tfrac12}{c^2}.
\end{equation}

\medskip
\noindent\textit{Step 3.} We identify the coefficients in \eqref{3.30}.

All quantities below are evaluated at $q$, with derivatives of
$R$ taken at $(q,q)$. By symmetry of $R$ and
$\Delta_0R(\cdot,p)=e^\phi$, we have
\[
 R_{z\bar z}=\tfrac14e^\phi,\qquad
 r_{g,z\bar z}=\tfrac12e^\phi+2R_{z\bar p}
 +\tfrac1{4\pi}\phi_{z\bar z},\qquad
 \phi_{z\bar z}=-\tfrac12K_ge^\phi.
\]
Thus
\[
\begin{aligned}
 \mathcal B
 &=\phi_{z\bar z}+(\rho/2-2\pi)e^\phi
   -2(\log h_1)_{z\bar z}-12\pi r_{g,z\bar z}
   +24\pi R_{z\bar p}\\
 &=-2(\log h_1)_{z\bar z}-2\phi_{z\bar z}
   +(\rho/2-8\pi)e^\phi.
\end{aligned}
\]
Using $4\partial_z\partial_{\bar z}=e^\phi\Delta_g$, we obtain
\begin{equation}\label{3.32}
 \mathcal B=e^\phi\left(\frac{\rho}{2}-8\pi+K_g
 -\frac12\Delta_g\log h_1\right)
 =-\frac12e^\phi\kappa_1.
\end{equation}
For $\mathcal A$, we use
\[
 r_{g,zz}=2R_{zz}+2R_{zp}+\frac{\phi_{zz}}{4\pi},
 \qquad U_0=\psi_2+\mathrm{const.},
 \qquad 2\psi_{1,z}+\psi_{2,z}=C_z=0.
\]
Hence
\begin{align*}
 \mathcal A
 &=\psi_{2,zz}-\tfrac12\psi_{2,z}^2-C_{zz}+24\pi R_{zp}\\
 &=-2(\log h_1)_{zz}-2\phi_{zz}-24\pi R_{zz}
   +w_{*,zz}-\tfrac12\psi_{2,z}^2\\
 &=-2(\psi_{1,zz}+\psi_{1,z}^2)=-2\mathcal T.
\end{align*}
Substitution in \eqref{3.30}, using \eqref{3.31}, gives \eqref{3.6}.

\medskip
\noindent\textit{Step 4.} We deduce \eqref{3.7} and exclude
collisions under \eqref{3.8}.

If $c_\nu\equiv-1/2$, divide \eqref{3.6} by $|\zeta_\nu|$
and take moduli in the limit to obtain \eqref{3.7}.
If $c_\nu\leq-1/2$ and
$\gamma=e^{\phi(q)}\kappa_1(q)-4|\mathcal T(q)|>0$,
the real part of \eqref{3.6} gives, for all large $\nu$,
\[
 0\leq\frac{-c_\nu-\tfrac12}{c_\nu^2}
 \leq-\frac{\gamma}{2}|\zeta_\nu|+o(|\zeta_\nu|)
 \leq-\frac{\gamma}{4}|\zeta_\nu|<0.
\]
\end{proof}

\begin{corollary}\label{co3.4}
Assume \eqref{3.8} for every residual scalar solution at parameter
$2(\rho-4\pi)$ and every $q\in\operatorname{Crit}(C)$.
For solutions of \eqref{3.1},
\[
 c_\nu\leq-1/2,\quad c_\nu\to-1/2
 \quad\Longrightarrow\quad
 \liminf_{\nu\to\infty}\operatorname{dist}_g(p_{+,\nu},p_{-,\nu})>0.
\]
Consequently, for some $\delta>0$, the solution set is compact
for $c\in[-1/2-\delta,-1/2]$.
\end{corollary}

\begin{proof}
If the center separation tends to zero along such a sequence,
then $p_{\pm,\nu}\to q$ along a subsequence. Theorem~\ref{th3.1}
gives $q\in\operatorname{Crit}(C)$, and \eqref{3.8} excludes the
collision. This gives uniform separation on a short closed interval.
If this failed, we could choose
$c_\nu\in[-1/2-1/\nu,-1/2]$ with center separation less than $1/\nu$.

For separated centers, every possible concentration point has
$m\leq1$. Quantization and \eqref{2.13} give
$n\in\mathbb Z_{\geq0}$ and $n<(m+1)/2\leq1$, so $n=0$.
The scalar a priori estimates and elliptic regularity give compactness.
\end{proof}

\subsection{Perturbing the weights to exclude collisions of two centers}\label{sec:convenient-shadow-degree}

Fix $\rho>4\pi$, $\rho\notin4\pi\mathbb N$, and put
$\sigma=\rho-4\pi$.  The regular residual scalar equation is
\begin{equation}\label{3.33}
 -\Delta_g z=2\sigma
 \left(\frac{h_2e^z}{\int_Mh_2e^z\dV}-1\right),
 \qquad \int_Mz\dV=0.
\end{equation}
Write $r_g$ for the intrinsic
Robin function and set
\begin{equation}\label{3.34}
 V=2\log h_1+\log h_2+12\pi r_g.
\end{equation}
One may first choose $h_2$ and a Morse function $V$, and then define
$h_1$ by \eqref{3.34}.

For $q\in\operatorname{Crit}(V)$, choose an isothermal coordinate
$\xi$ with $\xi(q)=0$, $g=e^\phi|\dif \xi|^2$, and $\phi(q)=0$. We write
$$
 G(\xi,q)=-\frac{1}{2\pi}\log|\xi|+R(\xi,q).
$$
For a residual solution $z_s$ define
\begin{align}
 \psi_{1,q,s}
 &=\log h_1+\phi+12\pi R(\cdot,q)-\frac12z_s,
 \label{3.35}\\
 \mathcal T_{q,s}
 &=\psi_{1,q,s,\xi\xi}(q)+\psi_{1,q,s,\xi}(q)^2,\qquad
 \kappa_q=\Delta_g\log h_1(q)+16\pi-\rho-2K_g(q).
 \label{3.36}
\end{align}
The coordinate normalization removes the factor $e^{\phi(q)}$ in
Theorem~\ref{th3.1}.

\begin{lemma}\label{le3.5}
Assume that $V$ is Morse and that every solution of \eqref{3.33}
is nondegenerate. Fix $0<\alpha<1$. For every $\eta,\varepsilon>0$
and every $C^1$ neighborhood $\mathcal U$ of $(h_1,h_2)$, there
exists a smooth positive pair
$(\widetilde h_1,\widetilde h_2)\in\mathcal U$ with the following
properties.
\begin{enumerate}
\item The center potential satisfies
\[
 2\log\widetilde h_1+\log\widetilde h_2+12\pi r_g=V
 \qquad\text{on }M.
\]

\item For each $z_s$, the equation with weight $\widetilde h_2$ has a unique solution $\widetilde z_s$ in a fixed $C^{2,\alpha}(M)$ neighborhood of $z_s$. Each $\widetilde z_s$ is nondegenerate and has the same local scalar Leray--Schauder degree as $z_s$.

\item For every $q\in\operatorname{Crit}(V)$ and every continuation
$\widetilde z_s$, the coefficients associated with
$(\widetilde h_1,\widetilde h_2,\widetilde z_s)$ satisfy
\begin{equation}\label{3.37}
 \widetilde\kappa_q>4|\widetilde{\mathcal T}_{q,s}|+\eta.
\end{equation}
\end{enumerate}
\end{lemma}

\begin{proof}
\emph{Step 1.}
Since $2\sigma\notin8\pi\mathbb N$, the solution set of
\eqref{3.33} is compact. By nondegeneracy, each solution is
isolated. Hence there are finitely many solutions, which we denote by $z_1,\ldots,z_N$.  The equation defines a smooth map
\[
\begin{aligned}
 \Phi:C^{2,\alpha}_0(M)\times C^{0,\alpha}_+(M)
       &\to C^{0,\alpha}_0(M)\\
 (z,h)&\mapsto-\Delta_g z
 -2\sigma\left(\frac{he^z}{\int_M he^z\dV}-1\right).
\end{aligned}
\]
Its derivative with respect to $z$ is invertible at each $z_s$.
The implicit function theorem therefore gives a unique nondegenerate
continuation near each $z_s$ for small $C^{0,\alpha}$ perturbations of
$h_2$, with convergence in $C^{2,\alpha}$.

For sufficiently small $C^1$ perturbations of $h_2$, there are no additional
solutions. Otherwise, there would be positive smooth $h_{2,n}\to h_2$ in $C^1$
and solutions outside those finitely many neighborhoods. If these solutions
were unbounded, the scalar blow-up theorem for positive coefficients
converging in $C^1$ would force
$2\sigma\in8\pi\mathbb N$. This contradicts the choice of $\rho$.
Since the solutions are bounded, elliptic estimates and the uniform
$C^1$ coefficient bounds give uniform $C^{2,\beta}$ estimates for any
$\alpha<\beta<1$. A subsequence then converges in $C^{2,\alpha}$ to an
original solution, contradicting the local uniqueness of its continuation.
This proves stability of the solution set. Nondegeneracy preserves each
local degree.

If $N=0$, the same compactness argument shows that the residual
solution set remains empty in a small $C^1$ neighborhood.  The
assertion about the positive gap is vacuous, and the original weights
already have all the required properties.  Hence assume $N\geq1$.

\emph{Step 2.}
Let $q_1,\ldots,q_m$ be the critical points of $V$. Choose pairwise disjoint
coordinate disks as above. Take a fixed smooth cutoff $\vartheta$ supported
in the unit disk and equal to one near zero. For numbers $L_j>0$ and
sufficiently small radii $a_j>0$, set
\begin{align}\label{def:br}
 b_j(\xi)=L_j|\xi|^2\vartheta(\xi/a_j),\qquad
 b=\sum_{j=1}^m b_j,\qquad
 \widetilde h_2=h_2e^{-b},\quad
 \widetilde h_1=h_1e^{b/2}.
\end{align}
For each fixed $L_j$ the estimates
$$
 \|b_j\|_{C^0}=O(L_ja_j^2),\qquad
 \|b_j\|_{C^1}=O(L_ja_j),\qquad
 \|b_j\|_{C^{0,\alpha}}=O(L_ja_j^{2-\alpha})
$$
hold as $a_j\downarrow0$.  At the center, however,
$$
 b(q_j)=0,\quad \dif b(q_j)=0,\quad
 b_{\xi\xi}(q_j)=0,\quad \Delta_g b(q_j)=4L_j.
$$
In particular,
$$
 2\log\widetilde h_1+\log\widetilde h_2
       =2\log h_1+\log h_2.
$$
For fixed $L_j$, the perturbation tends to zero in $C^1$ as
$a_j\to0$.

\emph{Step 3.}
Set $d_s=\widetilde z_s-z_s$. At $q_j$, we have
\begin{align}
 \widetilde\kappa_{q_j}
 &=\kappa_{q_j}+2L_j,\label{3.38}\\
 \widetilde{\mathcal T}_{q_j,s}-\mathcal T_{q_j,s}
 &=-\frac12d_{s,\xi\xi}(q_j)
   -\psi_{1,q_j,s,\xi}(q_j)d_{s,\xi}(q_j)
   +\frac14d_{s,\xi}(q_j)^2.
 \label{3.39}
\end{align}
Here $b_\xi(q_j)=b_{\xi\xi}(q_j)=0$. By Step~1, $d_s\to0$ in
$C^{2,\alpha}(M)$ as $b\to0$ in $C^{0,\alpha}(M)$. Equation~\eqref{3.39}
then gives $\max_{j,s}|\widetilde{\mathcal T}_{q_j,s}-\mathcal T_{q_j,s}|\to0$.

\emph{Step 4.}
Choose $L_j$ first, so large that
$$
 \kappa_{q_j}+2L_j>
 4\max_{1\leq s\leq N}|\mathcal T_{q_j,s}|+\eta+4.
$$
Shrink the radii so that $(h_1e^{tb/2},h_2e^{-tb})\in\mathcal U$
and Step~1 applies for all $t\in[0,1]$. At $t=1$, take
$\max_s\|d_s\|_{C^{2,\alpha}(M)}<\varepsilon$ and
$\max_{j,s}|\widetilde{\mathcal T}_{q_j,s}-\mathcal T_{q_j,s}|<1$.
Then \eqref{3.38} gives
\[
 4|\widetilde{\mathcal T}_{q_j,s}|+\eta
 <4|\mathcal T_{q_j,s}|+\eta+4
 <\kappa_{q_j}+2L_j=\widetilde\kappa_{q_j}.
\]
This proves \eqref{3.37}. Step~1 gives nondegeneracy along the path.
\end{proof}

\begin{remark}\label{re3.6}
Fix $p>2$ and an integer $k\geq2$. With $\ell=\log h_2$, consider
$(z,\ell)\mapsto\Phi(z,e^\ell)$ from
$W^{k+2,p}_0(M)\times W^{k,p}(M)$ to $W^{k,p}_0(M)$, where the
subscript $0$ denotes zero mean. Set
$\mu=h_2e^z/\int_Mh_2e^z\dV$. Then
\[
 D_\ell\bigl[\Phi(z,e^\ell)\bigr][\xi]
 =-2\sigma\mu\left(\xi-\int_M\mu\xi\dV\right).
\]
For any $f\in W^{k,p}_0(M)$, taking $\xi=-f/(2\sigma\mu)$ gives
$D_\ell[\Phi(z,e^\ell)][\xi]=f$.

The $z$-linearization $A$ is Fredholm of index zero. It is
$-\Delta_g$ plus a zeroth-order compact perturbation between the
indicated Sobolev spaces. Write $B=D_\ell\Phi$. The right inverse
of $B$ given above is bounded, so the zero set is a smooth
Banach submanifold with tangent space
$\{(\dot z,\dot\ell):A\dot z+B\dot\ell=0\}$.
For the projection to $\ell$, the kernel of its differential is
$\ker A\times\{0\}$ and the range is
$\{\dot\ell:B\dot\ell\in\operatorname{ran} A\}$.
Because $B$ is onto and $\operatorname{ran} A$ is closed,
this range has codimension $\dim\operatorname{coker}A$.
The projection is thus a smooth Fredholm map of index zero.
The two Sobolev spaces are separable, so the ambient product
and the zero set have countable bases. Theorem~\textup{(1.3)}
of \cite[p.~862]{Smale65} applies, component by component if
necessary. Its regular values are residual. At any such value,
the equation above shows that $A$ is onto and therefore invertible. The compactness
argument in Step~1 gives openness for smooth positive weights in
the $C^1$ topology. To obtain density in any prescribed $C^r$
neighborhood of a smooth weight, choose $k$ large enough that
$W^{k,p}$ embeds in $C^{\max(r,1)}$. First choose a regular
log-weight in a small $W^{k,p}$ neighborhood. To check that regularity
persists under smooth approximation of that weight in $W^{k,p}$,
consider a convergent sequence of log-weights in that norm. The
noncritical scalar estimate gives a uniform $C^0$ bound for the
solutions because the coefficients are uniformly positive and
Lipschitz and converge in $C^1$. The scalar equation
then gives boundedness first in $W^{2,p}$ and, by the Sobolev
algebra and elliptic estimates, in $W^{k+2,p}$. Compactness in
$W^{k,p}$ and convergence of the right-hand sides in $W^{k,p}$
give convergence of the solutions in $W^{k+2,p}$. In this topology the
linearizations converge in operator norm. Thus no sequence of
degenerate solutions for nearby weights can converge to a solution
at a regular weight. Smooth approximation in $W^{k,p}$ preserves
regularity in a sufficiently small neighborhood and keeps the
chosen $C^r$ bound. Thus the smooth positive weights with this property form an
open dense set in the $C^\infty$ topology.
For each such weight, the solution set is compact and discrete, hence
finite.
\end{remark}

\subsection{A compact homotopy for the shadow systems}\label{typed-sec:rotation}

Fix an integer $m\geq1$ and take $4\pi k<\rho<4\pi(k+1)$. A \emph{type} is a pair
of nonnegative integers $(a,b)$ with $a+2b=m$ and $b\leq k$. Put
$\sigma_b=\rho-4\pi b>0$.  Write $P=(p_1,\ldots,p_a)$ and
$Q=(q_1,\ldots,q_b)$, with all their coordinates distinct.  Permutations
within each type are identified.  Set
$$
 H_P=h_2e^{-4\pi\sum_iG(\cdot,p_i)},\qquad
 g_1=\log h_1+4\pi r_g,\qquad
 V=2\log h_1+\log h_2+12\pi r_g.
$$
Recall that for each type $(a,b)$, the shadow system \eqref{3.40} reads
\begin{equation}\tag{2.15}\label{e:2.155}
\left\{
\begin{aligned}
 F_0&=-\Delta_gw-2\sigma_b
 \left(\frac{H_Pe^w}{\int_MH_Pe^w\dV}-1\right)=0,
 \qquad \int_Mw\dV=0,\\
 F_i^A&=\nabla g_1(p_i)-\frac12\nabla w(p_i)
 +8\pi\sum_{l\ne i}\nabla_xG(p_i,p_l)
 +12\pi\sum_j\nabla_xG(p_i,q_j)=0,\\
 F_j^B&=\nabla V(q_j)+12\pi\sum_i\nabla_xG(q_j,p_i)
 +24\pi\sum_{l\ne j}\nabla_xG(q_j,q_l)=0.
\end{aligned}
\right.
\end{equation}

\subsubsection{The collision criterion at mixed centers}

Use $\psi_{1,j}$, $\mathcal T_j$ and $\kappa_m$ from the mixed-center
condition \eqref{3.41}. We check its coordinate invariance before
using it to control the compact homotopy.

Under a holomorphic change of coordinates
$z=f(\xi)=\lambda\xi+\beta\xi^2+O(\xi^3)$ with $\lambda\ne0$,
the regular term transforms as
\[
 \widetilde\psi_{1,j}
 =\psi_{1,j}\circ f+2\log|f'|-6\log|f(\xi)/\xi|.
\]
At $q_j$, this gives
\[
\begin{aligned}
 \widetilde\psi_{1,j,\xi}
 &=\lambda\psi_{1,j,z}-\beta/\lambda,\\
 \widetilde\psi_{1,j,\xi\xi}
 &=\lambda^2\psi_{1,j,zz}
   +2\beta\psi_{1,j,z}-\beta^2/\lambda^2.
\end{aligned}
\]
Thus $\widetilde{\mathcal T}_j=\lambda^2\mathcal T_j$ and
$e^{\widetilde\phi(q_j)}=|\lambda|^2e^{\phi(q_j)}$.
It follows that \eqref{3.41} is independent of the isothermal
coordinate, and we may take $\phi(q_j)=0$.

Suppose two ordinary centers converge to $q$, and let $I$ index
the remaining ordinary centers. Near $q$, the weights in the
local two-center equations are given by
\[
\begin{aligned}
 \log h_{1,\mathrm{ext}}
 &=\log h_1+8\pi\sum_{i\in I}G(\cdot,p_i)      +12\pi\sum_{\ell=1}^bG(\cdot,q_\ell),\\
 \log h_{2,\mathrm{ext}}
 &=\log h_2-4\pi\sum_{i\in I}G(\cdot,p_i).
\end{aligned}
\]
Since $\Delta_gG(\cdot,p)=1$ away from $p$, the coefficient
in \eqref{3.5} for these weights and parameter $\sigma_b$ is
\[
\begin{aligned}
 (\Delta_g \log h_{1,\mathrm{ext}})(q)
      &+16\pi-\sigma_b-2K_g(q)\\
 &=(\Delta_g\log h_1)(q)+8\pi(a-2)+12\pi b+16\pi-(\rho-4\pi b)-2K_g(q)\\
 &=\kappa_m(q).
\end{aligned}
\]
After the replacement, the regular term in \eqref{3.4} is
$\psi_{1,j}$ at the new mixed center $q_j=q$. For each other
pair collapsing at $r$, the two ordinary-center terms and the
corresponding term from the scalar limit combine as
$16\pi G(\cdot,r)-4\pi G(\cdot,r)=12\pi G(\cdot,r)$,
which is the contribution of the new mixed center at $r$.
\subsubsection{A rotational perturbation of the center equations}

Let $J$ denote rotation through $\pi/2$ in the positive direction
on each tangent plane of $M$. For $-1/2\leq c\leq0$ and
$\eta\in\mathbb R$, we replace the ordinary-center equations in
\eqref{3.40} by
\begin{equation}\label{3.42}
\begin{aligned}
0={}&\nabla g_1(p_i)+c\nabla w(p_i)
 +8\pi\sum_{l\ne i}\nabla_xG(p_i,p_l)
 +12\pi\sum_j\nabla_xG(p_i,q_j)\\
 &+8\pi\eta\sum_{l\ne i}J\nabla_xG(p_i,p_l),
 \qquad i=1,\ldots,a.
\end{aligned}
\end{equation}

Consider a sequence of solutions of \eqref{3.42} with $c_\nu\to c$
and $\eta_\nu$ bounded. Suppose that $u$ ordinary centers and $v$
mixed centers converge to $q$, and write $4\pi n$ for the local scalar
mass at $q$. As in Step~1 of Lemma~\ref{lem:successive-pair-collisions},
\[
 n\in\mathbb Z_{\geq0},\qquad n<\frac{u+1}{2},\qquad
 S_\nu=\sum_{i=1}^u p_{i,\nu}\cdot\nabla w_\nu(p_{i,\nu})
 \longrightarrow 4n(n-u-1),
\]
where $S_\nu$ is computed in an isothermal coordinate centered at $q$.
Pairing the center equations as in the same step gives
$c_\nu S_\nu=2Q(u,v)+o(1)$, where
$Q(u,v)=u(u-1)+3uv+3v(v-1)$. Indeed, the rotated logarithmic terms
cancel in pairs because $(p_i-p_l)\cdot J(p_i-p_l)=0$. The regular
parts contribute $O(|\eta_\nu|\sum_i|p_{i,\nu}|)=o(1)$. Hence
\begin{equation}\label{3.43}
 2c\,n(n-u-1)=Q(u,v).
\end{equation}

Since $0\leq-2c\leq1$, we have $Q(u,v)\leq n(u+1-n)$ and thus
\[
 (2n-u-1)^2\leq (u+1)^2-4Q(u,v)=4-3(u+2v-1)^2.
\]
For a collision, $u+v\geq2$, so $u+2v\leq2$ forces $(u,v)=(2,0)$.
Then $n<3/2$ and \eqref{3.43} give $n=1$ and $c=-1/2$.

If $(c_\nu,\eta_\nu)\to(-1/2,0)$, the limiting equations follow
as in Step~2 of Lemma~\ref{lem:successive-pair-collisions}.
Indeed, when we add the two center equations for each colliding pair in
Euclidean coordinates, the unrotated and rotated singular logarithmic
terms cancel separately. The remaining rotated terms are $O(|\eta_\nu|)$.
The resulting identity bounds the sum of the two gradients of $w_\nu$
because $c_\nu$ is bounded away from zero. It has the same limit as the
unperturbed pair identity when $c_\nu\to-1/2$ and $\eta_\nu\to0$.
The scalar translation identity is unchanged.
Replacing the $\ell$ colliding pairs by $\ell$ mixed centers
gives a solution of \eqref{3.40} of type $(a-2\ell,b+\ell)$,
with $\sigma_{b+\ell}=\sigma_b-4\pi\ell>0$.

\begin{lemma}\label{le3.7}
Suppose that two ordinary centers $p_{\pm,\nu}$ of a sequence
of solutions of \eqref{3.42} converge to $q$ and that
$(c_\nu,\eta_\nu)\to(-1/2,0)$. In an isothermal coordinate
with $\phi(q)=0$, write $p_{\pm,\nu}=s_\nu\pm d_\nu$ and
$\zeta_\nu=d_\nu^2$. In the limiting solution of \eqref{3.40}, let $q_j=q$ be the
mixed center formed by the colliding pair and set
$\mathcal T(q)=\mathcal T_j$. Then
\begin{equation}\label{3.44}
 -2\mathcal T(q)\zeta_\nu-\frac12\kappa_m(q)|\zeta_\nu|
 =\frac{-c_\nu-\tfrac12+i(c_\nu+1)\eta_\nu+\eta_\nu^2/2}
 {c_\nu^2}+o(|\zeta_\nu|).
\end{equation}
\end{lemma}

\begin{proof}
Choose a fixed isothermal disk $D$ about $q$ whose closure stays
away from all exterior limiting centers. Cancel the common conformal
factor in the center equations. In these coordinates $J$ is the positive
Euclidean quarter-turn. Write $G(x,y)=-(2\pi)^{-1}\log|x-y|+R(x,y)$ on
$D\times D$. Absorb the unrotated exterior potentials into the smooth
scalar weights. For the rotation of the exterior ordinary-center terms set
\[
 X_\nu(x)=8\pi\sum_{\ell\notin\{+,-\}} J\nabla_xG(x,p_{\ell,\nu}).
\]
Separation gives $\|X_\nu\|_{C^1(D)}\leq C$ uniformly. The regular
rotated pair terms are
$Y_{\pm,\nu}=8\pi J\nabla_xR(p_{\pm,\nu},p_{\mp,\nu})$. Smoothness
of $R$ in both variables gives $|Y_{\pm,\nu}|\leq C$ and
$|Y_{+,\nu}-Y_{-,\nu}|\leq C|p_{+,\nu}-p_{-,\nu}|$. The regular
part of the rotation therefore contributes the vectors
$\eta_\nu[X_\nu(p_{\pm,\nu})+Y_{\pm,\nu}]$, bounded by $C|\eta_\nu|$
at each center, with difference bounded by $C|\eta_\nu|\,|d_\nu|$.
The singular term at $p_{\pm,\nu}$ is
\[
 -4(I+\eta_\nu J)
       \frac{p_{\pm,\nu}-p_{\mp,\nu}}{|p_{\pm,\nu}-p_{\mp,\nu}|^2}.
\]
Multiplication by $|d_\nu|$ gives the displayed singular terms in the
positive-center equation below. The scaled regular rotation is
$O(|\eta_\nu||d_\nu|)$.
By \eqref{3.43}, the local scalar mass at $q$ is $4\pi$.
With $e_\nu=d_\nu/|d_\nu|$, the rescaled center equation is
\[
 c_\nu\nabla V_\nu(e_\nu)-2e_\nu-2\eta_\nu Je_\nu+O(|d_\nu|)=0.
\]
The scalar equation and its rescalings in the proof of
Lemma~\ref{lem:pair-midpoint} are unchanged by the rotation. The
one-bubble conclusion there uses the mass $4\pi$, separation of the
two zeros after rescaling, and the scalar mass bounds. In the step
that locates the center, the equation above has the same limit as
the unrotated equation because $c_\nu\to-1/2$ and $\eta_\nu\to0$.
The location and the profile in \eqref{3.2} therefore follow with
the same normalization.

Under the convention $u_z=(u_x-iu_y)/2$, the rotation $J$ corresponds
to multiplication by $-i$. Thus \eqref{3.21} takes the form
\[
 a_i=\frac{\alpha(c,\eta)}{p_i-p_j}-L_i^{c,\eta},
 \qquad
 \alpha(c,\eta)=-1-\frac2c+\frac{2i\eta}{c},
 \qquad \{i,j\}=\{+,-\}.
\]
Compare the coefficient expressions at the same centers and with
the same scalar data. Since $|c_\nu|$ is bounded away from zero,
division by $c_\nu$ and passage to the complex representation of
the vector terms preserve both bounds just proved. Thus
\[
 L_\pm^{c,\eta}-L_\pm^c=O(|\eta|),\qquad
 (L_+^{c,\eta}-L_-^{c,\eta})-(L_+^c-L_-^c)
 =O(|\eta|\,|d|).
\]
For clarity, the local scalar coefficient in the calculation is
\[
 \kappa_\nu=\Delta_0(\log h_{2,\mathrm{ext},\nu}+\phi)
                   +(2\sigma_b-8\pi)e^\phi
            =\Delta_0(\log h_2+\phi)+(2\sigma_b-4\pi a)e^\phi .
\]
If $\ell$ pairs collapse, the lower type has $a'=a-2\ell$
and $\sigma'=\sigma_b-4\pi\ell$, so $2\sigma'-4\pi a'=
2\sigma_b-4\pi a$. The limit of the second regular field, up
to a constant, is $\log h_2+\phi-4\pi\sum_{p\in P'}G(\cdot,p)+w'$.
At the other collapsing pairs, the extra terms $-8\pi G$ in the
external second weight cancel the corresponding $8\pi G$ terms in the
scalar limit. The limit of the second regular field satisfies the
same smooth log-density equation at $q$.
In the first regular field, each other collapsing pair contributes
$16\pi G-4\pi G=12\pi G$, as in the limiting type. In the trace
coefficient one must use the pre-collision $\sigma_b$ and the
pre-collision external first weight, as in the calculation
of the two-center coefficient using \eqref{3.5} for one colliding pair. This gives $\kappa_m(q)$.
The first regular field itself is $\psi_{1,j}$ of the lower
type at $q$, up to a constant, giving $\mathcal T(q)$.

With these identifications, Lemma~\ref{le3.3} and Steps~2 and~3
of the proof of Theorem~\ref{th3.1} apply to the two local sources.
The difference of the regular rotation terms is $O(|\eta_\nu||d_\nu|)$.
It becomes $o(1)$ after division by $d_\nu$ in the average.
We obtain \eqref{3.30} with $\alpha=\alpha(c_\nu,\eta_\nu)$,
$\mathcal A(q)=-2\mathcal T(q)$, and
$\mathcal B(q)=-\kappa_m(q)/2$. Substituting
\[
 \frac{3+2\alpha(c,\eta)-\alpha(c,\eta)^2}{8}
 =\frac{-c-\tfrac12+i(c+1)\eta+\eta^2/2}{c^2}
\]
gives \eqref{3.44}.
\end{proof}

\begin{proposition}\label{pr3.8}
Suppose that the solution sets of all lower types obtained by pair
replacements are compact and that \eqref{3.41} holds at their
mixed centers. Then \eqref{3.40} is compactly homotopic to the
system obtained by setting $(c,\eta)=(0,0)$ in \eqref{3.42}.
In particular, the two systems have the same degree.
\end{proposition}

\begin{proof}
By compactness and \eqref{3.41}, in coordinates with $\phi(q)=0$
there are $\gamma>0$ and $T_*\geq0$ such that
\[
 \kappa_m(q)-4|\mathcal T(q)|\geq\gamma,\qquad
 |\mathcal T(q)|\leq T_*
\]
at every mixed center of the lower solutions. Set
\[
 L(\zeta)=-2\mathcal T(q)\zeta-\tfrac12\kappa_m(q)|\zeta|,
 \qquad \lambda>\frac{4T_*}{\gamma},\qquad
 \beta=\frac{\lambda\gamma}{2}-2T_*>0.
\]
Then
\[
 \operatorname{Im}L(\zeta)+\lambda\operatorname{Re}L(\zeta)
 \leq-\beta|\zeta|.
\]
Take $(c,\eta)=(-1/2+s,\tau s)$ with $\tau=2\lambda+1$ and $s\geq0$.
Along a colliding sequence with $s_\nu\to0$, \eqref{3.43} and
Lemma~\ref{le3.7} give
\[
 (u,v,n)=(2,0,1),\qquad
 L(\zeta_\nu)=(-4+2i\tau)s_\nu+O(s_\nu^2)+o(|\zeta_\nu|).
\]
Applying $\operatorname{Im}+\lambda\operatorname{Re}$ gives, for large $\nu$,
\[
 2s_\nu+\beta|\zeta_\nu|
 \leq Cs_\nu^2+o(|\zeta_\nu|)
 \leq s_\nu+\tfrac{\beta}{2}|\zeta_\nu|
 <2s_\nu+\beta|\zeta_\nu|.
\]
Thus, for some $s_0\in(0,1/2)$, the centers are uniformly separated
along the first segment of the following path.
\[
 (-1/2,0)\longrightarrow(-1/2+s_0,\tau s_0)
 \longrightarrow(-1/2+s_0,0)\longrightarrow(0,0).
\]
On the last two segments, $c\geq-1/2+s_0>-1/2$, so
\eqref{3.43} gives uniform separation. At any candidate scalar
concentration point, separation yields
\[
 u\leq1,\qquad n\in\mathbb Z_{\geq0},\qquad
 n<\frac{u+1}{2}\leq1 \quad\Longrightarrow\quad n=0.
\]
As in Corollary~\ref{co3.4}, this gives scalar compactness.
Before introducing the tubular extension, we use the product of the
$X_s$ norm topology and the topology of the center manifold.
In this homotopy, $Z=(P,Q)$ denotes the center configuration.
For $Z\in\mathcal C_{a,b}(M)$ of Section~\ref{sec:shadow-setting},
define the scalar compact term
\[
 K_P(w)=2\sigma_bG\left(\frac{H_Pe^w}{\int_MH_Pe^w\dV}-1\right),
 \qquad \mathcal S_{\rm sc}(w,Z)=w-K_P(w).
\]
For separated centers, the weight is jointly smooth in $(x,P)$.
Indeed, each factor $e^{-4\pi G(x,p_i)}$ equals
$|z(x)-z(p_i)|^2e^{-4\pi R(x,p_i)}$ near its center. The derivatives
of the weight of each fixed order are bounded on compact sets of
separated configurations. The denominator satisfies
\begin{equation}\tag{3.45a}\label{eq:shadow-weight-lower-bound}
 \int_M H_Pe^w\dV
 \ge \exp\!\left(\int_M(\log H_P+w)\dV\right)
 = \exp\!\left(\int_M\log h_2\dV\right)>0 .
\end{equation}
The logarithms of the quadratic zeros are integrable. The equality
uses area one, the zero mean of $w$, and the zero mean of each
Green function. Jensen's inequality therefore applies even when
$P$ is nonempty. The Sobolev algebra and smooth composition in
$H^s$ give a uniform $H^s$ bound for the normalized density on
bounded $X_s$ sets and compact separated center sets. The Green
inverse gives a uniform $H^{s+2}$ bound there for $K_P(w)$.
The same multiplication and composition estimates give
continuity in $(w,P)$ in $H^s$. In particular the scalar compact
term is jointly compact on bounded products with compact center
closure.

Let $X_\tau(w,Z)\in T_Z\mathcal C_{a,b}(M)$ denote the center field
along the piecewise linear path above, parameterized by
$\tau\in[0,1]$. The Green function gradients in this field are smooth on separated
configurations. Since $X_s\hookrightarrow C^2(M)$, the evaluations
$\nabla w(p_i)$ are defined on the entire domain, and in fixed local
frames the map $(w,p_i)\mapsto\nabla w(p_i)$ is continuous. Indeed,
convergence in $X_s$ gives uniform convergence of the gradients,
and the limiting gradient is continuous in $p_i$. The field is
also jointly continuous in $\tau$ because the coefficients $c,\eta$
in \eqref{3.42} are continuous. The field is bounded on bounded scalar
sets and compact separated center sets. On each segment the maps are
$C^1$ in the variables $(w,Z)$. The $C^2$ embedding gives the derivative
with respect to a moving center.

By the separation proved above, the center projections of the
zeros lie in a compact separated set. Their uniform $C^0$ scalar
bound and \eqref{eq:shadow-weight-lower-bound} give uniform
$W^{2,p}$ bounds for every fixed finite $p$, by the scalar equation.
Starting with $H^2$, the Sobolev algebra estimates and the scalar
equation then give uniform bounds through $H^{s+2}$. Rellich
compactness in $X_s$ and continuity show that the union of zeros
over the closed path is compact in these variables. Choose a
relatively compact open center domain $\mathcal O$ containing
their center projections, and an open ball $B_s\subset X_s$
strictly containing their scalar projections.

Embed $\mathcal C_{a,b}(M)$ smoothly into $\R^N$ by $\iota$. Choose
a uniform tubular radius $r_{\rm tub}>0$ over a neighborhood of
$\overline{\mathcal O}$ so that $x=\iota(Z)+n$ with $n$ normal
gives unique smooth tubular coordinates there. The tube
$\mathcal T=\{x:Z\in\mathcal O,\ |n|<r_{\rm tub}\}$ is open and
bounded in $\R^N$. On $\overline{B_s\times\mathcal T}$ use the map
into $X_s\times\R^N$
\begin{equation}\tag{3.45b}\label{eq:shadow-tubular-homotopy}
 \mathfrak S_\tau(w,x)
 =\bigl(w-K_P(w),D\iota_Z X_\tau(w,Z)+n\bigr).
\end{equation}
Its compact term has scalar component $K_P(w)$ and
finite-dimensional component $\iota(Z)-D\iota_Z X_\tau(w,Z)$.
The latter is bounded and continuous, so this is a continuous
compact homotopy in the product norm. The tangent and normal
components of the second equation are orthogonal. Every zero
therefore has $n=0$ and $X_\tau=0$. Together with the scalar
equation, these identities put each zero strictly inside the scalar
ball and the center domain and away from the normal boundary. Thus the same bounded
open domain is admissible throughout the path.

At a nondegenerate zero, use the same tangent and normal
identifications of $\R^N$ in the domain and range. Derivatives of the
moving frames and of $D\iota_Z$ multiplying the center field or
normal variable vanish at the zero. The derivative is the
intrinsic coupled derivative with the normal identity block.
This proves agreement with the local orientation in
\eqref{eq:shadow-local-index} whenever its scalar block is
invertible, and fixes the normal factor in the product
calculation at the decoupled endpoint. Homotopy invariance of
Leray--Schauder degree now gives the asserted equality.
\end{proof}

\subsection{Nondegeneracy of shadow solutions and separation of center positions}

Let $f\geq0$ be smooth on a connected closed Riemannian surface $M$,
with $\sigma=\int_M f\,\dV>0$. Define
\[
 \avg{\xi}{f}=\frac1{\sigma}\int_M f\xi\,\dV,\qquad
 L_f=-\Delta_g+\mathcal K_f,\qquad
 \mathcal K_f\xi=-2f\bigl(\xi-\avg{\xi}{f}\bigr).
\]
Since $f$ is smooth and $\int_M\mathcal K_f\xi\,\dV=0$,
$\mathcal K_f$ is bounded on $H^r_0(M)$ for $r\geq0$. By the Rellich
embedding, it is compact from $H^{r+2}_0(M)$ to $H^r_0(M)$.
Thus $L_f$ is Fredholm of index zero on these spaces. Its realization
on $L^2_0(M)$ with domain $H^2_0(M)$ is self-adjoint, with bilinear form
\[
 Q_f(\xi,\zeta)=\int_M\nabla\xi\cdot\nabla\zeta\,\dV
 -2\int_M f\bigl(\xi-\avg{\xi}{f}\bigr)
             \bigl(\zeta-\avg{\zeta}{f}\bigr)\,\dV.
\]

\begin{lemma}\label{le3.9}
Let $U$ be a connected nonempty open set on which $f>0$. There exists
$v\in C_c^\infty(U)$ with $\int_M v\,\dV=0$ such that $L_{f+tv}$ is
invertible for all sufficiently small nonzero real $t$. For such $t$,
$f+tv\geq0$ and $\int_M(f+tv)\,\dV=\sigma$.
\end{lemma}

\begin{proof}
If $\ker L_f=\{0\}$, take $v=0$. Otherwise, let
$\xi_1,\ldots,\xi_d$ be an $L^2$-orthonormal basis of $\ker L_f$ and put
$\widehat\xi_i=\xi_i-\avg{\xi_i}{f}$. For $v\in C_c^\infty(U)$ with
$\int_U v\,\dV=0$, set
\[
 S(v)=\left(\int_U v\widehat\xi_i\widehat\xi_j\,\dV\right)_{i,j=1}^d.
\]
Suppose that all these matrices are singular. Choose a matrix
$S_0$ of maximal rank $r<d$ and a unit vector $z\in\ker S_0$.
By linearity in $v$, every $S_0+tS(v)$ has rank at most $r$. Its
compression to $\operatorname{ran}S_0\oplus\R z$ therefore has zero
determinant. The coefficient of $t$ in that determinant is
\[
 \det\bigl(S_0|_{\operatorname{ran}S_0}\bigr)\,z^TS(v)z=0,
\]
where the determinant is taken to be $1$ when $r=0$. Hence
\[
 \phi=\sum_i z_i\widehat\xi_i,\qquad
 \int_U v\phi^2\,\dV=0
 \quad\text{for all }v\in C_c^\infty(U)\text{ with }\int_Uv\,\dV=0.
\]
Fix $\chi\in C_c^\infty(U)$ with $\int_U\chi\,\dV=1$. Testing with
$v=\eta-(\int_U\eta\,\dV)\chi$ for arbitrary $\eta\in C_c^\infty(U)$
shows that $\phi^2$ is constant on $U$. By connectedness, $\phi$ is
constant there. Since $-\Delta_g\phi=2f\phi$ and $f>0$ on $U$,
we have $\phi=0$ on $U$. Unique continuation \cite{Aro57} gives
$\phi=0$ on $M$. Thus $\sum_i z_i\xi_i$ is constant and has zero
mean, contradicting $|z|=1$.

Choose $v$ with $S(v)$ invertible. Positivity of $f$ on the compact
support of $v$ gives $f+tv\geq0$ for small $|t|$, and
$\int_M(f+tv)\,\dV=\sigma$. For $\widehat\xi=\xi-\avg{\xi}{f}$
and $\widehat\zeta=\zeta-\avg{\zeta}{f}$, expansion of the bilinear form
gives
\[
 Q_{f+tv}(\xi,\zeta)-Q_f(\xi,\zeta)
 =-2t\int_M v\widehat\xi\widehat\zeta\,\dV
 +\frac{2t^2}{\sigma}
   \left(\int_Mv\xi\,\dV\right)\left(\int_Mv\zeta\,\dV\right).
\]
If $L_{f+tv}\xi=0$, write $\xi=\sum_i a_i\xi_i+\zeta$ with
$\zeta\perp\ker L_f$ in $L^2$. The restriction of $L_f$ to this
orthogonal complement is invertible onto its closed range with a bounded inverse.
Since $L_f\zeta=-(L_{f+tv}-L_f)\xi$, we have
\[
 \|\zeta\|_{H^2}\leq C\|L_f\zeta\|_{L^2}
 \leq C|t|\bigl(|a|+\|\zeta\|_{H^2}\bigr),
 \qquad
 \|\zeta\|_{H^2}\leq C|t|\,|a|.
\]
Testing the equation against each $\xi_j$ and using the form expansion,
we obtain
\[
 -2tS(v)a+O(t^2|a|)=0.
\]
Since $S(v)$ is invertible, $a=0$ and then $\zeta=0$ for sufficiently
small $t\ne0$. The index-zero property gives invertibility.
\end{proof}

\begin{corollary}\label{co3.10}
Let $(w,P,Q)$ be a solution of \eqref{3.40} with weights $(h_1,h_2)$.
Write $H_P=h_2W_P$ for the weight in \eqref{2.5} and set
\[
 \sigma=\sigma_b,\qquad
 f=\frac{\sigma H_Pe^w}{\int_M H_Pe^w\,\dV}.
\]
Choose a disk $U\Subset M\setminus(P\cup Q)$ and take $v$ as in
Lemma~\ref{le3.9}. Let $z$ be the mean-zero solution of
$-\Delta_g z=2v$, and set
\[
 w_t=w+tz,\qquad
 h_{1,t}=h_1e^{tz/2},\qquad
 h_{2,t}=h_2(1+tv/f)e^{-tz},
\]
where $v/f$ is extended by zero outside $U$. For all sufficiently small
$t\ne0$, $(w_t,P,Q)$ solves \eqref{3.40} with weights $(h_{1,t},h_{2,t})$,
and its scalar linearization at fixed centers is invertible. Moreover,
\[
 2\log h_{1,t}+\log h_{2,t}
 =2\log h_1+\log h_2
 \qquad\text{near }P\cup Q.
\]
\end{corollary}

\begin{proof}
Put $H_{P,t}=h_{2,t}W_P$ and $I=\int_M H_Pe^w\,\dV$. The definitions give
\[
 H_{P,t}e^{w_t}=H_Pe^w(1+tv/f)=H_Pe^w+\frac{tI}{\sigma}v.
\]
Since $\int_Mv\,\dV=0$, it follows that
\[
 \int_M H_{P,t}e^{w_t}\,\dV=I,\qquad
 \frac{\sigma H_{P,t}e^{w_t}}{I}=f+tv,\qquad
 -\Delta_gw_t=2(f+tv-\sigma).
\]
Also,
\[
\begin{aligned}
 \log h_{1,t}-\tfrac12w_t
 &=\log h_1-\tfrac12w,\\
 2\log h_{1,t}+\log h_{2,t}
 &=2\log h_1+\log h_2+\log(1+tv/f).
\end{aligned}
\]
The first identity gives the ordinary-center equations in \eqref{3.40}.
The last term in the second identity vanishes near $P\cup Q$, which
gives the mixed-center equations and the asserted identity there.
The scalar linearization is $L_{f+tv}$, so Lemma~\ref{le3.9} applies.
\end{proof}

Throughout this subsection, we work in an open set of smooth positive
weights for which the solution sets of \eqref{3.40} are uniformly
compact for a fixed finite collection of types $(a,b)$ with $a+2b=m$.
The joint transversality calculation uses independent weight jets
at the centers and variations supported away from them, as follows.

\begin{lemma}\label{le3.11}
There is a $C^\infty$-dense subset of the above open set such that,
for every pair of weights in this subset and every type $(a,b)$
under consideration, the shadow system \eqref{3.40} has finitely
many solutions, all of which are nondegenerate.
\end{lemma}

\begin{proof}
Fix a solution $Z=(w,P,Q)$ of \eqref{3.40} with weights
$h=(h_1,h_2)$ and density $f$ as in Corollary~\ref{co3.10}.
Elliptic regularity gives $w,f\in C^\infty(M)$. In ordered center
charts, for $s\geq6$,
\[
\begin{gathered}
 D_ZF:H^s_0(M)\times\R^{2(a+b)}
       \longrightarrow H^{s-2}_0(M)\times\R^{2(a+b)},\\
 D_ZF=\bigl((-\Delta_g)\oplus I\bigr)+\mathcal K,\qquad
 \operatorname{ind}D_ZF=0.
\end{gathered}
\]
Indeed, $D_wF_0=L_f$, so the scalar diagonal block of $\mathcal K$
is the compact operator $\mathcal K_f$. All other blocks of
$\mathcal K$ have finite rank. The equation map is $C^1$ in these
charts. Joint smoothness of the weight and continuity of the
moving point evaluations follow as in the proof of
Proposition~\ref{pr3.8}. Composing the first row with $G$ gives
the derivative of the intrinsic fixed-point and center equations
used in \eqref{eq:shadow-tubular-homotopy}. Thus one linearization is
invertible exactly when the other is.

Let $D_{\eta_i}$ denote differentiation along
$h_i\mapsto h_i e^{t\eta_i}$. In the normalized density,
differentiating the denominator contributes $\avg{\eta_2}{f}$.
Thus \eqref{3.40} gives
\begin{equation}\label{3.45}
\begin{aligned}
 D_{\eta_1}F_0&=0,
 &D_{\eta_2}F_0&=-2f\bigl(\eta_2-\avg{\eta_2}{f}\bigr),\\
 D_{\eta_1}F_i^A&=\nabla\eta_1(p_i),
 &D_{\eta_2}F_i^A&=0,\\
 D_{\eta_1}F_j^B&=2\nabla\eta_1(q_j),
 &D_{\eta_2}F_j^B&=\nabla\eta_2(q_j).
\end{aligned}
\end{equation}
The range of $D_{Z,h}F$ is closed because it contains the closed
finite-codimensional range of $D_ZF$. Let
$(\varphi,\zeta_A,\zeta_B)\in\ker\bigl((D_{Z,h}F)^*\bigr)$, where
$^*$ denotes the dual map and $\varphi$ is represented by a
zero-mean distribution. By \eqref{3.45},
\[
 \sum_i\zeta_{A,i}\cdot\nabla\eta_1(p_i)
 +2\sum_j\zeta_{B,j}\cdot\nabla\eta_1(q_j)=0
 \qquad\text{for every }\eta_1\in C^\infty(M).
\]
Using smooth functions supported in disjoint disks about the centers,
we can prescribe these gradients independently. Hence
$\zeta_A=\zeta_B=0$. The $w$-variations give $L_f\varphi=0$ in
$\mathcal D'(M)$ by formal self-adjointness and $L_f1=0$.
The $h_2$-variations then give
\[
 \left\langle
 f\left(\varphi-\frac{\langle\varphi,f\rangle}{\sigma_b}\right),
 \eta_2\right\rangle=0
 \qquad\text{for every }\eta_2\in C^\infty(M).
\]
Substituting this identity into $L_f\varphi=0$ yields
\[
 -\Delta_g\varphi=0,\qquad \langle\varphi,1\rangle=0,
 \qquad\text{hence }\varphi=0.
\]
Thus $\ker\bigl((D_{Z,h}F)^*\bigr)=\{0\}$, and the closed-range
property gives surjectivity.

At each solution, choose finitely many smooth log-weight directions
whose images span $\operatorname{coker}D_ZF$. With these directions
included, the joint derivative remains surjective in a neighborhood
of the solution. By compactness, a single finite family
$\{\eta^{(\nu)}\}_{\nu=1}^N$ suffices for all the solutions of the types
under consideration. Set
\[
 h_i(\lambda)=h_i\exp\!\left(
   \sum_{\nu=1}^N\lambda_\nu\eta_i^{(\nu)}\right),\qquad
 \mathcal S(Z,\lambda)=F(Z;h(\lambda)).
\]
If zeros for a sequence $\lambda_n\to0$ lay outside the neighborhoods
where the joint derivative is surjective, uniform compactness would
give a subsequence converging to an original solution in one of those
neighborhoods. Thus $D_{Z,\lambda}\mathcal S$ is surjective at every
zero for sufficiently small $\lambda$. Since
\[
 \operatorname{ind}D_{Z,\lambda}\mathcal S
 =\operatorname{ind}D_Z\mathcal S+N=N,
\]
the zero set $\mathfrak Z=\mathcal S^{-1}(0)$ for each type is a
$C^1$ manifold of dimension $N$.

Sard's theorem \cite{Sard42} gives arbitrarily small $\lambda$ that are
regular values of all the projections
$\pi:\mathfrak Z\to\R^N$, $\pi(Z,\lambda)=\lambda$, for the finitely
many types. At a zero over such a value, $D\pi$ is an isomorphism,
since $\dim\mathfrak Z=N$. Therefore
\[
 \ker D_Z\mathcal S\times\{0\}=\ker D\pi=\{0\}.
\]
Since $\operatorname{ind}D_Z\mathcal S=0$, the operator
$D_Z\mathcal S$ is invertible. Each solution set is compact and
discrete, hence finite.
\end{proof}

\begin{lemma}\label{le3.12}
For every pair of weights in a $C^\infty$-dense subset of the above
open set, the solutions of \eqref{3.40} for the types under
consideration satisfy the following properties.
\begin{enumerate}
\item Every solution is nondegenerate and its scalar linearization
with the centers fixed is invertible.
\item For any two solutions $(w,P,Q)$ and
$(\widetilde w,\widetilde P,\widetilde Q)$,
$P\cap\widetilde Q=\varnothing$.
\item If $(w,P,Q)$ and $(\widetilde w,\widetilde P,\widetilde Q)$
have different types, then $Q\cap\widetilde Q=\varnothing$.
\end{enumerate}
Solutions of the same type may share mixed center positions.
\end{lemma}

\begin{proof}
\noindent\emph{Step 1. Invertibility of the linearizations.}
By Lemma~\ref{le3.11}, an arbitrarily small weight perturbation gives
finitely many solutions $Z^1,\ldots,Z^N$ with $D_ZF(Z^\ell)$ invertible.
As in the proof of Lemma~\ref{le3.11},
\[
 D_ZF=\bigl((-\Delta_g)\oplus I\bigr)+\mathcal K,\qquad
 L_f=-\Delta_g+\mathcal K_f,\qquad
 \operatorname{ind}D_ZF=\operatorname{ind}L_f=0,
\]
where $\mathcal K_f$ is the scalar diagonal block of $\mathcal K$.
The implicit function theorem gives a unique local continuation of
each $Z^\ell$. By the compactness argument in the proof of
Lemma~\ref{le3.11}, every solution for nearby weights is one of these
continuations.

Choose a disk whose closure avoids the centers of all these solutions
for weights sufficiently close to the current weights. Apply
Lemma~\ref{le3.9} and Corollary~\ref{co3.10} successively to the
continued solutions. At each step, take the perturbation small enough
to preserve the invertibility of every full linearization and of
the scalar linearizations already treated. Thus, at every solution,
\[
 \ker D_ZF=\{0\},\qquad \ker L_f=\{0\}.
\]

\medskip
\noindent\emph{Step 2. The reduced center equations.}
For a solution of type $(a,b)$, set
\[
\begin{aligned}
 E_c(P,Q)
 &=\sum_i g_1(p_i)+\sum_j V(q_j)
   +8\pi\sum_{i<l}G(p_i,p_l)
   +12\pi\sum_{i,j}G(p_i,q_j)
   +24\pi\sum_{j<l}G(q_j,q_l),\\
 E_s(w,P)
 &=\frac12\int_M|\nabla w|^2\,\dV
   -2\sigma_b\log\int_M H_Pe^w\,\dV,
 \qquad E=E_c-\frac{1}{8\pi}E_s.
\end{aligned}
\]
For $(G*u)(x)=\int_M G(x,y)u(y)\,\dif V_y$, the formula for $H_P$
and symmetry of $G$ give
\[
 \nabla_{p_i}E_s
 =8\pi\int_M f(x)\nabla_{p_i}G(x,p_i)\,\dif V_x
 =8\pi\nabla(G*f)(p_i).
\]
The Green normalization gives $G*F_0=w-2G*f$. It follows that,
for every mean-zero $\varphi$,
\[
 D_wE[\varphi]=-\frac1{8\pi}\int_M F_0\varphi\,\dV,\qquad
 \nabla_{p_i}E=F_i^A+\tfrac12\nabla(G*F_0)(p_i),\qquad
 \nabla_{q_j}E=F_j^B.
\]
By Step~1 and the implicit function theorem, we solve
$F_0(w(P),P)=0$ locally. Along this graph, $D_wE=0$ and $G*F_0=0$,
so the center equations are the critical point equations of
$E(w(P),P,Q)$. Their linearization is the symmetric Hessian
\[
 S=D^2_{(P,Q)}E(w(P),P,Q),\qquad
 S_{q_j,p_i}=12\pi\nabla_{q_j}\nabla_{p_i}G(q_j,p_i).
\]
Differentiating the scalar and reduced center equations gives
\[
 \dot w=-L_f^{-1}D_PF_0[\dot P],\qquad
 D_ZF\begin{pmatrix}\dot w\\\dot P\\\dot Q\end{pmatrix}
 =\begin{pmatrix}0\\ S(\dot P,\dot Q)\end{pmatrix}.
\]
Thus $S$ is the Schur complement of $L_f$ in $D_ZF$, and is invertible
by Step~1. In orthonormal frames at the centers, a variation $\xi$
of the gradients of $\log h_1$ gives, by \eqref{3.45},
\[
 D_\xi(P,Q)=-S^{-1}D,\qquad
 D=\operatorname{diag}(I_{2a},2I_{2b}).
\]

\medskip
\noindent\emph{Step 3. Invertibility of the required blocks of $S^{-1}$.}
For fixed $q$, put $\gamma(p)=(\partial_{q_l}G(q,p))_{l=1}^2$
in local coordinates at $q$. Since $\Delta_{g,p}G(q,p)=1$ off the
diagonal, the components of $\gamma$ are harmonic on $M\setminus\{q\}$.
They are therefore real analytic in isothermal coordinates. In a
common isothermal chart near $q$, the logarithmic term of $G$ gives
\[
 D_p\gamma(p)=\frac{I_2-2e\otimes e}{2\pi|p-q|^2}+O(1),
 \qquad e=\frac{p-q}{|p-q|},
 \qquad
 \det D_p\gamma(p)=-\frac1{4\pi^2|p-q|^4}+O(|p-q|^{-2}).
\]
The zero set of the determinant is independent of the coordinate
chart. By analytic continuation on the connected punctured surface,
this zero set has empty interior in $M\setminus\{q\}$.

Choose the gradient variations in \eqref{3.45} independently at all
distinct center positions of the finite solution set. Since
$D_\xi(P,Q)=-S^{-1}D$, the parameter map to each center configuration
is a submersion. For a mixed--ordinary pair $(q,p)$ in a solution,
the condition $\det\nabla_q\nabla_pG(q,p)\ne0$ is open and dense in
the configuration space, by the preceding argument with $q$ fixed.
Pulling back these finitely many conditions under the submersions gives
a common open dense set of nearby parameters. We may thus take
every required block $S_{q,p}$ to be invertible.

Choose $\psi\in C^\infty(M)$ with $\nabla\psi=0$ and
$\nabla^2\psi=g$ at every current center position, using local
quadratic functions with disjoint supports. Under the perturbation
$\log h_1\mapsto\log h_1+t\psi$, each current $Z^\ell$ still solves
the equations by the vanishing of $\nabla\psi$ at its centers.
The reduced Hessian satisfies $S_t=S+tD$. As $t\to+\infty$,
\[
 \begin{gathered}
 S_t^{-1}=t^{-1}D^{-1}-t^{-2}D^{-1}SD^{-1}+O(t^{-3}),\\
 (S_t^{-1})_{q,q}=\frac1{2t}I_2+O(t^{-2}),\qquad
 (S_t^{-1})_{q,p}=-\frac1{2t^2}S_{q,p}+O(t^{-3}).
 \end{gathered}
\]
The determinants of these blocks are rational functions of $t$ and
the expansions show that each is nonzero for large $t$. Avoiding
their finitely many zeros and poles, we choose a single arbitrarily
small $t$ such that, with $R$ denoting the inverse of the current
reduced Hessian,
\[
 \det R_{q,q}\ne0,\qquad \det R_{q,p}\ne0
\]
for every mixed center $q$ and every mixed--ordinary pair $(q,p)$
in the current solutions.

\medskip
\noindent\emph{Step 4. Separation of center positions.}
Let $X=(w^X,P^X,Q^X)$ and $Y=(w^Y,P^Y,Q^Y)$ be two solutions of
\eqref{3.40}. Suppose that $p_i^X=q_j^Y=x$. In normal coordinates
centered at $x$, vary $\log h_2$ by
\[
 \eta_2(z)=\chi(z/r)\,\xi\cdot z,\qquad
 \|\eta_2\|_{C^{0,\alpha}}\leq Cr^{1-\alpha}|\xi|,
 \qquad 0<\alpha<1,
\]
where $\chi\in C_c^\infty(B_1)$ is equal to one near zero. Take $r$
small enough that the support avoids the other centers of $X$ and
$Y$. At fixed centers, the scalar variation is
\[
 \dot w=-L_f^{-1}D_{\eta_2}F_0,\qquad
 \|\dot w\|_{C^{2,\alpha}}\leq C r^{1-\alpha}|\xi|.
\]
Here we used \eqref{3.45} and the bounded inverse
$L_f^{-1}:C^{0,\alpha}_0(M)\to C^{2,\alpha}_0(M)$, which follows
from scalar invertibility and elliptic regularity.

In the reduced center equations, the contribution through $\dot w$
is $O(r^{1-\alpha}|\xi|)$. By \eqref{3.45}, the direct contribution
is $\xi$ in the mixed-center equation at $q_j^Y$. Multiplying the
reduced variations by the negative inverse Hessians gives
\[
 D_\xi p_i^X=O(r^{1-\alpha}),\qquad
 D_\xi q_j^Y=-R^Y_{q_j,q_j}+O(r^{1-\alpha}).
\]
Thus, in the common chart about $x$,
\[
 D_\xi(p_i^X-q_j^Y)=R^Y_{q_j,q_j}+O(r^{1-\alpha}),\qquad
 \operatorname{rank}D_\xi(p_i^X-q_j^Y)=2
\]
for fixed sufficiently small $r$. Combine the directions for all
current ordinary--mixed coincidences in a single finite-dimensional
family. Each difference map has rank two at the initial parameter,
since its derivative already has rank two in the directions chosen
for that pair. The implicit function theorem gives a coincidence set
of codimension two for each pair near that parameter. Pairs that
are already separated remain separated under small perturbations.
Choose parameters outside the finite union of coincidence sets.
This separates all ordinary--mixed pairs and preserves the
invertibility established in Steps~1 and~3.

Now suppose that $q_i^X=q_j^Y$ for solutions of different types.
Since $a+2b=m$ is fixed, we may assume $a_X>a_Y$ and choose
$p\in P^X\setminus P^Y$. The ordinary--mixed separation gives
$p\notin Q^Y$. A local gradient variation of $\log h_1$ at $p$,
supported away from the other centers of $X$ and $Y$, contributes
$\xi$ only to the ordinary-center equation at $p$ in $X$. Hence
\[
 D_\xi q_j^Y=0,\qquad
 D_\xi(q_i^X-q_j^Y)=-R^X_{q_i,p},\qquad
 \operatorname{rank}D_\xi(q_i^X-q_j^Y)=2.
\]
The same finite-dimensional argument removes all mixed-center
coincidences between different types. At each stage the smooth
perturbation directions are fixed before their coefficients are
chosen. Hence the total perturbation can be made arbitrarily small
in the $C^\infty$ topology.
\end{proof}

\subsection{Quantization for varying coefficients}
Let $\vartheta\in C_c^\infty(\R^2)$ and set $b_r(z)=L|z|^2\vartheta(z/r)$
with fixed $L$, as in (\ref{def:br}). Direct differentiation gives
\begin{equation}\label{3.46}
 \begin{split}
 \|b_r\|_\infty&=O(|L|r^2),\qquad \|\nabla_z b_r\|_\infty=O(|L|r),\\
 \|D_z^2b_r\|_\infty&=O(|L|),\qquad \|D_z^3b_r\|_\infty=O(|L|/r).
 \end{split}
\end{equation}
Using the local reduction in \eqref{2.12} for the scalar equation in \eqref{3.40}, we obtain
$h_\nu=A_\nu h_2e^{-b_{r_\nu}}$, where $A_\nu=\sigma_b e^{\phi+\psi}$ and $\Delta_0\psi=(2\sigma_b-4\pi a)e^\phi$.
On a smaller coordinate disk, \eqref{3.46} gives a uniform positive lower bound and a uniform $C^2$ bound for $h_\nu$ as $r_\nu\downarrow0$. Uniform positivity and a uniform $C^2$ bound
for $h_\nu$ give
$$
 \|\log h_\nu(q_\nu+r_\nu\,\cdot)-\log h_\nu(q_\nu)\|_{C^1(K)}
       =O(r_\nu),\qquad
 \|D_y^2\log h_\nu(q_\nu+r_\nu\,\cdot)\|_\infty=O(r_\nu^2).
$$
\begin{proposition}\label{pr3.13}
Let $\beta_j\in\mathbb N$ be fixed and let the distinct points
$p_{j,\nu}\to0$ lie in a disk $D$. Suppose that $0$ is the only
limiting blow-up point and that
\[
\begin{gathered}
 \Delta\widehat u_\nu+2h_\nu e^{\widehat u_\nu}
 =4\pi\sum_{j=1}^J\beta_j\delta_{p_{j,\nu}}
 \qquad\text{in }D,\\
 h_\nu\geq c>0,\qquad \|h_\nu\|_{W^{1,\infty}(D)}\leq C,\qquad
 \int_Dh_\nu e^{\widehat u_\nu}\dif A_x
 +\osc_{\partial D}\widehat u_\nu\leq C,
 \qquad \widehat u_\nu\text{ is continuous near }\partial D.
\end{gathered}
\]
Every positive local mass at $0$ of the sequence
$h_\nu e^{\widehat u_\nu}\dif A_x$ belongs to $4\pi\mathbb N$.
\end{proposition}

\begin{proof}
We apply the local theorem to the equation in the sense of distributions.
For fixed $\nu$, let
$\xi_\nu=\widehat u_\nu-2\sum_j\beta_j\log|x-p_{j,\nu}|$
with the sum empty for $J=0$. Since $\Delta\log|x-p|=2\pi\delta_p$,
its distributional equation is
\[
 -\Delta\xi_\nu=W_\nu e^{\xi_\nu}=2h_\nu e^{\widehat u_\nu},
 \qquad W_\nu=2h_\nu\prod_j|x-p_{j,\nu}|^{2\beta_j}.
\]
For each fixed $\nu$, the right-hand side is integrable and has no point masses.
On a small disk, write $\xi_\nu$ as its local logarithmic potential
plus a smooth harmonic function. For any fixed $p>2$, shrink that
disk so that the mass of the source is less than $4\pi/p$.
Jensen's inequality for the logarithmic potential then gives
$e^{\xi_\nu}\in L^p$ on smaller disks. Since $W_\nu$ is locally
Lipschitz and bounded, local $W^{2,p}$ and Schauder estimates give
$\xi_\nu\in C^2_{\rm loc}(D)$. Removing the pole logarithms changes the
boundary oscillation on a fixed disk by at most a bounded error.
As in the estimate of the regular weight term following \eqref{2.15},
the bounds $|x-q_\nu|\leq r$ on $B_r(q_\nu)\Subset D$ and
$\|\nabla\log h_\nu\|_{L^\infty(D)}\leq C/c$ give
\[
 \left|\int_{B_r(q_\nu)}(x-q_\nu)\cdot\nabla\log h_\nu\,
 h_\nu e^{\widehat u_\nu}\dif A_x\right|
 \leq \frac{C}{c}\,r\int_{B_r(q_\nu)}
 h_\nu e^{\widehat u_\nu}\dif A_x.
\]
The same gradient bound gives, under the rescaling $x=q_\nu+r_\nu y$,
\[
 \widetilde h_\nu(y)=\frac{h_\nu(q_\nu+r_\nu y)}{h_\nu(q_\nu)},
 \qquad \|\nabla_y\log\widetilde h_\nu\|_{L^\infty(K)}\leq Cr_\nu.
\]
For $J=0$, the local result of \cite[Proposition~2.1 and
Theorem~A(i)]{Tar23} applies to
$-\Delta\widehat u_\nu=2h_\nu e^{\widehat u_\nu}$.
The coefficient $2h_\nu$ is uniformly positive and Lipschitz. Thus it
has a uniformly bounded inverse moment and converges locally uniformly
after passage to a subsequence. The mass and boundary oscillation
are bounded by hypothesis. The local mass of
$2h_\nu e^{\widehat u_\nu}\dif A$ is $8\pi$, which gives the conclusion.
For $J\geq1$, use the equation for $\xi_\nu$ just established.
The factors $|x-p_{j,\nu}|^{2\beta_j}$ and the factor $2h_\nu$,
which is uniformly bounded and Lipschitz, give a uniformly bounded
and uniformly Lipschitz $W_\nu$ on $D$.
The logarithmic terms have uniformly bounded oscillation on $\partial D$
and are uniformly bounded on compact subsets of $D\setminus\{0\}$.
Hence the boundary oscillation hypothesis and the assumption that $0$
is the only limiting blow-up point give
\[
 \osc_{\partial D}\xi_\nu\leq C,\qquad
 \sup_K\xi_\nu\leq C_K
 \quad(K\Subset D\setminus\{0\}).
\]
A positive local mass at $0$ implies
$\sup_{B_r(0)}\xi_\nu\to+\infty$ for every $B_r(0)\Subset D$.
For $B=\sum_j\beta_j$ and $0<\epsilon<1/B$, H\"older's inequality gives
\[
 \int_D W_\nu^{-\epsilon}\dif A_x
 \leq C\prod_j
 \left(\int_D|x-p_{j,\nu}|^{-2\epsilon B}\dif A_x\right)^{\beta_j/B}
 \leq C_\epsilon.
\]
The full $W^{1,\infty}$ bound on $h_\nu$ gives local uniform
convergence of $h_\nu$ and hence of $W_\nu$ along a subsequence.
The mass, oscillation, and inverse moment estimates above give
the hypotheses of \cite[Proposition~2.1]{Tar23} for $\xi_\nu$.
A positive local mass places the sequence in the blow-up alternative of that result.
For $J\geq2$, the distinct collapsing zeros and their positive
integer orders meet \cite[Theorem~1 and \textup{(2.12)--(2.14)}]{Tar23}.
For $J=1$, \cite[Theorem~A(ii)]{Tar23} applies with
$\alpha=\beta_1$. In both cases the source masses satisfy
\[
 2\lim_{r\downarrow0}\lim_{\nu\to\infty}
 \int_{B_r(0)}h_\nu e^{\widehat u_\nu}\dif A_x
 =
 \lim_{r\downarrow0}\lim_{\nu\to\infty}
 \int_{B_r(0)}W_\nu e^{\xi_\nu}\dif A_x
 \in8\pi\mathbb N.
\]
\end{proof}

\subsection{Choice of weights and the shadow degree}\label{subsection:ProofofTheorem2.6}
We first choose a single pair of smooth positive weights for which
the generic properties established in Lemma~\ref{le3.12} hold
simultaneously for all admissible types and the hypotheses of
Proposition~\ref{pr3.8} are satisfied. The resulting a priori estimates
make the homotopy compact, so we can compute the shadow degrees
at its endpoint.
\begin{proof}[Proof of Theorem~\ref{th3.14}]
Write $h=(h_1,h_2)$ and let $\mathcal Z_{a,b}(h)$ denote the
solution set of \eqref{3.40}. We proceed by decreasing induction on
$b$, starting at $b_{\max}=\min\{k,\lfloor m/2\rfloor\}$.

At the step for $b$, assume that a single pair of weights satisfies
\textup{(i)}--\textup{(iii)} for all admissible types with second
index greater than $b$. We establish the three conclusions for
$(a,b)=(m-2b,b)$ while preserving them for those types.

\begin{enumerate}[label=\textup{(\roman*)},leftmargin=*]
\item \emph{Compactness.}
Consider solutions of a fixed type for weights $h^\nu\to h$ in
$C^1(M)$. On each disk isolating a possible concentration point,
the coefficient in \eqref{2.12} is uniformly positive and bounded
in $W^{1,\infty}$.

On an annulus containing the fixed isolating circle, we use the
Green representation and Jensen's inequality as in the proof of
Theorem~\ref{th2.2}. On that annulus, the mass on small disks is
uniformly small, the weights have uniform upper bounds, and the
denominators have uniform positive lower bounds. The estimate
gives a uniform bound for $w_\nu$ there, with the lower bound
following from that of $G$. The logarithms of the poles have
bounded oscillation on the circle, and the normalization is
a spatial constant. Hence Proposition~\ref{pr3.13} applies
with bounded boundary oscillation and gives the same mass
quantization as in Lemma~\ref{lem:successive-pair-collisions}.

The normalization, integrability and Pohozaev calculations in that
lemma use only $C^1$ convergence of the weights. In particular, a
collision gives $\ell\geq1$ disjoint pairs of ordinary centers and
a limit
\[
 Z_*\in\mathcal Z_{a-2\ell,b+\ell}(h),
 \qquad \sigma_{b+\ell}=\sigma_b-4\pi\ell>0.
\]
Each pair carries scalar mass $4\pi$. If the centers stay separated,
the bound $0\leq n<(u+1)/2$ with $u\leq1$ gives $n=0$ at every
point. The scalar estimates then give uniform $C^{2,\beta}$ bounds
for every $\beta<1$.

Suppose a pair converges to $q$ and the weights converge in $C^3$
on a fixed disk about $q$. The other centers contribute fields
converging smoothly on a smaller disk. We apply the calculation
of Theorem~\ref{th3.1} to these local coefficients.

Translate to the midpoint and rescale by the half-separation. By
\cite[Proposition~3.1]{Tar23}, there is a unique regular blow-up point
carrying the local mass $4\pi$. The rescaled
center equation locates it at the origin as in
Lemma~\ref{lem:pair-midpoint}, giving \eqref{3.2}.

On the disk, the background coefficient in \eqref{3.9} is
\[
 \kappa_\nu
 =\Delta_0(\log h_2^\nu+\phi)+(2\sigma_b-4\pi a)e^\phi.
\]
The convergence of this coefficient in $C^1$ gives convergence in the
Cauchy transform calculation \eqref{3.13}--\eqref{3.18}. The difference in
\eqref{3.28} is expanded at the midpoint, so $C^1$ convergence
of the regular terms gives the same $o(|d_\nu|)$ remainder.

Thus \eqref{3.7}, with the coefficients identified as in
\eqref{3.41}, gives
\[
 |\kappa_m(q)|=4e^{-\phi(q)}|\mathcal T(Z_*,q)|.
\]
At $b=b_{\max}$, either $a\leq1$ or $\sigma_b<4\pi$, so the
collision classification already excludes every pair. At any
later step, $b+\ell>b$ and \eqref{3.41} holds at $q$ by induction.
This contradicts the displayed equality.

The same argument applies to every previously treated type and
to weights converging to $h$ in $C^{s+2}$, where $s\geq6$ is
fixed. Elliptic estimates then give uniform $H^{s+1}$ bounds
along each such sequence.

Consequently, for weights in a sufficiently small $C^{s+2}$ neighborhood of $h$,
the solutions of all types considered so far have uniformly separated centers
and uniform $H^{s+1}$ bounds. Otherwise, a sequence of weights
converging to $h$ would violate one of these two bounds.

\item \emph{Nondegeneracy and separation of center positions.}
Apply Lemma~\ref{le3.12} in this neighborhood to all types with
second index at least $b$. Their solution sets are now finite,
and both $D_ZF$ and $L_f$ are invertible at every solution.

By compactness and continuity of the coefficients in \eqref{3.41},
a sufficiently small smooth perturbation preserves the strict
inequalities for the previously treated types. We continue to
denote the resulting weights by $h$.

\item \emph{The inequality at the mixed centers.}
If the current type has no mixed-center positions, this step is
complete. Otherwise, let $x_1,\ldots,x_J$ be its distinct
mixed-center positions. By Lemma~\ref{le3.12}, choose disjoint
coordinate disks centered at $x_j$ whose closures avoid all
ordinary positions and all mixed positions of the earlier types.

Take coordinates $z_j$ with $z_j(x_j)=0$ and $\phi_j(x_j)=0$.
Choose $L_j>0$ so that
\[
 \kappa_m(x_j)+2L_j
 >4\max_{\substack{Z\in\mathcal Z_{a,b}(h)\\x_j\in Q(Z)}}
      |\mathcal T(Z,x_j)|+2.
\]
Let $\vartheta\in C_c^\infty(B_1)$ satisfy $\vartheta=1$ on
$B_{1/2}$. With each summand extended by zero outside its disk,
set
\[
 b_r=\sum_{j=1}^J L_j|z_j|^2\vartheta(z_j/r),\qquad
 h_{1,r}=h_1e^{b_r/2},\qquad h_{2,r}=h_2e^{-b_r}.
\]
For fixed $L_j$ and $0<\alpha<1$, \eqref{3.46} gives
\[
 \|b_r\|_{C^1}=O(r),\qquad
 \|b_r\|_{C^{0,\alpha}}=O(r^{2-\alpha}),\qquad
 2\log h_{1,r}+\log h_{2,r}=2\log h_1+\log h_2.
\]

For any sequence $r_\nu\downarrow0$, the collision classification
in \textup{(i)} applies to solutions for $h_{r_\nu}$ of the current
and earlier types. Every collision limit belongs to a strictly
earlier type for the weights $h$.

A collision point $q$ would therefore be an earlier mixed-center
position. A fixed disk about $q$ avoids all the supports of
$b_{r_\nu}$ for large $\nu$. The local equality in \textup{(i)}
then contradicts \eqref{3.41} at $q$.

Hence the centers stay separated and the scalar estimates in
\textup{(i)}, with $\alpha<\beta<1$, give precompactness in
$C^{2,\alpha}$. In particular, every perturbed solution lies
near an original solution for all sufficiently small $r$.

For each type, let $\mathcal S_r$ be the map in \eqref{3.40}
with weights $h_r$. In local ordered center coordinates, take
\[
 X=C^{2,\alpha}_0(M)\times\R^{2(a+b)},\qquad
 Y=C^{0,\alpha}_0(M)\times\R^{2(a+b)},
\]
Here $a,b$ are the indices of that type. For an element of the kernel
of the H\"older linearization, the scalar row is a Poisson
equation with smooth coefficients, smooth terms from the
center variations, and a finite-rank scalar term. Elliptic
bootstrapping lifts its scalar component to $X_s$.
Thus the H\"older linearization has trivial kernel by \textup{(ii)}.
It has index zero by the same compact-perturbation argument as
in Lemma~\ref{le3.11}, now using Schauder estimates. Both
linearizations in \textup{(ii)} have bounded inverses on
these H\"older spaces.

Choose small disjoint neighborhoods $\mathcal U_Z$ of the
original solutions, with their ordinary coordinates away from
the perturbation disks. Since $V$ is unchanged and
$(x,p)\mapsto e^{-4\pi G(x,p)}$ is smooth, we have
\[
 \|\mathcal S_r-\mathcal S_0\|_{C^1(\mathcal U_Z;Y)}
 \leq C\|b_r\|_{C^{0,\alpha}}
 =O(r^{2-\alpha}).
\]
The implicit function theorem gives a unique solution $Z_r$ in
each $\mathcal U_Z$, with
\[
 \|Z_r-Z\|_X=O(r^{2-\alpha}).
\]
By the precompactness just proved, these are all the perturbed
solutions. Both linearizations remain invertible by
operator-norm continuity.

If a mixed center $q_r$ continues from $x_j$, then
$|z_j(q_r)|=O(r^{2-\alpha})=o(r)$, so $b_r=L_j|z_j|^2$ near $q_r$.
In particular,
\[
 (b_r)_{z_jz_j}(q_r)=0,\qquad
 (b_r)_{z_j}(q_r)=O(r^{2-\alpha}),\qquad
 \Delta_g b_r(q_r)=4L_j+o(1).
\]
The definitions of $\kappa_m$ and $\mathcal T$ give
\[
 \kappa_{m,r}(q_r)=\kappa_m(x_j)+2L_j+o(1),
 \qquad
 \mathcal T_r(Z_r,q_r)=\mathcal T(Z,x_j)+o(1).
\]
The choice of $L_j$ gives \eqref{3.41} at every current mixed
center for small $r$. At earlier mixed centers the weights agree
with $h$ on fixed neighborhoods, so their inequalities persist
by $Z_r\to Z$. All solutions at each treated type lie on the finitely
many continuing branches. Their centers move by $O(r^{2-\alpha})$.
The minimum of the finitely many positive distances between centers
of distinct solutions, when present, stays positive for small $r$.
\end{enumerate}

Choose one sufficiently small $r$ at each step. The induction
gives a single pair of weights for all admissible types. The
bounds are uniform over the finite collection, and a common
positive lower bound in \textup{(iii)} follows by taking the
minimum over all mixed centers of all solutions. If there are
no mixed centers, \textup{(iii)} is vacuous.
\end{proof}

\begin{proof}[Proof of Theorem~\ref{th3.15}]
\textit{Step 1.} We compute the degree of the center equations.
By Proposition~\ref{pr3.8}, we may compute the shadow degree at
the endpoint
$$
 \mathcal S_0(w,Z)=\bigl(\mathcal S_{\rm sc}(w,Z),\nabla E_c(Z)\bigr),
 \qquad Z=(P,Q)\in\mathcal C_{a,b}(M).
$$
Since the interaction coefficients are positive, $E_c\to+\infty$
as any two centers collide. Suppose a sequence of critical points
has centers indexed by $J$ tending to the same point $q$, with
$|J|\geq2$. In an isothermal coordinate centered at $q$, write $x_i$
for these centers and $c_{ij}\in\{8\pi,12\pi,24\pi\}$ for their
interaction coefficients. Pairing the logarithmic terms in
$x_i\cdot\nabla_{x_i}E_c$ as in
Lemma~\ref{lem:successive-pair-collisions} gives
\[
 0=\sum_{i\in J}x_i\cdot\nabla_{x_i}E_c
   =-\frac1{2\pi}\sum_{\substack{i<j\\ i,j\in J}}c_{ij}+o(1).
\]
The terms from centers outside $J$ and the regular parts of the
Green functions are $o(1)$. The sum is strictly negative, so
$\operatorname{Crit}(E_c)\Subset\mathcal C_{a,b}(M)$.

Perturb $E_c$ to a Morse function by a small perturbation supported
in a compact neighborhood of its critical set. Along the homotopy,
the critical points stay in that neighborhood and the scalar
equation is uniformly compact there. For the perturbed $E_c$,
choose a regular value $R$ above all its critical values and put
$\Omega_R=E_c^{-1}(( -\infty,R))$.
Since
$$
 \nabla E_c\cdot{\boldsymbol n}=|\nabla E_c|>0\quad\text{on }\partial\Omega_R,
$$
the Poincar\'e--Hopf theorem and deformation along the gradient
above level $R$ give
$$
 \sum_{Z\in\operatorname{Crit}(E_c)}\operatorname{sgn}\det D^2E_c(Z)
       =\chi(\Omega_R)=\chi(\mathcal C_{a,b}).
$$
The fibrations obtained by forgetting one center at a time give
the Euler characteristic of the ordered configuration space.
The free, orientation-preserving action of $S_a\times S_b$
then gives
$$
 \chi(\mathcal C_{a,b})
 =\frac{\chi(\chi-1)\cdots(\chi-a-b+1)}{a!b!}
 =\binom\chi{a+b}\binom{a+b}a.
$$

\textit{Step 2.} We compute the scalar degree at each critical point.
For fixed $Z_0=(P,Q)$, put $\lambda_0=2\sigma_b$ and let
\[
 K(w)=\lambda_0G\left(\frac{H_Pe^w}{\int_MH_Pe^w\dV}-1\right).
\]
This is the scalar compact term of
\eqref{eq:shadow-tubular-homotopy} after freezing the center at
$Z_0$. We identify this operator with the one in Chen--Lin
\cite[pp.~888--890, \textup{(1.2)--(1.8)} and Theorem~1.1]{CL15}.
Their kernel satisfies
$\Delta_g G_{\rm CL}(\cdot,p)=1-\delta_p$ and
$\int_M G_{\rm CL}(\cdot,p)\dV=0$. Our kernel has the same two
normalizations. The difference of the kernels is harmonic and hence
constant on the connected closed surface. Its mean is zero, so
$G_{\rm CL}=G$. On mean-zero data, the operators satisfy $\Delta_g^{-1}=-G$.

In their singular equation take the base weight $h^*=h_2$ and put
a mark of order $\alpha_i=1$ at each $p_i$. With
$v=w-4\pi\sum_iG(\cdot,p_i)$ one has
$\Delta_g v=\Delta_g w+4\pi\sum_i(\delta_{p_i}-1)$ and
$h^*e^v=H_Pe^w$, where $H_P$ is precisely \eqref{2.5}.
The substitution cancels the singular sources in their equation.
At $p_i$, in an isothermal chart centered there,
\[
 G(z,p_i)=-\frac1{2\pi}\log|z|+R_i(z),\qquad
 H_P(z)=|z|^2 h_2(z)e^{-4\pi R_i(z)}
            \prod_{l\ne i}e^{-4\pi G(z,p_l)}.
\]
The factor multiplying $|z|^2$ is smooth and strictly positive
because all centers are distinct and $h_2$ is smooth and positive.
The weight is strictly positive away from the marks. This is the
local weight condition in their \textup{(1.6)} with $\alpha_i=1$.
The surface is closed, oriented and of area one, as required.
The parameter in their equation is $\lambda_0=2(\rho_2-4\pi b)>0$.
The exceptional set in their definition is the set of positive
numbers of the form $8\pi n+16\pi|J|$ with $n\geq0$ integral and
$J$ a subset of the marks. The empty subset and arbitrary $n\geq1$
give all of $8\pi\mathbb N$. Conversely, every positive number of
the indicated form is one of these multiples. Since $0\leq b\leq k$, we have
$8\pi(k-b)<\lambda_0<8\pi(k-b+1)$.

In their \textup{(1.7)}, $T_{\lambda_0}$ is $\lambda_0\Delta_g^{-1}$
applied to the normalized density minus one. This difference has mean zero.
It follows that $I+T_{\lambda_0}=I-K$ as maps on $Y=H^1_0(M)$.
Here the subscript denotes zero mean on a closed surface. This
identification of the weight and inverse precedes the comparison
of the two function-space degrees.

We compare the scalar degrees on $X_s$ and $Y$ before using the
formula. The fixed points of $K$ in $Y$ are smooth. On bounded
$Y$ sets, Moser--Trudinger and \eqref{eq:shadow-weight-lower-bound}
give a bound on the density in $L^2$ and hence on $K(w)$ in $H^2$.
At a fixed point, the Sobolev algebra estimate and the elliptic
inverse can then be iterated from $H^2$ to every integer order.
The a priori bound for the fixed-center equation in the cited
noncritical range, together with the scalar equation and these
estimates, bounds the fixed points in $X_s$. Choose a ball
$B_X$ in $X_s$ with all these fixed points in its interior. Choose
a ball $B_Y$ in $Y$ with the image of $\overline{B_X}$ under inclusion
in its interior and large enough to compute the degree in (1.7) of
\cite{CL15}. A larger fixed-center scalar ball in
\eqref{eq:shadow-tubular-homotopy} gives the same scalar degree
by excision, since it contains the same fixed points.

Use real mean-zero Laplace eigenfunctions, ordered by nondecreasing
eigenvalue and orthonormal in $L^2$, and let $P_N$ be the projection onto
the span $V_N$ of the first $N$. Use the spectral Sobolev
norms with weights $(1+\lambda_j)^r$ for the squares of the
$H^r$ norms, where $\lambda_j>0$ are the Laplace eigenvalues.
In these equivalent norms $P_N$ is a contraction in every
$H^r$. The Green inverse and the normalization in $K$ are
the same in both spaces. Thus $P_NK$ is the same finite-rank
map on their common domain. On $\overline{B_Y}$, the approximation
error in $Y$ tends uniformly to zero because the image of $K$ there
is bounded in $H^2$. On $\overline{B_X}$, the error in $X_s$ tends
uniformly to zero because the image of $K$ there is bounded in
$H^{s+2}$. In each case the bound on the spectral tail from
$H^{r+d}$ to $H^r$ is $(1+\lambda_{N+1})^{-d/2}$.
Continuity of $K$ on $Y$ follows from strong $H^1$ convergence
and the uniform integrability of $e^w$ in each fixed $L^p$ given by
Moser--Trudinger. Compactness follows from the $H^2$ bound.
On $X_s$, it follows from the $H^{s+2}$ bound.

On the boundary of either ball, $\|w-K(w)\|$ in that space
has positive infimum. Otherwise, a sequence with vanishing
residual would have a strongly convergent subsequence by
compactness of $K$ on the closed ball. Its limit would be a
boundary fixed point. For sufficiently large $N$, the
homotopy $w-[(1-\theta)K(w)+\theta P_NK(w)]$,
$0\le\theta\le1$, thus has no boundary zero in either space.
Its finite-rank endpoint reduces the degree in each ball to
the Brouwer degree of $I-P_NK$ on the intersection of that
ball with $V_N$.

Suppose $w_N=P_NK(w_N)$ and $w_N\in\overline{B_Y}$. From the
$H^2$ bound of $K$ on $B_Y$ and the projection contraction,
$\|w_N\|_{H^2}$ is bounded independently of $N$. If $\|w_N\|_{H^r}$
is bounded for an integer $r\ge2$, the Sobolev algebra estimate
for the exponential, the smooth weight and
\eqref{eq:shadow-weight-lower-bound} give
$\|K(w_N)\|_{H^{r+2}}\le C_r$, so $\|w_N\|_{H^{r+2}}\le C_r$.
Iterating past order $s$ gives precompactness in $X_s$.
Any limit along $N\to\infty$ is a fixed point of $K$, by the
same spectral tail bound and continuity. Since all such
points are strictly inside $B_X$, for all sufficiently
large $N$ every zero of $I-P_NK$ in $\overline{B_Y}$
belongs to $B_X$. Brouwer excision from $B_Y\cap V_N$ to
$B_X\cap V_N$ now uses the same map and the same orientation
of $V_N$. We have proved
\begin{equation}\tag{3.48a}\label{eq:scalar-degree-space-comparison}
 \deg_{\rm LS}^{\,X_s}(I-K,B_X,0)
   =\deg_{\rm LS}^{\,Y}(I-K,B_Y,0).
\end{equation}

For the right side of \eqref{eq:scalar-degree-space-comparison},
apply \cite[\textup{(1.8)}, Theorem~1.1 and Section~5,
pp.~918--919]{CL15} to the operator and weight just identified.
The approximation in Section~5 of that paper covers integer orders.
If $a=0$, one may include a mark of order zero at any point.
The mark contributes neither a source nor a weight factor.
Its $1-y$ factor cancels the change in the exponent of
$(1-y)^{-1}$. Formula \textup{(1.8)} gives the increment series
$(1-y)^{\chi-a}(1-y^2)^a=(1-y)^\chi(1+y)^a$.
The interval index is $k-b$, and the cumulative sum in their theorem
multiplies this series by $(1-y)^{-1}$. Thus the scalar degree
in the shadow product is
\begin{equation}\tag{3.48b}\label{eq:frozen-scalar-degree}
 d_{\rm sc}(Z_0)=[y^{k-b}](1-y)^{\chi-1}(1+y)^a
                  =\sum_{j=0}^a\binom aj b_{k-b-j}.
\end{equation}

Choose disjoint small coordinate balls about the Morse critical
points and use the scalar $X_s$ ball and tube from
\eqref{eq:shadow-tubular-homotopy}, enlarging the balls if
needed to contain the separated critical configurations and
all the corresponding scalar fixed points strictly. The
separation and scalar compactness estimates used there also
apply during the compactly supported Morse perturbation.
Excision restricts the degree to these products. On the ball
of a critical point $Z_0$, use the homotopy
$\mathcal S_{\rm sc}(w,Z)\mapsto\mathcal S_{\rm sc}(w,\kappa_\theta(Z))$,
where $\kappa_\theta$ contracts the center ball to $Z_0$ in its
coordinates. Leave the center and normal components fixed.
A zero during this homotopy has $n=0$ and $Z=Z_0$. Its
scalar equation is therefore exactly $w=K(w)$ for the frozen weight.
No boundary zero is introduced. In the frozen product the
center and normal components can be expressed in a smooth
frame over the coordinate ball. Contracting the base point
of this frame to $Z_0$ in the range identification is a
homotopy through invertible matrices. It keeps the zero set
fixed. Use the same frame at $Z_0$ for the domain and range.
The center block then has local degree
$\operatorname{sgn}\det D^2E_c(Z_0)$, and the normal block
contributes $+1$. The product property
therefore uses \eqref{eq:frozen-scalar-degree} as the scalar
factor. Summing and using the configuration Euler characteristic gives
$$
 \begin{aligned}
 d^S_{a,b,k}
 &=\sum_{Z\in\operatorname{Crit}(E_c)}
       \operatorname{sgn}\det D^2E_c(Z)\,d_{\rm sc}(Z)\\
 &=\binom\chi{a+b}\binom{a+b}a
                    \sum_{j=0}^a\binom aj b_{k-b-j}.
 \end{aligned}
$$
\end{proof}

\section{Local analysis at a mixed blow-up point}\label{rev:sec:core}

We use the coordinate densities of Section~\ref{sec:notation}
and the single-point widths of Section~\ref{sec:scale-conventions}.
The application of Proposition~\ref{pr4.19} to the shadow
configurations uses the mixed-center inequality \eqref{3.41}.

\subsection{Scalar profile classification}
The rank-one case of \cite[Theorems~1.1 and~1.3(i)]{LWY} supplies
the scalar classification needed below. Suppose that on $\C$,
\[
 \Delta W+2e^W=4\pi\gamma\delta_0,\qquad
 \int_{\C}e^W\dif A<\infty,\qquad \gamma\in\{0,1\}.
\]
Here $W$ is smooth off the origin and, for $\gamma=1$,
$W-2\log|z|$ is regular there. The $A_1$ Cartan matrix is $(2)$.
In the notation for the inverse-Cartan variable in that theorem,
$W/2$ is $U_1$. The formula in Theorem~1.1 becomes
\[
 e^{-W/2}=|z|^{-\gamma}
        (\lambda_0+\lambda_1|z^{1+\gamma}-p|^2),\qquad
 \lambda_0\lambda_1=\frac1{4(1+\gamma)^2}.
\]
Thus, for some $s=\lambda_0/\lambda_1>0$ and $p\in\C$,
\begin{equation*}\tag{L}\label{eq:scalar-liouville}
 e^W(z)=\frac{4(1+\gamma)^2s\,|z|^{2\gamma}}
                  {(s+|z^{1+\gamma}-p|^2)^2},\qquad
 \int_{\C}e^W\dif A=4\pi(1+\gamma).
\end{equation*}
The last integral also follows by the change of variables
$\zeta=z^{1+\gamma}$, with degree $1+\gamma$ and Jacobian
$(1+\gamma)^2|z|^{2\gamma}$.

\subsection{An explicit family of entire Toda solutions}
We construct an explicit family with two scales by a Gram
determinant calculation.

Fix $\beta>0$, $\eps>0$, and $b\in\C$, and put
\begin{align}
 F(z)&=\eps^2+\beta\eps^2|z|^2+|z^2-b|^2,\label{4.1}\\
 Q(z)&=\beta\eps^4+4\eps^2|z|^2
          +\beta\eps^2|z^2+b|^2,\label{4.2}\\
 D&=4\beta\eps^4,\qquad
 n_1=\frac{4Q}{F^2},\qquad n_2=\frac{4DF}{Q^2}.
 \label{4.3}
\end{align}
We use $\Delta=4\partial_z\partial_{\bar z}$.

\begin{lemma}\label{le4.1}
The functions $n_i$ are positive and smooth on $\C$.  Their logarithms
$U_i=\log n_i$ satisfy
$$
 -\Delta U_1=2n_1-n_2,\qquad
 -\Delta U_2=2n_2-n_1,
$$
and
$$
 \int_{\C}n_1\dif A=\int_{\C}n_2\dif A=8\pi.
$$
\end{lemma}
\begin{proof}
Let $f(z)=(\eps,\sqrt\beta\eps z,z^2-b)$.  Then
$F=|f|^2$, $Q=|f\wedge f'|^2$, and
$D=|\det(f,f',f'')|^2$. The minors $\sqrt\beta\eps^2$, $2\eps z$,
and $\sqrt\beta\eps(z^2+b)$ give \eqref{4.1}--\eqref{4.3}.
By the Gram determinant identity,
$$
 \partial_z\partial_{\bar z}\log F=\frac Q{F^2},\qquad
 \partial_z\partial_{\bar z}\log Q=\frac{DF}{Q^2}.
$$
Thus $n_1=\Delta\log F$ and $n_2=\Delta\log Q$.  Since
$U_1=\log4+\log Q-2\log F$ and
$U_2=\log(4D)+\log F-2\log Q$, the Toda equations follow.
The bounds $F\geq\eps^2$ and $Q\geq\beta\eps^4$ give positivity
and smoothness. At infinity, $F=r^4(1+O(r^{-2}))$ and
$Q=\beta\eps^2r^4(1+O(r^{-2}))$. Differentiating these polynomial
expansions and applying the divergence theorem, we obtain
$$
 \int_{B_R}n_i\dif A=\int_{\partial B_R}\partial_r\log(F\text{ or }Q)\dif s
       =8\pi+O(R^{-2}).
$$
Letting $R\to\infty$ gives the two masses.
\end{proof}

\begin{lemma}\label{le4.2}
Fix $B\in\C$ and set $b=\eps B$. As $\eps\downarrow0$, the
following convergences hold in the sense of finite measures on $\C$.
\begin{align*}
 n_1\dif A&\rightharpoonup8\pi\delta_0,\\
 n_2\dif A&\rightharpoonup4\pi\delta_0+
       n_*\dif A,\qquad
 n_*(z)=\frac{\beta}{(1+\beta|z|^2/4)^2}.
\end{align*}
In particular, $\int_\C n_*\dif A=4\pi$. The same conclusions hold
when $B$ varies in a fixed compact set, after passing to a subsequence on
which $B$ converges.
\end{lemma}
\begin{proof}
Locally in $L^1(\C)$,
$$
 \log F\longrightarrow4\log|z|,\qquad
 \log Q-2\log\eps\longrightarrow
          2\log|z|+\log(4+\beta|z|^2).
$$
Away from zero the convergence is smooth. Near the origin, the logarithm
of $|z^2-\eps B|^2$ is a sum of logarithms with singularities tending to
the origin. Since $F\geq|z^2-\eps B|^2$, local $L^1$ convergence follows
from the uniform integrability of these logarithms on compact sets. The same argument applies to
$Q/\eps^2=\beta\eps^2+4|z|^2+\beta|z^2+\eps B|^2$.
Taking the distributional Laplacians gives the limits, since
$\Delta\log|z|=2\pi\delta_0$ and
$$
 \Delta\log(4+\beta|z|^2)
       =\frac{16\beta}{(4+\beta|z|^2)^2}=n_*.
$$
Both limiting measures have total mass $8\pi$, which equals the
mass of each nonnegative measure before taking the limit. Thus convergence against compactly supported test functions extends
to convergence against bounded continuous test functions.
\end{proof}

\begin{lemma}\label{le4.3}
Use the single-point widths in \eqref{7.7}, denoted here by
$d=\sqrt{\eps}$ for the outer width and $a$ for the inner width. For fixed $B$, on compact subsets of $\C$ avoiding zero,
$$
 d^2n_1(dy)\longrightarrow
       \frac{16|y|^2}{(1+|y^2-B|^2)^2},
 \qquad
 a^2n_2(ay)\longrightarrow\frac4{(1+|y|^2)^2}.
$$
The second convergence holds locally uniformly on all of $\C$.
The limiting densities have masses $8\pi$ and $4\pi$, respectively,
and $a/d\to0$.
\end{lemma}
\begin{proof}
Substitution gives
$$
 F(dy)=\eps^2(1+|y^2-B|^2+\beta\eps|y|^2),
$$
$$
 Q(dy)=\eps^3(4|y|^2+
                \beta\eps(1+|y^2+B|^2)).
$$
This proves the first limit.  At $z=ay$, uniformly on compact sets,
$$
 F(ay)=\eps^2(1+|B|^2)(1+o(1)),\quad
 Q(ay)=\beta\eps^4(1+|B|^2)(1+|y|^2)(1+o(1)).
$$
The second limit follows from \eqref{4.3}. To compute the first mass, use $w=y^2$, and for the second use polar coordinates.
\end{proof}

For the ordinary Liouville profile, the height and scale satisfy
$$
 U_a(x)=\log\frac{4a^2}{(a^2+|x|^2)^2}
       =\lambda_i-2\log\left(1+\frac{e^{\lambda_i}}4|x|^2\right),
 \qquad \lambda_i=\log\frac4{a^2}.
$$
Thus $-\Delta U_a=2e^{U_a}$, $a=2e^{-\lambda_i/2}$, and the
inner density is $e^{U_a}$. Put $\delta=\sqrt\eps$ and $b=\delta^2B$.
For the outer singular bubble set
$$
 U_{\delta,b}(x)
 =\log\frac{32\delta^4}{(\delta^4+|x^2-b|^2)^2}
 =\lambda_o-2\log\left(1+\frac{e^{\lambda_o}}{32}|x^2-b|^2\right),
 \qquad \lambda_o=\log\frac{32}{\delta^4}.
$$
Here $-\Delta U_{\delta,b}=|x|^2e^{U_{\delta,b}}$ and the outer
density is $\frac12|x|^2e^{U_{\delta,b}}$. In particular,
$a,\delta\to0$ corresponds to $\lambda_i,\lambda_o\to+\infty$,
with $\delta=32^{1/4}e^{-\lambda_o/4}$. We use the scales
$a$ and $\delta$ below.

\begin{remark}\label{re4.4}
The function
$$
 W_B(y)=\log\frac{32}{(1+|y^2-B|^2)^2}
$$
satisfies
$$
 -\Delta W_B=|y|^2e^{W_B},\qquad
 \int_\C |y|^2e^{W_B}\dif A=16\pi.
$$
As $|B|\to\infty$, the first-component density
$m_B=\tfrac12|y|^2e^{W_B}$ splits into two ordinary bubbles of
mass $4\pi$ each. In the original coordinate we may use the
reference centers and widths
\[
 z_\pm^0=\pm\sqrt b,\qquad r_\pm\sim\frac{\eps}{2\sqrt{|b|}}.
\]
The maxima of the outer density differ from these reference centers
by $o(r_\pm)$. Indeed, write $R=|B|$ and first maximize with respect to the angle.
The squared radius in the $y$ coordinate at a maximum is
$R+s_R$, where $3s_R^2+4Rs_R=1$ and $s_R=O(R^{-1})$. The resulting
radial displacement is $O(R^{-3/2})$, whereas the width in the
$y$ coordinate is of order $R^{-1/2}$.
Corollary~\ref{co4.11} excludes the analogous splitting of an
actual solution under the strict inequality on the weights.
\end{remark}

Differentiating with respect to $B_1=\Re B$ and $B_2=\Im B$ gives the functions
$\Phi_j=\partial_{B_j}(U_1,U_2)$, which lie in the kernel of the
linearized Toda operator. At $B=0$, their limits at the outer scale are
$4(Z_j,-Z_j/2)$, where
$$
 Z_1(y)=\frac{|y|^2\cos2\theta}{1+|y|^4},\qquad
 Z_2(y)=\frac{|y|^2\sin2\theta}{1+|y|^4}.
$$
See \cite{LWY} for the parametrization of the entire family.

\begin{remark}\label{re4.5}
Replacing $f$ by
$$
 f_c(z)=(\eps,\sqrt\beta\eps(z+c),z^2-b),\qquad c\in\C,
$$
gives another family of entire solutions, with
$$
 F_c=\eps^2+\beta\eps^2|z+c|^2+|z^2-b|^2,
$$
$$
 Q_c=\beta\eps^4+4\eps^2|z|^2
                   +\beta\eps^2|z^2+2cz+b|^2,
$$
and $D$ is unchanged. For fixed $\beta$ and $c$, as $\eps\downarrow0$
and $b\to0$, the second density converges to
$\beta/(1+\beta|z+2c|^2/4)^2$ away from zero. The gradient of this limit depends
on $c$. The outer scale and the parameter $B$ are given by
$$
 C_0=1+\beta|c|^2,\qquad \delta^4=\eps^2 C_0,
 \qquad B=b/\delta^2.
$$

Suppose the limiting second density at zero is prescribed as
$\beta_0>0$ and the $z$ derivative of the regular part of the first field at
zero is prescribed as $s$. The parameters matching these data are
$$
 C_0=1+\frac{4|s|^2}{\beta_0},\qquad
 \beta=\beta_0 C_0^2,\qquad
 c=\frac{2\overline s}{\beta_0 C_0}.
$$
Indeed, up to an additive constant, the regular part of the first field is
$\log(4+\beta|z+2c|^2)$. Its $z$ derivative at zero is
$\beta\overline c/(2C_0)=s$.
\end{remark}

\subsection{Second moments of the outer profile}

For $d>0$ and $B\in\C$, define the limiting outer density
$$
 m_{d,B}(z)=\frac{16d^4|z|^2}{(d^4+|z^2-d^2B|^2)^2},
$$
and use the moment $\mathcal M(B)$ from \eqref{7.21}, which is
defined for every $B\in\C$. The change of variables $w=v+B$ gives
$$
 \mathcal M(B)=\int_\C\frac{|w|}{(1+|w-B|^2)^2}\dif A_w.
$$
\begin{lemma}\label{le4.6}
For every $B\in\C$,
$$
 \int_{\C} m_{d,B}\dif A=8\pi,\qquad
 \int_{\C} z^2m_{d,B}\dif A=8\pi d^2B,\qquad
 \int_{\C} |z|^2m_{d,B}\dif A=8d^2\mathcal M(B).
$$
If $\varphi\in C_c^3(\C;\R)$, then, uniformly for $B$ in compact sets,
$$
 \int_{\C}\varphi(z)m_{d,B}(z)\dif A_z
 =8\pi\varphi(0)+d^2\left[
 2\Delta\varphi(0)\mathcal M(B)
 +8\pi\operatorname{Re}(\varphi_{zz}(0)B)\right]+o(d^2).
$$
\end{lemma}
\begin{proof}
The map $w=z^2$ has two preimages for $w\ne0$. Summing the
Jacobian factors of the two inverse branches gives $\dif A_z=\dif A_w/(2|w|)$ and hence
$$
 m_{d,B}\dif A_z\mapsto
       \frac{8d^4}{(d^4+|w-d^2B|^2)^2}\dif A_w.
$$
Integration against $1,w,|w|$ gives the three moment identities.
For the last assertion, use the Taylor expansion
$$
 \varphi(z)=\varphi(0)+2\operatorname{Re}(\varphi_z(0)z)
  +\tfrac14\Delta\varphi(0)|z|^2
  +\operatorname{Re}(\varphi_{zz}(0)z^2)+O(|z|^3).
$$
The odd term integrates to zero by $z\mapsto-z$. For bounded $B$,
the rescaled third moment is finite, so the local error is $O(d^3)$.
The integral outside any fixed neighborhood contributes $O(d^4)$.
\end{proof}

\begin{lemma}\label{le4.7}
$\mathcal M$ is smooth, radial, and strictly convex. It satisfies
$$
 \mathcal M(0)=\frac{\pi^2}{2},\qquad
 D^2\mathcal M(0)=\frac{\pi^2}{4}I,
$$
and $\nabla\mathcal M$ is a diffeomorphism from $\R^2$ onto
the open disk of radius $\pi$.  In particular,
$$
 \lim_{R\to\infty}\nabla\mathcal M(Re^{i\theta})
       =\pi e^{i\theta}
$$
uniformly in $\theta$.
\end{lemma}
\begin{proof}
Write $K(w)=(1+|w|^2)^{-2}$.  Then
$\mathcal M(B)=\int_{\C}|w+B|K(w)\dif A_w$.  The gradient is
$$
 \nabla\mathcal M(B)=\int_{\C}\frac{w+B}{|w+B|}K(w)\dif A_w.
$$
Its norm is strictly less than $\int_{\C}K\dif A=\pi$.
Differentiating again, for $v\ne0$ we obtain the Hessian form
$$
 \int_{\C}\frac{|v|^2-[v\cdot(w+B)/|w+B|]^2}{|w+B|}
             K(w)\dif A_w>0.
$$
The integrand is locally integrable. We can transfer higher derivatives
to $K$, so $\mathcal M$ is smooth. Rotational invariance makes it radial,
and polar integration gives the values at zero.
Dominated convergence gives the limit of the gradient at infinity.
Together with the positive Hessian, this yields the radial
diffeomorphism.
\end{proof}

\begin{corollary}
\label{co4.8}
Let $\kappa\ne0$, $\tau\in\R^2$ and $\mathcal F(B)=\kappa\nabla\mathcal M(B)+\tau$.
Then
$$
 \det D\mathcal F=\kappa^2\det D^2\mathcal M>0,\qquad
 \#\mathcal F^{-1}(0)=
 \begin{cases}1,&|\tau|<\pi|\kappa|,\\0,&|\tau|\geq\pi|\kappa|.\end{cases}
$$
For sufficiently large $R$ and $|\tau|\ne\pi|\kappa|$,
$$
 \deg(\mathcal F,B_R,0)=
 \begin{cases}1,&|\tau|<\pi|\kappa|,\\0,&|\tau|>\pi|\kappa|.\end{cases}
$$
At $|\tau|=\pi|\kappa|$, $\min_{|B|=R}|\mathcal F(B)|\to0$ as $R\to\infty$.
\end{corollary}
\begin{proof}
By Lemma~\ref{le4.7}, the image of $\nabla\mathcal M$ is the open
disk of radius $\pi$, and its limit along each ray is the corresponding
point on the boundary. The degree at a zero follows from
$\det(\kappa D^2\mathcal M)=\kappa^2\det D^2\mathcal M>0$.
\end{proof}

\subsection{Matching with the outer solution}

In an isothermal coordinate $g=e^\phi|dz|^2$, normalized by
$\phi(p)=0$, set
$$
 w_i=u_i+\log h_i+\phi+\log\frac{\rho_i}{\int_Mh_ie^{u_i}\dV}.
$$
Then \eqref{1.1} takes the form
$$
 -\Delta w_i=\sum_j a_{ij}e^{w_j}+f_i,
 \qquad
 f_i=-e^\phi\sum_j a_{ij}\rho_j-\Delta(\log h_i+\phi).
$$
Consequently, at $p$,
$$
 f_i(p)=-(\Delta_g\log h_i(p)+(A\rho)_i-2K_g(p)).
$$
At $\rho_1=8\pi$,
$$
 -f_1(p)=\Delta_g\log h_1(p)+16\pi-\rho_2-2K_g(p).
$$

\subsubsection{The center equation and the quadratic coefficient}

Suppose $p$ is the only blow-up point, with local masses
$(8\pi,4\pi)$, and the second density has a nonzero limit away from
$p$. Write $\rho_2=r>4\pi$. At noncritical $r$, the formulas below also follow from Lemma~\ref{le2.1}. In general the net
first-component atom is $12\pi$, with logarithmic coefficient $6$.
If the first normalization denominator remained bounded, the
circular-mean and Jensen argument in Section~\ref{sec:noncritical-degree}
would contradict integrability. Thus the first density has no
residual part. The net atom in the second component is zero.
The second normalization denominator stays bounded. If it diverged along a
subsequence, the local bounds would force the second density to vanish off $p$,
contrary to the positive residual mass. The residual density $f_*$ has no
point mass and has mass $r-4\pi>0$. Define the mean-zero
field $z_0=2G*f_*$, where $G*f_*(x)=\int_MG(x,y)f_*(y)\dif V_g(y)$.
The net point-mass contribution to $u_2$ is zero, so $u_2\to z_0$ off $p$.
Since the limiting density is integrable, its mass on a sufficiently
small disk about $p$ is small. The Green representation and Jensen's
inequality then give an $L^q$ bound for the density on a smaller disk
for some $q>1$, as in the proof of Theorem~\ref{th2.2}. Elliptic
bootstrapping gives smoothness across $p$. Since the residual density is proportional to
$h_2e^{z_0}$ and has mass $r-4\pi$, Green's representation gives
$$
 -\Delta_g z_0=2(r-4\pi)\left(
       \frac{h_2e^{z_0}}{\int_Mh_2e^{z_0}\dV}-1\right),
 \quad u_1\longrightarrow12\pi G(\cdot,p)-\tfrac12z_0
$$
away from $p$.  In an isothermal chart centered at $p$, with
$G(x,p)=-(2\pi)^{-1}\log|x|+H(x,p)$, define the regular fields
$$
 \psi_1=\log h_1+\phi+12\pi H(\cdot,p)-\tfrac12 z_0,
 \qquad \psi_2=\log h_2+\phi+z_0.
$$
The translation Pohozaev identity gives the center equation
$$
 2\nabla\psi_1(p)+\nabla\psi_2(p)=0.
$$
Equivalently, with $R_g(p)$ the intrinsic Robin function,
$$
 \nabla_g(2\log h_1+\log h_2+12\pi R_g)(p)=0.
$$
For $u_i=-a_i\log|x|+f_i$, the gradient terms in the boundary
integral give $-2\pi\sum_{ij}(A^{-1})_{ij}a_i\nabla f_j$.
Here $a=(6,0)$ and $2\pi A^{-1}a=(8\pi,4\pi)$. Including the
coefficient terms gives the displayed center equation. The $z_0$
terms cancel in $2\nabla\psi_1+\nabla\psi_2$.

Define
$$
 T=\psi_{1,zz}(p)+\psi_{1,z}(p)^2,
 \qquad
 \kappa_1=\Delta_g\log h_1(p)+16\pi-r-2K_g(p).
$$
The regular part of the first field of the entire solution satisfies
$\psi_{1,\mathrm{core},zz}=-\psi_{1,\mathrm{core},z}^2$. After matching
the first derivative, $T$ is the harmonic quadratic coefficient of the
difference. To check invariance under a coordinate change with derivative
one, write
$z=\zeta+a\zeta^2+b\zeta^3+O(\zeta^4)$. The regular field associated
with the singular term $-6\log|z|$ transforms as
$$
 \psi_{1,\mathrm{new}}(\zeta)=\psi_1(z(\zeta))
       +2\log|z'(\zeta)|-6\log|z(\zeta)/\zeta|.
$$
Consequently
$$
 \psi_{1,\mathrm{new},z}=\psi_{1,z}-a,\qquad
 \psi_{1,\mathrm{new},zz}=\psi_{1,zz}+2a\psi_{1,z}-a^2,
$$
so $T$ is unchanged. The center equation allows a quadratic
coordinate change for which both first derivatives vanish.

\subsubsection{The angular projection}

We test the correction equation against the derivatives with respect to $B$.
Let $U_\eps$ be the entire solution from \eqref{4.1}--\eqref{4.3},
let $\Phi_j=\partial_{B_j}U_\eps$, and put $C=A^{-1}$. Since the
additive constants in $CU_\eps$ are independent of $B$,
$$
 C\Phi_j=(-\partial_{B_j}\log F,-\partial_{B_j}\log Q).
$$
On a fixed matching circle of radius $r_0$, as $\delta\to0$,
$$
 (C\Phi_1)_1=\frac{2\delta^2\cos2\theta}{r_0^2}
       +o(\delta^2),\qquad
 (C\Phi_2)_1=\frac{2\delta^2\sin2\theta}{r_0^2}
       +o(\delta^2).
$$
For the family of Remark~\ref{re4.5}, the second component can
be $O(\delta^2/r_0)$. Its boundary contribution is estimated below
after matching the first derivatives.

For the constant part of the forcing at $p$, the first component
gives
$$
 \lim_{r_0\downarrow0}\lim_{\delta\downarrow0}
 \delta^{-2}\int_{B_{r_0}}\Phi^TC f\dif A
       =-2f_1(p)\nabla\mathcal M(B).
$$
The coefficient follows from the logarithmic moment identity
$$
 \lim_{\delta\downarrow0}\delta^{-2}
 \int_{B_{r_0}}[\log F_B-\log F_0]\dif A
        =2[\mathcal M(B)-\mathcal M(0)].
$$
The terms containing $f_i(x)-f_i(p)$ tend to zero in the same iterated limit.

Write the solution as $w_\eps=U_\eps+\eta_\eps$. After matching
the first derivatives, suppose $\eta_\eps$ and its ambient first
spatial derivatives converge uniformly on every fixed matching circle.
In particular, both tangential and normal derivatives converge.
Assume that the limit $\eta$ is $C^2$ across zero and satisfies
$\eta(0)=\nabla\eta(0)=0$.
The harmonic quadratic term of $\eta_1$ is
$\operatorname{Re}(Tz^2)$.  Green's identity for
$L_\eps=-\Delta-A\operatorname{diag}(n_1,n_2)$ gives
$$
 \int_{B_{r_0}}\Phi^TC L_\eps\eta_\eps\dif A
 =\int_{\partial B_{r_0}}
       [\partial_n(C\Phi)\cdot\eta_\eps
                     -(C\Phi)\cdot\partial_n\eta_\eps]\dif s.
$$
For the first component, the cosine mode gives
$$
 \int_0^{2\pi}
 \left[\left(-\frac{4\delta^2\cos2\theta}{r_0^3}\right)
                r_0^2\operatorname{Re}(Te^{2i\theta})
       -\frac{2\delta^2\cos2\theta}{r_0^2}
                2r_0\operatorname{Re}(Te^{2i\theta})\right]
                     r_0\dif\theta
       =-8\pi\delta^2\operatorname{Re}T.
$$
The sine mode gives $+8\pi\delta^2\operatorname{Im}T$. The radial
term $\frac14(\Delta\eta_1)(0)r_0^2$ integrates to zero since
$\int_0^{2\pi}e^{2i\theta}\dif\theta=0$. Matching the first
derivatives gives $o_{r_0}(\delta^2)$ for the second component.

The exact correction equation is
$$
 L_\eps\eta_\eps=f+
 A\operatorname{diag}(n_i)
             (e^{\eta_\eps}-1-\eta_\eps).
$$
If the nonlinear term satisfies
$$
 \lim_{r_0\downarrow0}\limsup_{\delta\downarrow0}
 \delta^{-2}\left|\int_{B_{r_0}}\sum_i
       \Phi_{ij}n_i(e^{\eta_{i,\eps}}-1-\eta_{i,\eps})\dif A\right|=0,
 \quad j=1,2,
$$
the leading angular equation is
$$
 \quad
 \kappa_1\nabla\mathcal M(B)
       +4\pi(\operatorname{Re}T,-\operatorname{Im}T)=0.
 \quad
$$
For the projected family, use the quadratic correction estimate
\eqref{7.27}. For actual sequences, first determine the limiting angular
parameter from \eqref{4.26}--\eqref{4.33}, independently of this
projection. The profile estimates in Section~\ref{sec:strong-chapter}
and the modulation in Section~\ref{sec:typedactual-chapter}
give $\delta_j\asymp D$ and the global mean-zero correction
$\|\eta\|_{H^1}=O(D^{4/3}\ell)$ by \eqref{7.53} and \eqref{7.62}.
The correction remains uniformly small in $L^\infty$ by
\eqref{7.54}. By \eqref{7.26} and the weighted estimates preceding
\eqref{7.61}, the normalization constants are
$O(D^2+\sqrt\ell\|\eta\|_{H^1})$. Thus the contribution of
the correction and its normalization to the absolute weighted
quadratic remainder is bounded by
$C(D^4+\ell\|\eta\|_{H^1}^2)=o(D^2)$.
The weights $n_i$ and the completed weights are comparable on the
fixed testing disks, and the angular kernel functions are uniformly bounded.

In both applications the smooth matching terms satisfy
$J_1=O(r^2)$ and $J_2=O(r)$ in the normalized mixed chart.
On $r\leq\delta_j$, the second moment of the inner peak is
$O(\delta_j^4|\log\delta_j|)$ and the bounded part contributes
$O(\delta_j^4)$. On $\delta_j\leq r\leq r_0$,
$n_2\leq C$ and $|\Phi_{2l}|\leq C\delta_j^2/r^2$ by the
explicit profile. The resulting integral of
$|\Phi_{2l}|n_2 J_2^2$ on this annulus is
$O(\delta_j^2r_0^2)$. The first component is bounded by its
fourth moment, $O_{r_0}(\delta_j^4|\log\delta_j|)$.
These bounds yield the stated angular remainder by first taking
$\delta_j\downarrow0$ and then $r_0\downarrow0$.
Under a uniform $C^0$ reduction, Corollary~\ref{co4.8} gives angular
degree $+1$ for $4|T|<|\kappa_1|$ and degree zero for
$4|T|>|\kappa_1|$. At a nonzero equality $4|T|=|\kappa_1|>0$,
the model has no finite zero and approaches zero along a direction at
infinity. If $\kappa_1=0$, the model is the constant $4\pi\overline T$,
which vanishes identically when $T=0$.
Section~\ref{sec:tward-chapter} treats the corresponding splitting
using an exact integral identity. Section~\ref{sec:inv-chapter}
obtains the local mass identity by testing against the scale
derivative, and Section~\ref{article:families} computes the local degree.

The entire family has masses $(8\pi,8\pi)$. To use it near a
mixed point on $M$, we match its limiting second density to the
outer solution at the prescribed noncritical second parameter.

\subsection{Estimates for the nonlinear remainder}
\label{sec:energy-remainder}

\begin{lemma}
\label{le4.9}
Let $M$ be a fixed closed surface, let $0<\delta<e^{-2}$, and
suppose $\mu_\delta\dV$ is a probability measure with
$$
 0\le\mu_\delta\le C\delta^{-4}.
$$
For a real $v\in H^1(M)$ with zero average, put
$E=\|v\|_{H^1}$ and $\ell_\delta=\log(1/\delta)$. For each fixed
$a>0$, one has
\begin{equation}\label{4.4}
 \int_M\mu_\delta v^2e^{a|v|}\dV
 \le C_a\ell_\delta E^2
             \exp(C_a\ell_\delta E^2).
\end{equation}
If $\ell_\delta E^2\to0$ and
$$
 c_\delta=\log\int_M\mu_\delta e^v\dV,
 \qquad \widehat v=v-c_\delta,
$$
then
\begin{equation}\label{4.5}
 |c_\delta|\le C\sqrt{\ell_\delta}\,E,
 \qquad
 \int_M\mu_\delta\widehat v^{\,2}e^{|\widehat v|}\dV
       \le C\ell_\delta E^2
\end{equation}
for all sufficiently small $\delta$.
\end{lemma}

\begin{proof}
Put $q=\ell_\delta$ and $p=q/(q-1)$. Then
$\|\mu_\delta\|_p\leq(C\delta^{-4})^{1/q}\leq C$.
H\"older's inequality, the two-dimensional Sobolev estimate, and
the Moser--Trudinger inequality \cite{Moser71} give
$$
 \int_{M}\mu_\delta v^2e^{a|v|}\dV
 \leq C\|v\|_{4q}^2
          \left(\int_{M} e^{2qa|v|}\dV\right)^{1/(2q)}
 \leq C_aqE^2e^{C_aqE^2}.
$$
Put $m=\int_{M}\mu_\delta v\dV$. By H\"older's inequality and
\eqref{4.4}, $|m|\leq C\sqrt qE$ and
$\int_{M}\mu_\delta e^v\dV=1+m+O(qE^2)$. If $qE^2=o(1)$, then
$c_\delta=m+O(qE^2)$. The same H\"older and Moser--Trudinger
estimates give $\int_M\mu_\delta e^{|v|}\dV\leq C e^{CqE^2}$.
Consequently,
\[
 \int_M\mu_\delta\widehat v^2e^{|\widehat v|}\dV
 \leq 2e^{|c_\delta|}\int_M\mu_\delta
              (v^2+c_\delta^2)e^{|v|}\dV
 \leq CqE^2,
\]
by \eqref{4.4}, $|c_\delta|\leq C\sqrt q E$, and $qE^2=o(1)$.
\end{proof}

\begin{corollary}
\label{co4.10}
Let $B,\beta,c$ vary in compact sets with $\beta\geq\beta_0>0$, and put
$n_i=e^{U_i}$, $\Phi_{ij}=\partial_{B_j}U_i$. Suppose
$$
 \eta_{i,\delta}=a_{i,\delta}+v_{i,\delta}-c_{i,\delta},\qquad
 |a_{i,\delta}(x)|\leq C|x|^2,\qquad \int_Mv_{i,\delta}\dV=0.
$$
Suppose the probability measures and normalization constants satisfy
the hypotheses of Lemma~\ref{le4.9}, and suppose each $n_i\dif A$ is bounded
by a fixed multiple of the corresponding probability measure. If
\begin{equation}\label{4.6}
 \|v_{i,\delta}\|_{H^1}
       =o\!\left(\frac\delta{\sqrt{\log(1/\delta)}}\right),\qquad i=1,2,
\end{equation}
then, for $j=1,2$,
$$
 \lim_{r_0\downarrow0}\limsup_{\delta\downarrow0}\delta^{-2}
 \left|\int_{B_{r_0}}\sum_i\Phi_{ij}n_i
             (e^{\eta_{i,\delta}}-1-\eta_{i,\delta})\dif A\right|=0.
$$
If $c_{i,\delta}=0$, it suffices to use the mass and $L^\infty$
bounds with a dominating probability measure of the same type.
\end{corollary}

\begin{proof}
Uniformly on compact parameter sets, $|\Phi_{ij}|\le C$,
$n_i\le C\delta^{-4}$, and
$$
 \int_{B_{r_0}}\sum_i|\Phi_{ij}|n_i|x|^4\dif A
       \le C(\delta^4+\delta^2r_0^4).
$$
For real $a,b$, with $a$ uniformly bounded,
$$
 |e^{a+b}-1-a-b|
 \le C a^2+C b^2e^{|b|}.
$$
Apply this inequality with $a=a_{i,\delta}$ and
$b=v_{i,\delta}-c_{i,\delta}$. By \eqref{4.5}, the absolute value
of the projected remainder, divided by $\delta^2$, is bounded by
$$
 C\left(\delta^2+r_0^4+
    \frac{\log(1/\delta)}{\delta^2}
       \sum_i\|v_{i,\delta}\|_{H^1}^2\right).
$$
Assumption \eqref{4.6} gives the required limit.
\end{proof}

\subsection{Local integral identities}\label{rev:sec:stress}

\label{sec:tward-chapter}
\subsubsection{The local equation and the limiting regular parts}

Let $M$ be a smooth closed surface of area one, let $h_i>0$ be smooth,
and let $u_{i,\nu}$ solve the regular SU(3) Toda equations with
$\rho_{1,\nu}\to8\pi$ and $\rho_{2,\nu}\to\rho>4\pi$.
Use $-\Delta_gG(x,p)=\delta_p-1$ and $\int_MG(x,p)\dV=0$.
Set
$$
 f_{i,\nu}=\rho_{i,\nu}
 \frac{h_i e^{u_{i,\nu}}}{\int_Mh_i e^{u_{i,\nu}}\dV}.
$$
Suppose $q$ is a mixed concentration point and, in a neighborhood of $q$,
\begin{equation}\label{4.7}
 f_{1,\nu}\dV\rightharpoonup8\pi\delta_q,
 \qquad f_{2,\nu}\dV\rightharpoonup4\pi\delta_q+f_*\dV,
\end{equation}
where $f_*$ is smooth near $q$. If $\rho\notin4\pi\mathbb N$,
Lemma~\ref{le2.1} applies to the converging parameter sequence.
Since $q$ is mixed and the first limiting mass is $8\pi$, the identity
$a+2b=2$ in \eqref{2.4} forces $a=0$ and $b=1$. Thus $q$ is the
only center and the global residual weight has no ordinary marked points.
In this case \eqref{2.5} and \eqref{2.7} give
$$
 f_*=(\rho-4\pi)\frac{h_2e^{w_*}}{\int_Mh_2e^{w_*}\dV},
 \qquad
 \Delta_gw_*+2(\rho-4\pi)
 \left(\frac{h_2e^{w_*}}{\int_Mh_2e^{w_*}\dV}-1\right)=0.
$$
The local identities below use the stated local data. They also
apply with $f_*=0$ when the regular derivative limits in
\eqref{4.11} are supplied independently.

In an isothermal coordinate $z$ with $q=0$ and $g=e^\phi|dz|^2$, put
$$
 F_{i,\nu}=e^\phi f_{i,\nu},\qquad U_{i,\nu}=\log F_{i,\nu}.
$$
The functions $U_{i,\nu}$ are smooth and satisfy the following
system with $\Delta_0=4\partial_z\partial_{\bar z}$.
\begin{equation}\label{4.8}
 \Delta_0U_{1,\nu}+2e^{U_{1,\nu}}-e^{U_{2,\nu}}=\kappa_{1,\nu}^c,
 \qquad
 \Delta_0U_{2,\nu}-e^{U_{1,\nu}}+2e^{U_{2,\nu}}=\kappa_{2,\nu}^c,
\end{equation}
where
$$
 \kappa_{i,\nu}^c
 =\Delta_0(\log h_i+\phi)+\sum_{j=1}^2 A_{ij}\rho_{j,\nu}e^\phi,
 \qquad A=\begin{pmatrix}2&-1\\-1&2\end{pmatrix}.
$$
The functions $\kappa_{i,\nu}^c$ converge smoothly. At $q$,
\begin{equation}\label{4.9}
 \kappa_1^c(q)=e^{\phi(q)}\kappa_1(q),\qquad
 \kappa_1=\Delta_g\log h_1+16\pi-\rho-2K_g.
\end{equation}

Write
$$
 G(z,q)=-\frac{1}{2\pi}\log|z|+R(z,q),\qquad
 r_g(q)=R(q,q)+\frac{\phi(q)}{4\pi},
$$
and, for the global single-center situation just described,
define the smooth parts of the limiting fields by
\begin{equation}\label{4.10}
 \psi_1=\log h_1+\phi+12\pi R(\cdot,q)-\frac12w_* ,
 \qquad \psi_2=\log h_2+\phi+w_*.
\end{equation}
The following convergences are smooth on compact sets away from $q$.
\begin{equation}\label{4.11}
 U_{1,\nu,z}\longrightarrow-\frac3z+\psi_{1,z},
 \qquad U_{2,\nu,z}\longrightarrow\psi_{2,z}
\end{equation}
The stress identities below also apply to the local system \eqref{4.8} with
smoothly convergent forcings, the measure limits \eqref{4.7},
and the derivative limits \eqref{4.11}. For the global fields above set
\begin{equation}\label{4.12}
 C=2\log h_1+\log h_2+12\pi r_g,
 \qquad \mathcal T=\psi_{1,zz}+\psi_{1,z}^{\,2}.
\end{equation}
The local stress identities use only the quantity $\mathcal T$ in this formula.

At a mixed center $q=q_j$ in a configuration with several centers,
form $\psi_i$ in \eqref{4.11} from the regular parts in \eqref{2.6}.
The corresponding $C$ in \eqref{4.12} includes the additional terms
$12\pi\sum_i G(\cdot,p_i)+24\pi\sum_{\ell\ne j}G(\cdot,q_\ell)$.
Equation \eqref{4.8} gives $8\pi m-\rho_2$ in place of
$16\pi-\rho$ in \eqref{4.9}.

\subsubsection{A necessary condition for splitting}

Suppose three ordinary Liouville bubbles carry the full local
masses $(8\pi,4\pi)$ at $q$.
\begin{enumerate}
\item The first component has two ordinary Liouville bubbles of
mass $4\pi$ each, centered at
$p_{\pm,\nu}=s_\nu\pm d_\nu+o(|d_\nu|)$, where $s_\nu\to q$,
$0\ne d_\nu\to0$.
\item The full local mass $4\pi$ of the second component is carried
by one ordinary Liouville bubble centered at $s_\nu+o(|d_\nu|)$.
\item The three scales are $o(|d_\nu|)$. After centering and
rescaling at its own scale, each measure converges narrowly to a
radial ordinary Liouville measure of mass $4\pi$.
\end{enumerate}
We may restrict each bubble measure to a disk whose radius is much
larger than the bubble scale and is $o(|d_\nu|)$. The omitted mass tends to zero.

\begin{corollary}\label{co4.11}
Along a subsequence with $\bar d_\nu/d_\nu\to\tau$, $|\tau|=1$,
we have
\begin{equation}\label{4.13}
 \mathcal T(q)=-\frac14\kappa_1^c(q)\tau.
\end{equation}
When \eqref{4.9} holds, the splitting condition reads
\begin{equation}\label{4.14}
 e^{\phi(q)}|\kappa_1(q)|=4|\mathcal T(q)|.
\end{equation}
Thus $e^{\phi(q)}\kappa_1(q)>4|\mathcal T(q)|$ excludes splitting.
\end{corollary}

We postpone the proof of Corollary~\ref{co4.11} until after
Theorem~\ref{th4.16}.

\subsubsection{Differential identities for $W_2$ and $W_3$}

Write the components of $A^{-1}U_\nu$ as
$$
 a_\nu=\frac{2U_{1,\nu}+U_{2,\nu}}3,
 \qquad b_\nu=\frac{U_{1,\nu}+2U_{2,\nu}}3,
$$
and
$$
 K_{1,\nu}=\frac{2\kappa_{1,\nu}^c+\kappa_{2,\nu}^c}{12},
 \qquad K_{2,\nu}=\frac{\kappa_{1,\nu}^c+2\kappa_{2,\nu}^c}{12}.
$$
Thus $\partial_{\bar z}a_{\nu,z}=K_{1,\nu}-F_{1,\nu}/4$ and
$\partial_{\bar z}b_{\nu,z}=K_{2,\nu}-F_{2,\nu}/4$.
Primes denote derivatives with respect to $z$. Set
\begin{align}
 W_{2,\nu}&=a_\nu''+b_\nu''-(a_\nu')^2+a_\nu'b_\nu'-(b_\nu')^2,
 \label{4.15}\\
 W_{3,\nu}&=a_\nu'''-2a_\nu'a_\nu''+a_\nu'b_\nu''
                +(a_\nu')^2b_\nu'-a_\nu'(b_\nu')^2.
 \label{4.16}
\end{align}
These functions are the coefficients of the Miura operator
$(\partial-b_\nu')(\partial+b_\nu'-a_\nu')(\partial+a_\nu')$.

\begin{lemma}\label{le4.12}
For every $\nu$, the following identities hold on the coordinate disk.
\begin{equation}\label{4.17}
 \partial_{\bar z}W_{2,\nu}
 =K_{1,\nu}'+K_{2,\nu}'
 -(2a_\nu'-b_\nu')K_{1,\nu}
 -(2b_\nu'-a_\nu')K_{2,\nu},
\end{equation}
and
\begin{equation}\label{4.18}
 \begin{split}
 \partial_{\bar z}W_{3,\nu}={}&K_{1,\nu}''-2a_\nu'K_{1,\nu}'
                              +a_\nu'K_{2,\nu}'\\
 &+K_{1,\nu}(-2a_\nu''+b_\nu''+2a_\nu'b_\nu'-(b_\nu')^2)\\
 &+K_{2,\nu}((a_\nu')^2-2a_\nu'b_\nu').
 \end{split}
\end{equation}
The nonlinear density terms cancel in both identities.
\end{lemma}
\begin{proof}
Suppress the index $\nu$ and put $P=a'$, $Q=b'$,
$X=e^{2a-b}/4$ and $Y=e^{2b-a}/4$. By \eqref{4.8},
$$
 \bar\partial P=K_1-X,\qquad \bar\partial Q=K_2-Y,
 \qquad X'=(2P-Q)X,\quad Y'=(2Q-P)Y.
$$
For $W_2=P'+Q'-P^2+PQ-Q^2$, the density terms in
$\bar\partial W_2$ are exactly
$$
 -X'-Y'+(2P-Q)X+(2Q-P)Y=0.
$$
The remaining terms are $K_1'+K_2'-(2P-Q)K_1-(2Q-P)K_2$,
which gives \eqref{4.17}.

For the cubic identity write
$W_3=P''-2PP'+PQ'+P^2Q-PQ^2$. The terms containing $X$ are
$$
 -X''+2XP'+2PX'-XQ'-2PQX+Q^2X.
$$
Substitution of $X''=(2P'-Q'+(2P-Q)^2)X$ cancels the derivative
and quadratic terms. The terms containing $Y$ are
$$
 -PY'-P^2Y+2PQY
 =(-P(2Q-P)-P^2+2PQ)Y=0.
$$
The terms containing $K_1,K_2$ are
$$
 K_1''-2PK_1'+PK_2'
 +K_1(-2P'+Q'+2PQ-Q^2)+K_2(P^2-2PQ),
$$
which is \eqref{4.18} for each smooth solution.
\end{proof}

Let $A_0=(2\psi_{1,z}+\psi_{2,z})/3$ and
$B_0=(\psi_{1,z}+2\psi_{2,z})/3$ as smooth functions of $z$.
By \eqref{4.11}, $a_\nu'\to-2/z+A_0$ and
$b_\nu'\to-1/z+B_0$. Substituting into \eqref{4.15} gives
\begin{equation}\label{4.19}
 W_{2,*}=\frac{3A_0}{z}+A_0'+B_0'-A_0^2+A_0B_0-B_0^2.
\end{equation}
Substituting into \eqref{4.16} gives
\begin{equation}\label{4.20}
 \begin{split}
 W_{3,*}={}&\frac{4A_0'-2B_0'-2A_0B_0-A_0^2+2B_0^2}{z}\\
 &+A_0''-2A_0A_0'+A_0B_0'+A_0^2B_0-A_0B_0^2.
 \end{split}
\end{equation}
The terms of order $z^{-2}$ and $z^{-3}$ cancel identically in
\eqref{4.20}, independently of the center equation.

\begin{lemma}\label{le4.13}
The assumptions above imply
$$
  2\psi_{1,z}(q)+\psi_{2,z}(q)=0.
$$
Consequently $A_0(q)=0$, $B_0(q)=-\psi_{1,z}(q)$, and
\begin{equation}\label{4.21}
 W_{3,*}(z)=\frac{2\mathcal T(q)}z+O(1)
 \quad\hbox{as }z\to0.
\end{equation}
The bounded remainder may contain terms proportional to $\bar z/z$.
\end{lemma}
\begin{proof}
Fix $1<p<2$. For the smoothly truncated Cauchy kernel $K_\varepsilon$,
$$
 K(z)=z^{-1},\qquad
 \|K-K_\varepsilon\|_{L^p(B_{2\varepsilon})}
       \leq C\varepsilon^{2/p-1},\qquad
 K_\varepsilon*\mu_\nu\longrightarrow K_\varepsilon*\mu
       \quad\text{locally uniformly}.
$$
Green representation and bounded mass therefore give
$$
 (a_\nu',b_\nu')\longrightarrow(-2/z+A_0,-1/z+B_0)
              \quad\text{in }L^p_{\rm loc}.
$$
Let $\chi_r(z)=\chi(|z|/r)$, with $\chi=1$ near the origin and
$\operatorname{supp}\chi_r\Subset B_r$. By \eqref{4.17},
$$
 \left|\int_{B_r}\chi_r\bar\partial W_{2,\nu}\dif A\right|
 \leq Cr^2+Cr^{2-2/p}
              (\|a_\nu'\|_p+\|b_\nu'\|_p)\longrightarrow0
              \quad(r\downarrow0),
$$
uniformly in $\nu$. The support of $\bar\partial\chi_r$ stays away from the origin. Hence
\eqref{4.19} and integration by parts give
$$
 \begin{aligned}
 3\pi A_0(q)
 &=\lim_{r\downarrow0}-\int_{B_r} W_{2,*}\bar\partial\chi_r\dif A\\
 &=\lim_{r\downarrow0}\lim_{\nu\to\infty}
                 \int_{B_r}\chi_r\bar\partial W_{2,\nu}\dif A=0,
 \qquad \bar\partial(1/z)=\pi\delta_0.
 \end{aligned}
$$
Thus
$$
  3A_0(q)=0,\quad B_0(q)=-\psi_{1,z}(q),\quad
 4A_0'(q)-2B_0'(q)+2B_0(q)^2=2\mathcal T(q).
$$
Substitution in \eqref{4.20} proves \eqref{4.21}.
\end{proof}
For the global fields, \eqref{4.10} and the Robin identity in the proof of Lemma~\ref{le2.1} give $2\psi_{1,z}(q)+\psi_{2,z}(q)=C_z(q)$. For a configuration with several centers, \eqref{2.6} gives the same equality with the additional Green terms above.
\subsubsection{Limits of the products of derivatives}

To pass to the limit in the quadratic terms of \eqref{4.18}, first
integrate the product of the two Cauchy kernels. The resulting kernel
is bounded, and we can compute its limit from the rescaled measures. Set $\Theta(s)=\bar s/s$ for $s\ne0$, so $|\Theta(s)|=1$.
Its value at zero is immaterial for integrals against smooth densities.

\begin{lemma}\label{le4.14}
For $k\in C_c^\infty(\C)$ and $x\ne y$, put
$$
 \cC k(x)=\frac1\pi\int_\C\frac{k(z)}{x-z}\dif A_z.
$$
Then
\begin{equation}\label{4.22}
 \int_\C\frac{k(z)}{(z-x)(z-y)}\dif A_z
 =-\pi\frac{\cC k(x)-\cC k(y)}{x-y}.
\end{equation}
As $(x,y)\to(q,q)$ off the diagonal,
\begin{equation}\label{4.23}
 \frac{\cC k(x)-\cC k(y)}{x-y}
 =(\cC k)_z(q)+k(q)\Theta(x-y)+o(1),
\end{equation}
with a uniform remainder.  For $k_r=\chi_r K$, $K$ smooth,
\begin{equation}\label{4.24}
 \sup_{x\ne y}\left|
 \frac{\cC k_r(x)-\cC k_r(y)}{x-y}\right|\le C,
 \qquad (\cC k_r)_z(q)=O(r).
\end{equation}
The constants can be chosen uniformly for $K$ in a bounded smooth family.
\end{lemma}
\begin{proof}
Partial fractions give \eqref{4.22}. Since $\bar\partial\cC k=k$
and $\cC k$ is $C^1$, we integrate its derivative along the segment
from $y$ to $x$ to obtain \eqref{4.23} with a uniform remainder, even
when $|x-y|$ is much smaller than $|x-q|+|y-q|$.

For \eqref{4.24}, rescale $z=r\xi$. The functions
$\chi(|\xi|)K(r\xi)$ have fixed compact support and uniform smooth
bounds. The corresponding Cauchy transforms have uniformly bounded
first derivatives on a fixed disk, and their derivatives decay outside
it. This proves the first bound in \eqref{4.24}. For the second, use
$$
 (\cC k_r)_z(0)=-\frac1\pi\operatorname{pv}
       \int_{B_r}\frac{\chi_r(z)K(z)}{z^2}\dif A_z.
$$
The contribution of the constant part $K(0)$ vanishes by angular
integration. The remaining absolute integral is at most $C\int_0^r ds=Cr$.
\end{proof}

Choose $r_\nu\downarrow0$ sufficiently slowly that the measures
$$
 \mu_{i,\nu}=\boldsymbol1_{B_{r_\nu}(q)}F_{i,\nu}\dif A
$$
have total masses tending to $M_1=8\pi$ and $M_2=4\pi$, respectively.
After passing to a subsequence, the following three limits exist.
\begin{equation}\label{4.25}
 S_{ij}=\lim_{\nu\to\infty}
    \iint_{\C\times\C}\Theta(x-y)\dif\mu_{i,\nu}(x)\dif\mu_{j,\nu}(y)
 \qquad (ij=11,12,22)
\end{equation}
These limits are unchanged under perturbations of the measures with
total variation $o(1)$.

\begin{lemma}\label{le4.15}
For $K$ smooth near $q$,
\begin{align*}
 \lim_{r\downarrow0}\lim_{\nu\to\infty}
 \int_{B_r}\chi_r K(a_\nu')^2\dif A&=-\frac{K(q)}{16\pi}S_{11},\\
 \lim_{r\downarrow0}\lim_{\nu\to\infty}
 \int_{B_r}\chi_r Ka_\nu'b_\nu'\dif A&=-\frac{K(q)}{16\pi}S_{12},\\
 \lim_{r\downarrow0}\lim_{\nu\to\infty}
 \int_{B_r}\chi_r K(b_\nu')^2\dif A&=-\frac{K(q)}{16\pi}S_{22}.
\end{align*}
The same formulas hold for smoothly convergent $K_\nu$.
\end{lemma}
\begin{proof}
On a smaller concentric disk, Green's representation gives
$$
 a_\nu'(z)=A_\nu^{\rm sm}(z)-\frac1{4\pi}
             \int_{D}\frac{F_{1,\nu}(x)}{z-x}\dif A_x,
 \quad
 b_\nu'(z)=B_\nu^{\rm sm}(z)-\frac1{4\pi}
             \int_{D}\frac{F_{2,\nu}(y)}{z-y}\dif A_y,
$$
The smooth terms converge with all derivatives and satisfy uniform
bounds. By the local $L^p$ bounds in the proof of Lemma~\ref{le4.13},
terms containing at most one Cauchy integral vanish as $r\downarrow0$.

For fixed $\nu$, the $L^\infty(D)$ bounds on the densities imply that
the triple integral of
$|z-x|^{-1}|z-y|^{-1}F_{i,\nu}(x)F_{j,\nu}(y)$ is finite. Fubini's
theorem and \eqref{4.22} with $k_r=\chi_r K$ give the coefficient
$-(4\pi)^{-2}\pi=-1/(16\pi)$ in all three formulas. The divided
difference is uniformly bounded by \eqref{4.24}. For the term with
$\mu_{i,\nu}$ and $\mu_{j,\nu}$, \eqref{4.23} gives the limit
$$
 -\frac1{16\pi}\left[
       (\cC k_r)_z(q)M_iM_j+K(q)S_{ij}\right].
$$
The first term vanishes as $r\downarrow0$ by \eqref{4.24}.

For the remaining terms, set
$\lambda_{i,\nu}=\boldsymbol1_D F_{i,\nu}\dif A-\mu_{i,\nu}$.
Choose $D$ so that its boundary has zero mass for the limiting measures.
Then these measures converge weakly to the smooth parts in \eqref{4.7} restricted to $D$.

Any product with at least one $\lambda_{i,\nu}$ factor has a weak
limit giving zero mass to the diagonal. Since the kernel in \eqref{4.24}
is bounded and continuous off the diagonal, the integrals converge
for fixed $r$. For $x\ne y$, the kernel tends to zero as $r\downarrow0$
by $\|\cC k_r\|_\infty=O(r)$. Dominated convergence gives zero in the
iterated limit. The bounds are uniform for smoothly convergent
$K_\nu\to K$, so the same limits hold in that case.
\end{proof}

\begin{theorem}\label{th4.16}
Assume \eqref{4.7}, \eqref{4.8}, and \eqref{4.11}, and suppose the
forcing terms converge smoothly. Then
\begin{equation}\label{4.26}
 32\pi^2\mathcal T(q)
 =-K_1(q)(2S_{12}-S_{22})
   -K_2(q)(S_{11}-2S_{12}).
\end{equation}
If $S_{11}=S_{12}=S$ and $S_{22}=0$, this reduces to
\begin{equation}\label{4.27}
 \mathcal T(q)=-\frac{\kappa_1^c(q)}{128\pi^2}S.
\end{equation}
\end{theorem}
\begin{proof}
We test \eqref{4.18} against the radial cutoff $\chi_r$. First let
$\nu\to\infty$ at fixed $r$, and then let $r\downarrow0$.

\emph{The boundary term.} Integration by parts gives
$-\int_{B_r} W_{3,\nu}\bar\partial\chi_r\dif A$. For fixed $r$,
\eqref{4.11} gives smooth convergence on the support of
$\bar\partial\chi_r$. By
\eqref{4.21} and $\bar\partial(1/z)=\pi\delta_0$,
$$
 \lim_{r\downarrow0}\lim_{\nu\to\infty}
       -\int_{B_r} W_{3,\nu}\bar\partial\chi_r\dif A
       =2\pi\mathcal T(q).
$$
The calculation uses only the convergence on that annulus.

\emph{Terms with at most one gradient.} The coefficients and their
derivatives are uniformly bounded. The local $L^p$ bounds on
$a_\nu',b_\nu'$ for $p<2$ imply that the integrals of
$K_1''$, $a'K_1'$ and $a'K_2'$ over $B_r$ vanish as $r\downarrow0$.

\emph{Second derivatives.} For example,
$$
 \int_{B_r}\chi_r K_1a_\nu''\dif A
      =-\int_{B_r}(\chi_r K_1)_z a_\nu'\dif A.
$$
For fixed $r$, use the strong $L^p$ limit $a_\nu'\to-2/z+A_0$.
The integral with $K_1(q)(\chi_r)_z/z$ vanishes by angular
integration of $e^{-2i\theta}$. Replacing $K_1$ by $K_1(q)$ on the
annulus incurs an $O(r)$ error. The term with $(K_1)_z\chi_r/z$
is also $O(r)$, and the contribution of $A_0$ tends to zero.
The same calculation for $b_\nu''$ gives zero in the iterated limit.

\emph{Products of gradients.} Lemma~\ref{le4.15} gives
$$
 \lim_{r\downarrow0}\lim_{\nu\to\infty}
 \int_{B_r}\chi_r(K_1(2a_\nu'b_\nu'-(b_\nu')^2)
            +K_2((a_\nu')^2-2a_\nu'b_\nu'))\dif A
 =-\frac{K_1(2S_{12}-S_{22})+K_2(S_{11}-2S_{12})}{16\pi}.
$$
Equating the two limits gives \eqref{4.26}. If $S_{11}=S_{12}=S$
and $S_{22}=0$, then $(2K_1-K_2)S=\kappa_1^cS/4$, which gives
\eqref{4.27}.
\end{proof}

\subsection{Exclusion of splitting and the mixed profile}\label{sec:mixed-formations}

\begin{proof}[Proof of Corollary~\ref{co4.11}]
We evaluate the integrals in \eqref{4.26}. Write $m=4\pi$ for the
mass of each ordinary bubble and $\tau=\lim\bar d_\nu/d_\nu$.
Since $\Theta(-z)=\Theta(z)$, the two orderings of the outer pair give
the same factor. For a finite radial measure $\lambda$ with
$\lambda(\{x\})=0$ for every $x\in\C$, simultaneous rotation of
$x$ and $y$ gives
$\iint_{\C\times\C}\Theta(x-y)\dif\lambda(x)\dif\lambda(y)=0$.
The kernel is bounded and continuous off the diagonal, which has
zero mass for the limiting product measure. Narrow convergence
therefore gives zero for the contribution of each bubble paired
with itself.

The distances between distinct bubbles are comparable to $|d_\nu|$,
and their scales are $o(|d_\nu|)$. Each pair of distinct bubbles
therefore contributes the factor $\tau$. There are two ordered
pairs from the first component alone and two pairs between the
components, each of product mass $m^2$. The second component has
only one bubble. Hence
$$
 \begin{aligned}
 S_{11}&=\underbrace{0+0}_{\text{self pairs}}+
                  \underbrace{m^2\tau+m^2\tau}_{\text{ordered cross pairs}}
                  =32\pi^2\tau,\\
 S_{12}&=m^2\tau+m^2\tau=32\pi^2\tau,\qquad S_{22}=0.
 \end{aligned}
$$
The omitted mass contributes $o(1)$ since $|\Theta|=1$.
Substitution into \eqref{4.27} gives
$$
 \mathcal T(q)
 =-\frac{\kappa_1^c(q)}{128\pi^2}\,32\pi^2\tau
 =-\frac14\kappa_1^c(q)\tau.
$$
Taking absolute values and using \eqref{4.9} gives \eqref{4.14},
which is incompatible with $e^{\phi(q)}\kappa_1(q)>4|\mathcal T(q)|$.
\end{proof}

For distinct bubbles, replacing $\Theta(x-y)$ by $\bar d_\nu/d_\nu$
is valid uniformly on the truncated cores in the stated scale regime.
Uniform tightness and $|\Theta|\leq1$ then control the discarded integrals
uniformly.

\subsubsection{The equation for the parameter $B$}

Suppose that for $s_\nu\to q$ and $\delta_\nu\downarrow0$ the
rescalings of $\mu_{1,\nu}$ under $y=(x-s_\nu)/\delta_\nu$ converge
narrowly to
\begin{equation}\label{4.28}
 m_B(y)dA_y=\frac{16|y|^2}{(1+|y^2-B|^2)^2}\dif A_y,
 \qquad B\in\C,
\end{equation}
and the full local mass $4\pi$ of the second component is carried by
a radial ordinary Liouville bubble centered at $s_\nu+o(\delta_\nu)$
with scale $o(\delta_\nu)$. We identify the real gradient of
$\mathcal M$ with the complex number
$\nabla\mathcal M=\mathcal M_{B_1}+i\mathcal M_{B_2}$.

\begin{proposition}
\label{pr4.17}
Assume the hypotheses of Theorem~\ref{th4.16}, the narrow convergence to the profile
in \eqref{4.28}, and the stated inner convergence. Then
\begin{equation}\label{4.29}
 4\pi\mathcal T(q)+\kappa_1^c(q)\overline{\nabla\mathcal M(B)}=0.
\end{equation}
The map $\nabla\mathcal M:\C\to B_\pi(0)$ is an
orientation-preserving diffeomorphism. If $|\kappa_1^c(q)|>4|\mathcal T(q)|$,
then \eqref{4.29} has a unique solution $B$, and
$$
 F_B(B)=\kappa_1^c(q)\nabla\mathcal M(B)+4\pi\overline{\mathcal T(q)},\qquad
 \det D_BF_B=(\kappa_1^c(q))^2\det D^2\mathcal M(B)>0.
$$
\end{proposition}
\begin{proof}
At the outer scale, the second measure converges to $4\pi\delta_0$,
and its product with $m_B\dif A$ gives zero mass to the diagonal.
At the inner scale, the second profile is radial. Thus
$$
 S_{12}=4\pi\int_\C m_B(y)\frac{\bar y}{y}\dif A_y,
 \qquad S_{22}=0.
$$
The substitution $w=y^2$ gives
$$
 \begin{aligned}
 \int_\C m_B(y)\frac{\bar y}{y}\dif A_y
 &=8\int_\C\frac{\bar w/|w|}{(1+|w-B|^2)^2}\dif A_w\\
 &=8\,\overline{\nabla\mathcal M(B)},
 \qquad
 S_{12}=32\pi\,\overline{\nabla\mathcal M(B)}.
 \end{aligned}
$$
The second equality follows by differentiating
$\mathcal M(B)=\int_\C|B+v|(1+|v|^2)^{-2}\dif A_v$.

Set $A_B=-\partial_y\log(1+|y^2-B|^2)$ and $b'=-1/y$.
Then
$$
 A_B'=A_B/y+A_B^2=A_B^2-A_Bb',\qquad
 A_B(y)=-2/y+O(|y|^{-3})\quad (|y|\to\infty).
$$
For a radial cutoff $\chi_R(y)=\chi(|y|/R)$ equal to one near zero,
\[
 \lim_{R\to\infty}\int_\C\chi_R(A_B^2-A_Bb')\dif A
 =-\lim_{R\to\infty}\int_\C(\chi_R)_y A_B\dif A=0.
\]
The term $-2/y$ integrates to zero against $(\chi_R)_y$ by angular
integration, and the remainder contributes $O(R^{-2})$.

The Cauchy representations of $A_B$ and $b'$ use the factor
$-1/(4\pi)$ with $m_B\dif A$ and $4\pi\delta_0$, respectively.
For $A_B$, this follows from $\partial_{\bar y}A_B=-m_B/4$
and the decay of the field at infinity.

By radial symmetry and $\partial_{\bar y}\cC\chi_R=\chi_R$,
$\cC\chi_R(y)=\bar y$ on the disk where $\chi_R=1$. The divided
difference therefore tends to $\Theta(x-y)$ for $x\ne y$ and is
uniformly bounded by scaling.

By \eqref{4.22} and dominated convergence, the limits of the two
product integrals are $-S_{11}/(16\pi)$ and $-S_{12}/(16\pi)$,
respectively. Thus $S_{11}=S_{12}$, and \eqref{4.27} gives
\eqref{4.29}.

The remaining assertions follow from
Lemma~\ref{le4.7} and
Corollary~\ref{co4.8}, with
$$
 \kappa=\kappa_1^c(q),\qquad
 \tau=4\pi\overline{\mathcal T(q)}.
$$
\end{proof}

\subsubsection{The local blow-up profiles}

We give the local classification used here. It yields the two
blow-up patterns described in \cite[Section~3, Remark~3.1]{WWZ} from
the local hypotheses below. Let $q$ be isolated and suppose the
local masses are $(8\pi,4\pi)$, with smooth bounded forcing,
bounded total masses, fixed boundary oscillation, and local
upper bounds for the logarithmic densities off $q$. After
subtracting a smooth Dirichlet Poisson potential as in the
proof of Proposition~\ref{pr4.19}, the regular local selection
theorem \cite[Proposition~3A]{LWYZ18} applies: the selected
rescaled disks have scalar or Toda limits, the core widths
are negligible compared to distances between their centers,
and the log-density plus twice the log-distance to the
selected set is bounded above. A fully bubbling regular
Toda limit has masses $(8\pi,8\pi)$ by the regular $A_2$ case
of \cite[Theorem~1.3(i)]{LWY}, so every selected core here is
scalar, of mass $4\pi$ by \eqref{eq:scalar-liouville}.
We first record the continuation argument for an isolated core.
Let $x_\nu$ be a selected center and let $R_\nu$ be an isolation
radius with $B_{2R_\nu}(x_\nu)$ in the working chart and with no other
selected point in that disk. Choose a core radius
$l_\nu$, enlarged slowly, at which the mass pair tends to that of the ordinary bubble
and both components have fast decay on the circle. This is possible
by the rescaled scalar convergence and a diagonal choice. Write
\[
 M_i(r)=\int_{B_r(x_\nu)}F_{i,\nu}\dif A,\qquad
 H_i(r)=\overline{U_{i,\nu}}(r)+2\log r ,
\]
where the bar denotes the circular mean about $x_\nu$. The selection
bound gives $F_{i,\nu}\leq C/r^2$ on balls of radius $r/4$ about
points on each circle of radius $r\leq R_\nu$. Split the Green
gradient integral into the ball and its complement. The pointwise
bound on the ball and the total mass off it give an $O(1/r)$
gradient bound on the circle. Thus each $U_{i,\nu}$ has uniformly
bounded oscillation on these circles.

Suppose a fixed positive mass increment occurs at a radius
$r_\nu=o(R_\nu)$ beyond the ordinary core. Choose $l_\nu=o(r_\nu)$
as above, and let $s_\nu\leq r_\nu$ be the first radius at which
the total increment beyond $B_{l_\nu}(x_\nu)$ reaches a small
fixed $\varepsilon_0>0$. The flux identity is
\[
 \frac{\dif H_i}{\dif\log r}
 =2-\frac{2M_i(r)-M_{3-i}(r)}{2\pi}+O(r^2).
\]
For $l_\nu\leq r\leq s_\nu$, this slope is at most $-1$ in the
component $i$ carrying the ordinary bubble, and at least $3$
in the other component $j$, for small $\varepsilon_0$. Both
$H_i(l_\nu)$ and $H_j(l_\nu)$ tend to $-\infty$. Bounded slopes
and bounded oscillation imply $s_\nu/l_\nu\to\infty$, since a
bounded ratio would give $o(1)$ increment in both components.
Integration of the upper bound on $H_i(r)$ now gives a mass
increment of $o(1)$ in component $i$.
The lower bound on the slope of $H_j$ and the oscillation bound
give
\[
 M_j(s_\nu)-M_j(l_\nu)\leq C e^{H_j(s_\nu)},\qquad
 M_j(\eta s_\nu)-M_j(l_\nu)\leq C\eta^3e^{H_j(s_\nu)}
 \quad (0<\eta<1)
\]
for each fixed $\eta$ and all large $\nu$. Thus $H_j(s_\nu)$ is
bounded below, and the selection bound gives an upper bound.

Rescale at $s_\nu$. The slopes are uniformly bounded on each fixed
annulus in the rescaled plane by the total mass bound, and the
oscillation bound holds there. The logarithmic density of component $i$ tends to $-\infty$
on these annuli, while component $j$ has a finite limit
$W$ by local elliptic estimates. The limit is defined on the
punctured plane since $R_\nu/s_\nu\to\infty$.
The measure in component $i$ tends to $4\pi\delta_0$, and the
$\eta^3$ bound rules out a point mass in component $j$ at zero.
The limiting flux gives $r\overline W'(r)=2+O(r^3)$ near zero,
and bounded oscillation gives $W=2\log r+O(1)$. Subtracting
$2\log r$ and a bounded local Poisson potential leaves
a bounded harmonic function with a removable singularity. Consequently
$\Delta W+2e^W=4\pi\delta_0$ on $\C$. The mass bound persists,
and the integer-source Liouville classification recalled below
gives $\int_\C e^W\dif A=8\pi$.

Suppose there are at least two selected centers tending to $q$.
Choose a closest pair and rescale about one of its points by the
minimum separation $d_\nu$. After passing to a subsequence, at least
two selected centers remain at finite distance and have distinct limits.
Any other finite limit of a selected center is separated from
them by at least one. The selected widths divided by $d_\nu$
tend to zero for the centers at finite distance. The local
upper bounds away from those limits follow from the selection
estimate. Green's formula gives uniformly bounded rescaled
gradients on compact sets away from them. To see this, split
the singular kernel into a small disk on which the rescaled
density is bounded, a compact region away from its
singularity, and a tail bounded by the total mass divided
by the tail radius. In particular, after passing to a subsequence, each rescaled
log-density either tends locally to $-\infty$ off the limiting centers
or has a finite limit there. In the finite case, local elliptic
estimates give convergence to a solution of the limiting equation.

At each of the finite centers, apply the regular local
quantization leading to \eqref{2.3} on a small fixed rescaled
disk. Its smooth forcing is $O(d_\nu^2)$, its local mass is
bounded, its boundary oscillation is bounded by the gradient
estimate, and there are no other blow-up points in that disk
by the selection bound. Each atom contains at least the mass
of its selected core. The total mass on fixed rescaled disks
is at most $(8\pi,4\pi)$, since their original radii shrink
to $q$. There is at least one other nonzero atom. The list
\eqref{2.3} therefore forces each atom to be pure, of mass
$(4\pi,0)$ or $(0,4\pi)$.

A component with a pure point mass cannot have a finite rescaled
limit off the centers. Indeed, its circular mean about the
concentration point is $4\log(1/r)+o(\log(1/r))$ up to bounded terms.
The other component has no point mass there. Its potential has
lower logarithmic order by the circle-average argument following
Theorem~\ref{th1.1}. Jensen's inequality rules
out integrability of the component with the point mass. If only first-component
atoms are present, there must be two of them and they exhaust
its mass. A finite second-component limit then satisfies the
equation in Remark~\ref{re4.18}, which contradicts its mass
bound. An arrangement with only second-component atoms is
impossible by its mass bound. Thus both densities vanish
locally off the finite selected centers.

Apply the radial Pohozaev identity on a fixed rescaled circle
enclosing these centers. Take the sequence limit first and then
let the radius tend to infinity. The density boundary term tends
to zero for each fixed circle. The forcing tends to zero on
fixed disks. The limiting gradient is the gradient of the
logarithmic potentials of the finite atoms, by the same
kernel splitting and the vanishing of both densities off
the centers. If $M=(M_1,M_2)$ is the sum of their masses,
the boundary gradients are $-(AM)_i/(2\pi r)+O(r^{-2})$ in
the radial direction and $O(r^{-2})$ in the tangential
direction. In the identity obtained by pairing the equations
against $A^{-1}(y\cdot\nabla U)$, the limiting flux is
$M^TAM/(4\pi)$ and the density term is $2(M_1+M_2)$.
Thus $M_1^2-M_1M_2+M_2^2=4\pi(M_1+M_2)$. The only sums of
at least two pure atoms within the mass cap are
$(8\pi,0)$, $(4\pi,4\pi)$, and $(8\pi,4\pi)$. Only the last
satisfies this equation. There are exactly three centers. Their distances are comparable,
and they carry the full local masses.

At their common separation scale, write $y_+,y_-$ for the two
limiting first-component centers and $y_0$ for the second.
The logarithmic coefficient vector at $y_+$ is
$A(4\pi,0)^T/(2\pi)$. In the translation Pohozaev identity,
the boundary integral is therefore $-4\pi$ times the gradient
of the regular first field. This gives
\[
 -4\frac{y_+-y_-}{|y_+-y_-|^2}
 +2\frac{y_+-y_0}{|y_+-y_0|^2}=0.
\]
Thus $y_0=(y_++y_-)/2$, as required in Corollary~\ref{co4.11}.

If there is only one selected disk, use a fixed small isolation
disk. The full mass pair $(8\pi,4\pi)$ gives a fixed positive mass
increment beyond the ordinary core at radii tending to zero.
For this use of the continuation argument, choose the core
radius $l_\nu$ by a diagonal argument. Take a larger radius, still
tending to zero, that captures a fixed positive increment, as
ensured by the local mass pair. The selection bound holds
out to the fixed isolation radius. The continuation argument
therefore applies.
The ordinary core cannot belong to the first component, since the
outer singular bubble would require $8\pi$ in the second
component. Hence the second component has the ordinary inner
bubble, and the first component at the larger scale has a
finite-mass limit $W$ satisfying $\Delta W+2e^W=4\pi\delta_0$
on $\C$. The profiles have separated widths since the ordinary
width is $o(l_\nu)$ and $l_\nu=o(s_\nu)$ in the continuation
argument. They account for the full local masses.

\begin{remark}\label{re4.18}
Suppose the first rescaled measure converges to
$4\pi(\delta_{p_+}+\delta_{p_-})$ with $p_+\ne p_-$ and the second
converges to $e^V\dif A$. Then
$$
 \Delta V+2e^V=4\pi(\delta_{p_+}+\delta_{p_-})\quad\text{in }\C,
 \qquad \int_\C e^V\dif A\leq4\pi.
$$
For $R\geq R_0>\max|p_\pm|$, integration gives
$$
 R\overline V'(R)=4-\frac1\pi\int_{B_R}e^V\dif A\geq0.
$$
By Jensen's inequality,
$$
 \int_{\C\setminus B_{R_0}}e^V\dif A
 \geq2\pi\int_{R_0}^\infty R e^{\overline V(R)}\,dR
 \geq2\pi e^{\overline V(R_0)}\int_{R_0}^\infty R\,dR=\infty,
$$
a contradiction.
\end{remark}

\paragraph{Local selection with a prescribed reference center.}
\begin{lemma*}
Let $D$ be a disk. Let $F_{i,\nu}=e^{U_{i,\nu}}$, $i=1,2$,
be smooth positive coordinate densities on $D$ with uniformly
bounded integrals. Assume
$\Delta_0U_{i,\nu}+\sum_j a_{ij}F_{j,\nu}=\kappa_{i,\nu}^c$
there, with $\kappa_{i,\nu}^c$ uniformly bounded and $A$ the
regular $A_2$ Cartan matrix. Suppose $x_\nu,q_\nu\to q\in D$,
$r_\nu=|x_\nu-q_\nu|>0$, and
$H_\nu=r_\nu^2\max_iF_{i,\nu}(x_\nu)\to\infty$.
Then, on a subsequence, there are
$\widehat x_\nu=x_\nu+o(r_\nu)$ and $s_\nu=o(r_\nu)$ such that
$B_{s_\nu}(\widehat x_\nu)\subset D$ and
\begin{equation}\tag{4.29a}\label{eq:prescribed-center-selection}
 \liminf_{\nu\to\infty}
 \int_{B_{s_\nu}(\widehat x_\nu)}
       (F_{1,\nu}+F_{2,\nu})\dif A\geq4\pi .
\end{equation}
The reference center $q_\nu$ is retained.
\end{lemma*}
\begin{proof}
Put $h_\nu=r_\nu H_\nu^{-1/4}$ and
$g_\nu=\max_iF_{i,\nu}$. On $\overline{B_{h_\nu}(x_\nu)}$, let
$(h_\nu-|x-x_\nu|)^2g_\nu(x)$ attain its maximum at $\widehat x_\nu$.
The value at $x_\nu$ is $H_\nu^{1/2}$, so the maximizer is
interior. Write $d_\nu=h_\nu-|\widehat x_\nu-x_\nu|$ and
$\lambda_\nu=g_\nu(\widehat x_\nu)^{-1/2}$.
Then $d_\nu/\lambda_\nu\ge H_\nu^{1/4}\to\infty$ and
$\lambda_\nu\to0$. On $B_{d_\nu/2}(\widehat x_\nu)$ the
distance to the boundary of $B_{h_\nu}(x_\nu)$ is at least
$d_\nu/2$, so the maximality gives $g_\nu\le4\lambda_\nu^{-2}$.
Set
\[
 W_{i,\nu}(y)=U_{i,\nu}(\widehat x_\nu+\lambda_\nu y)
                                    +2\log\lambda_\nu.
\]
On expanding disks $W_{i,\nu}\le\log 4$, and at least one
component is zero at the origin, with its index fixed along a
subsequence. The scaled equation is
\[
 \Delta_0W_{i,\nu}+\sum_j a_{ij}e^{W_{j,\nu}}
   =\lambda_\nu^2\kappa_{i,\nu}^c(\widehat x_\nu+\lambda_\nu y).
\]
On any fixed ball $B_{2R}$, the right-hand side tends to zero in
$L^\infty$ and the exponentials are uniformly bounded. The
zero-Dirichlet potential of $\Delta_0W_{i,\nu}$ on this ball is
uniformly bounded by comparison with the quadratic disk barrier.
Subtract it to get a harmonic function bounded above. Applying
Harnack's inequality to a constant minus that harmonic function shows on $B_R$
that a component with $W_{i,\nu}(0)$ bounded below is bounded on
both sides. If $W_{i,\nu}(0)\to-\infty$, the same Harnack
inequality in the other direction gives $W_{i,\nu}\to-\infty$
uniformly on $B_R$. Pass to a subsequence to make this choice
for each component. Interior $W^{2,p}$ estimates for any finite
$p$ and a diagonal subsequence give local $C^1$ convergence for
the bounded components. The limiting equations are smooth by
elliptic regularity. Change of variables and Fatou's lemma give
finite mass on the plane for every surviving component.

With one surviving component, the equation is
$\Delta_0W+2e^W=0$, so \eqref{eq:scalar-liouville} with
$\gamma=0$ gives mass $4\pi$. With two surviving components,
\cite[Theorem~1.3(i)]{LWY} applies with $n=2$ and
$\gamma_1=\gamma_2=0$. The entire finite-mass limit satisfies
$A(\int e^{W_1},\int e^{W_2})^T=(8\pi,8\pi)^T$, so each
mass is $8\pi$. Choose $R_\nu\to\infty$ sufficiently slowly that
$R_\nu\lambda_\nu<d_\nu/2$. For each fixed $R$, the enlarged
ball contains the scaled $B_R$ for large $\nu$, and the mass on
the latter converges to the limiting mass there. Taking the
liminf and then $R\to\infty$ proves
\eqref{eq:prescribed-center-selection} with
$s_\nu=R_\nu\lambda_\nu$. Both the displacement and $s_\nu$
are bounded by $h_\nu=o(r_\nu)$.
\end{proof}

\begin{proposition}
\label{pr4.19}
Assume \eqref{4.7}, \eqref{4.8} and
\eqref{4.11} at an isolated mixed point $q$, and suppose
\begin{equation}\label{4.30}
 \kappa_1^c(q)>4|\mathcal T(q)|.
\end{equation}
After passing to a subsequence, there are $q_\nu\to q$ and scales
$a_\nu=o(\delta_\nu)$ with $\delta_\nu\to0$ such that
\begin{align}
 a_\nu^2 F_{2,\nu}(q_\nu+a_\nu y)&\longrightarrow
                \frac4{(1+|y|^2)^2},\label{4.31}\\
 \delta_\nu^2F_{1,\nu}(q_\nu+\delta_\nu y)&\longrightarrow
                \frac{16|y|^2}{(1+|y^2-B_*|^2)^2}.
                \label{4.32}
\end{align}
The convergence in \eqref{4.31} is locally smooth on $\C$ and
also holds for the logarithms. In \eqref{4.32}, the convergence
is locally smooth on $\C\setminus\{0\}$. Both profiles carry the
full local masses at $q$, and $B_*$ is uniquely determined by
\begin{equation}\label{4.33}
 \kappa_1^c(q)\nabla\mathcal M(B_*)
                       +4\pi\overline{\mathcal T(q)}=0.
\end{equation}
For some fixed $r_0>0$,
\begin{equation}\label{4.34}
 r^2F_{i,\nu}(q_\nu+re^{i\theta})\leq C,\qquad
 \osc_{|z-q_\nu|=r}\log F_{i,\nu}\leq C,
 \quad 2a_\nu\leq r\leq r_0.
\end{equation}
The bounds are uniform near any fixed finite set of limiting
solutions satisfying \eqref{4.30}.
\end{proposition}

\begin{proof}
To remove the smooth forcing, take a fixed disk $D$ about $q$ and set
$$
 \Delta\psi_{i,\nu}=\kappa_{i,\nu}^c\quad\text{in }D,
 \qquad \psi_{i,\nu}=0\quad\text{on }\partial D,
 \qquad V_{i,\nu}=U_{i,\nu}-\psi_{i,\nu}.
$$
Then
$$
 \begin{gathered}
 \Delta V_{i,\nu}+\sum_j a_{ij}e^{\psi_{j,\nu}}e^{V_{j,\nu}}=0,
 \qquad c\leq e^{\psi_{i,\nu}}\leq C,\quad
 \|e^{\psi_{i,\nu}}\|_{C^3}\leq C,\\
 \sum_i\int_D e^{\psi_{i,\nu}}e^{V_{i,\nu}}\dif A\leq C,
 \qquad \osc_{\partial D}V_{i,\nu}\leq C.
 \end{gathered}
$$
On compact annuli in $D\setminus\{q\}$, the oscillation bound from
\eqref{4.11} and the mass bound give a uniform upper bound for
$V_{i,\nu}$. These are the hypotheses of the local classification
just proved. The smooth coefficient bounds and boundary
oscillation are uniform on the fixed coordinate disk.

In the split case, selection gives three regular radial profiles
of mass $4\pi$, which account for the full local masses in
\eqref{4.7}. By Corollary~\ref{co4.11},
$|\kappa_1^c(q)|=4|\mathcal T(q)|$, contrary to \eqref{4.30}.

In the remaining blow-up pattern, the inner component has a regular
finite-mass scalar limit. The classification
\eqref{eq:scalar-liouville} at $\gamma=0$ gives its radial
profile. The strict maximum of the limiting profile gives a
local maximum $q_\nu$ of $F_{2,\nu}$ within $o(a_\nu)$ of the
initial inner center. Setting $a_\nu=2F_{2,\nu}(q_\nu)^{-1/2}$
fixes the rescaled maximum to $4$, and \eqref{4.31} follows.

For the outer finite-mass limit $W$, the inner point mass gives
$\Delta W+2e^W=4\pi\delta_0$. Use
\eqref{eq:scalar-liouville} with $\gamma=1$. In the notation
$s=\lambda^2>0$ and $p=\gamma_0\in\C$, it gives
\[
 e^W(y)=\frac{16\lambda^2|y|^2}
                   {(\lambda^2+|y^2-\gamma_0|^2)^2}.
\]
Rescaling the initial outer scale by $\sqrt\lambda$ gives
\eqref{4.32} with $B_*=\gamma_0/\lambda$. At the inner scale,
the selection also gives $a_\nu^2 F_{1,\nu}(q_\nu+a_\nu y)\to0$
locally uniformly. Integrating the two profiles yields
$$
 \lim_{R\to\infty}\lim_{\nu\to\infty}
       \int_{B_{R\delta_\nu}(q_\nu)}F_{1,\nu}\dif A=8\pi,
 \qquad
 \lim_{R\to\infty}\lim_{\nu\to\infty}
       \int_{B_{Ra_\nu}(q_\nu)}F_{2,\nu}\dif A=4\pi.
$$
These limits equal the full local masses at $q$. Hence, if we first restrict the measures to a sufficiently slowly
shrinking disk and then rescale them, they converge narrowly. Proposition~\ref{pr4.17} gives
\eqref{4.33}. By Lemma~\ref{le4.7}, $g(r)=|\nabla\mathcal M(r)|$
is strictly increasing from $0$ to $\pi$. Thus, over a finite set of
limiting solutions,
$$
 \vartheta:=\max_q\frac{4|\mathcal T(q)|}{\kappa_1^c(q)}<1,
 \qquad |B_*|\leq g^{-1}(\pi\vartheta).
$$

Choose $r_0$ so that the disks about $q_\nu$ of radius $2r_0$
are in the isolating chart. Suppose the pointwise bound fails
on those disks at $x_\nu$, with $r_\nu=|x_\nu-q_\nu|$ and
$r_\nu^2\max_iF_{i,\nu}(x_\nu)\to\infty$. The off-center
bounds from \eqref{4.11} and the mass bound force $r_\nu\to0$.
Apply the local selection lemma \eqref{eq:prescribed-center-selection}.
Its equation and bounds follow from \eqref{4.8} and the local
mass and forcing bounds. The lemma gives a small ball $E_\nu$ of radius $s_\nu=o(r_\nu)$.
Its center differs from $x_\nu$ by $o(r_\nu)$, and its mass has
limit inferior at least $4\pi$.

For the exhaustion argument, let $P_1,P_2$ be the density limits in \eqref{4.32}, \eqref{4.31}. For fixed $R>1$, define the cores
\[
 C_{1,\nu}(R)=\{\delta_\nu/R<|z-q_\nu|<R\delta_\nu\},
 \qquad C_{2,\nu}(R)=B_{Ra_\nu}(q_\nu).
\]
Their masses tend to the integrals of $P_i$ on the corresponding
fixed scaled regions, by the two local convergences. As $R\to\infty$
these integrals tend to $M_1=8\pi$ and $M_2=4\pi$, respectively.
Whenever $E_\nu$ avoids $C_{i,\nu}(R)$ for every fixed $R$ and
all sufficiently large $\nu$, fix a small disk $B_\eta(q)$ in the
chart. Both sets lie inside it for large $\nu$. Positivity and
\eqref{4.7} give
\[
 \limsup_{\nu\to\infty}\int_{E_\nu}F_{i,\nu}\dif A
 \leq M_i+\int_{B_\eta(q)}f_{i,*}\dV
      -\lim_{\nu\to\infty}\int_{C_{i,\nu}(R)}F_{i,\nu}\dif A ,
 \qquad f_{1,*}=0,\quad f_{2,*}=f_* .
\]
The limiting measure assigns no mass to the boundary of the fixed disk.
First let $R\to\infty$ and then let $\eta\downarrow0$. The contribution
of the residual density vanishes. Thus the $i$th mass on $E_\nu$ tends to
zero in this case.

Pass to a subsequence on which $r_\nu/a_\nu$ and $r_\nu/\delta_\nu$
have limits in $[0,\infty]$. Since $a_\nu=o(\delta_\nu)$, the
following cases are exhaustive.
\begin{enumerate}[label=\textup{(\roman*)}]
\item If $r_\nu=O(a_\nu)$, both inner-scale densities are bounded
on a fixed compact set containing the scaled $E_\nu$. This follows
from \eqref{4.31} and the locally uniform convergence of the first
density to zero on that scale. Its area is $O((s_\nu/a_\nu)^2)=o(1)$, so
both masses vanish.
\item If $a_\nu\ll r_\nu\ll\delta_\nu$, the ball avoids the inner
core of the second density and lies inside the hole of every fixed
outer annular core of the first density. Both masses vanish by
the exhaustion argument.
\item If $r_\nu\asymp\delta_\nu$, the scaled $E_\nu$ stays on
a fixed compact annulus away from zero. The first density is
bounded there by \eqref{4.32}. Its mass is
$O((s_\nu/\delta_\nu)^2)=o(1)$. The second mass vanishes by
inner-core exhaustion, since $r_\nu/a_\nu\to\infty$.
\item If $\delta_\nu\ll r_\nu\to0$, the ball is outside both
cores for every fixed $R$. Again both masses vanish by exhaustion.
\end{enumerate}
This contradicts \eqref{eq:prescribed-center-selection}, proving
the pointwise bound up to radius $2r_0$.
The estimate on circles from Green's representation gives the
oscillation bound. For a point $z$ on the circle of radius $r$, the
disk of radius $r/4$ about $z$ stays in the region
where $F_{i,\nu}\le C/r^2$. Its contribution to the gradient
is $O(r^{-1})$, by integration of the $1/|z-x|$ kernel.
Outside this disk the gradient is bounded by $C/r$ times the
total mass. The smooth kernel contributes $O(1)$. Integration
along the circle gives a uniform oscillation bound.
The contradiction argument and $|B_*|\leq g^{-1}(\pi\vartheta)$
give uniform constants over the finite set of limiting solutions.
\end{proof}

\section{Local completion and matching estimates}
\label{rev:sec:linear}\label{sec:inv-chapter}\label{rev:sec:matching}

We use the coordinate measures of Section~\ref{sec:notation}
and the single-point width notation of Section~\ref{sec:scale-conventions}.

The mixed profiles proved in Section~\ref{rev:sec:core} provide
the local inputs for the estimates below. We retain the annular
bounds and the normalization needed to apply these estimates to
actual configurations in Section~\ref{sec:typedactual-chapter}.

\subsection{Local profiles and density estimates}
Let $M$ be a closed orientable surface of area one. In an isothermal
coordinate disk about $p$, write $g=e^\phi|dx|^2$ with $\phi(p)=0$.
For a Euclidean density $m$, write
$$
 m_g=e^{-\phi}m,\qquad m_g\dV=m\dif A_x.
$$
We use fixed cutoffs supported in this disk.

Let $B$ and $\beta$ range in compact subsets of $\C$ and $(0,\infty)$,
respectively. Use the single-point widths $\delta$ and $a$ from
\eqref{7.7}, with $\delta\to0$. Set
\[
 m_o(x)=\frac{16\delta^4|x|^2}
 {(\delta^4+|x^2-\delta^2B|^2)^2},\qquad
 m_i(x)=\frac{4a^2}{(a^2+|x|^2)^2}.
\]
The masses of $m_o$ and $m_i$ on $\C$ are $8\pi$ and $4\pi$, respectively.
Denote the profiles in the variables $x=\delta y$ and $x=ay$ by
$m_{o,B}(y)$ and $m_i^0(y)$. The finite-energy kernels of
$-\Delta-2m_{o,B}$ and $-\Delta-2m_i^0$ each have dimension three.
A basis for the inner kernel is
$$
 Z_{i0}=\frac{1-|y|^2}{1+|y|^2},\quad
 Z_{i1}=\frac{y_1}{1+|y|^2},\quad
 Z_{i2}=\frac{y_2}{1+|y|^2}.
$$
We obtain a basis for the outer kernel by differentiating with
respect to $\lambda,\Re B,\Im B$ in
$$
 \log\frac{32\lambda^2}
 {(\lambda^2+|y^2-B|^2)^2},\qquad \lambda=1.
$$
Each basis function has zero mean with respect to its profile density.
For the bounded outer kernel, use \cite[Theorem~1]{DEM} with
$N=1$. In their notation the entire density solving the
Liouville equation with unit source is $2m_{o,B}$, which
has finite mass and produces the linearization
$\Delta+2m_{o,B}$. The theorem identifies its entire bounded
solutions with the three parameter derivatives above, for
each $B\in\C$.
For the inner kernel, use the spherical calculation in the proof
of Lemma~\ref{le5.12}. The coefficients
extend smoothly across zero after inversion, with a quadratic zero
for the outer coefficient. Finite energy gives an $H^1$ extension
across zero, and elliptic regularity gives boundedness.

\begin{lemma}\label{le5.1}
Take $\eps=\delta^2$ and $b=\delta^2B$ in \eqref{4.1}--\eqref{4.3}, so that
$$
 F=\delta^4+\beta\delta^4r^2+|z^2-\delta^2B|^2,\quad
 Q=\delta^4(\beta\delta^4+4r^2+\beta|z^2+\delta^2B|^2),
$$
$$
 n_1=4Q/F^2,\qquad n_2=16\beta\delta^8 F/Q^2.
$$
On a fixed small disk, for each $1<p<2$ the following estimates hold
uniformly on the stated compact parameter sets.
$$
 \|n_1-m_o\|_{L^p}=O(\delta^{2/p}),
$$
$$
 \left\|n_2-m_i-
       \frac{16\beta}{(4+\beta r^2)^2}\right\|_{L^p}
 =O(\delta^{4/p-2}+\delta^2).
$$
Both errors are $o(|\log\delta|^{-1})$. Moreover,
$n_1\ge c m_o$, $n_2\ge c m_i$, $\|n_1\|_\infty\le C\delta^{-2}$,
and $\|n_2\|_\infty\le C\delta^{-4}$. On each fixed annulus,
$n_2$ is bounded below by a positive constant.
\end{lemma}
\begin{proof}
Set $F_0=\delta^4+|z^2-\delta^2B|^2$. Uniformly for bounded $B$,
$F_0\asymp\delta^4+r^4$ and
$\beta\delta^4r^2/F_0=O(\delta^2)$. Thus
$$
 |n_1-m_o|\le C\delta^2m_o+
                C\frac{\delta^4}{\delta^4+r^4}.
$$
Both terms have $L^p$ norm $O(\delta^{2/p})$.
For $n_2$, write $A_B=1+|B|^2$,
$$
 D_0=\beta\delta^4 A_B+4r^2,\qquad
 D=D_0+\beta r^4+2\beta\delta^2\Re(z^2\overline B).
$$
Then
$$
 n_2=16\beta\frac{\delta^4A_B+r^4
        -2\delta^2\Re(z^2\overline B)+\beta\delta^4r^2}{D^2},
 \qquad m_i=\frac{16\beta\delta^4 A_B}{D_0^2}.
$$
For small $\delta$, $D\asymp\delta^4+r^2$ on the fixed disk.
The absolute value of $n_2-m_i-16\beta/(4+\beta r^2)^2$ is bounded by
$$
 C\left((\delta^2+r^2)m_i+
 \min(1,\delta^4/r^2)+\delta^2+
 \frac{\delta^2r^2}{(\delta^4+r^2)^2}\right).
$$
The first and last terms are $O(\delta^{4/p-2})$ in $L^p$, and the
remaining terms are of smaller order. The bounds on $n_1,n_2$
follow from $F_0\asymp\delta^4+r^4$ and $D\asymp\delta^4+r^2$.
\end{proof}

\begin{remark}\label{re5.2}
Multiplication by a smooth positive factor $1+O(r^2)$ preserves the
profile estimates after completing $16\beta/(4+\beta r^2)^2$ to a smooth
global density. Indeed, $\|r^2m_o\|_p=O(\delta^{2/p})$ and
$\|r^2m_i\|_p=O(\delta^{4/p})$. Fixed cutoffs and the change to metric
densities satisfy the same bounds. Proposition~\ref{pr5.3} estimates
the nonlinear residual of the global approximation.
\end{remark}

\subsection{The global approximation and its residual}

Fix $r>4\pi$, put $\sigma=r-4\pi$, and suppose the mean-zero
residual scalar solution $z_*$ satisfies
$$
 -\Delta_g z_*=2\sigma\left(\frac{h_2e^{z_*}}{I}-1\right),
 \qquad I=\int_Mh_2e^{z_*}\dV.
$$
Write $n_*^g=\sigma h_2e^{z_*}/I$ for the density with respect to $\dV$.
Assume that the scalar linearization at $z_*$ has trivial kernel.
For $q$ in a fixed coordinate neighborhood, let $H$ denote the
regular part from Section~\ref{sec:notation} in the chart centered
at $q$ and set
$$
 G_q(x)=-\frac1{2\pi}\log|x|+H(x,q),\qquad
 \psi_1=\log h_1+\phi+12\pi H(\,\cdot\,,q)-\frac12z_*,
 \quad \psi_2=\log h_2+\phi+z_* .
$$
Choose the isothermal coordinate smoothly in $q$, with $x(q)=0$,
$\phi(0)=0$, and $\partial_z\psi_1(0)=0$. As in the coordinate
change in Section~\ref{typed-sec:rotation},
$z=\zeta+a\zeta^2+O(\zeta^3)$ gives
$\psi_{1,\zeta}=\psi_{1,z}-a$ and $\psi_{2,\zeta}=\psi_{2,z}+2a$
at the origin. In the resulting coordinate,
$$
 \nabla\psi_2(0)=\nabla V(q),\qquad
 V=2\log h_1+\log h_2+12\pi r_g .
$$
Indeed, $r_g(q)=H(q,q)+\phi(q)/(4\pi)$ and symmetry of $G$ give
$\nabla V=2\nabla\psi_1+\nabla\psi_2$ at $q$. We construct the
approximation for all $q$ in this neighborhood.

Set $\eps=\delta^2$ and $\beta=e^{\phi(0)}n_*^g(q)$, and take
$U_i=\log n_i$ as in Lemma~\ref{le5.1}, with $b=\eps B$. In the
chosen chart, put
$$
 L_1=2\log\eps+\log16-6\log|x|
                  +\log(1+\beta|x|^2/4),\qquad
 L_2=\log\beta-2\log(1+\beta|x|^2/4).
$$
Let $\chi$ be radial, with $\chi=1$ on $|x|\le r_0$ and $\chi=0$
on $|x|\ge2r_0$. Define the global approximation by
\begin{equation}
 v^\delta=
 \binom{12\pi G_q-z_*/2}{z_*}
       +\chi\binom{U_1-L_1}{U_2-L_2},
 \qquad
 u_i^{\rm app}=v_i^\delta-\int_Mv_i^\delta\dV .
 \label{5.1}
\end{equation}
The logarithmic terms cancel at $q$, so $v^\delta$ and $u^{\rm app}$
are smooth on $M$.

\begin{proposition}\label{pr5.3}
For $B$ in a fixed compact set, $q$ in the chosen neighborhood,
and $|\mu|\le C_0\delta^2$, set
$\rho=(8\pi+\mu,r)$ and
$$
 V_i^\delta=\rho_i
       \frac{h_i e^{u_i^{\rm app}}}
            {\int_Mh_i e^{u_i^{\rm app}}\dV}.
$$
Let $m_{o,g}=e^{-\phi}m_o$ and $m_{i,g}=e^{-\phi}m_i$, with the
fixed cutoffs understood.
Write $V_1^\delta=m_{o,g}+e_{1,\delta}$ and
$V_2^\delta=m_{i,g}+n_*^g+e_{2,\delta}$. Then
\begin{equation}
 \|e_{1,\delta}\|_{L^p}=O(\delta^{2/p}),\qquad
 \|e_{2,\delta}\|_{L^p}
             =O(\delta^{4/p-2}),\qquad 1<p<2 .
 \label{5.2}
\end{equation}
The normalized residual satisfies
\begin{equation}
 \left\|-\Delta_g u^{\rm app}
          -A\binom{V_1^\delta-\rho_1}{V_2^\delta-\rho_2}
       \right\|_{H^{-1}}
       \le C_p\delta^{4/p-2}.
 \label{5.3}
\end{equation}
For $p=6/5$, the bound is $O(\delta^{4/3})$.
\end{proposition}
\begin{proof}
\emph{Step 1. The residual before normalization.} Set
$$
 \lambda_1=2\log\eps+\log16-\psi_1(0),\qquad
 \lambda_2=\log(\sigma/I),\qquad
 N_i^g=e^{\lambda_i}h_i e^{v_i^\delta}.
$$
On $B_{r_0}$, the corresponding coordinate densities are
$$
 N_i^{\rm loc}=e^\phi N_i^g=n_i e^{J_i},
$$
where
$$
 J_1=\psi_1-\psi_1(0)-\log(1+\beta|x|^2/4)=O(|x|^2),
$$
\begin{equation}
 J_2=\log(e^\phi n_*^g)-\log\beta
                       +2\log(1+\beta|x|^2/4)
       =\nabla V(q)\cdot x+O(|x|^2).
 \label{5.4}
\end{equation}
Thus $e^{J_2}=n_*^{\rm loc}/n_{\rm diff}$, where
$n_{\rm diff}=16\beta/(4+\beta|x|^2)^2$.
At $\rho^0=(8\pi,r)$, define
$f_i=-\Delta(\log h_i+\phi)-(A\rho^0)_i e^\phi$.
The equations for the regular terms give
$$
 -\Delta J_1=f_1-n_*^{\rm loc}+n_{\rm diff},\qquad
 -\Delta J_2=f_2+2n_*^{\rm loc}-2n_{\rm diff}.
$$
Hence the residual on $B_{r_0}$ is
\begin{equation}
 -\Delta(U+J)-A\binom{n_1e^{J_1}}{n_2e^{J_2}}-f
 =-A\binom{n_1(e^{J_1}-1)}
                 {(n_2-n_{\rm diff})(e^{J_2}-1)}.
 \label{5.5}
\end{equation}
The subtraction in the second component of \eqref{5.5} leaves
$n_2-n_{\rm diff}$. By Lemma~\ref{le5.1} and its pointwise estimates,
\begin{equation}
 \||x|^2n_1\|_{L^p}=O(\delta^{2/p}),\qquad
 \||x|(n_2-n_{\rm diff})\|_{L^p}
          =O(\delta^{4/p-2}),\quad 1<p<2.
 \label{5.6}
\end{equation}
For the second estimate, $\||x|m_i\|_p=O(a^{2/p-1})=O(\delta^{4/p-2})$.
Multiplying the pointwise bound for $n_2-m_i-n_{\rm diff}$ in the
proof of Lemma~\ref{le5.1} by $|x|$ gives the same bound.

On the cutoff annulus, expansion of $F,Q$ gives
$U_i-L_i=O(\eps)$ in $C^2$. Outside $B_{2r_0}$,
$N_1^g=O(\eps^2)$ and $N_2^g=n_*^g$. Away from $q$ the outer pair
in \eqref{5.1} satisfies
$$
 -\Delta_g\binom{12\pi G_q-z_*/2}{z_*}
 =A\left[\binom0{n_*^g}-\rho^0\right].
$$
Together with \eqref{5.5}--\eqref{5.6}, this yields
\begin{equation}
 \|R^{\rm pre}\|_{L^p}
 :=\left\|-\Delta_gv^\delta
          -A\left[\binom{N_1^g}{N_2^g}-\rho^0\right]\right\|_{L^p}
 =O(\delta^{4/p-2}).
 \label{5.7}
\end{equation}

\emph{Step 2. The masses before normalization.} We prove
\begin{equation}
 M_1:=\int_MN_1^g\dV=8\pi+O(\eps),\qquad
 M_2:=\int_MN_2^g\dV=r+O(\eps).
 \label{5.8}
\end{equation}
For $M_1$, use $\int_\C n_1\dif A=8\pi$, the $O(\eps^2)$ tail outside
$B_{r_0}$, and $\int_{B_{r_0}}|x|^2 n_1\dif A=O(\eps)$ in
$N_1^{\rm loc}=n_1e^{J_1}$. The cutoff annulus contributes $O(\eps)$ to
$\int_{M}(N_2^g-n_*^g)\dV$. The integral over $B_{r_0}$ is
$\int_{B_{r_0}}(n_2-n_{\rm diff})e^{J_2}\dif A$. On $\partial B_{r_0}$,
$$
 \log Q=2\log\eps+\log(|x|^2(4+\beta|x|^2))
 +\frac{2\beta\eps\Re(e^{2i\theta}\overline B)}
           {4+\beta|x|^2}
 +O(\eps^2)
$$
in $C^1$. The term linear in $\eps$ and its radial derivative have zero
circular means. Since $n_2=\Delta\log Q$, we obtain
$$
 \int_{B_{r_0}}(n_2-n_{\rm diff})\dif A=4\pi+O(\eps^2).
$$
Since $n_2-n_{\rm diff}$ is even, the linear term of $e^{J_2}-1$
integrates to zero. By the pointwise estimates in the proof of
Lemma~\ref{le5.1}, the remaining error is bounded by
$$
 C\int_{B_{r_0}}|x|^2|n_2-n_{\rm diff}|\dif A=O(\eps).
$$
Adding $\int_M n_*^g\dV=r-4\pi$ gives the estimate for $M_2$.

\emph{Step 3. Normalization.} The mean subtraction in \eqref{5.1}
gives $V_i^\delta=(\rho_i/M_i)N_i^g$, with
$\log(\rho_i/M_i)=O(\eps)$ by \eqref{5.8} and $|\mu|\le C_0\eps$.
The pointwise bounds follow from Lemma~\ref{le5.1}. Its density
estimates, $\||x|^2m_o\|_p=O(\delta^{2/p})$,
$\||x|m_i\|_p=O(\delta^{4/p-2})$, and the expansion on the cutoff
annulus give \eqref{5.2} after normalization. Finally,
$$
 R^{\rm post}=R^{\rm pre}
   -A\binom{V_1^\delta-N_1^g}{V_2^\delta-N_2^g}
   +A(\rho-\rho^0),
$$
and
$$
 \|n_1\|_p=O(\delta^{2/p-2}),\qquad
 \|n_2\|_p=O(\delta^{4/p-4})+O(1).
$$
Since $\rho_i/M_i-1=O(\delta^2)$, these bounds and
$L^p\hookrightarrow H^{-1}$ give \eqref{5.3}.
\end{proof}

\subsection{Green's identity for the scale derivative}

Use the coordinates of the angular projection in
Section~\ref{rev:sec:core}. Write $w_\delta=U_\delta+\eta_\delta$
for a solution of $-\Delta w=Ae^w+f$. Assume that $\eta_\delta$ and its first spatial derivatives converge
uniformly on each fixed matching circle to the corresponding quantities
for $\eta$. We assume that $\eta$ is $C^2$ across the origin, with
$\eta(0)=0$ and $\nabla\eta(0)=0$ in both components. Set
$$
 \kappa_1=-f_1(p),\qquad
  \eta_1(z)=\frac14(\Delta\eta_1)(0)|z|^2+\Re(Tz^2)+o(|z|^2).
$$
Assume that the first Toda density has mass $o(\delta^2)$ outside
$B_{r_0}$ for each fixed $r_0$. Assume also that the nonlinear integral
in \eqref{5.10}, divided by $\delta^2$, tends to zero as $\delta\to0$ and then $r_0\downarrow0$. On compact
sets of $B$, we assume that the matching convergence, the $C^2$
bounds, and the small-error conditions all hold uniformly. We also assume that the Taylor remainders of $\eta_1$
and its gradient are uniformly $o(r^2)$ and $o(r)$, respectively.

\begin{proposition}\label{pr5.4}
Under these matching assumptions, uniformly on compact $B$ sets,
$$
 \rho_1-8\pi
 =\delta^2(\kappa_1\mathcal M(B)+4\pi\Re(TB))
      +o(\delta^2),
$$
where
$\mathcal M(B)=\int_\C |w|(1+|w-B|^2)^{-2}\dif A_w$.
If $\kappa_1\nabla\mathcal M(B)+4\pi(\Re T,-\Im T)=0$, then
$$
 \quad
 \rho_1-8\pi
 =\delta^2\kappa_1
       (\mathcal M(B)-B\cdot\nabla\mathcal M(B))
       +o(\delta^2).
 \quad
$$
\end{proposition}
\begin{proof}
Set $\eps=\delta^2$ and $E=\eps\partial_\eps U$, with $B,\beta$ fixed.
The constant limit of $E$ away from the center is $E_\infty=(2,0)$.
Since
$$
 CU=\left(\log4+\tfrac13\log D-\log F,\
          \log4+\tfrac23\log D-\log Q\right),\quad
 D=4\beta\eps^4,
$$
we have
\begin{equation}
 K:=CE-CE_\infty
     =(-\eps\partial_\eps\log F,\
             2-\eps\partial_\eps\log Q).
 \label{5.9}
\end{equation}
Apply Green's identity to the correction equation with $CE$. Write
$N_i=e^{\eta_{\delta,i}}-1-\eta_{\delta,i}$. Then
\[
 \begin{aligned}
 &\int_{B_{r_0}}(CE)\cdot f\dif A+
   \int_{B_{r_0}}\sum_i E_i n_i N_i\dif A\\
 &\qquad =
   \int_{\partial B_{r_0}}
    \bigl(\partial_{\boldsymbol n}(CE)\cdot\eta_\delta
                 -(CE)\cdot\partial_{\boldsymbol n}\eta_\delta\bigr)\dif s.
 \end{aligned}
\]
Integrate $-\Delta\eta_\delta=A(n_i(e^{\eta_{\delta,i}}-1))+f$
against the constant $CE_\infty$. Subtracting this identity gives
\begin{equation}
 \begin{aligned}
 2\int_{B_{r_0}}n_1(e^{\eta_{\delta,1}}-1)\dif A
 ={}&\int_{B_{r_0}}K\cdot f\dif A+
   \int_{B_{r_0}}\sum_i E_i n_i N_i\dif A\\
 &-\int_{\partial B_{r_0}}
    \bigl(\partial_{\boldsymbol n}K\cdot\eta_\delta
                    -K\cdot\partial_{\boldsymbol n}\eta_\delta\bigr)\dif s .
 \end{aligned}
 \label{5.10}
\end{equation}
For fixed $r_0$,
\begin{equation}
 \int_{B_{r_0}}(\log F-4\log r)\dif A
      =2\eps\mathcal M(B)+O_{r_0}(\eps^2|\log\eps|),
 \label{5.11}
\end{equation}
The expansion also holds after applying $\eps\partial_\eps$.
For $F_0=\eps^2+|z^2-\eps B|^2$, use $\Delta\log F_0=m_o$
in Green's identity with test function $r^2$ on $B_R$ and let
$R\to\infty$. By Lemma~\ref{le4.6},
$$
 4\lim_{R\to\infty}\int_{B_R}(\log F_0-4\log r)\dif A
      =\int_\C r^2m_o\dif A=8\eps\mathcal M(B).
$$
The $r^{-2}$ term has zero circular mean, and the remaining
boundary contribution tends to zero. Restricting the integral to
$B_{r_0}$ introduces an error of $O_{r_0}(\eps^2)$. The term $\beta\eps^2r^2$ in $F$
contributes $O_{r_0}(\eps^2|\log\eps|)$, also after differentiation.

Similarly,
\begin{equation}
 \int_{B_{r_0}}
 \left[\log(Q/\eps^2)-\log(r^2(4+\beta r^2))\right]\dif A
       =O_{r_0}(\eps^2|\log\eps|),
 \label{5.12}
\end{equation}
also after applying $\eps\partial_\eps$. The term linear in $\eps B$
has zero circular mean, and the radial correction contributes
$O_{r_0}(\eps^2|\log\eps|)$.
Equations \eqref{5.9}--\eqref{5.12} imply
\begin{equation}
 \lim_{r_0\downarrow0}\lim_{\eps\downarrow0}
 \eps^{-1}\int_{B_{r_0}}K\cdot f\dif A
       =-2f_1(p)\mathcal M(B)=2\kappa_1\mathcal M(B).
 \label{5.13}
\end{equation}
The contribution of $f-f(p)$ to this limit vanishes as $r_0\downarrow0$.

On $\partial B_{r_0}$, the leading term of $K_1$ is
$2\eps\Re(Be^{-2i\theta})/r_0^2$. For the second component,
\eqref{4.2} and \eqref{5.9} give, on every fixed matching annulus,
\[
 \eps^{-1}K_2\longrightarrow
 -\frac{2\beta\Re(Be^{-2i\theta})}{4+\beta r^2}
 \qquad\text{in } C^1.
\]
The limiting function is bounded, and its radial derivative is
$O(r)$. Since $\eta_2=O(r^2)$ and $\partial_r\eta_2=O(r)$, the
limiting second-component boundary pairing divided by $\eps$ is
$O(r_0^2)$. Pairing the first component with the harmonic quadratic
term of $\eta_1$ thus gives
\begin{equation}
 \lim_{r_0\downarrow0}\lim_{\eps\downarrow0}
 \eps^{-1}\int_{\partial B_{r_0}}
   \bigl(\partial_{\boldsymbol n}K\cdot\eta_\delta
                -K\cdot\partial_{\boldsymbol n}\eta_\delta\bigr)\dif s
          =-8\pi\Re(TB).
 \label{5.14}
\end{equation}
The $|z|^2$ term integrates to zero since
$\int_0^{2\pi}\Re(Be^{-2i\theta})\dif\theta=0$.

By the matching assumptions, the nonlinear integral is $o(\eps)$
in the stated order of limits. Since the first profile has mass
$O_{r_0}(\eps^2)$ outside $B_{r_0}$, the left side of \eqref{5.10}
is $2(\rho_1-8\pi)+o(\eps)$. Substituting \eqref{5.13}--\eqref{5.14}
gives the first formula. Taking the scalar product of the angular
equation with $(\Re B,\Im B)$ gives the second.
\end{proof}

\begin{remark}\label{re5.5}
The $E_i$ are uniformly bounded. For $r\ge R\delta$ and $B$ in a
compact set, $E_1=O(1)$ and $E_2=O(\delta^2/r^2+\delta^4/r^4)$.
Thus a smooth correction $\eta_{\rm app}=O(r^2)$ satisfies
$$
 \int_{B_{r_0}}\sum_i |E_i|n_i|\eta_{{\rm app},i}|^2\dif A
 \le C(\delta^4|\log\delta|+\delta^2r_0^4).
$$
If the weighted exponential quadratic integral of the remaining
normalized correction is $o(\delta^2)$, the exponential Taylor
bound and Cauchy--Schwarz give an $o(\delta^2)$ nonlinear integral
in \eqref{5.10} in the stated order of limits.
\end{remark}

\subsection{Estimates for blowing-up solutions}\label{sec:apriori-chapter}
\subsubsection{Local hypotheses}

Let $M$ be a smooth closed surface of area one and let $h_i>0$ be
smooth. Consider mean-zero solutions $u_{i,\nu}$ of
\begin{equation}\label{5.15}
 -\Delta_g u_{i,\nu}=\sum_{j=1}^2 A_{ij}(f_{j,\nu}-\rho_{j,\nu}),
 \qquad
 f_{j,\nu}=\rho_{j,\nu}
 \frac{h_j e^{u_{j,\nu}}}{\int_Mh_j e^{u_{j,\nu}}\dV},
 \qquad
 A=\begin{pmatrix}2&-1\\-1&2\end{pmatrix},
\end{equation}
with bounded positive parameters. The single-point case is described
by (\ref{apriori:formation-measures}). The rescaled limits and local
bounds are given in (\ref{apriori:formation-profiles}) and
(\ref{apriori:formation-harnack}).

\begin{enumerate}[label=\textnormal{(F\arabic*)},ref=F\arabic*]
\item\label{apriori:formation-measures} Suppose $\rho_{1,\nu}\to8\pi$
and $\rho_{2,\nu}=\rho>4\pi$. The only concentration point is $q$, and
\begin{equation}\label{5.16}
 f_{1,\nu}\dV\rightharpoonup8\pi\delta_q,
 \qquad
 f_{2,\nu}\dV\rightharpoonup4\pi\delta_q+f_*\dV,
\end{equation}
where $f_*>0$ is smooth and
\begin{equation}\label{5.17}
 f_*=(\rho-4\pi)\frac{h_2e^{w_*}}{\int_Mh_2e^{w_*}\dV},
 \quad \int_Mw_*\dV=0,
 \quad -\Delta_g w_*=2(f_*-(\rho-4\pi)).
\end{equation}
On compact subsets of $M\setminus\{q\}$, we have
$(f_{2,\nu},u_{2,\nu})\to(f_*,w_*)$ smoothly and $|\nabla_g u_{1,\nu}|\leq C$.

\item\label{apriori:formation-profiles} In a fixed isothermal coordinate
$x$ about $q$, let $z=x-x(q_\nu)$ with $q_\nu\to q$ and
$g=e^{\phi_\nu(z)}|dz|^2$. Set $F_{i,\nu}=e^{\phi_\nu}f_{i,\nu}$
and $U_{i,\nu}=\log F_{i,\nu}$. For scales $a_\nu=o(\delta_\nu)$,
$\delta_\nu\to0$, and a complex parameter $B_\nu\to B$, we have
\begin{align}
 a_\nu^2F_{2,\nu}(a_\nu y)&\longrightarrow
       m_i(y):=\frac4{(1+|y|^2)^2},\label{5.18}\\
 \delta_\nu^2F_{1,\nu}(\delta_\nu y)&\longrightarrow
       m_B(y):=\frac{16|y|^2}{(1+|y^2-B|^2)^2}.
       \label{5.19}
\end{align}
The convergence in \eqref{5.18} is locally smooth on $\C$, and
$\log(a_\nu^2F_{2,\nu}(a_\nu\,\cdot))\to\log m_i$ in $C^2_{\rm loc}(\C)$.
The convergence in \eqref{5.19} is locally smooth on $\C\setminus\{0\}$.

We first establish the smallness of the first density on the inner
scale. The argument uses the local first mass in \eqref{5.16},
the two profiles and \eqref{5.20}. It works on any fixed
isolating disk, including the local hypotheses specified below.
For a fixed $S>2$, set $V_\nu(y)=\log(a_\nu^2F_{1,\nu}(a_\nu y))$
on $B_S$. Given $\gamma>0$, choose a fixed annulus
$c\le|y|\le C$ on which $m_B$ has mass at least $8\pi-\gamma$.
Its image at scale $\delta_\nu$ is disjoint from $B_{S a_\nu}$
for large $\nu$ because $a_\nu/\delta_\nu\to0$. By
\eqref{5.19} and the first mass $8\pi+o(1)$ on the isolating
disk, this gives
\[
 m_\nu:=\int_{B_S}e^{V_\nu}\dif A\longrightarrow0 .
\]
The bounds in \eqref{5.20} give $V_\nu\le C_S$ on $\partial B_S$.
By \eqref{2.2},
\[
 -\Delta_y V_\nu=2e^{V_\nu}+g_\nu,\qquad
 g_\nu=-a_\nu^2F_{2,\nu}(a_\nu y)+a_\nu^2 f_1^{\rho_\nu}(a_\nu y),
 \qquad \|g_\nu\|_{L^\infty(B_S)}\le C_S.
\]
The last bound follows from \eqref{5.18} and the smooth forcing.
With $G_S$ the positive Dirichlet Green operator on $B_S$,
the maximum principle gives $V_\nu\le C_S+2G_S(e^{V_\nu})$.
For any fixed $p>1$, Jensen's inequality with the probability
density $e^{V_\nu}/m_\nu$ and
$G_S(y,x)\le (2\pi)^{-1}\log(C_S/|x-y|)$ gives
\[
 \int_{B_S} e^{2pG_S(e^{V_\nu})}\dif A
 \le \sup_{x\in B_S}\int_{B_S}
                 e^{2p m_\nu G_S(y,x)}\dif A_y\le C_{S,p}.
\]
Here $p m_\nu/\pi$ is uniformly below $2$ for large $\nu$.
Thus the zero-boundary potential
$P_\nu=G_S(2e^{V_\nu}+g_\nu)$ is uniformly bounded by the
$L^p$ estimate and the Dirichlet inverse. The function
$V_\nu-P_\nu$ is harmonic. Its exponential is subharmonic,
so for each $0<R<S$ the interior mean-value inequality gives
\begin{equation}\tag{5.19a}\label{eq:inner-first-density-small}
 \sup_{B_R}a_\nu^2F_{1,\nu}(a_\nu\,\cdot)
   =\sup_{B_R}e^{V_\nu}\le C_{R,S}\,m_\nu\longrightarrow0 .
\end{equation}
The constant uses only the fixed radii, the annular bound,
the uniform forcing and the bound for the rescaled second density
on $B_S$. The argument applies in the moving charts with these
uniform bounds. For a fixed finite set of limiting configurations, we can choose
the constants uniformly over the corresponding disks. The mass
convergence used here holds on each profile subsequence.

The masses of $m_i$ and $m_B$ are $4\pi$ and $8\pi$. Replacing
$a_\nu$ by $a_\nu'=2F_{2,\nu}(0)^{-1/2}=a_\nu(1+o(1))$ gives
$a_\nu'^2F_{2,\nu}(0)=4$ and preserves both profile limits and
\eqref{eq:inner-first-density-small}. Only $a'_\nu<a_\nu$ adds
radii to \eqref{5.20}. The right-hand sides of the equations
\eqref{2.2} for $U_{i,\nu}(a_\nu y)$ are uniformly bounded on
$B_3$ by \eqref{eq:inner-first-density-small}, \eqref{5.18}
and the rescaled forcing bound. Subtract the mean of the boundary
values on $\partial B_3$ for each equation. The original
\eqref{5.20} bounds the remaining boundary values. The Dirichlet
estimate then bounds the oscillations on the additional band
$2a'_\nu/a_\nu\le|y|\le2$. The pointwise upper bounds for the densities there
follow from the two inner density bounds. Thus the asymptotic
normalization suffices for all the local hypotheses.

\item\label{apriori:formation-harnack} There are $r_0>0$ and $C$
independent of $r$ and $\nu$ such that, for $2a_\nu\le r\le r_0$,
\begin{equation}\label{5.20}
 \max_{|z|=r}U_{i,\nu}-\min_{|z|=r}U_{i,\nu}\le C,
 \qquad r^2 F_{i,\nu}(z)\le C.
\end{equation}
\end{enumerate}

For the local estimates at a mixed point $q$, let $D$ be a disk of
fixed radius in the moving coordinate. Assume \eqref{5.16} on a
neighborhood of $\overline D$ containing no other concentration
point, with $f_*$ smooth and positive there. On compact subsets of
this neighborhood away from $q$, assume $f_2\to f_*$ smoothly and
$|\nabla_g u_1|\leq C$. We also assume
(\ref{apriori:formation-profiles}) and (\ref{apriori:formation-harnack}).
These are the local hypotheses at $q$. They may be imposed in a
smoothly varying family of isothermal charts with uniform coordinate
bounds. The local forcing in \eqref{2.2} uses the full parameter
vector $\rho_\nu$. The Green potentials of $\mathbf1_{M\setminus D}f_{i,\nu}$ have uniformly
bounded derivatives on smaller disks by separation from their
supports and bounded total mass.

In applications to actual Toda solutions, \eqref{2.6} gives the
regular field for testing \eqref{4.30}, and
Proposition~\ref{pr4.19} supplies
(\ref{apriori:formation-profiles}) and
(\ref{apriori:formation-harnack}). The local estimates below use
these two conditions directly. With only the local density data above, we obtain the same regular
part of the first field from the local Green representation by subtracting
the net point mass $2(8\pi)-4\pi=12\pi$. The first density converges uniformly
to zero on compact sets off the center. A positive value, together with
the log-gradient bound, would give a fixed positive mass on a smaller ball. On a
fixed subdisk $D'$ the first potential can be written
\[
 u_{1,\nu}(z)=
 -\frac1{2\pi}\int_{D'}\log|z-y|\,
      (2F_{1,\nu}(y)-F_{2,\nu}(y))\dif A_y+\Pi_\nu(z)+H_\nu(z),
\]
where $\Delta_0\Pi_\nu=(2\rho_{1,\nu}-\rho_{2,\nu})e^\phi$
with zero Dirichlet data, and $H_\nu$ is harmonic. Indeed,
applying $\Delta_0$ to the first term gives
$-2F_{1,\nu}+F_{2,\nu}$ and the Poisson term restores the
smooth background. The Poisson potentials converge smoothly
on smaller disks. The logarithmic potential converges in first
derivatives on compact annuli: near an evaluation point off the
center, the densities are uniformly bounded and the $1/|z-y|$
kernel is integrable and has small integral on a small ball. The
remaining truncated kernel is bounded and continuous, so we can
test it against the weakly convergent measures. This argument
is uniform on any fixed smaller annulus. On a fixed isolating
circle, the logarithmic potential and $\Pi_\nu$ have bounded
oscillation. The same holds for $u_{1,\nu}$ by the local gradient bound.
After subtracting the mean on that circle, Poisson's formula
bounds $H_\nu$ in each $C^j$ on smaller disks and gives smooth
subsequential convergence. After subtracting the point mass
$12\pi$ in the potential, the remaining logarithmic potential
has smooth density near zero. Together with $\Pi_\nu$ and
$H_\nu$, this gives the smooth field used in the
first-derivative normalization. We work along the profile
subsequence and often omit $\nu$. Constants may depend on a fixed
small $\epsilon>0$ and a compact set of $B$.

\subsubsection{Pointwise estimates and comparison of scales}

\begin{lemma}\label{le5.6}
Under the local hypotheses at $q$, for every sufficiently small
$\epsilon>0$ there are $R_\epsilon>1$,
$r_\epsilon>0$ and $C_\epsilon>0$ such that, for all large $\nu$,
\begin{align}
 F_1(z)&\le C_\epsilon\delta^{-4+\epsilon}|z|^{2-\epsilon},
 &F_2(z)&\le C_\epsilon a^{2-\epsilon}|z|^{-4+\epsilon}
 &&(R_\epsilon a\le|z|\le\delta/R_\epsilon),
 \label{5.21}\\
 F_1(z)&\le C_\epsilon\delta^{4-\epsilon}|z|^{-6+\epsilon},
 &F_2(z)&\le C_\epsilon|z|^{-\epsilon}
 &&(R_\epsilon\delta\le|z|\le r_\epsilon).
 \label{5.22}
\end{align}
On the transition annulus $\delta/R_\epsilon\le|z|\le
R_\epsilon\delta$, we have
\begin{equation}\label{5.23}
 F_1\le C_\epsilon\delta^{-2},\qquad
 F_2\le C_\epsilon\delta^{-\epsilon}.
\end{equation}
On $|z|\le R_\epsilon a$, we have
\begin{equation}\label{5.24}
 F_2\le C_\epsilon a^{-2},\qquad
 F_1\le C_\epsilon\delta^{-4+\epsilon}a^{2-\epsilon}.
\end{equation}
On $D\setminus B_{r_\epsilon}$, $f_1\le C_\epsilon\delta^{4-\epsilon}$
and $f_2$ is uniformly bounded. Under (\ref{apriori:formation-measures})
these bounds hold on $M\setminus B_{r_\epsilon}$. In both cases,
\begin{equation}\label{5.25}
 a=\delta^{2+o(1)}.
\end{equation}
\end{lemma}

\begin{proof}
By \eqref{2.2}, the local equations are
\begin{equation}\label{5.26}
 \Delta_0U_1+2F_1-F_2=\kappa_1^c,
 \qquad
 \Delta_0U_2-F_1+2F_2=\kappa_2^c,
\end{equation}
where the smooth terms are uniformly bounded and given by
$$
 \kappa_i^c=\Delta_0(\log h_i+\phi_\nu)
             +\sum_j A_{ij}\rho_j e^{\phi_\nu}.
$$
Set $m_j(r)=(2\pi)^{-1}\int_{B_r}F_j\dif A$ and let $\bar U_j(r)$
be the circular mean. Integration of \eqref{5.26} gives
\begin{equation}\label{5.27}
 r\bar U_1'(r)=-2m_1(r)+m_2(r)+O(r^2),\qquad
 r\bar U_2'(r)=m_1(r)-2m_2(r)+O(r^2).
\end{equation}

Fix $\eta>0$. By \eqref{5.16}, \eqref{5.18}, and \eqref{5.19},
we can choose $R$ large and then $r_\epsilon$ small so that
\begin{align*}
 &0\le m_1(r)\le\eta,\qquad |m_2(r)-2|\le\eta
             &&(Ra\le r\le\delta/R),\\
 &|m_1(r)-4|+|m_2(r)-2|\le 2\eta
             &&(R\delta\le r\le r_\epsilon).
\end{align*}
Indeed, the profile $m_B$ accounts for the full first mass at $q$.
Thus $\int_{B_{\delta/R}}F_1\dif A=\int_{B_{1/R}}m_B\dif A+o(1)$.
For $F_2$, the inner profile has mass $4\pi$ and
$\int_{B_{r_\epsilon}}F_2\dif A=4\pi+O(r_\epsilon^2)+o(1)$. The
bounds at intermediate radii follow by monotonicity of $m_j(r)$.

The slopes $r\bar U_1',r\bar U_2'$ in \eqref{5.27} are
$2+O(\eta),-4+O(\eta)$ on $Ra\le r\le\delta/R$ and
$-6+O(\eta),O(\eta)$ on $R\delta\le r\le r_\epsilon$. Choose
$\eta$ sufficiently small in terms of $\epsilon$. The profiles give
\[
 \bar U_2(Ra)=-2\log a+O_R(1),\qquad
 \begin{aligned}
 \bar U_1(\delta/R)&=-2\log\delta+O_R(1),\\
 \bar U_1(R\delta)&=-2\log\delta+O_R(1).
 \end{aligned}
\]
Also $\bar U_2(r_\epsilon)=O(1)$ by outer convergence. Integrating
\eqref{5.27} and using \eqref{5.20} gives \eqref{5.21}--\eqref{5.22}.
The first bound in \eqref{5.23} follows from \eqref{5.19} on the
fixed rescaled annulus. For the second, integrate the bounded
slope $r\bar U_2'(r)$ from $R\delta$ to $r\in[\delta/R,R\delta]$
and use \eqref{5.20} again.

By \eqref{5.18}, \eqref{eq:inner-first-density-small}, and \eqref{5.26},
$$
 \sup_{B_{2R}}a^2F_2(a\,\cdot)\leq C_R,\qquad
 \sup_{B_{2R}}a^2F_1(a\,\cdot)=o(1),\qquad
 \|\Delta_yU_1(a\,\cdot)\|_{L^\infty(B_{2R})}\leq C_R.
$$
The Dirichlet formula on $B_{2R}$ and \eqref{5.21} imply
$$
 \sup_{B_R}U_1(a\,\cdot)
    \leq\sup_{\partial B_{2R}}U_1(a\,\cdot)+C_R,
 \qquad
 \sup_{B_{Ra}}F_1\leq C_R\delta^{-4+\epsilon}a^{2-\epsilon}.
$$
This proves \eqref{5.24}, with a change in $C_R$ if $R$ is enlarged
by a fixed factor. On $D\setminus B_{r_\epsilon}$, the first bound
follows from \eqref{5.22} on $\partial B_{r_\epsilon}$ and the bound
on $\nabla_g\log f_1$. The second follows from the smooth convergence
of $f_2$. Under (\ref{apriori:formation-measures}), the same argument
applies on $M\setminus B_{r_\epsilon}$.

Finally, integrate $\bar U_2'$ from $Ra$ to $r_\epsilon$. The two
bounds from \eqref{5.27} for $r\bar U_2'$ apply off the annulus
$\delta/R\le r\le R\delta$, whose contribution is $O_R(1)$. Thus
\begin{equation}\label{5.28}
 |2\log a-4\log\delta|
 \le C\eta( |\log a|+|\log\delta|)+C_\eta.
\end{equation}
For fixed small $\eta$, \eqref{5.28} bounds $|\log a|/|\log\delta|$
above and below. Dividing by $|\log\delta|$ and letting $\nu\to\infty$
and then $\eta\downarrow0$ gives \eqref{5.25}.
\end{proof}

\subsubsection{Local mass estimates}

Set $P=(2U_1+U_2)/3$, $Q=(U_1+2U_2)/3$, and
$$
 K_1=(2\kappa_1^c+\kappa_2^c)/12,\qquad
 K_2=(\kappa_1^c+2\kappa_2^c)/12.
$$
Then $\bar\partial P'=K_1-F_1/4$ and
$\bar\partial Q'=K_2-F_2/4$. As in \eqref{4.15}, put
\begin{equation}\label{5.29}
 W_2=P''+Q''-(P')^2+P'Q'-(Q')^2,
\end{equation}
By \eqref{4.17},
\begin{equation}\label{5.30}
 \bar\partial W_2=K_1'+K_2'-(2P'-Q')K_1-(2Q'-P')K_2.
\end{equation}
\begin{lemma}
\label{le5.7}
Let $D=B_{r_0}$ be the fixed coordinate disk. Uniformly for
$0<r<r_0/4$ and all $\nu$,
\begin{equation}\label{5.31}
 \frac1{2i}\int_{\partial B_r}zW_2(z)\,dz
       =\int_{B_r}z\bar\partial W_2\dif A=O(r^2).
\end{equation}
On $|z|=r$, approximate $\int_D F(x)(z-x)^{-1}\dif A_x$ by
$z^{-1}\int_{B_r}F\dif A$. If $L=\int_{B_{r/2}}|x|F(x)\dif A_x$
and $F+r|\nabla F|\le D$ on $r/2<|x|<2r$, the errors in the
potential and its $z$ derivative are bounded, respectively, by
\begin{align}
 C\left(\frac L{r^2}+rD+
             \int_{D\setminus B_{2r}}\frac{F(x)}{|x|}\dif A_x\right),
 \label{5.32}\\
 C\left(\frac L{r^3}+D+
             \int_{D\setminus B_{2r}}\frac{F(x)}{|x|^2}\dif A_x\right).
 \label{5.33}
\end{align}
The regular part of $G$ and the densities outside $D$ contribute
uniformly bounded smooth terms on $B_{r_0/2}$.
\end{lemma}

\begin{proof}
Green's representation gives
$$
 P'=H_1-\frac1{4\pi}\int_{D}\frac{F_1(x)}{z-x}\dif A_x,\qquad
 Q'=H_2-\frac1{4\pi}\int_{D}\frac{F_2(x)}{z-x}\dif A_x,
 \qquad \|H_i\|_{C^1(B_{r_0/2})}\leq C.
$$
Since $\int_{D} F_i\dif A\leq C$,
$$
 \sup_x\int_{B_r}\frac{\dif A_z}{|z-x|}\leq Cr,
 \qquad \int_{B_r}(|P'|+|Q'|)\dif A\leq Cr.
$$
Equation \eqref{5.30} and Stokes' formula imply
$$
 \left|\frac1{2i}\int_{\partial B_r}zW_2\,dz\right|
 \leq\int_{B_r}|z\bar\partial W_2|\dif A
 \leq Cr\left(r^2+\int_{B_r}(|P'|+|Q'|)\dif A\right)
 \leq Cr^2.
$$
For $|z|=r$ and $j=0,1$,
$$
 \begin{aligned}
 \left|\partial_z^j\left(\frac1{z-x}-\frac1z\right)\right|
 &\leq\frac{C|x|}{r^{2+j}},&&|x|<r/2,\\
 |\partial_z^j(z-x)^{-1}|&\leq C|x|^{-1-j},&&|x|>2r.
 \end{aligned}
$$
Put $A_r=B_{2r}\setminus B_{r/2}$. The assumption $F+r|\nabla F|\leq D$ gives
$$
 \begin{aligned}
 \int_{A_r}F\dif A&\leq Cr^2D,
 &\int_{A_r}\frac{F(x)}{|z-x|}\dif A_x&\leq CrD,\\
 \operatorname{pv}\int_{B_{r/4}(z)}\frac{F(z)}{(z-x)^2}\dif A_x&=0,
 &\int_{B_{r/4}(z)}\frac{|F(x)-F(z)|}{|z-x|^2}\dif A_x
       &\leq\frac{CD}r\int_{B_{r/4}(z)}\frac{\dif A_x}{|z-x|}\leq CD,\\
 \int_{A_r\setminus B_{r/4}(z)}\frac{F(x)}{|z-x|^2}\dif A_x&\leq CD.
 \end{aligned}
$$
Thus, with $M_r=\int_{B_r}F\dif A$ and $L=\int_{B_{r/2}}|x|F(x)\dif A_x$,
$$
 \left|\partial_z^j\left(\int_{D}\frac{F(x)}{z-x}\dif A_x-\frac{M_r}z\right)\right|
 \leq C\left(\frac L{r^{2+j}}+r^{1-j}D
       +\int_{D\setminus B_{2r}}\frac{F(x)}{|x|^{1+j}}\dif A_x\right),
 \qquad j=0,1.
$$
For $F_i$, \eqref{5.20} and the mass bound give
$$
 F_i\leq Cr^{-2},\qquad |\nabla U_i|\leq C/r,
 \qquad |\nabla F_i|=F_i|\nabla U_i|\leq CF_i/r
$$
on the relevant annuli.
\end{proof}

\begin{proposition}\label{pr5.8}
Set $\tau=\delta^{3/2}$. Let $D$ be a disk of fixed radius in the
moving coordinate containing a single mixed concentration point $q$.

Under the local hypotheses of Lemma~\ref{le5.6},
\[
 \int_{B_\tau}F_2\dif A=4\pi+O(\delta^{1/4}),
 \qquad \int_D F_1\dif A=8\pi+O(\delta^{1/4}).
\]
Under the global condition~(\ref{apriori:formation-measures}),
\begin{equation}\label{5.34}
 \int_{B_\tau}F_2\dif A=4\pi+O(\delta^{1/4}),
 \qquad \rho_{1,\nu}-8\pi=O(\delta^{1/4}).
\end{equation}
\end{proposition}

\begin{proof}
The local hypotheses give the estimates of Lemma~\ref{le5.6} inside
$D$. In particular, $\int_D F_1\dif A\to8\pi=\int_\C m_B\dif A$.
The bound on $\nabla_g u_1$ on fixed annuli extends the estimate
for $F_1$ to $\partial D$.

Fix $0<\epsilon<1/10$ in Lemma~\ref{le5.6}. By \eqref{5.25},
for large $\nu$,
\begin{equation}\label{5.35}
 \delta^{9/4}\le a\le\delta^{7/4}.
\end{equation}
Constants may depend on $R_\epsilon$.

\paragraph{The inner contour.}
We first prove, uniformly on $|z|=\tau=\delta^{3/2}$,
\begin{equation}\label{5.36}
 |P'(z)|\le C\delta^{-1},\qquad |P''(z)|\le C\delta^{-2}.
\end{equation}
In Green's representation of Lemma~\ref{le5.7}, split the source
disk into
$$
 B_{\tau/2},\qquad
 A_\tau=B_{2\tau}\setminus\overline{B_{\tau/2}},\qquad
 D\setminus\overline{B_{2\tau}}.
$$
The kernel is singular at $x=z$ in $A_\tau$. The smooth terms in
$P'$ and $P''$ are uniformly bounded.

For $x\in B_{\tau/2}$, $|z-x|\ge\tau/2$. Since $a/\tau\to0$,
\eqref{5.21} and \eqref{5.24} give
\begin{align*}
 \int_{B_{\tau/2}}F_1\dif A
 &\le C\delta^{-4+\epsilon}
       \left(a^{4-\epsilon}+\int_{R_\epsilon a}^{\tau/2}
                                      r^{3-\epsilon}\,dr\right)\\
 &\le C\delta^{-4+\epsilon}\tau^{4-\epsilon}.
\end{align*}
Consequently, for $j=1,2$,
\begin{equation}\label{5.37}
 \int_{B_{\tau/2}}\frac{F_1(x)}{|z-x|^j}\dif A_x
 \le C\delta^{-4+\epsilon}\tau^{4-\epsilon-j}.
\end{equation}

On $A_\tau$, \eqref{5.21} and the gradient bound in
Lemma~\ref{le5.7} give
$$
 F_1+\tau|\nabla F_1|\le D_\tau,
 \qquad D_\tau=C\delta^{-4+\epsilon}\tau^{2-\epsilon}.
$$
The Cauchy integral over $A_\tau$ is bounded by $C\tau D_\tau$.
For the $z$ derivative, use $B_{\tau/4}(z)\subset A_\tau$ and
subtract $F_1(z)$ in the principal value integral. This gives
$$
 \left|\int_{B_{\tau/4}(z)}
       \frac{F_1(x)-F_1(z)}{(z-x)^2}\dif A_x\right|
 \le \frac{CD_\tau}{\tau}
          \int_{B_{\tau/4}(z)}\frac{\dif A_x}{|z-x|}
 \le CD_\tau.
$$
The constant term has zero principal value. On
$A_\tau\setminus B_{\tau/4}(z)$, the kernel of
$\partial_z(1/(z-x))=-\operatorname{pv}(z-x)^{-2}$ is bounded by
$C\tau^{-2}$, so its contribution is $O(D_\tau)$.

For $|x|>2\tau$, $|z-x|\ge |x|/2$. Splitting at
$|x|=\delta/R_\epsilon$, use \eqref{5.21} and the first mass bound
to obtain
\begin{align*}
 \int_{2\tau<|x|<\delta/R_\epsilon}
                     \frac{F_1(x)}{|z-x|^j}\dif A_x
 &\le C\delta^{-4+\epsilon}
              \int_{2\tau}^{\delta/R_\epsilon}r^{3-\epsilon-j}\,dr
 \le C\delta^{-j},\\
 \int_{D\setminus B_{\delta/R_\epsilon}}
                     \frac{F_1(x)}{|z-x|^j}\dif A_x
 &\le C\delta^{-j}\int_{D} F_1\dif A\le C\delta^{-j}.
\end{align*}

Combining the estimates on the three regions, we get $|P'|\le C(1+\delta^{-1}+
\delta^{-4+\epsilon}\tau^{3-\epsilon})$ and $|P''|\le C(1+\delta^{-2}+
\delta^{-4+\epsilon}\tau^{2-\epsilon})$. The last terms are
$\delta^{1/2-\epsilon/2}$ and $\delta^{-1-\epsilon/2}$, respectively,
so \eqref{5.36} follows.

Set $b_\nu=(4\pi)^{-1}\int_{B_\tau}F_2\dif A\to1$ by \eqref{5.16}
and \eqref{5.18}. The first moment of $F_2$ on $B_{\tau/2}$ is
$O(a)$ by \eqref{5.24} and
$$
 a^{2-\epsilon}\int_{Ra}^{\tau}r^{-2+\epsilon}\,dr\le Ca.
$$
On $A_\tau$, $F_2\le Ca^{2-\epsilon}\tau^{-4+\epsilon}$.
On $2\tau<|x|<\delta/R$, apply \eqref{5.21}. On $|x|\ge\delta/R$,
the mass bound gives Cauchy integrals of orders $O(\delta^{-1})$
and $O(\delta^{-2})$ for the kernel and its derivative, respectively.
By \eqref{5.32}--\eqref{5.33},
\begin{equation}\label{5.38}
 Q'=-\frac{b_\nu}{z}+O(E/\tau),\qquad
 Q''=\frac{b_\nu}{z^2}+O(E/\tau^2),
 \quad
 E\le C\left(\frac a\tau+
        \left(\frac a\tau\right)^{2-\epsilon}+
        \frac\tau\delta\right)=O(\delta^{1/4}).
\end{equation}
Substitute \eqref{5.36} and \eqref{5.38} into \eqref{5.29}. The
contour integral in \eqref{5.31} is $\pi(b_\nu-b_\nu^2)+O(\delta^{1/4})$.
The terms containing $P'$ or $P''$ contribute $O(\tau/\delta)$.
Hence $b_\nu-b_\nu^2=O(\delta^{1/4})$ and, since $b_\nu\to1$,
$b_\nu=1+O(\delta^{1/4})$.

\paragraph{The outer contour.}
Set $r=\delta^{1/2}$. By \eqref{5.21}--\eqref{5.23}, the mass of
$F_2$ on $B_r\setminus B_\tau$ is bounded by
\begin{equation}\label{5.39}
 C\left(\left(\frac a\tau\right)^{2-\epsilon}
                  +\delta^{2-\epsilon}+r^{2-\epsilon}\right)
       =O(\delta^{1/4}).
\end{equation}
Consequently
$\widetilde b_\nu=(4\pi)^{-1}\int_{B_r}F_2\dif A
=1+O(\delta^{1/4})$.

The first moment of $F_1$ in $D$ is $O(\delta)$, and its mass
in $D\setminus B_r$ is
$$
 O(\delta^{4-\epsilon}r^{-4+\epsilon})
                  =O(\delta^{2-\epsilon/2}).
$$
Put $\alpha_\nu=(4\pi)^{-1}\int_D F_1\dif A\to2$.
Applying \eqref{5.32}--\eqref{5.33} and replacing the first mass
in $B_r$ by the first mass in $D$, we obtain
\begin{align*}
 P'&=-\alpha_\nu/z+O(\delta^{1/2}/r),&
 P''&=\alpha_\nu/z^2+O(\delta^{1/2}/r^2),\\
 Q'&=-\widetilde b_\nu/z+O(\delta^{1/2}/r),&
 Q''&=\widetilde b_\nu/z^2+O(\delta^{1/2}/r^2).
\end{align*}
To estimate $Q'$ and $Q''$, use that the first moment of $F_2$ on
$B_r$ is bounded by $C(a+r^{3-\epsilon})$, and $F_2=O(r^{-\epsilon})$
on $A_r$. The bound $F_2(x)\le C|x|^{-\epsilon}$ yields
$\int_{D\setminus B_{2r}}F_2(x)|x|^{-1}\dif A_x=O(1)$ and
$\int_{D\setminus B_{2r}}F_2(x)|x|^{-2}\dif A_x=O(r^{-\epsilon})$.

Substitution in \eqref{5.31} gives
\begin{equation}\label{5.40}
 \alpha_\nu+\widetilde b_\nu-\alpha_\nu^2
       +\alpha_\nu\widetilde b_\nu-\widetilde b_\nu^2
       =O(\delta^{1/2}).
\end{equation}
Put $x=\alpha_\nu-2$ and $y=\widetilde b_\nu-1$. Then \eqref{5.40}
becomes
$$
 -2x+y-x^2+xy-y^2=O(\delta^{1/2}).
$$
Since $x,y\to0$, we absorb the terms in $x^2$ and $xy$ to obtain
$|x|\leq C(|y|+y^2+\delta^{1/2})=O(\delta^{1/4})$.
This proves the first mass estimate because
$\int_D F_1\dif A-8\pi=4\pi x$. Under
(\ref{apriori:formation-measures}), Lemma~\ref{le5.6} gives
$\int_{M\setminus D}f_1\dV=O(\delta^{4-\epsilon})$, proving \eqref{5.34}.
\end{proof}

\subsubsection{The limiting second density and the scale ratio}

\begin{lemma}\label{le5.9}
Under (\ref{apriori:formation-measures}), for some fixed $p>1$,
in particular for $p=6/5$,
\begin{equation}\label{5.41}
 \big\|f_{2,\nu}{\bf1}_{M\setminus B_\tau}-f_*\big\|_{L^p(M)}
       \longrightarrow0.
\end{equation}
Here and below $B_\tau$ is defined in the moving coordinate.
\end{lemma}
\begin{proof}
Take $p=6/5$ and $0<\epsilon<1/10$. By \eqref{5.21}--\eqref{5.23}
and \eqref{5.25},
$$
 \begin{aligned}
 \|F_2\|_{L^p(D_{\delta/R}\setminus D_\tau)}
 &\leq Ca^{2-\epsilon}\tau^{2/p-4+\epsilon}
   =C\delta^{3/p-2-\epsilon/2+o(1)}=o(1),\\
 \|F_2\|_{L^p(D_{R\delta}\setminus D_{\delta/R})}
 &\leq C\delta^{2/p-\epsilon}=o(1),\\
 F_2(z)&\leq C|z|^{-\epsilon},\qquad R\delta<|z|<d.
 \end{aligned}
$$
Smooth convergence on $M\setminus D_d$ gives
$$
 \limsup_{\nu\to\infty}
 \|f_{2,\nu}{\bf1}_{M\setminus B_\tau}-f_*\|_{L^p(M)}
 \leq Cd^{2/p-\epsilon}\longrightarrow0\qquad(d\downarrow0).
$$
\end{proof}

\begin{lemma}\label{le5.10}
Assume (\ref{apriori:formation-measures}). In the moving coordinate
centered at $q_\nu$, on a disk $B_{r_0}$ of fixed radius,
\begin{align}
 \int_{B_{r_0}} F_{1,\nu}(z)\log|z|\dif A
       &=8\pi\log\delta+2\pi\log(1+|B|^2)+o(1),
       \label{5.42}\\
 \int_{B_\tau}F_{2,\nu}(z)\log|z|\dif A
       &=4\pi\log a+o(1).
       \label{5.43}
\end{align}
The first-component mass outside $B_{r_0}$ is $O(\delta^{4-\epsilon})$.
\end{lemma}
\begin{proof}
Put $t=a/\delta$, $\rho_o(y)=\delta^2F_1(\delta y)$ and
$\rho_i(y)=a^2F_2(ay){\bf1}_{|ay|<\tau}$. By Lemma~\ref{le5.6},
$$
 \rho_o(y)\leq C
 \begin{cases}
 (|y|+t)^{2-\epsilon},&|y|\leq1,\\
 |y|^{-6+\epsilon},&|y|\geq1,
 \end{cases}
 \qquad \rho_i(y)\leq C(1+|y|)^{-4+\epsilon}.
$$
Extend $\rho_o$ by zero outside the rescaled coordinate disk. By
\eqref{5.18}--\eqref{5.19} and dominated convergence,
$$
 \begin{aligned}
 \int_{\C}\rho_o(y)\log|y|\dif A_y
     &\longrightarrow\int_{\C} m_B(y)\log|y|\dif A_y,\\
 \int_{\C}\rho_i(y)\log|y|\dif A_y
     &\longrightarrow\int_{\C} m_i(y)\log|y|\dif A_y.
 \end{aligned}
$$
Proposition~\ref{pr5.8} gives
$$
 \left|\int_{\C}\rho_o\dif A-8\pi\right||\log\delta|
 +\left|\int_{\C}\rho_i\dif A-4\pi\right||\log a|=o(1).
$$
Set $R=|B|$. Circular averaging, integration by parts, and the
substitution $w=y^2$ give
$$
 \begin{aligned}
 \int_\C\frac{\log|w|}{(1+|w-B|^2)^2}\dif A_w
 &=2\pi\int_0^\infty\frac{r\log\max(R,r)}{(1+r^2)^2}\,dr
 =\frac\pi2\log(1+R^2),\\
 \int_\C m_B(y)\log|y|\dif A_y&=2\pi\log(1+|B|^2),\qquad
 \int_\C m_i(y)\log|y|\dif A_y=0.
 \end{aligned}
$$
Thus
$$
 \begin{aligned}
 \int_{B_{r_0}} F_1(z)\log|z|\dif A
 &=\log\delta\int_{\C}\rho_o\dif A+\int_{\C}\rho_o\log|y|\dif A
   =8\pi\log\delta+2\pi\log(1+|B|^2)+o(1),\\
 \int_{B_\tau}F_2(z)\log|z|\dif A
 &=\log a\int_{\C}\rho_i\dif A+\int_{\C}\rho_i\log|y|\dif A
   =4\pi\log a+o(1).
 \end{aligned}
$$
\end{proof}

\begin{proposition}\label{pr5.11}
In the fixed isothermal coordinate, set $\beta_*=e^{\phi(q)}f_*(q)$.
Then
\begin{equation}\label{5.44}
 \frac{a_\nu^2}{\delta_\nu^4}
       \longrightarrow\frac{\beta_*(1+|B|^2)}4.
\end{equation}

\end{proposition}
\begin{proof}
Write $G(q_\nu,z)=-(2\pi)^{-1}\log|z|+R_\nu(z)$. By
Lemmas~\ref{le5.9} and \ref{le5.10},
$$
 \begin{aligned}
 |G(f_{2,\nu}{\bf1}_{M\setminus B_\tau}-f_*)(q_\nu)|
 &\leq\|G(q_\nu,\cdot)\|_{p'}
        \|f_{2,\nu}{\bf1}_{M\setminus B_\tau}-f_*\|_p=o(1),\\
 Gf_{1,\nu}(q_\nu)
 &=-4\log\delta-\log(1+|B|^2)+8\pi R_\nu(0)+o(1),\\
 Gf_{2,\nu}(q_\nu)
 &=-2\log a+4\pi R_\nu(0)+\tfrac12w_*(q_\nu)+o(1).
 \end{aligned}
$$
By \eqref{5.17} and the outer convergence of $f_{2,\nu}$ and $u_{2,\nu}$,
$$
 c_{2,\nu}=\log\rho-\log\int_Mh_2e^{u_{2,\nu}}\dV
 \longrightarrow\log\frac{\rho-4\pi}{\int_Mh_2e^{w_*}\dV}.
$$
Since $u_{2,\nu}=-Gf_{1,\nu}+2Gf_{2,\nu}$,
$$
 \begin{aligned}
 \log4-2\log a+o(1)
 &=U_{2,\nu}(0)\\
 &=c_{2,\nu}+\log h_2(q_\nu)+\phi(q_\nu)+u_{2,\nu}(q_\nu)\\
 &=\log\beta_*+4\log\delta-4\log a+\log(1+|B|^2)+o(1).
 \end{aligned}
$$
Therefore $4a^2=\beta_*\delta^4(1+|B|^2)(1+o(1))$.
\end{proof}

\subsection{Profile estimates and scale matching}\label{sec:strong-chapter}
\subsubsection{A uniform estimate for the scalar linearizations}

Set
$$
 L(y)=\log(2+|y|),\quad
 \|v\|_{\log,R}=\sup_{|y|<R}\frac{|v(y)|}{L(y)},\quad
 \|f\|_{\sigma,R}=\sup_{|y|<R}(1+|y|)^{2+\sigma}|f(y)|.
$$
For the inner and outer profiles, respectively, set
$$
 m_i(y)=\frac4{(1+|y|^2)^2},\qquad
 m_B(y)=\frac{16|y|^2}{(1+|y^2-B|^2)^2},\qquad |B|\le K.
$$
The corresponding operators are $\mathcal L_i=-\Delta-2m_i$ and
$\mathcal L_B=-\Delta-2m_B$. For $\mathcal L_i$, impose
\begin{equation}
 v(0)=0,\qquad \nabla v(0)=0.                       \label{5.45}
\end{equation}
For $\mathcal L_B$, let $Z_{\alpha,B}$, $\alpha=0,1,2$, be the derivatives
with respect to $\log\lambda,\Re B,\Im B$ of
$$
 v_{\lambda,B}(y)=
 \log(16\lambda^2)-2\log(1+\lambda^2|y^2-B|^2)
$$
at $\lambda=1$. Fix a nonnegative $\chi\in C_c^\infty(\C)$ positive on
a sufficiently large fixed disk and impose
\begin{equation}
 \Lambda_{\alpha,B}(v):=\int_{\C}\chi vZ_{\alpha,B}\dif A=0,
 \qquad\alpha=0,1,2.                               \label{5.46}
\end{equation}
We take $D_R$ to contain $\operatorname{supp}\chi$. The three kernel functions are
linearly independent. The Gram matrix $\bigl(\int_{\C}\chi Z_{\alpha,B}Z_{\gamma,B}\dif A\bigr)$ is
uniformly positive definite for $|B|\le K$.

\begin{lemma}\label{le5.12}
Fix $0<\sigma<2$. There exist $R_0,C>0$, uniform for $|B|\le K$, such
that for $R\ge R_0$,
\begin{equation}
 \|v\|_{\log,R}\le
 C\bigl(\|\mathcal L v\|_{\sigma,R}
                 +\osc_{\partial D_R}v\bigr).      \label{5.47}
\end{equation}
Here $\mathcal L=\mathcal L_i$ under \eqref{5.45}, or
$\mathcal L=\mathcal L_B$ under \eqref{5.46}. The estimate holds for
arbitrary boundary mean. For general $v$, add $CN_i(v)$ or $CN_B(v)$
to the right-hand side, respectively, with
\begin{equation}
 N_i(v)=|v(0)|+|\nabla v(0)|,
 \qquad\hbox{or}\qquad
 N_B(v)=\sum_{\alpha=0}^2|\Lambda_{\alpha,B}(v)|.
 \label{5.48}
\end{equation}
The constants are uniform in $R$ and $B$.
\end{lemma}

\begin{proof}
For the outer operator, the bounded kernel is given by
\cite{DEM} with $N=1$. For the inner operator, stereographic
projection to the unit sphere gives $-\Delta_{S^2}v=2v$.
A bounded solution extends across the missing point by subtracting
a bounded local Poisson potential and removing the bounded
harmonic singularity. The corresponding eigenspace consists of the first
spherical harmonics, which give the three inner modes.
Thus in either case
\[
 \ker_{L^\infty}\mathcal L=\operatorname{span}(Z_0,Z_1,Z_2),
 \qquad \det(N_\alpha(Z_\beta))\ne0 ,
\]
where $N_\alpha$ are the functionals in \eqref{5.45} or \eqref{5.46}.
For the inner modes, the determinant follows from their values
and first derivatives at zero. For the outer modes, it follows
from the Gram matrix.
Writing $D_R(y,z)=G_R(y,z)-G_R(0,z)$ for the difference of the
Dirichlet Green kernels, we have
\begin{align}
 \sup_{y\in D_R}L(y)^{-1}
 \int_{D_R}|D_R(y,z)|(1+|z|)^{-2-\sigma}\dif A_z
 &\leq C_\sigma,\label{5.49}\\
 \sup_{R,y}L(y)^{-1}
 \int_{D_R\setminus\overline{D_A}}|D_R(y,z)|m(z)L(z)\dif A_z
 &=:\varepsilon(A)\longrightarrow0.\label{5.50}
\end{align}
The decay $m_i=O(r^{-4})$ and $m_B=O(r^{-6})$ makes these estimates uniform for
$|B|\leq K$.

If \eqref{5.47} fails, choose $R_n\to\infty$ and $v_n$ with
$$
 \|v_n\|_{\log,R_n}=1,\qquad
 \|\mathcal Lv_n\|_{\sigma,R_n}
       +\osc_{\partial D_{R_n}}v_n\longrightarrow0,
 \qquad N_\alpha(v_n)=0.
$$
By elliptic estimates and the logarithmic expansion of the Newton
potential, along a subsequence,
$$
 v_n\longrightarrow\varphi\quad\text{in }C^1_{\rm loc},\qquad
 \mathcal L\varphi=0,\quad
 \varphi=c\log|y|+O(1),\quad N_\alpha(\varphi)=0.
$$
Green's identity with the dilation kernel function $Z_0$ yields
$$
 0=\lim_{r\to\infty}\int_{\partial D_r}
   (Z_0\partial_r\varphi-\varphi\partial_rZ_0)\dif s
   =2\pi cZ_0(\infty),\qquad Z_0(\infty)\ne0.
$$
Hence $c=0$, and the conditions that eliminate the kernel give $\varphi=0$. Let $H_Rv$
be the harmonic extension of $v|_{\partial D_R}$. Then
\begin{equation}\label{5.51}
 v(y)=v(0)+H_Rv(y)-H_Rv(0)
       +\int_{D_R}D_R(y,z)(2m(z)v(z)+\mathcal Lv(z))\dif A_z.
\end{equation}
Using \eqref{5.49}--\eqref{5.50} and the local convergence to zero, we obtain
$$
 \begin{aligned}
 1=\|v_n\|_{\log,R_n}
 &\leq C\bigl(|v_n(0)|+\osc_{\partial D_{R_n}}v_n
               +\|\mathcal Lv_n\|_{\sigma,R_n}\bigr)
       +o_n(1)+C\varepsilon(A)\\
 &=o_n(1)+C\varepsilon(A)\longrightarrow0
       \qquad(n\to\infty,\ A\to\infty),
 \end{aligned}
$$
a contradiction. For general $v$, subtract its kernel component,
$$
 v=v_0+\sum_{\beta=0}^2c_\beta Z_\beta,
 \qquad N_\alpha(v_0)=0,\qquad
 |c|\leq C\sum_\alpha|N_\alpha(v)|,
$$
and apply \eqref{5.47} to $v_0$.
\end{proof}

\begin{remark}\label{re5.13}
Let $N=\|v\|_{\log,R}=o(1)$ and $0<\sigma<2$. Since
$|v(y)|\leq NL(y)$,
$$
 |e^v-1-v|\leq Cv^2e^{|v|}
       \leq CN^2L(y)^2(1+|y|)^N.
$$
Using $m(y)\leq C(1+|y|)^{-4}$, we obtain
\begin{equation}
 \sup_{D_R}(1+|y|)^{2+\sigma}
 m(y)|e^v-1-v|\leq C_\sigma N^2.\label{5.52}
\end{equation}
\end{remark}

\subsubsection{Estimate at the inner scale}

Write $U_i=\log F_i$ and center the chart at the local maximum of
$F_2$ chosen in Proposition~\ref{pr4.19}. Set
$$
 a^2F_2(0)=4,\qquad
 \eta_i(y)=\log(a^2F_2(ay))-\log m_i(y).
$$
Then \eqref{5.45} holds. Under a smooth change of coordinates fixing
the center, normalize $a$ again by $a^2F_2(0)=4$. The Jacobian factor
gives $|\nabla\eta_i(0)|=O(a)$ since $D_y=aD_x$, contributing
$O(a)=o(\delta^{1/4})$ in \eqref{5.48}. Set
$$
 \tau=\delta^{3/2},\qquad R_i=\tau/a.
$$

By \eqref{5.36} and \eqref{5.38}, on $|z|=\tau$,
$$
 P'=O(\delta^{-1}),\qquad
 Q'=-b/z+O(\delta^{1/4}/\tau),\qquad
 b\in\mathbb R,
$$
where $U_2=2Q-P$. Since $b$ is real, the term $-b/z$ makes no
contribution to the tangential derivative of $Q$. Integrating the
tangential derivative of $U_2$ gives
\begin{equation}
 \osc_{|z|=\tau}U_2\le C\delta^{1/4}.              \label{5.53}
\end{equation}
The second equation in \eqref{5.26} gives
\begin{equation}
 \mathcal L_i\eta_i
 =-a^2(F_1(ay)+\kappa_2^c(ay))
       +2m_i(e^{\eta_i}-1-\eta_i).                 \label{5.54}
\end{equation}

\begin{proposition}\label{pr5.14}
Under the local hypotheses of Lemma~\ref{le5.6}, uniformly for
$|y|\le R_i$,
\begin{equation}
 |\eta_i(y)|\le C\delta^{1/4}\log(2+|y|).          \label{5.55}
\end{equation}
\end{proposition}

\begin{proof}
The circular mean identity \eqref{5.27} and the oscillation bound
\eqref{5.20} give, for every $\epsilon>0$,
$$
 \begin{gathered}
 |\bar\eta_i(r)-\bar\eta_i(R_\epsilon)|
       \leq\epsilon\log(r/R_\epsilon),\qquad
 \osc_{|y|=r}\eta_i\leq C,\qquad R_\epsilon\leq r\leq R_i,\\
 \limsup_{\delta\to0}\sup_{R\leq|y|\leq R_i}
       \frac{|\eta_i(y)|}{\log(2+|y|)}
       \leq\epsilon+\frac C{\log(2+R)},\qquad
 \sup_{|y|\leq R}|\eta_i(y)|\longrightarrow0.
 \end{gathered}
$$
Letting $\delta\to0$, then $R\to\infty$ and $\epsilon\downarrow0$
gives $N_i:=\|\eta_i\|_{\log,R_i}=o(1)$. For the quantitative forcing
bound, \eqref{eq:inner-first-density-small} first justifies the bounded
Dirichlet right-hand side used to obtain \eqref{5.24}.
Use \eqref{5.21} and \eqref{5.24} of Lemma~\ref{le5.6} with
$0<\epsilon<1/10$ to get
$$
 a^2F_1(ay)\leq C\delta^{-4+\epsilon}a^{4-\epsilon}
                          (1+|y|)^{2-\epsilon},
 \quad a=\delta^{2+o(1)},\quad R_i=\delta^{-1/2+o(1)}.
$$
Consequently,
\begin{equation}\label{5.56}
 \sup_{|y|<R_i}(1+|y|)^3a^2
       (F_1(ay)+|\kappa_2^c(ay)|)
 \leq C\delta^{3/2-\epsilon/2+o(1)}\leq C\delta^{4/3}.
\end{equation}
Lemma~\ref{le5.12} and \eqref{5.52}--\eqref{5.54} imply
$$
 N_i\leq C(\delta^{1/4}+\delta^{4/3}+a)+CN_i^2
     \leq C\delta^{1/4}+\tfrac12N_i,
 \qquad N_i\leq C\delta^{1/4}.
$$
\end{proof}

\subsubsection{Estimate at the outer scale}

Set $s=\delta^{1/2}$ and let $\Phi$ solve the zero Dirichlet problem on
$D_s$ for
\begin{equation}
 \Delta\Phi=F_2-m_i^a+\kappa_1^c,\qquad
 m_i^a(z)=\frac{4a^2}{(a^2+|z|^2)^2}.             \label{5.57}
\end{equation}

\begin{lemma}\label{le5.15}
Under the local hypotheses of Lemma~\ref{le5.6}, with
$\ell=\log(1/\delta)$,
\begin{equation}
 \|\Phi\|_{L^\infty(D_s)}
       \le C\delta^{1/4}\ell^2.                   \label{5.58}
\end{equation}
\end{lemma}
\begin{proof}
On $D_\tau$, \eqref{5.55} gives
$$
 |F_2-m_i^a|\leq C\delta^{1/4}\log(2+|z|/a)m_i^a,
 \qquad \delta^{1/4}\log(2+\tau/a)=o(1).
$$
For the positive Dirichlet Green kernel $G_s$, put
$\mathcal P_E(g)=\sup_{x\in D_s}\int_EG_s(x,z)|g(z)|\dif A_z$.
The bounds of Lemma~\ref{le5.6} and the finite logarithmic moments of
$m_i$ give
$$
 \begin{aligned}
 \mathcal P_{D_\tau}(F_2-m_i^a)&\leq C\delta^{1/4}\ell^2,\\
 \mathcal P_{D_{\delta/R}\setminus D_\tau}(F_2)
     &\leq C(a/\tau)^{2-\epsilon}\ell,\\
 \mathcal P_{D_{R\delta}\setminus D_{\delta/R}}(F_2)
     &\leq C\delta^{2-\epsilon}\ell,\\
 \mathcal P_{D_s\setminus D_{R\delta}}(F_2)
     &\leq Cs^{2-\epsilon}\ell,\\
 \mathcal P_{D_s\setminus D_\tau}(m_i^a)
       +\mathcal P_{D_s}(\kappa_1^c)
     &\leq C((a/\tau)^2+s^2)\ell.
 \end{aligned}
$$
Near the logarithmic singularity at $z=x$, use the pointwise density
bounds. Away from it, use the mass estimates on the indicated annuli. Since
$$
 (a/\tau)^{2-\epsilon}=\delta^{1-\epsilon/2+o(1)},
 \qquad s^{2-\epsilon}=\delta^{1-\epsilon/2},
$$
Green's representation yields
$$
 \|\Phi\|_\infty
 \leq C(\delta^{1/4}\ell^2+(a/\tau)^{2-\epsilon}\ell
                   +\delta^{2-\epsilon}\ell+s^{2-\epsilon}\ell)
 \leq C\delta^{1/4}\ell^2.
$$
\end{proof}

Put $t=a/\delta$ and
$$
 v(y)=U_1(\delta y)+2\log\delta
       -\log(|y|^2+t^2)-\Phi(\delta y),
 \qquad |y|<R_o:=s/\delta .
$$
Using $\Delta\log(a^2+|z|^2)=m_i^a$ in \eqref{5.57} and the first
equation of \eqref{5.26}, we obtain
\begin{equation}
 -\Delta v=2(|y|^2+t^2)e^{\Phi(\delta y)}e^v.       \label{5.59}
\end{equation}
Write $v_B=\log16-2\log(1+|y^2-B|^2)$.
On a fixed small disk $D_d$, the local profile limit
\eqref{5.19} and the bound $\|\Phi\|_\infty=o(1)$ from
Lemma~\ref{le5.15} combine with Lemma~\ref{le5.6} and \eqref{5.59} to give
$$
 -\Delta v=2\delta^2F_1(\delta y)
        \leq C(|y|+t)^{2-\epsilon},\qquad
 \|v\|_{L^\infty(\partial D_d)}\leq C .
$$
Interior estimates for the Dirichlet problem give uniform $C^{1,\alpha}$
bounds on $D_{d/2}$ for $0<\alpha<1$. Together with the convergence
away from zero, this gives $v\to v_B$ in $C^1_{\mathrm{loc}}(\C)$.

We adjust $p=(\log\delta,B_1,B_2)$ to impose \eqref{5.46}. The limiting
derivatives with respect to these parameters are
$$
 \begin{gathered}
 T_0=4+y\cdot\nabla v_B,\qquad
 T_1=\partial_{B_1}v_B,\qquad T_2=\partial_{B_2}v_B,\\
 \operatorname{span}(T_0,T_1,T_2)
     =\operatorname{span}(Z_{0,B},Z_{1,B},Z_{2,B}).
 \end{gathered}
$$
On a fixed disk, set $\mathcal F_\alpha(p)
=\Lambda_{\alpha,B}(v-v_B)$. Local convergence gives
$$
 \mathcal F(p_0)=o(1),\qquad
 D_p\mathcal F(p_0)=M_B+o(1),\qquad
 \det M_B\ne0.
$$
Here $M_B$ is the matrix of the functionals $\Lambda_{\alpha,B}$
applied to $T_0,-T_1,-T_2$. With $a$ and $\Phi$ fixed,
$$
 v_{e^\vartheta\delta}(y)=v_\delta(e^\vartheta y)+4\vartheta .
$$
The local $C^1$ convergence therefore gives uniform convergence of the
parameter derivatives on a fixed neighborhood of $p_0$ in the
coordinates $(\vartheta,B'-B)$. The implicit function theorem gives
$p'=p_0+o(1)$ such that
$$
 \Lambda_{\alpha,B'}(v'-v_{B'})=0\quad(\alpha=0,1,2),
 \qquad \delta'/\delta\longrightarrow1,\quad B'-B=o(1).
$$
We keep $s$ and $\tau$ fixed during this adjustment, so that
$s\asymp (\delta')^{1/2}$ and $\tau\asymp (\delta')^{3/2}$.

\begin{proposition}\label{pr5.16}
Under the local hypotheses of Lemma~\ref{le5.6}, the adjusted
parameters satisfy
\begin{equation}
 |v(y)-v_B(y)|
 \le C\delta^{1/4}\ell^2\log(2+|y|),
 \qquad |y|<R_o.                                  \label{5.60}
\end{equation}
\end{proposition}
\begin{proof}
The adjustment preserves local convergence. As in the proof of
Proposition~\ref{pr5.14}, we use \eqref{5.27} and \eqref{5.20} to obtain
$$
 \limsup_{\delta\to0}\sup_{R\leq|y|<R_o}
 \frac{|v-v_B|}{\log(2+|y|)}
 \leq\epsilon+\frac C{\log(2+R)},
 \qquad N_o:=\|v-v_B\|_{\log,R_o}=o(1).
$$
The estimates on the outer circle in the proof of
Proposition~\ref{pr5.8} give
$$
 \begin{gathered}
 P'=-\alpha/z+O(\delta^{1/2}/s),\qquad
 Q'=-b/z+O(\delta^{1/2}/s),\qquad \alpha,b\in\R,\\
 \Phi=0,\qquad \osc_{\partial D_s}U_1\leq C\delta^{1/2},\qquad
 \osc_{\partial D_{R_o}}v_B\leq CR_o^{-2},\\
 \osc_{\partial D_{R_o}}(v-v_B)\leq C\delta^{1/2}.
 \end{gathered}
$$
For $\eta_o=v-v_B$, equations \eqref{5.52}, \eqref{5.58},
and \eqref{5.59} yield
$$
 \begin{aligned}
 \mathcal L_B\eta_o
 &=2e^{v_B}\bigl[(|y|^2+t^2)e^{\Phi(\delta y)+\eta_o}
                       -|y|^2(1+\eta_o)\bigr],\\
 \|\mathcal L_B\eta_o\|_{1,R_o}
 &\leq C(\delta^{1/4}\ell^2+\delta^{2-o(1)}+N_o^2).
 \end{aligned}
$$
Since $\Lambda_{\alpha,B}(\eta_o)=0$, Lemma~\ref{le5.12} gives
$$
 N_o\leq C(\delta^{1/4}\ell^2+\delta^{2-o(1)}+\delta^{1/2})
          +CN_o^2
     \leq C\delta^{1/4}\ell^2+\tfrac12N_o.
$$
Thus $N_o\leq C\delta^{1/4}\ell^2$.
\end{proof}

\subsubsection{Green potential estimates and scale matching}

Multiply $m_i^a$ and
$$
 m_o^{\delta,B}(z)
 =\frac{16\delta^4|z|^2}
        {(\delta^4+|z^2-\delta^2B|^2)^2}
$$
by fixed radial cutoffs supported in $D$ and equal to one on a
smaller disk. In the Green integrals, $m_i^a$ and $m_o^{\delta,B}$
denote the corresponding metric densities, so
$m_g\dV=m\dif A$ as above. We write $Gf$ for the mean-zero Green
potential. The tails beyond the fixed cutoff radii have masses bounded
by positive powers of $\delta$.

\begin{lemma}\label{le5.17}
Under the local hypotheses of Lemma~\ref{le5.6}, for some fixed
integer $N$,
\begin{align}
 \|G(1_{D_\tau}f_2-m_i^a)\|_\infty
       &\le C\delta^{1/4}\ell^N,\label{5.61}\\
 \|G(\mathbf1_Df_1-m_o^{\delta,B})\|_\infty
       &\le C\delta^{1/4}\ell^N .                 \label{5.62}
\end{align}
The estimate \eqref{5.61} also holds with $\tau$ replaced by any
radius comparable to $\delta^{3/2}$.
\end{lemma}
\begin{proof}
Applying the logarithmic kernel bound used in the proof of
Lemma~\ref{le5.15} to \eqref{5.55} and to the inner profile outside
$D_\tau$ gives
$$
 \|G({\bf1}_{D_\tau}f_2-m_i^a)\|_\infty
 \leq C(\delta^{1/4}\ell^2+(a/\tau)^2\ell).
$$
By \eqref{5.21}, changing $\tau$ by a bounded factor changes the
potential by at most $C(a/\tau)^{2-\epsilon}\ell$. For the outer
profile, put $\varepsilon_\delta=\delta^{1/4}\ell^2$.
On $D_s$, equations \eqref{5.59}--\eqref{5.60} give
$$
 \begin{aligned}
 F_1(\delta y)&=\delta^{-2}(|y|^2+t^2)e^{v_B(y)}
                         e^{\Phi(\delta y)+\eta_o(y)},\\
 \|\Phi\|_\infty+\|\eta_o\|_\infty
       &\leq C\delta^{1/4}\ell^3=o(1),\\
 |F_1(\delta y)-\delta^{-2}m_B(y)|
       &\leq C\delta^{-2}
          (\varepsilon_\delta m_B(y)L(y)+t^2e^{v_B(y)}).
 \end{aligned}
$$
For $h=m_BL$ or $h=e^{v_B}$, the finite logarithmic moments give
$$
 \sup_x\int_{D_s}|G(x,z)|\delta^{-2}h(z/\delta)\dif A_z
 \leq C\ell.
$$
On $D\setminus D_s$, the pointwise bounds of Lemma~\ref{le5.6} give
$$
 \|G({\bf1}_{D\setminus D_s}(f_1-m_o^{\delta,B}))\|_\infty
 \leq C\delta^{2-\epsilon/2}\ell.
$$
Hence
$$
 \|G(\mathbf1_Df_1-m_o^{\delta,B})\|_\infty
 \leq C(\varepsilon_\delta\ell+t^2\ell
                       +\delta^{2-\epsilon/2}\ell)
 \leq C\delta^{1/4}\ell^N.
$$
\end{proof}

For the rest of this subsection, assume the single-point condition
(\ref{apriori:formation-measures}). The exterior bound of
Lemma~\ref{le5.6} gives
$\|G(\mathbf1_{M\setminus D}f_1)\|_\infty=O(\delta^{4-\epsilon})$,
so \eqref{5.62} holds with $f_1$ in place of $\mathbf1_Df_1$. Set
$$
 w_\delta=2G(1_{M\setminus D_\tau}f_2),\quad
 T_\delta=2G(m_i^a)-G(m_o^{\delta,B}).
$$
By the second component of \eqref{2.1} and \eqref{5.61}--\eqref{5.62},
\begin{equation}
 u_2=w_\delta+T_\delta+E_\delta,\qquad
 \|E_\delta\|_\infty\le C\delta^{1/4}\ell^N.       \label{5.63}
\end{equation}
By Lemma~\ref{le5.9} and elliptic regularity,
$$
 \|w_\delta-w_*\|_{W^{2,6/5}}
 \leq C\|{\bf1}_{M\setminus D_\tau}f_2-f_*\|_{6/5}=o(1).
$$

The Green potentials of the planar profiles give, on a fixed disk,
\begin{equation}
 T_\delta(z)=
 \log\frac{\delta^4+|z^2-\delta^2B|^2}
                {(a^2+|z|^2)^2}
       +O(\delta^2\ell).                          \label{5.64}
\end{equation}
The regular Green terms proportional to the masses cancel by
$2(4\pi)-8\pi=0$. The vanishing first moments and the second moment bounds give
errors of $O(a^2\ell)$ and $O(\delta^2)$.
Outside the disk, $T_\delta=O(\delta^2\ell)$. Thus, for $p=6/5$,
\begin{equation}
 \|1_{M\setminus D_\tau}(e^{T_\delta}-1)\|_{L^p}
 \le C(\delta^2\ell+
        \delta^2\tau^{2/p-2}+
        \delta^4\tau^{2/p-4})
 \le C\delta^{1/2}\ell.                           \label{5.65}
\end{equation}

Write $A_\delta=\rho_2/\int_Mh_2e^{u_2}\dV$.
The smooth convergence away from $q$ gives
$A_\delta\to\sigma/I_*$, where
$\sigma=\rho_2-4\pi$ and $I_*=\int_{M} h_2e^{w_*}\dV$.
Since $w_\delta$ is uniformly bounded, \eqref{5.63}--\eqref{5.65} yield
\begin{equation}
 \|1_{M\setminus D_\tau}f_2
               -A_\delta h_2e^{w_\delta}\|_{L^{6/5}}
       \le C\delta^{1/4}\ell^N.                   \label{5.66}
\end{equation}
Estimating the contribution on $D_\tau$ and integrating \eqref{5.66}, we obtain
$$
 \begin{aligned}
 \|{\bf1}_{D_\tau}A_\delta h_2e^{w_\delta}\|_{6/5}
     &\leq C\tau^{5/3},\\
 A_\delta\int_Mh_2e^{w_\delta}\dV
     &=\sigma+O(\delta^{1/4}\ell^N),\\
 \left\|-\Delta_gw_\delta
       -2\sigma\left(\frac{h_2e^{w_\delta}}
                    {\int_Mh_2e^{w_\delta}\dV}-1\right)\right\|_{6/5}
     &\leq C\delta^{1/4}\ell^N.
 \end{aligned}
$$
Assume that the scalar linearization at $w_*$ is invertible on the
mean-zero space. Set $e_\delta=w_\delta-w_*$. Apply the inverse and
use $\|e_\delta\|_\infty=o(1)$ to obtain
$$
 \|e_\delta\|_{W^{2,6/5}}
 \leq C\delta^{1/4}\ell^N
        +C\|e_\delta\|_\infty\|e_\delta\|_{W^{2,6/5}}
 \leq C\delta^{1/4}\ell^N+\tfrac12\|e_\delta\|_{W^{2,6/5}}.
$$
Consequently,
\begin{equation}\label{5.67}
 \|w_\delta-w_*\|_{W^{2,6/5}}
       +|A_\delta-\sigma/I_*|\leq C\delta^{1/4}\ell^N.
\end{equation}

Evaluating \eqref{5.63}--\eqref{5.64} at the inner maximum and using
$a^2F_2(0)=4$, we obtain
\begin{equation}
 \frac{4a^2}{\delta^4(1+|B|^2)}
   =\beta_{*,\nu}\bigl(1+O(\delta^{1/4}\ell^N)\bigr),
 \qquad \beta_{*,\nu}=e^{\phi_\nu(0)}f_*(q_\nu).
 \label{5.68}
\end{equation}

Set
$$
 \mathcal P_{\delta,a}=
 \binom{2G(m_o)-G(m_i^a)-w_*/2}{2G(m_i^a)-G(m_o)+w_*}.
$$
By \eqref{2.1}, \eqref{5.61}--\eqref{5.62}, and \eqref{5.67},
$$
 \|u-\mathcal P_{\delta,a}\|_\infty\leq C\delta^{1/4}\ell^N.
$$
Let $\zeta$ be the coordinate of Proposition~\ref{pr5.3} at $q_\nu$,
and write the change of coordinates as $z=g_\nu(\zeta)$ with
$g_\nu'(0)=\lambda_\nu>0$. Put
$$
 \bar\delta=\delta/\lambda_\nu,\qquad
 \bar a=a/\lambda_\nu,\qquad
 \bar\beta=\lambda_\nu^2\beta_{*,\nu}=f_*(q_\nu).
$$
Both $\lambda_\nu$ and $\lambda_\nu^{-1}$ are bounded. To express the
Green potentials in the new coordinate, set
$J_\nu(\zeta)=g_\nu(\zeta)/\zeta=\lambda_\nu+O(|\zeta|)$.
The regular parts satisfy
$R^{z}(g_\nu(\zeta),0)=R^{\zeta}(\zeta,0)+(2\pi)^{-1}\log|J_\nu(\zeta)|$.
The two denominators in the logarithmic terms become
$$
 |\zeta|^2+a^2/|J_\nu|^2,\qquad
 \delta^4/|J_\nu|^4+
       |\zeta^2-\delta^2B/J_\nu^2|^2.
$$
Replacing $J_\nu$ by $\lambda_\nu$ gives relative errors
$O(a)$ and $O(\delta)$, respectively. Indeed, with
$r=|\zeta|$, the denominator comparisons reduce to
$a^2r/(a^2+r^2)=O(a)$ and
$\delta^2r(\delta^2+r^2)/(\delta^4+r^4)=O(\delta)$.

The tails beyond the fixed cutoff radii and the regular Green terms
satisfy the same bounds as in \eqref{5.64}. Thus the $C^0$ estimate for
$u-\mathcal P$ also holds in the normalized chart with parameters
$(\bar\delta,\bar a,B)$.

Set $\bar a_*^2=\bar\beta\bar\delta^4(1+|B|^2)/4$.
We apply \eqref{5.1}--\eqref{5.8} at $q_\nu$ in this chart with
inner scale $\bar a_*$. For this comparison, we evaluate the normalized
densities $V_i^{\rm app}$ and residual $\mathcal R^{\rm app}$ at
$\rho=\rho^0=(8\pi,\rho_2)$, that is, at $\mu=0$. The approximation
$u^{\rm app}$ in \eqref{5.1} is independent of $\mu$. Let $\widehat V_i$
be the densities before normalization. Then
$$
 \begin{aligned}
 \|V_1^{\rm app}-m_{o,g}\|_{6/5}
       +\|V_2^{\rm app}-m_{i,g}^{\bar a_*}-n_*^g\|_{6/5}
      &\leq C\delta^{4/3},\\
 \|\mathcal R^{\rm app}\|_{6/5}
      &\leq C\delta^{4/3}
       +C\delta^2(\|\widehat V_1\|_{6/5}+\|\widehat V_2\|_{6/5})\\
      &\leq C\delta^{4/3}
       +C\delta^2(\delta^{-1/3}+\delta^{-2/3})
       \leq C\delta^{4/3}.
 \end{aligned}
$$
The cutoffs change the two planar masses by $O(\delta^4)$, so the
constant terms in $-\Delta_g\mathcal P_{\bar\delta,\bar a_*}$ agree
with those at $\rho^0$ up to $O(\delta^4)$. Since both fields have
zero mean, elliptic regularity gives
$$
 \|u^{\rm app}-\mathcal P_{\bar\delta,\bar a_*}\|_\infty
 \leq C\|\Delta_g(u^{\rm app}-\mathcal P_{\bar\delta,\bar a_*})\|_{6/5}
 \leq C\delta^{4/3}.
$$
Equations \eqref{5.64} and \eqref{5.68} imply
$$
 \begin{aligned}
 |\log(\bar a/\bar a_*)|&\leq C\delta^{1/4}\ell^N,\\
 \|G(m_{i,g}^{\bar a}-m_{i,g}^{\bar a_*})\|_\infty
 &\leq C|\log(\bar a/\bar a_*)|+C\delta^2\ell
 \leq C\delta^{1/4}\ell^N.
 \end{aligned}
$$
Combining the estimates, we may take $c=1/8$ in
\begin{equation}\label{5.69}
 \|u-u^{\rm app}_{\bar\delta,q_\nu,B,\bar\beta}\|_{C^0(M)}
 \leq C\delta^{1/4}\ell^N\leq C\delta^c.
\end{equation}

\begin{lemma}
\label{le5.18}
Let $0<r<e^{-1}$ and set $L=\log(1/r)$. Suppose $M_r,E_r\geq0$ satisfy,
for some $\alpha>0$ and constants independent of $r$,
$$
 M_r=O(r^\alpha),\qquad
 E_r\leq C(r^{4/3}L+M_rL^{3/2}),
$$
and
$$
 M_r\leq C(r^2+r\sqrt L\,E_r+LE_r^2+M_r^2).
$$
Then $M_r=O(r^2)$ and $E_r=O(r^{4/3}L)$.
The estimates are uniform when the constants in the hypotheses are
uniform and $\alpha$ is fixed.
\end{lemma}
\begin{proof}
Substituting the estimate for $E_r$ into the inequality for $M_r$ gives
$$
 M_r\leq C(r^2+r^{7/3}L^{3/2}+r^{8/3}L^3)
                +CrL^2M_r+CM_r^2L^4.
$$
Since $M_rL^4=O(r^\alpha L^4)=o(1)$ and $rL^2=o(1)$, we can absorb
the last two terms to obtain
$$
 M_r\leq Cr^2+\tfrac12M_r,\qquad M_r\leq Cr^2.
$$
The estimate for $E_r$ now gives
$$
 E_r\leq C(r^{4/3}L+r^2L^{3/2})\leq Cr^{4/3}L.
$$
\end{proof}

\section{The linearized operator with common normalization}
\label{rev:sec:multiple}\label{sec:mc-chapter}

We use the energy product and weighted averages of
Section~\ref{sec:degree-conventions}, with the scales of
Section~\ref{sec:scale-conventions}.

For the projected equations in both local constructions, we need the
inverse estimate and the negative index in \eqref{6.3}--\eqref{6.4}.
These estimates hold uniformly on the sets of separated centers and
scales specified below. Comparison with the completed weights also
determines the parameter ranges in which the estimates apply.

\subsection{The model operator}

Set $N=a+b\geq1$ and $L=\log(1/s)$, and assume that the centers stay
uniformly separated. The scales at the ordinary centers and the inner
scales at the mixed centers are comparable to $s$. The outer scales at
the mixed centers are comparable to $\sqrt s$. The scale ratios lie in
a compact subset of $(0,\infty)$ and $|B_k|\leq C$.
The Euclidean profile at an ordinary center and at the inner scale of a
mixed center is $p_{\rm in}(y)=4(1+|y|^2)^{-2}$. At the outer scale it is
\[
 p_B(y)=\frac{16|y|^2}{(1+|y^2-B|^2)^2}.
\]
At each center, choose a fixed cutoff $0\leq\chi\leq1$ that equals one
near the center and has support disjoint from the other cutoffs.
For a profile $p$ of width $r$, define
\[
 m\dV=\alpha_r\chi(x)r^{-2}p(x/r)\dif A_x,
 \qquad
 \alpha_r=\frac{\int_{\R^2}p\dif A}
                 {\int\chi(x)r^{-2}p(x/r)\dif A_x}.
\]
Here $x$ is the coordinate centered at the point in question. Let
$m_j$ denote this cutoff density for the ordinary profile at an
ordinary center and for the outer profile at a mixed center. At
each mixed center we use the same cutoff at both scales and write
$\widetilde m_k$ for the cutoff inner density. The density $n_*$
is the residual density of mass $\sigma$ for the given scalar
solution at fixed centers, with the normalization of \eqref{2.5}.
Set
$$
 V_1=\sum_{j=1}^Nm_j,\qquad
 V_2=\sum_{k\in\typeB}\widetilde m_k+n_*.
$$
The masses are
\[
 M_j=\int_Mm_j\dV=
 \begin{cases}4\pi&j\in\typeA,\\8\pi&j\in\typeB,\end{cases}
 \qquad
 \int_M\widetilde m_k\dV=4\pi,\qquad \int_Mn_*\dV=\sigma>0.
\]
The smooth density $n_*$ is positive on an open set and vanishes at the
ordinary centers prescribed by the shadow solution. It varies smoothly
over the given compact set of scalar solutions. The scalar form
$$
 \frac12\int_M|\nabla v|^2\dV
       -\int_Mn_*(v-\avg v{n_*})^2\dV
$$
is nondegenerate on mean-zero functions and has index $m_*$.
It is one half of the scalar form $Q_{n_*}$ introduced before
Lemma~\ref{le3.9}, and a weak null vector
satisfies $-\Delta_gv=2n_*(v-\avg v{n_*})$. The coefficient of this
scalar equation is smooth, so elliptic bootstrapping puts the vector
in the fixed-center kernel in $X_s$, which is trivial.
The sign in \eqref{eq:shadow-local-index} uses this same
scalar index.

For $\eta\in\mathcal H$ define
$$
 \begin{aligned}
 \mu_i(\eta)&=\avg{\eta_i}{V_i},\\
 I_s(\eta,\xi)&=\sum_i\int_MV_i(\eta_i-\mu_i(\eta))
                               (\xi_i-\mu_i(\xi))\dV,\\
 \qform_s(\eta,\xi)&=\ipA\eta\xi-I_s(\eta,\xi).
 \end{aligned}
$$
For a density $m$ of width $r$ and a kernel function $Z$, set
$g=\chi Z(x/r)$ and $m\widehat Z=m(g-\avg{g}{m})$. At an ordinary
center and at the inner scale of a mixed center, the translation
functions are $y_l/(1+|y|^2)$, $l=1,2$. At the outer scale, the angular
functions are the real and imaginary parts of
$(y^2-B)/(1+|y^2-B|^2)$. The radial functions are
\[
 Z_{\typeA}=Z_{\rm in}=\frac{1-|y|^2}{1+|y|^2},\qquad
 Z_{0,B}=\frac{1-|y^2-B|^2}{1+|y^2-B|^2}.
\]
The ordinary translation sources and the outer angular sources enter
the first component. The inner translation sources at the mixed centers
enter the second component. The common radial source in the first
component is $\sum_jm_j\widehat Z_j$.

Let $F_{\alpha,s}$ denote these vector sources, and let
$S_{\alpha,s}=GA F_{\alpha,s}$ be their Riesz representatives with respect
to the energy inner product. Define
\[
 E_s=\operatorname{span}(S_{1,s},\ldots,S_{2N+2b+1,s})^{\perp_A},
 \qquad P_s:\mathcal H\longrightarrow E_s
\]
where $P_s$ is the orthogonal projection for the energy inner product.
With
\begin{equation}\label{6.1}
 d_j=M_j/3,\qquad e=4\pi/3,\qquad t_k=(1+|B_k|^2)^{-1},
\end{equation}
the common radial constraint is
\begin{equation}\label{6.2}
 \sum_j\int_Mm_j\widehat Z_j\eta_1\dV=0.
\end{equation}
Define
$$
 \mathcal A_s\eta
 =\eta-GA
 \binom{V_1(\eta_1-\mu_1(\eta))}
       {V_2(\eta_2-\mu_2(\eta))},
 \qquad \mathcal L_s=P_s\mathcal A_s|_{E_s}.
$$
Equivalently,
$$
 \ipA{\mathcal L_s\eta}\xi=\qform_s(\eta,\xi),
 \qquad \eta,\xi\in E_s.
$$

\begin{theorem}
\label{th6.1}
For sufficiently small $s$, uniformly on the stated compact parameter sets,
\begin{equation}\label{6.3}
 \operatorname{ind}\mathcal L_s=m_*+a+4b-1,\qquad
 \operatorname{dist}(0,\operatorname{spec}\mathcal L_s)\geq c/L,
\end{equation}
and
\begin{equation}\label{6.4}
 \|\mathcal L_s^{-1}f\|_{H^1}\leq CL\|f\|_{H^{-1}}.
\end{equation}
The same bounds hold for self-adjoint perturbations of size $o(L^{-1})$
after the projected subspaces are identified isometrically.
\end{theorem}

For fixed $s>0$, the weights are bounded and the centered
forms factor through the compact embedding of $H^1$ into $L^2$.
Thus on the projected energy space we have $\mathcal L_s=I-K_s$,
where $K_s$ is compact and self-adjoint. We use its orthogonal spectral
decomposition. For the variational comparison below, see \cite[Chapter~V, Section~2.3, Theorem~2.10]{Kato}.
In particular, the spectrum below any $c<1$ consists of finitely
many eigenvalues counted with multiplicity. Any subspace on which
all Rayleigh quotients are strictly below $c$ projects injectively
onto the spectral space below $c$, since the quadratic form on its
orthogonal complement is at least $c$ times the norm squared.
Taking the spectral space itself gives the converse dimension
bound. This is the min--max comparison used below.

\paragraph{Green potentials of the centered sources.}
Write $p(y)\leq C(1+|y|)^{-2-\alpha}$, where $\alpha=2$ at the ordinary
and inner scales and $\alpha=4$ at the outer scales. Set
$Z_\infty=\lim_{|y|\to\infty}Z(y)$.
The profile kernel functions satisfy $-\Delta Z=2pZ$,
$\int pZ=0$, and $|Z-Z_\infty|+|y||\nabla Z|=O(|y|^{-1})$
at infinity. These estimates hold uniformly for bounded $B$.

For $\Psi_{m,Z}=2G(m\widehat Z)$ we have
\[
 \Psi_{m,Z}=\chi(x)(Z(x/r)-Z_\infty)+O(r+r^\alpha L)
 \qquad\text{in }L^\infty(M).
\]
Indeed, $\alpha_r=1+O(r^\alpha)$ and $\avg{g}{m}=O(r^\alpha)$
by the decay of $p$. Applying $-\Delta_g$ to
$\chi(Z(x/r)-Z_\infty)$ gives $2e^{-\phi}\chi r^{-2}p(x/r)Z(x/r)$
up to an error bounded by $Cr$ on the cutoff annulus.

Replacing this leading term by $2m\widehat Z$ gives a cutoff error
bounded by $Cr^\alpha$ on the same annulus. The normalization and
centering errors are bounded by $Cr^\alpha m$. Since
$\sup_x\int_M |G(x,y)|m(y)\dV(y)\leq CL$, their Green potentials are
$O(r^\alpha L)$. The mean of $\chi(Z-Z_\infty)$ is $O(r)$ by the
decay of $Z-Z_\infty$.

In particular, $S_{\alpha,s}$ is uniformly bounded in $L^\infty(M)$.
We compute its energy pairing with $S_{\beta,s}$ by integrating the
potential against the corresponding source. Since $a_{ii}=2$, a pairing
at the same scale converges to the weighted $L^2$ product of the
kernel functions.

For distinct centers, each potential is $o(1)$ on the support of the
other source. At a mixed center, the outer potential converges to a
constant at the inner scale, while the inner potential converges to
zero away from the origin at the outer scale. The sources have zero
mass, so bounded convergence gives cross terms of order $o(1)$.

At each scale, the translation, radial and angular functions are
orthogonal in the profile-weighted $L^2$ space and have positive
norms. At the outer scale, this follows by pulling back the spherical
density under the two-to-one map $y\mapsto y^2-B$. The common radial
entry is the sum of the contributions at the first-component scales.
This gives the Gram bound in \eqref{6.17}. The $L^\infty$ estimate and
the mass bounds give the variance bounds in the same equation.

\subsection{Negative eigenvalues}

The exact variance identity
\begin{equation}\label{6.5}
 \int_{M}\Big(\sum_\alpha w_\alpha\Big)(v-\bar c)^2\dV
 =\sum_\alpha\int_{M} w_\alpha(v-c_\alpha)^2\dV
      +\sum_\alpha r_\alpha(c_\alpha-\bar c)^2
\end{equation}
holds for $r_\alpha=\int_{M} w_\alpha\dV$, $c_\alpha=\avg v{w_\alpha}$,
and $\bar c=\sum r_\alpha c_\alpha/\sum r_\alpha$. The weighted
Poincar\'e estimate following \eqref{6.7} bounds the first sum uniformly
by the energy. The second sum defines a finite-rank form.

For $V_1$, the matrix of this form is
$$
 C_1=\operatorname{diag}(M_j)-\frac{(M_j)_j(M_j)_j^T}{\sum M_j},
 \qquad \operatorname{rank}C_1=N-1.
$$
For $V_2$, express the form in terms of the $b$ differences between
the inner profile means and the mean with respect to $n_*\dV$. Its
matrix is
$$
 C_2=4\pi I_b-\frac{16\pi^2}{4\pi b+\sigma}\mathbf1\mathbf1^T>0.
$$
Since $\sigma>0$,
$$
 \operatorname{rank}C_1+\operatorname{rank}C_2
       =(N-1)+b=:r.
$$

The energy Gram matrix of the local mean functionals is
\begin{equation}\label{6.6}
 H_s=LH+O(1),\qquad
 H=\bigoplus_{\typeA}\frac1\pi
       \ \oplus\!
       \bigoplus_{\typeB}\frac1{4\pi}
           \begin{pmatrix}2&-1\\-1&4\end{pmatrix}>0.
\end{equation}
For the means of components $i,j$ at the same center and at widths
$r_0,r_1$, Green's representation gives
$$
 H_s[(i,r_0),(j,r_1)]
 =-\frac{a_{ij}}{2\pi}\log\max(r_0,r_1)+O(1).
$$
The terms involving distinct centers or $n_*$ are $O(1)$. The bounded
source potentials and \eqref{6.17} give an $O(1)$ change in $H_s$ upon
projection onto $E_s$. Thus
$$
 \begin{gathered}
 \mathcal L_s=\mathcal B_s-\mathcal R_s,\qquad
 \|\mathcal B_s\|\leq C,\qquad \mathcal R_s\geq0,\\
 \operatorname{rank}\mathcal R_s=r,
 \qquad \operatorname{spec}(\mathcal R_s)\cap(0,\infty)
                 \subset[cL,CL].
 \end{gathered}
$$
Order the eigenvalues increasingly. The min--max principle gives
$$
 -CL\leq\tau_1\leq\cdots\leq\tau_r\leq-cL,
 \qquad \operatorname{spec}(\mathcal L_s)\cap(-cL,\infty)\subset[-C,C].
$$
Put $Q_s=\operatorname{proj}\operatorname{ran}\mathcal R_s$ and
$\Pi_s^-=\mathbf1_{(-\infty,-cL]}(\mathcal L_s)$. Then
$$
 \begin{aligned}
 \mathcal L_s\eta=\tau\eta,\quad |\tau|\geq cL
 &\ \Longrightarrow\ 
 \tau(I-Q_s)\eta=(I-Q_s)\mathcal B_s\eta,\\
 \|(I-Q_s)\eta\|_A&\leq CL^{-1}\|\eta\|_A,
 \qquad \|\Pi_s^--Q_s\|\leq CL^{-1}.
 \end{aligned}
$$
For the vector $\mathfrak m_s$ of local mean differences,
$$
 \|\mathfrak m_s\|_{\mathcal H^*}\leq C\sqrt L,
 \qquad \Pi_s^-\eta=0
 \ \Longrightarrow\ 
 |\mathfrak m_s(\eta)|\leq CL^{-1/2}\|\eta\|_A.
$$

\paragraph{Limiting scalar equations.}
For any model profile $p$, use
$\|u\|_{H_p}^2=\int_{\R^2}(|\nabla u|^2+pu^2)\dif A$.
Choose a circle in an annulus where $p$ is bounded below so that the
mean of $u$ on the circle is controlled by $\|u\|_{H_p}$. The
Cauchy--Schwarz inequality in the radial variable gives
$|\bar u(r)|^2\leq C(1+\log r)\|u\|_{H_p}^2$ for large $r$.
Poincar\'e's inequality on circles and $p(y)\leq C|y|^{-2-\alpha}$ give
the following bound for $q\geq0$ and $R\geq2$.
\begin{equation}\label{6.7}
 \int_{\R^2\setminus B_R}p(y)\log^q(2+|y|)u^2\dif A
 \leq C_q R^{-\alpha}(1+\log R)^{q+1}
       \left(\int_{\R^2}|\nabla u|^2\dif A+\int_{\R^2}pu^2\dif A\right).
\end{equation}
Here $\alpha=2$ for the ordinary and inner profiles and $\alpha=4$
for the outer profiles. With the radial upper bound for $p$, the
angular part is bounded by $\sup_{r>R}(r^2p(r)\log^q(2+r))$ times the
angular energy.

To estimate the centered variance at width $r_s$, subtract the mean of
$\eta_i$ on an annulus with radii comparable to $r_s$. After rescaling,
Poincar\'e's inequality bounds the resulting $L^2$ norm on a fixed disk
by the Dirichlet energy. Starting from the same annulus, the proof of
\eqref{6.7} bounds the weighted integral up to the edge of the chart.
Replacing the annular mean by the profile mean decreases the weighted
integral. For $n_*$, use Poincar\'e's inequality on $M$.

Suppose
\[
 \|\eta_s\|_A+I_s(\eta_s,\eta_s)\leq C,\qquad
 \mathfrak m_s(\eta_s)\longrightarrow0.
\]
The local mean bounds give $|\mu_1|\leq C\sqrt L$. The inner means of
the second component differ by $o(1)$ from the mean against $n_*$, so
$|\mu_2|\leq C$. At width $r_s$, apply \eqref{6.7} to
$\chi(r_sy)\eta_i(q+r_sy)-\mu_i$. Its Dirichlet energy is bounded.
The cutoff and the tail change the profile-weighted squared $L^2$
norm by $O(r_s^\alpha(1+L))$.

On fixed rescaled disks, the bound in the weighted norm on an annulus,
together with Poincar\'e's inequality, gives an $H^1$ bound. Rellich
compactness and
\eqref{6.7} give strong convergence in the profile-weighted $L^2$
spaces, including with fixed logarithmic weights. Denote the limits of
the first-component and inner second-component profiles by $v_j^o$
and $v_k^i$, respectively. Let $v_*$ be the weak $H^1$ limit of
$\eta_{2,s}$ on $M$. By the mean difference condition and Rellich
compactness, $\mu_2\to\avg{v_*}{n_*}$, where $n_*$ denotes the limiting
density if the parameters vary. Thus
\begin{equation}\label{6.8}
 I_s(\eta_s,\eta_s)\longrightarrow
 \sum_j\int_{\R^2}p_j(v_j^o)^2\dif A+\sum_k\int_{\R^2}p_{\rm in}(v_k^i)^2\dif A
 +\int_Mn_*|v_*-\avg{v_*}{n_*}|^2\dV.
\end{equation}
The profile means of $v_j^o$ and $v_k^i$ vanish by
$\mathfrak m_s(\eta_s)\to0$.

For an energy-normalized eigenvector with bounded eigenvalue
$\tau_s\leq 1/2$, write $\lambda_s=(1-\tau_s)^{-1}$ and
$v_{i,s}=\eta_{i,s}-\mu_i$. The eigenvector satisfies the following unprojected equation.
\[
 -\Delta_g\eta_s
 =\lambda_s A\left[
 \binom{V_1 v_{1,s}}{V_2v_{2,s}}
 +\sum_\alpha c_{\alpha,s}F_{\alpha,s}\right].
\]
Indeed, $(\mathcal A_s-\tau_s)\eta_s=\sum_\alpha c_{\alpha,s}
S_{\alpha,s}$. Test this identity with $S_{\beta,s}$ and use \eqref{6.17}
and $I_s(\eta_s,\eta_s)=1-\tau_s$ to bound the multipliers.

Let $\lambda_s\to\lambda$. At the outer scale of a mixed center, the
inner source has mass $o(1)$ by the bound on the mean differences.
Subtract the value at the origin from a fixed test function. The
weighted Cauchy--Schwarz inequality and the profile tail then show
that the remaining pairing tends to zero. The centered inner constraint
sources have zero mass. The same weighted inequality shows that the
outer sources have vanishing mass on each fixed disk at the inner scale.

Passing to the limit at each scale gives
\[
 -\Delta v=2\lambda p\left(v+\sum_\alpha c_\alpha Z_\alpha\right),
 \qquad \int_{\R^2}pv\dif A=0,\qquad
 \int_{\R^2}|\nabla v|^2\dif A<\infty,
\]
where the sum runs over the constraint sources at that scale. The
radial coefficient at each first-component scale is the same limit
$c_0$. Kelvin inversion and finite energy give a removable singularity
at infinity and allow us to test with the bounded kernel functions.

Testing with each constrained translation or angular function
$Z_\alpha$ gives
\[
 (1-\lambda)\int p v Z_\alpha\dif A
 =\lambda c_\alpha\int p Z_\alpha^2\dif A.
\]
The constraint makes the moment on the left vanish. Testing with
the first-component radial functions and summing gives
\[
 (1-\lambda)\sum_j\int p_jv_j^oZ_j\dif A
 =\lambda c_0\sum_j d_j.
\]
By \eqref{6.2}, the left-hand side also vanishes. Thus all multipliers
tend to zero, including when $\lambda=1$.

In the original variables, the centered local sources and the terms
$m_jv_{1,s}$, $\widetilde m_kv_{2,s}$ converge to zero in distributions.
The two equations and the zero-mean condition give $-v_*/2$ as the
limit of the first component. For the second component, we have
\[
 -\Delta_gv_*=2\lambda n_*(v_*-\avg{v_*}{n_*}).
\]
In particular, when $\lambda=1$, scalar nondegeneracy gives
$v_*=0$ and hence $\mu_2\to0$.

Denote the scalar forms for the ordinary profile, the outer
profile and $n_*$ by $Q_{\typeA}$, $Q_{\typeB,B}$ and $Q_*$.
The ordinary profile is the pullback of the round density on the
sphere. On zero-mean functions, the eigenvalues satisfy
$2\lambda=l(l+1)$ for $l\geq1$, so
$\operatorname{ind}Q_{\typeA}=0$. For the outer profile at $B=0$,
the substitution $\chi=2\arctan|y|^2$ gives the metric
$\dif\chi^2+4\sin^2\chi\,\dif\theta^2$.

For the angular mode $e^{in\theta}$, write the radial function as
$(\sin\chi)^{|n|/2}g(\cos\chi)$. With $x=\cos\chi$ and
$h=|n|/2$, the equation for $g$ is
\[
 -(1-x^2)g''+(2h+2)xg'+h(h+1)g=2\lambda g.
\]
Finite energy selects the regular solution at both endpoints.
The Gegenbauer polynomials form an orthogonal basis for the
weight $(1-x^2)^h$ and give
\[
 2\lambda=(j+|n|/2)(j+|n|/2+1),\qquad
 j\geq0,\quad n\in\mathbb Z .
\]
The zero-mean condition removes the constant mode. The modes with
$0<\lambda<1$ are $(j,n)=(0,\pm1)$. Along the path $B_\theta=\theta B$,
use the fixed energy space $H^1(S^2)/\R$ with a fixed round metric in
stereographic coordinates. The measures $p_{B_\theta}\dif A$
have bounded densities relative to that metric, including in the
Kelvin chart, where the densities are $O(|z|^2)$. Their densities
vary continuously in $L^\infty$ because the formula is smooth
in $B$ on the compact sphere and $B_\theta$ remains bounded.
Subtracting the $p_{B_\theta}$-mean is continuous in the energy
norm on representatives with round mean zero. Rellich compactness
therefore gives a norm-continuous family of centered compact
self-adjoint operators for the weighted forms.

The zero eigenspace of the identity-minus-compact operator
corresponds, by subtracting the weighted mean, to the solutions
of $-\Delta Z=2p_{B_\theta}Z$ in the finite-energy class.
In the Kelvin chart a logarithmic cutoff that vanishes in a
shrinking disk and has Dirichlet energy tending to zero allows
the weak equation to be tested across the missing point.
The density is bounded in that chart, so the right-hand side also
converges as the cutoff is removed. Elliptic regularity then gives
bounded representatives. Conversely, a bounded planar solution
extends across the missing point. Subtract the local Poisson potential
of the bounded right-hand side in the Kelvin chart, then remove the
singularity of the bounded harmonic remainder. Testing the extended equation
with constants gives zero weighted mean. The kernel result
\cite[Theorem~1]{DEM} applies with the exponent parameter $N=1$
and the complex parameter $B_\theta$ in the potential
$2p_{B_\theta}=32|z|^2/(1+|z^2-B_\theta|^2)^2$.
The kernel thus has dimension three for every $\theta$.

At a fixed $\theta$, take a small contour in the complex plane around
zero that excludes the rest of the spectrum. For nearby parameter
values, the Neumann series gives continuity of the resolvent on the
contour. Its contour integral then gives norm-continuous spectral
projections. Two projections at distance less than one have the same
finite rank because the restriction of each to the range of the
other is injective. Thus the kernel fills the spectral cluster of
rank three throughout the neighborhood. A contour around the negative
spectrum shows that the rank of its spectral projection is locally
constant, and hence constant along the path. At $B=0$, the explicit
eigenvalues above give rank two. The separate
residual form has index $m_*$ by definition. Hence
\[
 \operatorname{ind}Q_{\typeB,B}=2,\qquad \operatorname{ind}Q_*=m_*.
\]

For an energy-orthonormal eigenfamily in a bounded negative
spectral interval, $I_s(\eta_l,\eta_r)=(1-\tau_l)\delta_{lr}$.
Polarizing \eqref{6.8} shows that the limiting profiles are
linearly independent in the direct sum of the scalar negative
eigenspaces. Thus, for each $\epsilon>0$,
\[
 \limsup_{s\to0}\dim\mathbf1_{[-C,-\epsilon]}(\mathcal L_s)
 \leq m_*+2b.
\]

For the lower bound, use $(w,-w/2)$ at the outer scales of the mixed
centers and $(-v/2,v)$ on $M$, where $w$ and $v$ range over the negative
eigenspaces of the outer and background scalar forms, respectively.
At an outer scale, multiply by a logarithmic cutoff equal to one
for $|x|\leq r_s^{3/4}$ and zero for $|x|\geq r_s^{1/4}$, with
Dirichlet energy $O(L^{-1})$. The profile functions are bounded and
have finite Dirichlet energy, so their energies and centered profile
variances converge. Each outer test converges to a constant at the
inner scale, giving zero limiting centered variance there. Its
shrinking support makes its energy pairing with smooth background
functions tend to zero.

The source pairings tend to zero because the tests are orthogonal to
the profile kernels and converge to constants on the smaller scales.
By \eqref{6.17}, projecting onto $E_s$ changes the tests by $o(1)$ in
energy norm. Their local mean differences are bounded. To place the
tests in $\ker\mathcal R_s$, subtract linear combinations of the
representatives in $\operatorname{ran}\mathcal R_s$. The inverse Gram
matrix gives the coefficients. By \eqref{6.6}, the correction is $O(L^{-1/2})$
in energy norm and changes the $\mathcal B_s$ pairings by $o(1)$.

On the resulting trial space $E_s^-$,
\[
 \dim E_s^-=m_*+2b,\qquad
 \qform_s|_{E_s^-}=\mathcal B_s|_{E_s^-}\leq-c_0\|\cdot\|_A^2.
\]
The cross terms with $\operatorname{ran}\mathcal R_s$ are
bounded by $\|\mathcal B_s\|$, since $\mathcal R_s E_s^-=0$.
For $u\in\operatorname{ran}\mathcal R_s$ and $v\in E_s^-$,
\[
 \qform_s(u+v,u+v)
 \leq-cL\|u\|_A^2+2C\|u\|_A\|v\|_A-c_0\|v\|_A^2
 \leq-\tfrac c2L\|u\|_A^2-\tfrac{c_0}2\|v\|_A^2.
\]
The min--max principle therefore gives at least
\begin{equation}\label{6.9}
 r+(m_*+2b)=m_*+a+4b-1
\end{equation}
negative eigenvalues bounded away from zero.

\subsection{Eigenvalues near zero}

Suppose $\tau_s\to0$ and $\|\eta_s\|_A=1$. By the compactness and
kernel classification above, the limits are radial:
$$
 v_1\longrightarrow\gamma_jZ_j,\qquad
 v_2\longrightarrow\beta_kZ_{\rm in},\qquad v_i=\eta_i-\mu_i,
$$
at the corresponding scales, while $v_*=0$.
By \eqref{6.2} and \eqref{6.8},
\begin{equation}\label{6.10}
 \sum_jd_j\gamma_j=0,\qquad
 \mathcal D(\Gamma):=\sum_jd_j\gamma_j^2+e\sum_k\beta_k^2=1,
 \qquad \Gamma=(\gamma,\beta).
\end{equation}
Put $\vartheta=\mu_1$, $I_j=\int_{M} m_jv_1\dV$, $J_k=\int_{M}\widetilde m_kv_2\dV$.
The normalization of the first component gives $\sum_jI_j=0$.

\begin{lemma}\label{le6.2}
Every such eigenvector satisfies
\begin{equation}\label{6.11}
 \tau_sL\,\mathcal D(\Gamma)=\mathcal F(\Gamma)+o(1),
\end{equation}
where
\begin{align}
 \mathcal F(\Gamma)&=\min_{\vartheta\in\R}\mathcal F_\vartheta(\Gamma),
                                                        \label{6.12}\\
 \mathcal F_\vartheta(\Gamma)
 &=\pi\sum_{\typeA}x_j^2+
       \frac{4\pi}{7}\sum_{\typeB}(4x_k^2-2x_kz_k+2z_k^2),
                                                        \label{6.13}\\
 x_j&=\vartheta-\gamma_j,\qquad z_k=\beta_k-t_k\gamma_k.\notag
\end{align}
\end{lemma}

\begin{proof}
Set $\lambda_s=(1-\tau_s)^{-1}$. At width $r$, subtract from each
Green potential the product of $-(2\pi)^{-1}\log r$ and its source
mass. The remaining logarithmic convolution at the same scale
converges to the potential of the limiting profile by strong
weighted $L^2$ convergence, local square integrability of the
logarithmic kernel, and \eqref{6.7}.

At the outer scale of a mixed center, the remaining inner potential
tends to zero on fixed annuli. Indeed, for
$\varrho=r_{\rm in}/r_{\rm out}$, the kernel $\log|y-\varrho z|$
converges to $\log|y|$ in $L^2(p_{\rm in}(z)\dif A_z)$ uniformly in $y$
on each such annulus. Split the integral into a fixed disk in $z$
and its complement. In the part of the complement where
$|y-\varrho z|<|y|/2$, we have $p_{\rm in}(z)\leq C\varrho^4$, and
$\zeta=y-\varrho z$ gives an integral bounded by $C\varrho^2$.
On the rest, \eqref{6.7} controls the tail. The limiting kernel
$\log|y|$ pairs with the inner mass $J_k=o(1)$.

At the inner scale, after subtracting the logarithmic mass term,
the outer potential converges to its value at zero by the same
weighted $L^2$ estimate for $\log|\varrho y-z|$. The regular Green
terms from the other centers are $o(1)$ by the vanishing local
masses. Since $v_*=0$ and $\mu_2\to0$, the background potential
converges uniformly to zero. The centered-source potentials are
uniformly bounded, and their multipliers tend to zero.

For the radial functions the planar identities are
\[
 2G_{\R^2}(pZ)=Z+1,\qquad
 -\tfrac12\gamma_k(Z_{0,B_k}(0)+1)=-t_k\gamma_k.
\]
For the first identity, both sides have the same Laplacian and decay
at infinity. Substitute these identities into the Green representation
on fixed rescaled annuli. Comparing the resulting expressions with the
radial limits gives
\begin{equation}\label{6.14}
 \begin{aligned}
 x_j&=\lambda_s LI_j/\pi+o(1),&&j\in\typeA,\\
 x_k&=\lambda_s L(2I_k-J_k)/(4\pi)+o(1),&&k\in\typeB,\\
 z_k&=\lambda_s L(I_k-4J_k)/(4\pi)+o(1),&&k\in\typeB.
 \end{aligned}
\end{equation}
The logarithms of the ordinary and inner widths are
$-L+O(1)$, and those of the outer widths are $-L/2+O(1)$.
The $O(1)$ terms are multiplied by $I_j,J_k=o(1)$.
Since $\sum_jI_j=0$,
$$
 \vartheta=
 \frac{7\sum_{\typeA}\gamma_j+
       4\sum_{\typeB}((4-t_k)\gamma_k+\beta_k)}
      {7a+16b}+o(1)=O(1),
 \qquad |I_j|+|J_k|=O(L^{-1}).
$$
As $\lambda_s\to1$, inversion of \eqref{6.14} yields
\begin{equation}\label{6.15}
 I_j=\frac{\pi x_j}{L}+o(L^{-1})\quad(j\in\typeA),\qquad
 \binom{I_k}{J_k}=\frac{4\pi}{7L}
             \binom{4x_k-z_k}{x_k-2z_k}+o(L^{-1}).
\end{equation}
Set
$$
 \begin{aligned}
 Q_j&=\int_{M} m_j\widehat Z_jv_1\dV,&
 d_{j,s}&=\int_{M} m_j\widehat Z_j^2\dV=d_j+O(s^\theta),\\
 P_k&=\int_{M}\widetilde m_k\widehat Z_{\rm in}v_2\dV,&
 e_{k,s}&=\int_{M}\widetilde m_k\widehat Z_{\rm in}^2\dV=e+O(s^\theta).
 \end{aligned}
$$
Write $\Psi_j=2G(m_j\widehat Z_j)$ and
$\widetilde\Psi_k=2G(\widetilde m_k\widehat Z_{\rm in})$.
The estimate for the centered Green potentials gives
\[
 \|\Psi_j-(\widehat Z_j+1)\|_{L^2(m_j)}
   +\|\widetilde\Psi_k-(\widehat Z_{\rm in}+1)\|_{L^2(\widetilde m_k)}
 \leq Cs^\theta
\]
for some $\theta>0$. These potentials are also $O(s^\theta)$ in
$L^2(n_*)$ and in the weighted spaces of the other centers.
At a mixed center, with $r=a_k/\delta_k\asymp\sqrt s$, the
cross-scale estimates are
\[
 \begin{aligned}
 \|Z_{0,B}(r\,\cdot)-Z_{0,B}(0)\|_{L^2(p_{\rm in})}
 &\leq Cr,\\
 \|Z_{\rm in}(\cdot/r)+1\|_{L^2(p_B)}
 &\leq Cr^2\sqrt{1+|\log r|}.
 \end{aligned}
\]
For the first estimate, use
$|Z_{0,B}(ry)-Z_{0,B}(0)|\leq C\min(r^2|y|^2,1)$.
For the second, use $Z_{\rm in}(y/r)+1=2r^2/(r^2+|y|^2)$
and $p_B(y)\leq C|y|^2$ near zero.

Test the first component of the unprojected equation with $\Psi_j$.
The left-hand side is $2Q_j$ and the pairing with the density at
the same scale is $Q_j+I_j+O(s^\theta)$. At a mixed center, the
pairing of $\Psi_k$ with the inner density is $2t_kJ_k+O(s^\theta)$.
Testing the second component with $\widetilde\Psi_k$ gives $2P_k$
on the left. The inner pairing is $P_k+J_k+O(s^\theta)$, and the
outer pairing is $O(s^\theta L)$.

The source multipliers are bounded. By the cutoff estimate and the
weighted orthogonality on $\R^2$, the pairing of the radial test
potential with its radial source is $d_{j,s}+O(s^\theta)$, and its
pairings with the angular or translation sources are $O(s^\theta)$.
For pairings between the two scales at a mixed center, the constant
part gives zero because the sources have zero mass. The cross-scale
estimates above control the remaining pairings.

Dividing by $2\lambda_s$ now gives
\begin{equation}\label{6.16}
 \begin{aligned}
 I_j&=-\tau_s Q_j-c_s d_{j,s}+o(L^{-1}),&&j\in\typeA,\\
 I_k-t_kJ_k&=-\tau_s Q_k-c_s d_{k,s}+o(L^{-1}),&&k\in\typeB,\\
 J_k&=-\tau_s P_k+o(L^{-1}),&&k\in\typeB.
 \end{aligned}
\end{equation}
Here $c_s$ is the common radial multiplier. Weighted
Cauchy--Schwarz and the bound on $I_s(\eta_s,\eta_s)$ give
$O(s^\theta L^C)=o(L^{-1})$ for all the error terms.

With
$$
 \begin{gathered}
 \gamma_{j,s}=Q_j/d_{j,s},\qquad \beta_{k,s}=P_k/e_{k,s},\qquad
 \mathcal D_s=\sum_jd_{j,s}\gamma_{j,s}^2+
                    \sum_ke_{k,s}\beta_{k,s}^2\longrightarrow1,\\
 x_{j,s}=\vartheta-\gamma_{j,s},\qquad
 z_{k,s}=\beta_{k,s}-t_k\gamma_{k,s},\qquad
 \sum_jQ_j=\sum_jI_j=0,
 \end{gathered}
$$
equations \eqref{6.16} and \eqref{6.15} imply
$$
 \begin{aligned}
 \tau_s\mathcal D_s
 &=\sum_jx_{j,s}I_j-\sum_kz_{k,s}J_k+o(L^{-1})\\
 &=L^{-1}\mathcal F_\vartheta(\Gamma_s)+o(L^{-1}),\\
 \partial_\vartheta\mathcal F_\vartheta(\Gamma_s)
 &=2L\sum_jI_j+o(1)=o(1),\qquad
 \partial_\vartheta^2\mathcal F_\vartheta=2\pi a+32\pi b/7>0.
 \end{aligned}
$$
Thus
$$
 \tau_s=O(L^{-1}),\qquad
 \mathcal F_\vartheta(\Gamma_s)=\min_\theta
       \mathcal F_\theta(\Gamma_s)+o(1),
$$
which proves \eqref{6.11}.
\end{proof}

\begin{lemma}\label{le6.3}
On the subspace $\sum_jd_j\gamma_j=0$, the form $\mathcal F$ satisfies
$c\mathcal D\leq\mathcal F\leq C\mathcal D$ uniformly for
$0\leq t_k\leq1$.
\end{lemma}
\begin{proof}
For $D_0=\sum_jd_j$ and $\sum_jd_j\gamma_j=0$,
$$
 \vartheta D_0=\sum_jd_jx_j,\qquad
 \sum_jd_j\gamma_j^2
 =\sum_jd_jx_j^2-D_0\vartheta^2\leq\sum_jd_jx_j^2.
$$
Since $0\leq t_k\leq1$,
$$
 \beta_k^2=(z_k+t_k\gamma_k)^2\leq2z_k^2+2\gamma_k^2,
 \qquad 4x^2-2xz+2z^2=3x^2+(x-z)^2+z^2.
$$
Consequently, for every $\vartheta$,
$$
 \mathcal D(\Gamma)
 \leq C\left(\sum_jx_j^2+\sum_kz_k^2\right)
 \leq C\mathcal F_\vartheta(\Gamma),
 \qquad
 \mathcal F(\Gamma)\leq\mathcal F_0(\Gamma)\leq C\mathcal D(\Gamma).
$$
\end{proof}

At fixed $b_k=s_kB_k$ and $\beta_k^{\rm match}$,
$$
 a_k^2=\frac{\beta_k^{\rm match}}4(s_k^2+|b_k|^2)
 \quad\Longrightarrow\quad
 d\log a_k=t_k\,d\log s_k,
 \qquad t_k=\frac1{1+|B_k|^2}.
$$
For the full Toda core, the scale derivative at fixed $B_k$ equals
the derivative at fixed $b_k$ plus the result of applying
$\sum_{l=1}^2 B_{k,l}\partial_{B_{k,l}}$ to the core.
The added term has an inner radial component because
$d\log a_k=d\log s_k$ at fixed $B_k$. In \eqref{6.13}, the radial
coefficients correspond to the derivative at fixed $b_k$. Hence
$$
 \mathcal F_\vartheta(\Gamma)=0
 \quad\Longleftrightarrow\quad
 \gamma_j=\vartheta,\quad\beta_k=t_k\vartheta,
$$
which is the common dilation direction.

\begin{proof}[Proof of Theorem~\ref{th6.1}]
By Lemmas~\ref{le6.2} and \ref{le6.3},
$$
 \tau_s\to0\quad\Longrightarrow\quad
 c\leq\tau_sL\leq C.
$$
Combining this with \eqref{6.9} gives
$$
 \operatorname{ind}\mathcal L_s=m_*+a+4b-1,
 \qquad \operatorname{dist}(0,\operatorname{spec}\mathcal L_s)\geq c/L,
 \qquad \|\mathcal L_s^{-1}\|\leq CL.
$$
The constants are uniform on the stated compact parameter sets.
\end{proof}

\subsection{Stability under perturbations}\label{sec:range-stability}

\begin{lemma}\label{le6.4}
After isometric identification of the projected spaces, a
self-adjoint perturbation of size $o(L^{-1})$ preserves the index and
the $O(L)$ inverse bound of Theorem~\ref{th6.1}.
\end{lemma}
\begin{proof}
For $0\leq\theta\leq1$, the spectral gap of
$\mathcal L_s$ is at least $c/L$. The perturbation on this space
satisfies $\|\mathcal L_s^{-1}R_s\|=o(1)$, so
\[
 (\mathcal L_s+\theta R_s)^{-1}
   =(I+\theta\mathcal L_s^{-1}R_s)^{-1}\mathcal L_s^{-1}
\]
by the Neumann series, with norm $O(L)$ uniformly in $\theta$.
Write the energy space as the orthogonal sum of the negative
spectral space of $\mathcal L_s$ and its positive complement.
On the former, the unperturbed form is at most $-c/L$ times the
squared norm. On the latter, it is at least $c/L$ times the squared norm.
Since $\|R_s\|=o(L^{-1})$, the perturbed form is strictly negative
on the former and strictly positive on the latter. Any negative
subspace of larger dimension would intersect the positive complement
nontrivially, which is impossible. This proves equality of the indices
and the required uniform gap.
\end{proof}

The estimates for the centered Green potentials give
\begin{equation}\label{6.17}
 \|S_{\alpha,s}\|_A^2+I_s(S_{\alpha,s},S_{\alpha,s})\leq C,
 \qquad (\ipA{S_{\alpha,s}}{S_{\beta,s}})_{\alpha,\beta}\geq cI.
\end{equation}

\begin{corollary}
\label{co6.5}
For $f\in\mathcal H^{-1}$ there is a unique pair $(\eta,c)$ such that
$$
 \begin{cases}
 \qform_s(\eta,\xi)=\langle f,\xi\rangle
       +\displaystyle\sum_\alpha c_\alpha\ipA{S_\alpha}\xi,
       &\xi\in\mathcal H,\\
 \ipA\eta{S_\alpha}=0,&\text{for every }\alpha.
 \end{cases}
$$
Moreover,
$$
 \|\eta\|_{H^1}+|c|\leq CL\|f\|_{H^{-1}}.
$$
The same conclusion holds for the perturbations of the form and
centered sources considered below.
\end{corollary}
\begin{proof}
The projected inverse gives $\|\eta\|_A\leq CL\|f\|_{H^{-1}}$.
Testing with $\eta$ gives
\[
 I_s(\eta,\eta)=\|\eta\|_A^2-\langle f,\eta\rangle
 \leq CL^2\|f\|_{H^{-1}}^2 .
\]
Testing with $S_\alpha$ then yields
\[
 \sum_\beta\ipA{S_\beta}{S_\alpha}c_\beta
 =-I_s(\eta,S_\alpha)-\langle f,S_\alpha\rangle .
\]
By \eqref{6.17} and Cauchy--Schwarz, $|c|\leq CL\|f\|_{H^{-1}}$.
The argument uses only the projected inverse bound and \eqref{6.17}.
\end{proof}

For nonnegative weights, write $B_V(\eta,\xi)=\ipA\eta\xi-I_V(\eta,\xi)$.
The orthogonality functionals are
$\int_{M} V_i(z_\alpha-\avg{z_\alpha}{V_i})\eta_i\dV$.

\begin{lemma}\label{le6.6}
Suppose the masses of $V_i,\widetilde V_i$ lie in a compact
subset of $(0,\infty)$,
$\sum_i\int_{M} V_i\xi_i^2\dV\leq CL\|\xi\|_{H^1}^2$, and
$e_p=\sum_i\|\widetilde V_i-V_i\|_{L^p}\leq1$, $p>1$.
Then
\begin{equation}\label{6.18}
 \|B_{\widetilde V}-B_V\|\leq C_pe_pL.
\end{equation}
Suppose also that the Riesz representatives of the centered
functionals defined by the bounded tests $z_\alpha$ satisfy
\eqref{6.17}. If $e_p\sqrt L=o(1)$, there is an isometry $U$ between
the common kernels associated with $V$ and $\widetilde V$ such that
\begin{equation}\label{6.19}
 \|U-I\|\leq C_pe_p\sqrt L,\qquad
 \|B_{\widetilde V}(U\cdot,U\cdot)-B_V\|
                         \leq C_pe_pL^{3/2}.
\end{equation}
\end{lemma}
\begin{proof}
For each component put
$m_V=\int_{M} V\dV$, $\ell_V(\xi)=\int_{M} V\xi\dV$ and
$Q_V(\eta,\xi)=\int_{M} V\eta\xi\dV$. H\"older and Sobolev give
$$
 \begin{gathered}
 I_V=\sum_i(Q_{V_i}-m_{V_i}^{-1}\ell_{V_i}\otimes\ell_{V_i}),\\
 |m_{\widetilde V}-m_V|+
 \|\ell_{\widetilde V}-\ell_V\|+\|Q_{\widetilde V}-Q_V\|
 \leq C_pe_p,\qquad
 \|\ell_V\|+\|\ell_{\widetilde V}\|\leq C\sqrt L,\\
 \|I_{\widetilde V}-I_V\|\leq C_pe_pL.
 \end{gathered}
$$
For $J_{V,\alpha}(\xi)=\int_{M} V(z_\alpha-\avg{z_\alpha}V)\xi\dV$,
$$
 \|J_{\widetilde V,\alpha}-J_{V,\alpha}\|
 \leq C_pe_p\sqrt L.
$$
Let $S$ collect their Riesz representatives and let
$P=I-S(S^*S)^{-1}S^*$. By \eqref{6.17},
$$
 \|\widetilde P-P\|\leq C_pe_p\sqrt L=o(1),\qquad
 U=\widetilde PP(P\widetilde PP)^{-1/2}|_{\operatorname{ran}P},
 \qquad \|U-I\|\leq C_pe_p\sqrt L.
$$
Since $\|B_V\|+\|B_{\widetilde V}\|\leq CL$,
$$
 \|B_{\widetilde V}(U\cdot,U\cdot)-B_V\|
 \leq C_pe_pL+CL\|U-I\|\leq C_pe_pL^{3/2}.
$$
\end{proof}

\begin{corollary}
\label{co6.7}
Let $V$ be the model density with $s=t$, and put $\ell=\log(1/t)$.
Suppose $V_i^{\rm app}\geq0$ have masses in a compact subinterval
of $(0,\infty)$ and satisfy
\[
 \sum_i\|V_i^{\rm app}-V_i\|_{6/5}\leq Ct^{2/3}.
\]
Use the model cutoff test functions, including the sum of the
first-component radial tests in \eqref{6.2}. Multiply each test
by a fixed nonzero factor and center it with respect to the
corresponding $V_i^{\rm app}$. Here $\mathcal L_t$ is the model
operator $\mathcal L_s$ at $s=t$, and $\mathcal L_t^{\rm app}$
is the operator projected onto the completed space. Then there is
an isometry $U$ from the model space to the completed space such that
\[
 \|U^*\mathcal L_t^{\rm app}U-\mathcal L_t\|
 =O(t^{2/3}\ell^{3/2})=o(\ell^{-1}).
\]
The completed operator has the index and inverse bounds in
\eqref{6.3}--\eqref{6.4}, and its source representatives
satisfy \eqref{6.17}.
\end{corollary}
\begin{proof}
Divide each completed test by its fixed factor.
This preserves the common kernel and rescales the source bounds
by fixed constants. First compare the model sources with the sources
obtained by centering the same tests against $V_i$. In the first
component, $V_1g=m_jg$ for a test $g$ supported at one center.
The local mean of the test is $O(r^\alpha)$, with $r\asymp t$, $\alpha=2$
at an ordinary center and $r\asymp\sqrt t$, $\alpha=4$ at a mixed
center. The same estimate applies to each summand of the common
radial test.

In the second component, an inner translation test satisfies
$|g(x)|\leq C\min(1,t/|x|)$ and
$V_2g=\widetilde m_kg+n_*g$. Hence
$\|n_*g\|_{6/5}+\|n_*g\|_1\leq Ct$. The local and global
weighted means are $O(t)$.

Since $\|V_i\|_{6/5}\leq Ct^{-1/3}$, the change from local to global
centering is $O(t^{2/3})$ in $L^{6/5}$. By Sobolev's inequality,
the Riesz representatives change by the same order in energy norm.
The Gram bound \eqref{6.17} and the projection formula in the proof
of Lemma~\ref{le6.6} give an isometry at distance $O(t^{2/3})$ from
the identity.

For the second comparison, let $S_{\alpha,s}$ and
$\widehat S_{\alpha,s}$ be the locally and globally centered model
representatives. These representatives differ by $O(t^{2/3})$ in energy
norm. The norms of this fixed finite set of vectors are uniformly
bounded, so their Gram matrix changes by $O(t^{2/3})$. The model
form $I_s$ is a sum of nonnegative weighted variances. Therefore,
\[
 I_s(\widehat S_{\alpha,s},\widehat S_{\alpha,s})^{1/2}
 \le I_s(S_{\alpha,s},S_{\alpha,s})^{1/2}
     +C\sqrt\ell\,t^{2/3}\le C,
\]
using $\|I_s\|\leq C\ell$. This proves \eqref{6.17} for the
intermediate representatives. Next apply Lemma~\ref{le6.6} to the
completed sources. The form errors for the two isometries are
$O(t^{2/3}\ell)$ and $O(t^{2/3}\ell^{3/2})$, respectively. Composing
the isometries and applying Lemma~\ref{le6.4} gives the index and
inverse bounds.

The completed source representatives differ from the original
representatives by $O(t^{2/3}\sqrt\ell)$ in energy norm, so the Gram
lower bound persists. Equation \eqref{6.18}, the estimate $\|I_s\|\leq C\ell$,
and the original source bounds \eqref{6.17} give the centered variance bounds.
\end{proof}

\begin{lemma}\label{le6.8}
Suppose the completed profile is independent of $\mu=\rho_1-4\pi m$,
with $m\geq1$, and its projected linearized operator at $\mu=0$
satisfies \eqref{6.3}, \eqref{6.4}, and \eqref{6.17}. Suppose also
that its first weight satisfies the uncentered bound
\begin{equation}\tag{6.19a}\label{eq:completed-weight-capacity}
 \int_M V_1(t,0)\xi^2\dV\leq C\ell\|\xi\|_{H^1(M)}^2
 \quad\text{for all }\xi\in\mathring H^1(M),\qquad
 \ell=\log(1/t).
\end{equation}
The same spectral, inverse and source bounds hold whenever
$|\mu|\ell^2=o(1)$. The residual changes by at most
$C|\mu|\sqrt\ell$ in $H^{-1}$.
\end{lemma}
\begin{proof}
By the normalization,
$$
 V_1(t,\mu)=\left(1+\frac{\mu}{4\pi m}\right)V_1(t,0),
 \qquad V_2(t,\mu)=V_2(t,0).
$$
The weighted means and the space $E_t$ are unchanged. The centered
variance of the first component is bounded by the uncentered
quadratic form in \eqref{eq:completed-weight-capacity}. Thus
$$
 \mathcal L_{t,\mu}-\mathcal L_{t,0}
       =-\frac{\mu}{4\pi m}I_{V_1(t,0)},\qquad
 \|\mathcal L_{t,\mu}-\mathcal L_{t,0}\|\leq C|\mu|\ell,
$$
where the quadratic form is identified with its energy operator.
It follows that
$$
 \|\mathcal L_{t,0}^{-1}
       (\mathcal L_{t,\mu}-\mathcal L_{t,0})\|
       \leq C|\mu|\ell^2=o(1).
$$
The factors by which the centered sources change tend to one, so the
norm, variance and Gram estimates persist. Apply Cauchy--Schwarz with
respect to $V_1(t,0)\dV$. Since $\int_M V_1(t,0)\dV=4\pi m$,
\eqref{eq:completed-weight-capacity} also gives
$\|V_1(t,0)/(4\pi m)-1\|_{H^{-1}}\leq C\sqrt\ell$.
The change in the residual is the product of $\mu(V_1(t,0)/(4\pi m)-1)$ and a fixed vector, and
$$
 \left\|\mu\left(\frac{V_1(t,0)}{4\pi m}-1\right)\right\|_{H^{-1}}
       \leq C|\mu|\sqrt\ell.
$$
\end{proof}

We will use the fact that \eqref{eq:completed-weight-capacity} follows from
$\|V_1(t,0)\|_\infty\leq Ct^{-2}$ and bounded mass.
With $q=1+1/\ell$, interpolation gives
$\|V_1(t,0)\|_q\leq C$. The two-dimensional Sobolev estimate
$\|\xi\|_{2(\ell+1)}\leq C\sqrt{\ell+1}\|\xi\|_{H^1}$ then gives
the assertion by H\"older's inequality. The scale ranges in
Section~\ref{rev:sec:multiple} give the same supremum bound for the first
model weight. The same capacity bound holds for a completed weight
within $O(t^{2/3})$ of the model in $L^{6/5}$, since the difference
in the quadratic forms is bounded by
$Ct^{2/3}\|\xi\|_{12}^2$.

\begin{lemma}\label{le6.9}
If $V\geq0$, $\int_{M} V\dV$ stays in a compact subinterval of $(0,\infty)$,
$\|V\|_\infty\leq Ct^{-2}$, and $E=\|\eta\|_{H^1}$ with
$t^{-1/3}E=o(1)$, then
\begin{equation}\label{6.20}
 Z=(\int_{M} V\dV)^{-1}\int_{M} Ve^\eta\dV=1+O(t^{-1/3}E),\qquad
 \|Z^{-1}Ve^\eta-V\|_{12/11}\leq Ct^{-1/2}E.
\end{equation}
\end{lemma}
\begin{proof}
By the exponential Sobolev inequality and interpolation,
$$
 \|e^\eta-1\|_q\leq C_qE,\qquad
 \|V\|_{6/5}\leq Ct^{-1/3},\qquad
 \|V\|_{12/11}\leq Ct^{-1/6}.
$$
Hence, by H\"older's inequality,
$$
 \begin{aligned}
 |Z-1|&\leq C\|V\|_{6/5}\|e^\eta-1\|_6
                 \leq Ct^{-1/3}E,\\
 \|Z^{-1}Ve^\eta-V\|_{12/11}
 &\leq C\|V\|_{6/5}\|e^\eta-1\|_{12}
                +C|Z-1|\|V\|_{12/11}\\
 &\leq C(t^{-1/3}+t^{-1/2})E\leq Ct^{-1/2}E.
 \end{aligned}
$$
\end{proof}

\begin{corollary}
\label{co6.10}
Suppose a completed operator satisfies the index and inverse bounds
of Corollary~\ref{co6.7}, with $\|V_i^{\rm app}\|_\infty\leq Ct^{-2}$.
Let $\|\eta\|_{H^1}\leq Ct^{2/3}\ell$. For $0\leq\theta\leq1$, multiply
each density by $e^{\theta\eta_i}$ and normalize it to preserve its mass.
Write $\mathcal L_t(\theta\eta)$ for the resulting completed
operator with its projected space held fixed. Then
\[
 \sup_{0\leq\theta\leq1}
 \|\mathcal L_t(\theta\eta)-\mathcal L_t(0)\|
 =O(t^{1/6}\ell^2)=o(\ell^{-1}),\qquad
 \operatorname{ind}\mathcal L_t(\theta\eta)=m_*+a+4b-1.
\]
\end{corollary}
\begin{proof}
By \eqref{6.20}, the $L^{12/11}$ density difference is at most
\[
 Ct^{-1/2}\|\theta\eta\|_{H^1}\leq Ct^{1/6}\ell.
\]
Apply H\"older's inequality with exponent $1+1/\ell$ to obtain
the weighted $L^2$ bound for $V_i^{\rm app}$. Use the mass and
supremum bounds, together with $\|\xi\|_q\leq C\sqrt q\|\xi\|_{H^1}$.
Equation \eqref{6.18} gives the form bound, and Lemma~\ref{le6.4}
shows that the index is preserved.
\end{proof}

\section{Local construction and degree at mixed configurations}
\label{article:families}

We use the choice of weights in Section~\ref{sec:wall-setting}
and the wall and normalization conventions of
Section~\ref{sec:scale-conventions}.

Fix a nondegenerate shadow solution with $b\geq1$.

\subsection{Construction of the approximate family}
\label{multi-sec:construction}

For separated ordinary and mixed centers near the given shadow
solution, we seek a Toda solution of the form $u=u^{\rm app}+\eta$.
We construct $u^{\rm app}$ by matching the profiles of
Section~\ref{rev:sec:core} to the outer Green fields with a common
first-component normalization. The residual and density estimates
allow us to apply Theorem~\ref{th6.1} to solve the projected equation
for $\eta$. The remaining equations determine the centers, the
parameters $B_j$, and the common scale.

Write
$$
 m=a+2b,\qquad N=a+b,\qquad b\geq1,\qquad
 \rho_1^0=4\pi m,\qquad \rho_1=4\pi m+\mu,
 \qquad \sigma=\rho_2-4\pi b>0.
$$
We use $C=A^{-1}$ from \eqref{2.1}. Fix a nondegenerate
solution $(P^0,Q^0,w^0)$ of
\eqref{3.40} whose scalar linearization at fixed centers is invertible.
Use the scalar branch $w(P)$ and the density $n_*^g$ in
\eqref{7.1}, with $w(P^0)=w^0$.  The density $n_*^g$ has a zero of
order two at each ordinary center and is positive on $M\setminus P$.
We work in a neighborhood where all centers are uniformly separated.
Write $m_*$ for the scalar Morse index and $F=(F^A,F^B)$ for the
center field in \eqref{3.40} evaluated at $w=w(P)$.

\subsubsection{Common normalization and local scales}
Recall the outer pair $u^0$ in \eqref{7.2}.
Away from the centers it satisfies
\begin{equation}\label{7.3}
 -\Delta_gu^0=A\left(\binom0{n_*^g}
                              -\binom{4\pi m}{\rho_2}\right).
\end{equation}
The Dirac masses at $p_i$ and $q_j$ have coefficient vectors
$A(4\pi,0)^T$ and $A(8\pi,4\pi)^T$, respectively.

Use the common normalization \eqref{7.4} and the ordinary
charts and scales in \eqref{7.5}. By \eqref{2.8},
\[
 \nabla\psi_i^A(p_i)=F_i^A(P,Q).
\]
This derivative is invariant under the quadratic holomorphic
change of coordinates at a mixed center used in
\eqref{7.6}--\eqref{7.7}. With $D_i$ as in \eqref{7.9}, we add
$(D_i,-D_i/2)$ to the outer pair near $p_i$.

At each mixed center, use the chart and parameters in
\eqref{7.6}--\eqref{7.7} and the model densities in \eqref{7.8}.
The projected equation determines the corresponding dilation
and the relative dilations of the first component.
By Lemma~\ref{le4.1}, $-\Delta U_j=Ae^{U_j}$ on $\C$.
Use $L_j$ from \eqref{7.9}. Write $\chi_i$ and $\chi_j$ for
$\chi_i^A$ and $\chi_j^B$, respectively. The ansatz
\eqref{eq:ansatz} then reads
\begin{equation}\label{7.10}
 \begin{aligned}
 v^{\rm app}
 &=\binom{u_1^0}{u_2^0}+\sum_i\chi_i\binom{D_i}{-D_i/2}
       +\sum_j\chi_j(U_j-L_j),\\
 u^{\rm app}
 &=v^{\rm app}-\int_Mv^{\rm app}\dV.
 \end{aligned}
\end{equation}
The logarithmic singularities cancel in \eqref{7.10}, so
$u^{\rm app}$ is smooth. The scales in \eqref{7.5} and \eqref{7.7}
are determined by the common constant $c_1$ in \eqref{7.4}.

\begin{proposition}
\label{pr7.1}
Uniformly on the stated compact parameter sets, with $|\mu|\leq C_0t$,
put
$$
 N_i^g=e^{c_i}h_ie^{v_i^{\rm app}},\qquad
 M_i=\int_MN_i^g\dV,\qquad
 V_i^{\rm app}=\rho_i\frac{h_ie^{u_i^{\rm app}}}
                             {\int_Mh_ie^{u_i^{\rm app}}\dV}.
$$
We then have
\begin{equation}\label{7.11}
 M_1=4\pi m+O(t),\qquad M_2=\rho_2+O(t),\qquad
 V_i^{\rm app}=\frac{\rho_i}{M_i}N_i^g,
 \qquad \log(\rho_i/M_i)=O(t).
\end{equation}
For $1<p<2$, the residual of the normalized equation satisfies
\begin{equation}\label{7.12}
 \left\|-\Delta_gu^{\rm app}
 -A\binom{V_1^{\rm app}-\rho_1}{V_2^{\rm app}-\rho_2}\right\|_{H^{-1}}
 \leq C_p t^{2/p-1}.
\end{equation}
Let $V^{\rm mod}$ be the model of Section~\ref{rev:sec:multiple}
with the same scales and parameters. We use $m\dif A=e^{-\phi}m\dif V_g$
for each planar density $m$, with the cutoffs and mass corrections
of that section. Then
\[
 \sum_{i=1}^2\|V_i^{\rm app}-V_i^{\rm mod}\|_{L^p(M)}
 \leq C_pt^{2/p-1}.
\]
\end{proposition}

\begin{proof}
On each ordinary disk, the densities before normalization are
\begin{equation}\label{7.13}
 N_1^{\rm loc}=m_{A,i}e^{J_{1,i}},\quad
 J_{1,i}=F_i^A\cdot x+O(|x|^2),\quad
 N_2^{\rm loc}=\gamma_i(x)(s_i^2+|x|^2),
\end{equation}
where $\gamma_i=n_*^{\rm loc}/|x|^2$ is smooth.
Set $R^{\rm pre}=-\Delta_g v^{\rm app}-A(N^g-(4\pi m,\rho_2)^T)$.
Subtracting the background equation gives
\begin{equation}\label{7.14}
 e^\phi R^{\rm pre}
 =-A\binom{m_{A,i}(e^{J_{1,i}}-1)}{s_i^2\gamma_i}.
\end{equation}
Since $\|rm_{A,i}\|_p=O(s_i^{2/p-1})$, the local residual has
$L^p$ norm $O(t^{2/p-1})$. Parity and the second moment of $m_{A,i}$
give a first mass of $4\pi+O(t^2\ell)$ and an excess of $O(t^2)$
over the limiting second mass.

On each mixed disk, use the local calculation from the proof of
Proposition~\ref{pr5.3}. Separation ensures uniform smoothness of the
Green fields from the other centers on the disk. Here the forcing
corresponds to the full parameter $\rho^0$. We normalize the masses
globally below.
\begin{equation}\label{7.15}
 \begin{aligned}
 N_1^{\rm loc}&=n_{1,j}e^{J_{1,j}},&
 J_{1,j}&=\psi^B_j-\psi^B_j(q_j)
                        -\log(1+\beta_j|x|^2/4)=O(|x|^2),\\
 N_2^{\rm loc}&=n_{2,j}e^{J_{2,j}},&
 J_{2,j}&=\log(n_*^{\rm loc}/n_{{\rm diff},j})
                              =F_j^B\cdot x+O(|x|^2).
 \end{aligned}
\end{equation}
Write $W^0=u^0+(\log h_i+\phi+c_i)_{i=1,2}$ in this chart,
so $J=W^0-L$ and
$-\Delta_0W^0=A(0,n_*^{\rm loc})^T+f^{\rho^0}$ off the center,
where $f^{\rho^0}$ is the forcing in \eqref{2.2} at
$\rho^0=(4\pi m,\rho_2)^T$. Also
$-\Delta_0L=A(0,n_{\rm diff})^T$ off the center. Subtracting
these equations and using $n_*^{\rm loc}=n_{\rm diff}e^{J_2}$
gives the exact residual
\begin{equation}\label{7.16}
 e^\phi R^{\rm pre}
 =-A\binom{n_{1,j}(e^{J_{1,j}}-1)}
          {(n_{2,j}-n_{{\rm diff},j})(e^{J_{2,j}}-1)}.
\end{equation}
Separation gives uniform bounds on the background derivatives.
Thus \eqref{5.6}--\eqref{5.8} give an $L^p$ residual of order
$O(\delta_j^{4/p-2})=O(t^{2/p-1})$, a first mass of $8\pi+O(t)$, and
an excess second mass of $4\pi+O(t)$. All three estimates use the
common normalization \eqref{7.4}.

The cutoff errors are $O(t^2)$ on ordinary annuli and $O(t)$ on
mixed annuli.
Outside the disks, $N_1^g=O(t^2)$ and $N_2^g=n_*^g$.
Summing the mass contributions gives $M_1=4\pi(a+2b)+O(t)$ and
$M_2=4\pi b+\sigma+O(t)$.
For the comparison with the model, the ordinary density errors
are bounded in $L^p$ by
$C(\|r m_{A,i}\|_p+s_i^2)$. At a mixed center, write $m_o,m_i$
for the two planar models of Lemma~\ref{le5.1}. There
\[
 \begin{aligned}
 N_1^{\rm loc}-m_o
 &= (n_1-m_o)e^{J_1}+m_o(e^{J_1}-1),\\
 N_2^{\rm loc}-n_*^{\rm loc}-m_i
 &= (n_2-n_{\rm diff}-m_i)e^{J_2}+m_i(e^{J_2}-1).
 \end{aligned}
\]
By Lemma~\ref{le5.1} and the bounds $J_1=O(r^2)$ and $J_2=O(r)$,
these errors are bounded by $C_p(\delta_j^{2/p}+\delta_j^{4/p-2})$ in $L^p$.
On the cutoff annuli, $U_j-L_j=O(t)$ in $C^2$, and both
model tails are $O(t^2)$. The model mass corrections are $O(t^2)$
relative to the cut-off densities.

The normalized density is $V_i^{\rm app}=(\rho_i/M_i)N_i^g$.
Since $\rho_i/M_i-1=O(t)$ and $\|N_i^g\|_p=O(t^{2/p-2})$,
normalization changes both the density and the residual by
$O(t^{2/p-1})$ in $L^p$. The embedding
$L^p(M)\hookrightarrow H^{-1}(M)$ gives \eqref{7.12}.
\end{proof}

\subsubsection{The projected equation}

Let $Z_A(y)=(1-|y|^2)/(1+|y|^2)$,
$Z_{B,0}(y)=(1-|y^2-B|^2)/(1+|y^2-B|^2)$, and
\[
 Z_l(y)=\frac{y_l}{1+|y|^2},\qquad
 Z_{B,l}(y)=\frac{4(y^2-B)_l}{1+|y^2-B|^2},\qquad l=1,2 .
\]
At $p_i$, use $\chi_i Z_l(x/s_i)$ in the first component. At $q_j$,
use $\chi_j Z_l(x/a_j)$ in the second component and
$\chi_j Z_{B_j,l}(x/\delta_j)$ in the first. The common radial test
function in the first component is
\begin{equation}\label{7.17}
 \zeta_0=\sum_{i\in A}\chi_i(-2Z_A(x/s_i))
                  +\sum_{j\in B}\chi_j(-2Z_{B_j,0}(x/\delta_j)).
\end{equation}
For each test function $\zeta_\alpha$ in component $k(\alpha)$,
subtract its mean with respect to $V_{k(\alpha)}^{\rm app}\dif V_g$.
The resulting source is $V_{k(\alpha)}^{\rm app}\widehat\zeta_\alpha$,
with functional
$$
 \mathcal L_\alpha(\xi)
   =\int_MV_{k(\alpha)}^{\rm app}\widehat\zeta_\alpha
                                      \xi_{k(\alpha)}\dV.
$$
The number of orthogonality conditions is
$$
 d=2a+4b+1.
$$
The projected space is $E_t=\bigcap_\alpha\ker\mathcal L_\alpha$
in $\mathcal H$, with energy $\int_M\nabla\xi^TC\nabla\xi\dV$.
Write $\mathcal L_t$ for the linearization projected onto $E_t$.

\begin{proposition}
\label{pr7.2}
Under the hypotheses of Proposition~\ref{pr7.1}, for small $t$ and
$|\mu|\leq C_0t$, there is a unique correction
$\eta\in E_t\cap B_{Ct^{2/3}\ell}^{H^1}$ satisfying
\begin{equation}\label{7.18}
 \begin{aligned}
 -\Delta_g(u^{\rm app}+\eta)
 &=A\begin{pmatrix}
 \displaystyle\rho_1\frac{h_1e^{u_1^{\rm app}+\eta_1}}
      {\int_Mh_1e^{u_1^{\rm app}+\eta_1}\dV}-\rho_1\\[1.2ex]
 \displaystyle\rho_2\frac{h_2e^{u_2^{\rm app}+\eta_2}}
      {\int_Mh_2e^{u_2^{\rm app}+\eta_2}\dV}-\rho_2
 \end{pmatrix}\\
 &\quad+A\sum_\alpha e_{k(\alpha)}c_\alpha V_{k(\alpha)}^{\rm app}\widehat\zeta_\alpha.
 \end{aligned}
\end{equation}
The correction and multipliers depend continuously on the parameters
for $t>0$ and satisfy
\begin{equation}\label{7.19}
 \|\eta\|_{H^1}=O(t^{2/3}\ell),\qquad
 |c|=O(t^{2/3}\ell^2).
\end{equation}
The linearization projected onto $E_t$ has Morse index
$m_*+a+4b-1$.
\end{proposition}

\begin{proof}
Put $\delta=\sqrt t$ and $\nu_i=V_i^{\rm app}/\rho_i$.
The residual has $H^{-1}$ norm $O(\delta^{4/3})$ and the inverse
has operator norm $O(\ell)$. For fixed $K>0$ and
$\|v\|_{H^1},\|w\|_{H^1}\leq K\delta^{4/3}\ell$, the nonlinear
remainder satisfies
\begin{equation}\label{7.20}
 \|\mathcal N(v)-\mathcal N(w)\|_{H^{-1}}
 \leq C_\vartheta\delta^{-\vartheta}
       (\|v\|_{H^1}+\|w\|_{H^1})\|v-w\|_{H^1},
 \qquad\vartheta>0.
\end{equation}
To estimate the normalization, set
$m_i=\int_{M}\nu_iv_i\dV$, $n_i=\int_{M}\nu_i(e^{v_i}-1-v_i)\dV$ and
$D_i=m_i+n_i$. The exact remainder divided by $V_i^{\rm app}$ is
$$
 \frac{(e^{v_i}-1-v_i)-n_i-D_i(v_i-m_i)}{1+D_i}.
$$
On this ball, Lemma~\ref{le4.9} gives
$|m_i|\leq C\sqrt\ell\|v\|_{H^1}$ and
$|n_i|\leq C\ell\|v\|_{H^1}^2$, so $1+D_i=1+o(1)$.
H\"older's inequality and the exponential Sobolev inequality give
\eqref{7.20}, using
$\|V_i^{\rm app}\|_{1+\varepsilon}
\leq C_\varepsilon\delta^{-4\varepsilon/(1+\varepsilon)}$
with $\varepsilon>0$ sufficiently small for the prescribed $\vartheta$.
The fixed-point map
$\mathcal L_t^{-1}(-R^{\rm app}+\mathcal N(v))$
maps this ball into itself for sufficiently large fixed $K$. Its Lipschitz
constant is $O(\delta\ell^2)=o(1)$ with $\vartheta=1/3$.
The contraction mapping theorem gives existence and uniqueness.
The implicit function theorem gives parameter dependence, and
Corollary~\ref{co6.5} gives the multipliers and their bound in
\eqref{7.19}. By Corollary~\ref{co6.10}, the correction bound
preserves the projected Morse index.
\end{proof}

\subsubsection{The reduced equations}

Use the reduced quantities from \eqref{7.21} in the normalized
charts \eqref{7.6}--\eqref{7.7}.
In an arbitrary isothermal chart the coefficient is
$\psi^B_{j,zz}+(\psi^B_{j,z})^2$, with the transformation law
from the mixed-center calculation in Section~\ref{typed-sec:rotation}.
In the chosen charts the smooth corrections have jets
$$
 J_{1,j}=\frac{\kappa_j}4|z|^2+\Re(T_jz^2)+o(|z|^2),
 \qquad J_{2,j}=F_j^B\cdot x+O(|x|^2).
$$
The background equation and \eqref{7.15} give
\[
 \Delta_0 J_{1,j}(q_j)
 =\kappa_j+n_*^{\rm loc}(q_j)-n_{{\rm diff},j}(q_j)
 =\kappa_j.
\]

Write $c_i^A\in\mathbb R^2$ and $c_j^B\in\mathbb R^2$ for the
translation multipliers at $p_i$ and at the inner scale of $q_j$.
The multipliers for $B_j$ are $c_j^{\rm sh}\in\mathbb R^2$.
The common radial multiplier is $c_0$. We use the order
$(c^A,c^B,c^{\rm sh},c_0)$.

\begin{theorem}
\label{th7.3}
Uniformly on compact center and angular parameter sets, with
$|\mu|\leq C_0t$, the constructed family satisfies
\begin{align}
 c_i^A&=-3s_iF_i^A+o(t),\label{7.22}\\
 c_j^B&=-3a_jF_j^B+o(t),\label{7.23}\\
 c_j^{\rm sh}&=-\frac{3\epsilon_j}{16\pi}
       (\kappa_j\nabla\mathcal M(B_j)+4\pi\mathbf t_j)+o(t),
                    \label{7.24}\\
 c_0&=\frac3{8\pi m}
       \left(\mu-\sum_{j=1}^b\epsilon_j\kappa_j\mathcal H(B_j)\right)
       +o(t).\label{7.25}
\end{align}
The remainders are uniform in $C^0$.
\end{theorem}

\begin{proof}
Choose $r_0$ so that the cutoffs in \eqref{7.10} equal one on all
testing disks $B_{r_0}$. With the construction cutoffs fixed, we first let
$t\to0$ at fixed $r_0$ and then let $r_0\downarrow0$.
In the local pairings, $o(t)$ denotes an error that is uniform in
the parameters in this iterated limit.

In each core the corrected field is $U+J+\theta$, where
\begin{equation}\label{7.26}
 \theta_i=\log(\rho_i/M_i)+\eta_i
                     -\log\int_M\nu_i e^{\eta_i}\dV.
\end{equation}
H\"older's inequality with exponent $1+1/\ell$ and the exponential
Sobolev inequality give
\begin{equation}\label{7.27}
 \sum_{\rm cores}\sum_i\int_{B_{r_0}} n_i|\theta_i|^2e^{|\theta_i|}\dif A
 \leq C\ell\|\eta\|_{H^1}^2e^{C\ell\|\eta\|_{H^1}^2}+O(t^2)
 =O(t^{4/3}\ell^3)=o(t).
\end{equation}
At an ordinary center, $U_{A,1}=\log m_A$ and
$U_{A,2}=\log[\gamma_i(0)(s_i^2+r^2)]$. Thus
$-\Delta U_A=A(m_A,0)^T$ and
$J_{A,2}=\log(\gamma_i/\gamma_i(0))$.

On every core disk, let $q_\alpha^{\rm loc}
=e^\phi V_{k(\alpha)}^{\rm app}\widehat\zeta_\alpha$. Then
\[
 -\Delta(U+J+\theta)
 =A\binom{F_1}{F_2}+f^\rho
     +A\sum_\alpha e_{k(\alpha)}c_\alpha q_\alpha^{\rm loc},
 \qquad F_i=e^{U_i+J_i+\theta_i},
\]
where $f^\rho$ is the forcing in \eqref{2.2} at
$\rho=(4\pi m+\mu,\rho_2)^T$.
Write $v=J+\theta$ on each testing disk $D$. For a test vector $X$
and a boundary test vector $K$, set
\[
 \begin{aligned}
 \mathcal P_D(X)
   &=\sum_\alpha c_\alpha\int_D X_{k(\alpha)}q_\alpha^{\rm loc}\dif A,\\
 \mathcal B_D(K,v)
   &=\int_{\partial D}
      (\partial_{\boldsymbol n}K\cdot v-K\cdot\partial_{\boldsymbol n}v)\dif s .
 \end{aligned}
\]

On fixed annuli, the approximate densities are bounded. The residual
tends to zero in every finite $L^p$, and the multiplier terms also
tend to zero. By \eqref{7.19} and the exponential Sobolev inequality,
$e^{\eta_i}-1\to0$ in $L^p$ for every finite $p$. Local elliptic
estimates and the normalization give
\[
 \sup_{\partial B_{r_0}}(|\theta|+|\nabla\theta|)=o(1),\qquad
 \int_{M\setminus\bigcup B_{r_0}}f_1\dV=O_{r_0}(t^2).
\]
The last bound follows from $V_1^{\rm app}=O(t^2)$ on the complement,
$\int_M e^{\eta_1}\dV=O(1)$, and the denominator estimate in
Lemma~\ref{le6.9}.

At each mixed center, use the kernel functions
$$
 \Psi_j=-a_j\nabla_xU_j,\qquad
 \Phi_{j,l}=\partial_{B_{j,l}}U_j,\qquad
 E_j=\epsilon_j\partial_{\epsilon_j}U_j,
 \qquad E_\infty=(2,0),
$$
The last derivative holds $B_j,\beta_j$ fixed. For any planar
kernel function $X$, Green's identity gives
\[
 \begin{aligned}
 \mathcal P_D(X)
 ={}&\mathcal B_D(CX,v)-\int_D(CX)\cdot f^\rho\dif A\\
    &-\int_D\sum_i X_i n_i(e^{v_i}-1-v_i)\dif A .
 \end{aligned}
\]
At an ordinary center, the last integral uses $n=(m_A,0)^T$
and the right-hand side includes $-\int_{B_{r_0}}X_2 F_2\dif A$.

For a mixed translation,
$C\Psi_{j,l}=a_j(\partial_l\log F_j,\partial_l\log Q_j)^T$.
On each fixed circle, its leading part is
$a_j(4,2)^T x_l/r^2$. The smooth error is $O(a_jr)$
with gradient $O(a_j)$, and the remaining error is
$o(a_j)$ in $C^1$ for fixed $r_0$. The explicit polynomials also give
\[
 \int_{B_{r_0}}|C\Psi_{j,l}|\dif A\leq C a_jr_0.
\]
Indeed, $|\nabla\log F_j|\leq C/\delta_j$ for $r\leq C\delta_j$
and $C/r$ farther out, while
$|\nabla\log Q_j|\leq Cr/(r^2+a_j^2)$ on the disk.
The forcing integral is therefore $O(a_jr_0)$.

In the iterated limit, the translation boundary integral is
determined by the term $(F_j^B)_l x_l$ in $J_2$. For $l=1$ it gives
$-4a_j(F_j^B)_1\int_0^{2\pi}\cos^2\theta\dif\theta
=-4\pi a_j(F_j^B)_1$. For $l=2$ the same calculation uses
$\sin^2\theta$.

For $\Phi_{j,l}$, use the first-component boundary and forcing
pairings from the angular projection of Section~\ref{rev:sec:core}.
At the center, $-f_1^\rho(q_j)=\kappa_j+2\mu$, and the
$2\mu$ contribution to the forcing integral is $o(t)$ by the
bounded normalized first-component pairing.

To estimate the second-component forcing, put $\eps=\epsilon_j$,
$\beta=\beta_j$, $B=B_j$ and $r=|z|$. The identity from
\eqref{7.8} is $(C\Phi_{j,l})_2=-\partial_{B_l}\log Q_j$. With
$Q_j/\eps^2=\beta\eps^2+4r^2+\beta|z^2+\eps B|^2$, we have
\[
 Q_j/\eps^2\ge c(r^2+\eps^2),\qquad
 |\partial_{B_l}(Q_j/\eps^2)|\le C\eps(r^2+\eps).
\]
Here $\beta$ stays in a compact subset of $(0,\infty)$ and $B$ stays
in the fixed bounded parameter set. Hence the uniformly bounded
forcing $f_2^\rho$ satisfies
\[
 \left|\int_{B_{r_0}}(C\Phi_{j,l})_2 f_2^\rho\,\dif A\right|
 \le C\int_0^{r_0}\frac{\eps(r^2+\eps)}{r^2+\eps^2}r\,dr
 \le C\eps r_0^2+C\eps^2\log(1+r_0^2/\eps^2).
\]
Since $\eps\asymp t$, divide by $t$. Let $t\downarrow0$ at fixed $r_0$
and then let $r_0\downarrow0$. The result tends to zero uniformly
in the specified parameters.

For the boundary term, consider the limit of
$\eps^{-1}(C\Phi_{j,l})_2$ on each fixed circle. It is a bounded
angular function with radial derivative $O(r_0)$. In fact, it equals
the negative of $2\beta\,\operatorname{Re}(z^2\overline{e_l})/
(4r^2+\beta r^4)$ with $e_1=1$, $e_2=i$.
Since $J_2(q_j)=0$ and its first derivatives are bounded, both
terms of the Green boundary pairing, divided by $\eps$, tend to
zero as $r_0\downarrow0$.

The remainder estimates below give the following expressions
for the translation and angular pairings
$\mathcal R_j^B=\mathcal P_D(\Psi_j)$ and
$\mathcal R_j^{\rm sh}=\mathcal P_D(\Phi_j)$.
\[
 \mathcal R_j^B=-4\pi a_jF_j^B+o(t),\qquad
 \mathcal R_j^{\rm sh}
  =-2\epsilon_j(\kappa_j\nabla\mathcal M(B_j)+4\pi\mathbf t_j)+o(t).
\]
At an ordinary center, the translation and dilation kernel functions
satisfy
$$
 \begin{gathered}
 \Psi_l^A=-s_i\partial_{x_l}U_A,\qquad E_A=s_i\partial_{s_i}U_A,\\
 E_A=\left(2-\frac{4s_i^2}{s_i^2+r^2},
                 \frac{2s_i^2}{s_i^2+r^2}\right),\qquad
 K_A=C(E_A-E_\infty)=\left(-\frac{2s_i^2}{s_i^2+r^2},0\right),\\
 C\Psi_l^A=\left(\frac{2s_ix_l}{s_i^2+r^2},0\right),
 \qquad e^{U_{A,2}+J_{A,2}}=\gamma_i(x)(r^2+s_i^2).
 \end{gathered}
$$
The second-density term in the translation pairing is bounded by
\[
 s_i^{-1}\left|\int_D(\Psi_l^A)_2 F_2\dif A\right|
 \leq Cr_0^3+o(1).
\]
The dipole $C\Psi_l^A$ pairs with the linear part of $J_{1,A}$ on
the boundary to give $\mathcal R_i^A=\mathcal P_D(\Psi^A)
=-4\pi s_iF_i^A+o(t)$.

For the radial pairing at the $i$th ordinary center, put
$\Delta M_{1,i}^A=\int_D(F_1-m_A)\dif A$. Subtracting the model
equation gives
\[
 -\Delta v=A\binom{F_1-m_A}{F_2}+f^\rho+
                  A\sum_\alpha e_{k(\alpha)}c_\alpha q_\alpha^{\rm loc}.
\]
Pair this equation with the constant $CE_\infty$. Since
$E_\infty=(2,0)$, integration by parts gives
\[
 \mathcal B_D(CE_\infty,v)
   =2\Delta M_{1,i}^A+\int_D CE_\infty\cdot f^\rho\dif A
       +\mathcal P_D(E_\infty).
\]
Subtract this identity from the kernel identity for $E_A$, including
its second-density term. The resulting exact pairing is
\[
 \begin{aligned}
 \mathcal R_{i,\rm rad}^A
   :=\mathcal P_D(E_A-E_\infty)
   ={}&2\Delta M_{1,i}^A+\mathcal B_D(K_A,v)
      -\int_D K_A\cdot f^\rho\dif A\\
     &-\int_D E_{A,1}m_A(e^{v_1}-1-v_1)\dif A
      -\int_D E_{A,2}F_2\dif A .
 \end{aligned}
\]
On each fixed testing circle, $K_A$ and its normal derivative
are $O_{r_0}(s_i^2)$, and $v$ and its normal derivative are bounded.
The forcing term is $O_{r_0}(s_i^2\log(1/s_i))$ by radial
integration. The last integrand is
$2s_i^2\gamma_i(x)e^{\theta_2}$, which has integral $O(s_i^2)$
by \eqref{7.26} and the exponential Sobolev inequality.
For the nonlinear term, $E_{A,1}$ is bounded and
$J_{1,A}=O(r)$. The quadratic bound and \eqref{7.27} give
\[
 \int_D m_A|e^{v_1}-1-v_1|\dif A
 \leq C\int_D r^2m_A\dif A+O(t^{4/3}\ell^3)
 =O(s_i^2\log(1/s_i)+t^{4/3}\ell^3).
\]
Thus $\mathcal R_{i,\rm rad}^A=2\Delta M_{1,i}^A+o(t)$ for fixed
$r_0$. We will combine these mass differences with those on the mixed disks.

At a mixed center, for $\Xi=\Phi$ or $E-E_\infty$,
$$
 \begin{aligned}
 \sum_i\int_{B_{r_0}}|\Xi_i|n_i|J_i|^2\dif A
       &\leq C(t^2\ell+tr_0^2),\\
 \sum_i\int_{B_{r_0}}|\Psi_i|n_i|J_i|^2\dif A
       &\leq C(t^2\ell+tr_0^3).
 \end{aligned}
$$
We test the nonlinear term in \eqref{5.10} against $E$.
The contribution of $E_\infty=(2,0)$ is bounded by
$C\int_{B_{r_0}}n_1(J_1^2+\theta_1^2e^{|\theta_1|})\dif A$.
At an ordinary center, $n_1=m_A$ and the second moment gives
$\int_{B_{r_0}}n_1J_1^2\dif A=O(t^2\ell)$. At a mixed center,
$J_1=O(r^2)$, the first mass is bounded, and $n_1\leq Ct^2r^{-6}$
for $\delta_j\leq r\leq r_0$. Integrating $r^4n_1$ gives the same
bound. By \eqref{7.27}, the contribution of $E_\infty$ is
$O(t^2\ell+t^{4/3}\ell^3)=o(t)$.

By Cauchy--Schwarz, \eqref{7.27}, and the convergence on
$\partial B_{r_0}$, the remaining terms satisfy
$$
 \lim_{r_0\downarrow0}\limsup_{t\downarrow0}
 \frac{|\mathcal E_{\rm nonlinear}|+|\mathcal E_{\rm boundary}|}{t}=0.
$$

On the $j$th mixed disk $D$, set
$\Delta M_{1,j}^B=\int_D(F_1-n_{1,j})\dif A$ and
$K_j=C(E_j-E_\infty)$. Subtract the constant test as above to obtain
\[
 \mathcal P_D(E_j-E_\infty)
 =2\Delta M_{1,j}^B+\mathcal B_D(K_j,v)
  -\int_D K_j\cdot f^\rho\dif A
  -\int_D\sum_i E_{j,i}n_i(e^{v_i}-1-v_i)\dif A.
\]
Equations \eqref{5.11}--\eqref{5.13} give the leading forcing
pairing $2\epsilon_j\kappa_j\mathcal M(B_j)$. Here
$-f_1^\rho(q_j)=\kappa_j+2\mu$, and $\epsilon_j\mu=O(t^2)$.
The first boundary component gives $-8\pi\epsilon_j\Re(T_jB_j)$
as in \eqref{5.14}, using the harmonic quadratic term of $J_{1,j}$.
For the second component, $K_{j,2}/\epsilon_j$ has a bounded limit
on fixed matching circles with limiting radial derivative $O(r_0)$.
The bounds $J_{2,j}=O(r)$ and $\nabla J_{2,j}=O(1)$, together with
the convergence of $\theta$ on each fixed circle, give a boundary
pairing of $o(t)$ in the iterated limit. Consequently
\[
 \mathcal P_D(E_j-E_\infty)
 =2\Delta M_{1,j}^B
  -2\epsilon_j\bigl(\kappa_j\mathcal M(B_j)+4\pi B_j\cdot\mathbf t_j\bigr)
  +o(t).
\]

The source pairings are
$$
 \begin{gathered}
 \int_{\C} m_AZ_{A,l}Z_{A,r}\dif A=\frac\pi3\delta_{lr},\qquad
 \int_{\C} m_AZ_A^2\dif A=\frac{4\pi}3,\qquad (\Psi_l^A)_1=4Z_{A,l},\\
 \mathcal Z_B=(-2Z_{B,0},Z_{B,1},Z_{B,2}),\qquad
 \int_{\C} m_B\mathcal Z_{B,\alpha}\mathcal Z_{B,\beta}\dif A
       =\frac{32\pi}3\delta_{\alpha\beta},\\
 E_1\longrightarrow-2Z_{B,0}+B\cdot(Z_{B,1},Z_{B,2}).
 \end{gathered}
$$
The outer-profile identity follows from the change of variables
$w=y^2$. The other limiting cross terms vanish by centering. With
$\Gamma_B=32\pi/3$ and $N=a+b$,
$$
 \mathsf G_0=
 \begin{pmatrix}
 (4\pi/3)I_{2N}&0&0\\
 0&\Gamma_BI_{2b}&0\\
 0&\Gamma_B(B_1^T,\ldots,B_b^T)&\Gamma_0
 \end{pmatrix},\qquad \det\mathsf G_0>0,
$$
whose last row is
\begin{equation}\label{7.28}
 \Gamma_0c_0+\sum_{j=1}^b\Gamma_BB_j\cdot c_j^{\rm sh},\qquad
 \Gamma_0=\frac{16\pi}3a+\frac{32\pi}3b=\frac{16\pi m}3.
\end{equation}
The total first-component mass is $\rho_1$. The full model masses
sum to $4\pi a+8\pi b=4\pi m$. The first-component mass outside the
testing disks and the tails of the first model densities are
$O_{r_0}(t^2)$. Thus
\[
 \sum_{i=1}^a\Delta M_{1,i}^A+\sum_{j=1}^b\Delta M_{1,j}^B
     =\mu+O_{r_0}(t^2).
\]
Let $\mathsf G(t,r_0)$ be the exact source matrix with the
translation and angular rows above. Its last row is the sum of the
shifted radial pairings. Subtracting the constant $E_\infty$ on each
disk changes each entry of the summed row by $O_{r_0}(t^2)$.
Indeed, each first-component source has zero integral on $M$,
and $V_1^{\rm app}=O_{r_0}(t^2)$ outside the testing disks.
Therefore $\mathsf G(t,r_0)\to\mathsf G_0$ for fixed $r_0$,
uniformly on the parameter set, by the density comparison in
Proposition~\ref{pr7.1} and the profile pairings.

The right-hand sides of the exact system are $O(t)$ for fixed
sufficiently small $r_0$. For the radial row, this bound uses the
sum of the local mass differences, which contributes $2\mu$ to leading
order. Inverting the exact matrix first gives
$|c|=O(t)$. For each fixed $r_0$ we may now replace
$\mathsf G(t,r_0)c$ by $\mathsf G_0c+o(t)$. Taking the stated
iterated limit gives the radial equation
\begin{equation}\label{7.29}
 \Gamma_0c_0+\sum_j\Gamma_BB_j\cdot c_j^{\rm sh}
 =2\mu-2\sum_j\epsilon_j
       (\kappa_j\mathcal M(B_j)+4\pi B_j\cdot\mathbf t_j)+o(t).
\end{equation}

The correction, multipliers, and leading pairings are independent
of $r_0$. For each row the limsup as $t\downarrow0$ of the remainder
divided by $t$ tends to zero as $r_0\downarrow0$, uniformly on the
parameter set. Dividing the center and angular equations by their
limiting Gram factors gives \eqref{7.22}--\eqref{7.24}.
Finally,
$$
 \begin{aligned}
 \Gamma_Bc_j^{\rm sh}
 &=-2\epsilon_j(\kappa_j\nabla\mathcal M(B_j)+4\pi\mathbf t_j)+o(t),\\
 \Gamma_0c_0
 &=2\mu-2\sum_j\epsilon_j
       (\kappa_j\mathcal M(B_j)+4\pi B_j\cdot\mathbf t_j)
       -\sum_j\Gamma_BB_j\cdot c_j^{\rm sh}+o(t)\\
 &=2\mu-2\sum_j\epsilon_j\kappa_j
          (\mathcal M(B_j)-B_j\cdot\nabla\mathcal M(B_j))+o(t).
 \end{aligned}
$$
This is \eqref{7.25}.
\end{proof}

\subsection{Reduction to finite dimensions}\label{sec:fiber-degree}

We reduce the degree in the parameter-dependent projected spaces
by a homotopy in the correction bundle.
Fix the nonzero parameter $\mu$ at which we compute the degree. Let
$\Omega$ be the open parameter domain, $P_p$ the energy projection
onto the correction space $E_p$, and $S_{\alpha,p}$ the Riesz
representatives of the constraint functionals for the energy inner
product. Write $u=u^{\rm app}(p)+\eta$ in the
open image $\mathcal U$ of the correction bundle. Let $\eta_p$
be the projected solution and $u_p=u^{\rm app}(p)+\eta_p$.
Then $P_pF(u_p)=0$ and $F(u_p)=\sum_\alpha g_\alpha(p)S_{\alpha,p}$.
The graph is $C^1$ by the implicit function theorem because the
projected linearization at the corrected profile is invertible in a
local trivialization of $E_p$.

Let $L_p=P_pDF(u^{\rm app}(p))|_{E_p}$ and $w=\eta-\eta_p$.
For $0\leq\lambda\leq1$, use the following map on $\mathcal U$:
\[
\begin{aligned}
 H_\lambda(u)
 ={}&L_pw+\lambda\{P_p[F(u)-F(u_p)]-L_pw\}\\
 &+(I-P_p)\{(1-\lambda)F(u_p)+\lambda F(u)\}.
\end{aligned}
\]
This is an identity-minus-compact homotopy on the bounded domain.
Indeed $F=I-K_\rho$, the projections differ from the identity
by finite-rank terms with continuously varying representatives,
and the compact terms involving $K_\rho$ and its derivative
are uniform over the compact parameter closure. The remaining
translations depend continuously on that closure. Use the estimate
for the projected contraction on a larger ball of the same order if
needed. If the first line vanishes, we obtain
$\|w\|_{H^1}\leq o(1)\|w\|_{H^1}$ after applying $L_p^{-1}$ to the
difference of nonlinear remainders. Hence $w=0$. At such a point,
the second line is $\sum_\alpha g_\alpha(p)S_{\alpha,p}$.
The homotopy has no boundary zero if $g$ has none on
$\partial\Omega$ and $\eta_p$ lies strictly inside the correction
ball. Lemma~\ref{lem:parameter-modulation} identifies the boundary
with these two possibilities, namely the parameter boundary and
the correction boundary.

At $\lambda=0$ approximate $g$ uniformly by a smooth map with
regular zero, keeping the parameter boundary free of zeros. This
perturbs only the finite-dimensional source term. Near one of its
zeros, use $(p,w)$ coordinates about the graph and identify the
fiber at that point with $E_p$. The pairing matrix
$\bigl(\mathcal L_{\alpha,p}D_{p_\beta}u_p\bigr)$ has positive
determinant: in the scaled parameter coordinates it differs
from the pairing matrix for $u^{\rm app}$ by
$O(\sqrt\ell\|\eta_p\|_{H^1})=o(1)$, by differentiating
$\mathcal L_{\alpha,p}(\eta_p)=0$ and using the derivative
bound in Lemma~\ref{lem:parameter-modulation}. The pairing
matrix for $u^{\rm app}$ has the positive limit \eqref{7.34}
(or its $b=0$ version). In these coordinates the upper block of
the source identification is the inverse pairing matrix times
the positive Gram matrix. The derivative of $H_0$ is therefore
block triangular up to this positive change of finite-dimensional
coordinates, with blocks $Dg$ and $L_p$. The local degree is the
sign of $\det Dg$ times $(-1)^\nu$, where $\nu$ is the projected
Morse index. This follows by the finite-dimensional approximation
of the degree or by the linear product formula. Summing over the
regular zeros and removing the approximation gives
\begin{equation}\label{7.30}
 \deg_{\rm LS}(F,\mathcal U,0)
                  =(-1)^\nu\deg_{\rm B}(g,\Omega,0).
\end{equation}
The uniform $C^0$ estimate below gives an admissible homotopy
from $g$ to its leading map.

\subsection{Local degree at a mixed configuration}

\begin{theorem}
\label{th7.4}
Suppose $b\geq1$, $(P^0,Q^0,w^0)$ is a nondegenerate solution of
\eqref{3.40}, and its scalar linearization at fixed centers is invertible.
Assume $\kappa_j>4|T_j|$ at each mixed center. For sufficiently small
$\mu>0$, the Leray--Schauder degree in the bubbling neighborhood of
the family in Proposition~\ref{pr7.2} is
\begin{equation}\label{7.31}
 (-1)^ad_S(P^0,Q^0,w^0)
   =(-1)^{a+m_*}\operatorname{sgn}\det D_{(P,Q)}F(P^0,Q^0).
\end{equation}
For $\mu\leq0$, $|B_j|\leq C$ and $|\mu|\leq C_0t$,
this family has no zero for sufficiently small $t$.
\end{theorem}

\begin{proof}
Write $\epsilon_j=td_j(P,Q)$, with $d_j>0$, and set
$$
 \mathcal S(P,Q,B)=\sum_{j=1}^b d_j(P,Q)\kappa_j\mathcal H(B_j)>0.
$$
The positivity follows from
$\mathcal H(B)=\frac12\int_\C(1+|v|^2)^{-1}|v+B|^{-1}\dif A_v$.
For $\mu>0$, set $t=\mu\tau$. Divide each two-dimensional center
and angular block by its negative factor in \eqref{7.22}--\eqref{7.24},
and multiply the radial row by $8\pi m/(3\mu)>0$. These operations
preserve orientation. The reduced map converges uniformly to
\begin{equation}\label{7.32}
 (P,Q,B,\tau)\longmapsto
 \left(F(P,Q),\,
       (\kappa_j\nabla\mathcal M(B_j)+4\pi\mathbf t_j)_{j=1}^b,\,
       1-\tau\mathcal S(P,Q,B)\right).
\end{equation}
Evaluate $F$ on the scalar branch $w=w(P)$.

Choose a small center neighborhood $\Omega$ that isolates $(P^0,Q^0)$
and on whose closure all inequalities $\kappa_j>4|T_j|$ remain strict.
Since $\mathcal M$ is radial and $\mathcal M'(R)\to\pi$, a sufficiently large
fixed $R$ satisfies the following inequality for $|B_j|=R$ and $0\leq\lambda\leq1$.
$$
 \bigl(\kappa_j\nabla\mathcal M(B_j)
                  +4\pi\lambda\mathbf t_j\bigr)\cdot B_j
 \geq R(\kappa_j\mathcal M'(R)-4\pi|T_j|)>0.
$$
Choose
$$
 0<S_-\leq\mathcal S\leq S_+,\qquad
 0<\tau_-<S_+^{-1}\leq S_-^{-1}<\tau_+,
 \qquad S_0\in[S_-,S_+].
$$
Use $C_0\geq1/\tau_-$ in Proposition~\ref{pr7.2}, so that the
construction applies throughout this domain.
For $0\leq\lambda\leq1$, put
$$
 G_j^\lambda(B_j)=\kappa_j\nabla\mathcal M(B_j)
                              +4\pi\lambda\mathbf t_j,
 \qquad \mathcal S_\lambda=(1-\lambda)S_0+\lambda\mathcal S.
$$
On the boundary of $\Omega\times\prod_jB_R\times(\tau_-,\tau_+)$,
$$
 \begin{aligned}
 F&\ne0&&\text{on }\partial\Omega,\\
 G_j^\lambda(B_j)\cdot B_j&>0&&\text{if }|B_j|=R,\\
 1-\tau_-\mathcal S_\lambda&>0>1-\tau_+\mathcal S_\lambda.
 \end{aligned}
$$
The homotopy has no zeros on the boundary. At $\lambda=0$, each
angular degree is $+1$ and the radial degree is $-1$. Hence
\begin{equation}\label{7.33}
 \deg_{\rm B}=-\operatorname{ind}_{(P^0,Q^0)}F.
\end{equation}
The uniform error estimate preserves this degree. For $\mu\leq0$,
the radial equation $\mu-t\mathcal S+o(t)=0$ has no solution.

By Lemma~\ref{lem:parameter-modulation}, the map from parameters
and corrections to $u$ is injective on this domain. The pairings of
the source functionals with the parameter derivatives converge to
\begin{equation}\label{7.34}
 J_0=
 \begin{pmatrix}
  \frac{4\pi}{3}I_{2(a+b)}&0&0\\
  0&\frac{32\pi}{3}I_{2b}&\frac{32\pi}{3}\mathbf B\\
  0&0&\frac{16\pi m}{3}
 \end{pmatrix},\qquad
 \mathbf B=(B_1,\ldots,B_b)^T,\qquad \det J_0>0.
\end{equation}
Here the blocks are in real coordinates.

By Proposition~\ref{pr7.2}, the projected correction lies in the
interior of its ball and the projected Morse index is $m_*+a+4b-1$.
Combining \eqref{7.30} and \eqref{7.33} gives
$$
 (-1)^{m_*+a+4b-1}
       \bigl(-\operatorname{ind}F\bigr)
       =(-1)^a d_S,
$$
as asserted.
\end{proof}

\begin{remark}\label{re7.5}
For $b=0$, the first possible logarithmic mass term has order
$t^2\log(1/t)$ by the calculation in Section~\ref{rev:sec:ordinary}.
This is the leading term when its coefficient $L_0$ is nonzero.
Lemma~\ref{le8.4} ensures that the coefficient is nonzero for
the weights used to compute the degree jump across the wall.
\end{remark}

\subsection{Localization of mixed blowup solutions}\label{sec:typedactual-chapter}

Let $u_\nu$ be a sequence of Toda solutions concentrating at the fixed
shadow solution. We seek parameters for which
$$
 u_\nu=u^{\rm app}(t_\nu,P_\nu,Q_\nu,B_\nu)+\eta_\nu,
 \qquad \eta_\nu\in E_{t_\nu}.
$$
Set $E=\lVert\eta\rVert_{H^1}$ and use $D,\ell$ from
\eqref{7.36}. We prove
$$
 \begin{aligned}
 \lvert\mu\rvert&=O(D^\alpha),\\
 E&\leq C(D^{4/3}\ell+\lvert\mu\rvert\ell^{3/2}),\\
 \lvert\mu\rvert&\leq C(D^2+D\sqrt\ell E+\ell E^2+\mu^2).
 \end{aligned}
$$
Lemma~\ref{le5.18} then gives
$\lvert\mu\rvert=O(D^2)$ and $E=O(D^{4/3}\ell)$.

We use the densities and canonical fields in \eqref{2.1}--\eqref{2.2},
with the limiting masses and scalar equation of Lemma~\ref{le2.1}.
Assume that the limiting centers are separated and that the scalar
solution of \eqref{2.5} is nondegenerate. At a limiting mixed
center $q_j$, use the fixed chart of \eqref{7.6}--\eqref{7.7}, so
$\phi(q_j)=0$ and $\partial_z\psi_j^B(q_j)=0$. Set
$\kappa_B=\kappa_j$ and $T_B=T_j$, where $\kappa_j$ and $T_j$ are
defined in \eqref{7.21}. The trace coefficient of the matched
smooth correction $J_{1,j}$ is $\kappa_j$ by the cancellation
$n_*^{\rm loc}(q_j)=n_{{\rm diff},j}(q_j)$ in that computation.
The quadratic expansion of the regular field $\psi_j^B$ has
harmonic coefficient $T_j$. The first derivative of the regular
field vanishes at the normalized center, so $T_j=\mathcal T_j$ from
\eqref{3.41}. Equation \eqref{2.2} gives
$\kappa_j=-f_1^{(4\pi m,\rho_2)}(q_j)=\kappa_m(q_j)$. Thus
\eqref{3.41} in these charts gives
\begin{equation}\label{7.35}
 \kappa_B>4|T_B|+\eta_0\qquad(\eta_0>0),
\end{equation}
with a common $\eta_0$ for the finite set of limiting configurations.
Proposition~\ref{pr4.19} applies with the regular fields from
\eqref{7.2}. It gives the profile $m_{B_{k,*}}$, where $B_{k,*}$
solves \eqref{4.33} in that chart, and the bounds in \eqref{4.34}.

\subsubsection{Local mass and profile estimates}

Choose the moving ordinary centers $p_{j,\nu}$ to be local maxima of
$f_{1,\nu}$ as in Lemma~\ref{le8.6}. In the charts of \eqref{7.5},
$\phi(p_{j,\nu})=0$ and $\nabla\phi(p_{j,\nu})=0$, so
$\nabla F_{1,\nu}(p_{j,\nu})=0$. Define $s_{j,\nu}$ by
$s_{j,\nu}^2F_{1,\nu}(p_{j,\nu})=4$.

At a mixed center, let $q_{k,\nu}$ be the physical inner peak
given by Proposition~\ref{pr4.19} in a fixed chart. Use the chart
of \eqref{7.6}--\eqref{7.7} centered at this point, with the scalar
branch \eqref{7.1} at the moving ordinary centers. With the rotation
fixed, the old coordinate centered at $q_{k,\nu}$ is
$z_{\rm old}=g_\nu(z)=\lambda_\nu z+O(|z|^2)$, where
$\lambda_\nu>0$ and $\lambda_\nu,\lambda_\nu^{-1}$ are bounded.
The coordinate density is multiplied by $|g_\nu'|^2$. Divide the
preliminary inner and outer scales by $\lambda_\nu$.
Thus the profile limits persist in these charts. The pointwise bound
in \eqref{4.34} persists because $|g_\nu(z)|\asymp |z|$.
By Green's representation, this bound and the mass bound give
$|\nabla U_{i,\nu}|\leq C/r$ on $|z|=r$ and hence the oscillation
bound there. Reset $a_{k,\nu}$ by imposing
$a_{k,\nu}^2F_{2,\nu}(q_{k,\nu})=4$. The local first mass,
the transformed profiles and the annular bounds just verified
give \eqref{eq:inner-first-density-small} in these moving charts.
The reset preserves this estimate and the annular bounds as shown
there. The Jacobian contributes $O(a_{k,\nu})$ to the gradient condition
in \eqref{5.45}, as in the inner estimate of Section~\ref{sec:strong-chapter}.

Take disjoint coordinate disks $D_{A,j}$ and $D_{B,k}$ of fixed
radius about the moving centers. Apply Lemma~\ref{le8.6} on each
ordinary disk. On each mixed disk, \eqref{2.4}--\eqref{2.6} give
the local mass and outer convergence hypotheses of Lemma~\ref{le5.6},
Proposition~\ref{pr5.8}, and Section~\ref{sec:strong-chapter}.
The limiting second density is smooth and positive there. The forcing
is $f^\rho$ from \eqref{2.2}, and the Green potentials of sources
outside each disk have uniformly bounded derivatives on smaller disks.

Use \eqref{5.46} to modulate the outer scale and shape on each mixed
disk in these charts. Let $\widehat\delta_{k,\nu}$ be the preliminary
outer scale. The resulting pair $(\delta_{k,\nu},B_{k,\nu})$ satisfies
$\delta_{k,\nu}/\widehat\delta_{k,\nu}\to1$ and
$B_{k,\nu}\to B_{k,*}$. Use this pair in \eqref{5.60}, retaining
$\tau_{k,\nu}=\widehat\delta_{k,\nu}^{3/2}\asymp\delta_{k,\nu}^{3/2}$.

Suppress $\nu$ in the local notation and set
\begin{equation}
 D=\max(\delta_1,\ldots,\delta_b,\sqrt{s_1},\ldots,\sqrt{s_a}),
 \qquad \ell=|\log D|.
 \label{7.36}
\end{equation}
The maximum ranges over the centers present. We consider $b\ge1$
and show below that the scales are comparable.
The local bounds are
\begin{equation}\label{7.37}
 F_1\leq C_\eps s^{-2}(1+r/s)^{-4+\eps},\qquad
 F_2\leq C_\eps(s+r)^{2-\eps},\qquad
 |\nabla Gf_1|\leq C/(s+r)
\end{equation}
at an ordinary point, and
\begin{equation}\label{7.38}
 a=\delta^{2+o(1)},\qquad
 \begin{cases}
 F_1\leq C_\eps\delta^{-4+\eps}r^{2-\eps},\
 F_2\leq C_\eps a^{2-\eps}r^{-4+\eps},&Ra<r<\delta/R,\\
 F_1\leq C_\eps\delta^{4-\eps}r^{-6+\eps},\
 F_2\leq C_\eps r^{-\eps},&R\delta<r<r_0
 \end{cases}
\end{equation}
at a mixed point. At $r\asymp\delta$, the bounds are
$F_1\leq C\delta^{-2}$ and $F_2\leq C_\eps\delta^{-\eps}$. At
$r\lesssim a$, the corresponding bounds are
$F_1\leq C_\eps\delta^{-4+\eps}a^{2-\eps}$ and $F_2\leq Ca^{-2}$.
In particular, the second density is locally Lipschitz at ordinary
points.

\subsubsection{Profile and potential comparisons}
Write the three model densities as
\begin{equation}\label{7.39}
 \begin{aligned}
 m_A^s(z)&=\frac{4s^2}{(s^2+|z|^2)^2},&
 m_i^a(z)&=\frac{4a^2}{(a^2+|z|^2)^2},\\
 m_o^{\delta,B}(z)&=
 \frac{16\delta^4|z|^2}{(\delta^4+|z^2-\delta^2B|^2)^2}.
 \end{aligned}
\end{equation}
For the Green potentials, cut off each model in its fixed coordinate
disk and use $m\dif A=e^{-\phi}m\dV$. At an ordinary center,
Lemma~\ref{le8.6} gives \eqref{8.30} in the chart of \eqref{7.5}, namely
\begin{equation}\label{7.40}
 \left|\log\frac{s^2F_1(sy)}{m_A^1(y)}\right|
       \leq Cs^{1/4}\log(2+|y|),\qquad |y|\leq s^{-1/2}.
\end{equation}
At a mixed center, Propositions~\ref{pr5.14} and \ref{pr5.16}
give the inner and outer errors \eqref{5.55} and \eqref{5.60}.
In Lemma~\ref{le5.15}, $\Delta_0\Phi=F_2-m_i^a-f_1^\rho$ by
\eqref{2.2}.

Set $\Omega_{\rm ext}=M\setminus(\bigcup_j D_{A,j}\cup\bigcup_k D_{B,k})$.
The first Green potential splits as
\[
 \begin{aligned}
 Gf_1-\sum_A Gm_A-\sum_B Gm_o
  ={}&\sum_j G(\mathbf 1_{D_{A,j}}f_1-m_A^{s_j})\\
     &+\sum_k G(\mathbf 1_{D_{B,k}}f_1-m_o^{\delta_k,B_k})
       +G(\mathbf 1_{\Omega_{\rm ext}}f_1).
 \end{aligned}
\]
On an ordinary disk, use \eqref{7.40} and the tail bound in
\eqref{7.37}. On a mixed disk, use the local first-density
estimate \eqref{5.62} and the local second-density
estimate \eqref{5.61} with the cutoff $\tau_k$.
$f_1=O(e^{c_1})$ with $c_1=\log(\rho_1/\int_Mh_1e^{u_1}\dV)$.
The outer bound on any mixed disk gives $e^{c_1}\leq C_\eps D^{4-\eps}$,
so the exterior potential has the same $L^\infty$ bound.
Since $x^\gamma|\log x|^N$ increases for small $x$ when $\gamma>0$,
\eqref{7.36} gives, for some fixed integer $N$,
\begin{equation}\label{7.41}
 \begin{aligned}
 \left\|G\!\left(f_1-\sum_A m_A^{s_j}
                        -\sum_B m_o^{\delta_k,B_k}\right)\right\|_\infty
     &\leq CD^{1/4}\ell^N,\\
 \left\|G\!\left(\sum_B1_{D_{\tau_k}}f_2
                              -\sum_Bm_i^{a_k}\right)\right\|_\infty
     &\leq CD^{1/4}\ell^N.
 \end{aligned}
\end{equation}
The omitted Green tails have orders $O(s^{1-\eps/2}|\log s|)$,
$O(\delta^{2-\eps/2}|\log\delta|)$, and
$O(\delta^{1-\eps/2+o(1)}|\log\delta|)$, respectively.
The last estimate concerns the tail of the inner model and any
omitted portion of the actual second density inside its inner
comparison region. For a cutoff $r_*\asymp\delta_k^{3/2}$,
the estimate $Ca_k^{2-\eps}r^{-4+\eps}$ in that region gives a Green
potential of size at most $C(a_k/r_*)^{2-\eps}(1+|\log r_*|)$.
To obtain a bound uniform in the evaluation point $y$, split the integral
at $d(x,y)=r_*$. Bound the near part by the supremum of the density on
$r\ge r_*$ times $Cr_*^2(1+|\log r_*|)$. Bound the far part by the
mass of the tail times $C(1+|\log r_*|)$. For the model, the
same proof uses $\eps=0$. With $a_k=\delta_k^{2+o(1)}$ these
give the displayed order. The actual density on
$M\setminus\bigcup_BD_{\tau_k}$ belongs to $w_D$ below.
The first mass on a fixed mixed disk is $8\pi+O(\delta_k^{1/4})$
by Proposition~\ref{pr5.8}. At an ordinary center, integrating
\eqref{7.40} and using \eqref{7.37} on the rest of the disk gives
$4\pi+O(s_j^{1/4})$. Summing these masses and the exterior estimate yields
\begin{equation}\label{7.42}
 \mu:=\rho_1-\rho_1^*=O(D^{1/4}).
\end{equation}

\subsubsection{Estimates for the residual scalar solution}

With the radii $\tau_k\asymp\delta_k^{3/2}$ above, set
$$
 w_D=2G(1_{M\setminus\cup_BD_{\tau_k}}f_2),\quad
 H_{P,s}=h_2e^{-\sum_AGm_A^{s_j}},\quad
 T_k=2Gm_i^{a_k}-Gm_o^{\delta_k,B_k}.
$$
By Green's representation and \eqref{7.41},
\begin{equation}
 u_2=w_D-\sum_AGm_A^{s_j}+\sum_BT_k+E_D,\qquad
 \|E_D\|_\infty\le CD^{1/4}\ell^N .
 \label{7.43}
\end{equation}
The bounds in \eqref{7.37}--\eqref{7.38} give
$w_D\to w_*$ in $W^{2,6/5}$. On the annulus from $\tau$ to the
outer scale at a mixed center, the integral of $F_2^{6/5}$ is
$$
 O\!\left(a^{(2-\eps)6/5}
       \tau^{\,2-(4-\eps)6/5}\right)
 =O(\delta^{3/5-(3/5)\eps+o(1)}).
$$
Farther out, $r^{-\eps}$ belongs to $L^{6/5}$ for small $\eps$.
At an ordinary center the second density satisfies the stronger
bound in \eqref{7.37}.

The ordinary Green identity is
$$
 Gm_A^s=-\log(s^2+r^2)+4\pi R(x,p)+O(s^2|\log s|).
$$
It implies
\begin{equation}
 \|H_{P,s}-H_P\|_\infty
       \le C\sum_A s_j^2|\log s_j|\le CD^4\ell .
 \label{7.44}
\end{equation}
Near $p_j$, write $H_P=r^2k_j(x)$ with $k_j$ smooth and bounded.
Multiplication by $1+s_j^2/r^2$ replaces $r^2k_j$ with
$(r^2+s_j^2)k_j$, with uniform error $O(s_j^2)$. The factors from
the other centers are smooth and uniformly bounded there.

The mixed Green identities give
$$
 T_k(z)=\log\frac{\delta_k^4+|z^2-\delta_k^2B_k|^2}
                       {(a_k^2+|z|^2)^2}
                      +O(\delta_k^2|\log\delta_k|),
$$
and, outside $D_{\tau_k}$,
\begin{equation}
 \|\boldsymbol 1_{M\setminus D_{\tau_k}}(e^{T_k}-1)\|_{6/5}
 \le C(\delta_k^2|\log\delta_k|
    +\delta_k^2\tau_k^{-1/3}
    +\delta_k^4\tau_k^{-7/3})
 \le C\delta_k^{1/2}|\log\delta_k|.
 \label{7.45}
\end{equation}
Outside the $k$-th disk, $e^{T_k}=1+O(\delta_k^2|\log\delta_k|)$.
Since the disks are disjoint, \eqref{7.45} also controls the product
of these factors.

Write $c_2^D=\rho_2/\int_{M} h_2e^{u_2}\dV$. Equations \eqref{7.43}--\eqref{7.45}
yield
$$
 \|1_{M\setminus\cup_BD_{\tau_k}}f_2
          -c_2^D H_Pe^{w_D}\|_{6/5}
       \le CD^{1/4}\ell^N .
$$
Integrating and using the inner second-component masses of
Proposition~\ref{pr5.8} gives
\[
 c_2^D\int_{M} H_Pe^{w_D}\dV=\sigma+O(D^{1/4}\ell^N).
\]
Write $w_*(P)=w(P)$ for the scalar branch of \eqref{7.1} evaluated
at the moving ordinary centers. The equation for $w_D$ yields
\[
 \left\|-\Delta_gw_D
 -2\sigma\left(\frac{H_Pe^{w_D}}{\int_M H_Pe^{w_D}\dV}-1\right)
 \right\|_{6/5}\leq CD^{1/4}\ell^N .
\]
Along this branch, the fixed-center differential operator is invertible
as a map $W^{2,6/5}_0\to L^{6/5}_0$. It is a compact perturbation of the
Poisson isomorphism. By Sobolev embedding, a kernel element first belongs
to $C^0$. The smooth coefficient in the scalar row and elliptic
regularity then place the element in $X_s$, where fixed-center
nondegeneracy forces it to vanish. Continuity of the coefficients
gives a uniform inverse bound near the limit. Since
$w_D-w_*(P)=o(1)$ in $W^{2,6/5}$, we absorb the quadratic
remainder and obtain
\begin{equation}
 \|w_D-w_*(P)\|_{W^{2,6/5}}
   +\left|c_2^D-\frac{\sigma}{\int_{M} H_Pe^{w_*(P)}\dV}\right|
       \le CD^{1/4}\ell^N .
 \label{7.46}
\end{equation}
This estimate uses the branch evaluated at the moving centers.

\subsubsection{Comparison of the scales}

Green's representation, \eqref{7.41}, and \eqref{7.46} give the following
global $L^\infty$ approximation with error $O(D^{1/4}\ell^N)$.
\begin{equation}
 \mathcal{U}^{G}=
 \begin{pmatrix}
  2\sum_AGm_A+2\sum_BGm_o-\sum_BGm_i-w_*(P)/2\\
  -\sum_AGm_A+2\sum_BGm_i-\sum_BGm_o+w_*(P)
 \end{pmatrix}.
 \label{7.47}
\end{equation}
The values of the regular first-component fields are
\begin{align}
 S_A(p)&=8\pi R(p,p)+8\pi\!\sum_{A'\ne A}G(p,p')
                    +12\pi\sum_BG(p,q)-w_*(P)(p)/2,\label{7.48}\\
 S_B(q)&=12\pi R(q,q)+8\pi\sum_AG(q,p)
                    +12\pi\!\sum_{B'\ne B}G(q,q')
                    -w_*(P)(q)/2 .
 \label{7.49}
\end{align}
In the normalized moving charts, \eqref{7.47} and $s^2F_1(p)=4$
at an ordinary peak give
\begin{equation}
 s_A^2=\frac{e^{c_1}}4\,h_1(p)e^{S_A(p)}
                  (1+O(D^{1/4}\ell^N)).
 \label{7.50}
\end{equation}
At a moving mixed peak $q=q_{k,\nu}$, let
$\beta_k=n_*^g(P_\nu)(q_{k,\nu})>0$, where $n_*^g(P)$ is the
density on the branch \eqref{7.1}. By \eqref{7.43}--\eqref{7.46}
and $\phi(q)=0$,
\[
 F_2(q)e^{-T_k(q)}
 =c_2^D H_{P,s}(q)e^{w_D(q)+E_D(q)}
       \prod_{l\ne k}e^{T_l(q)}
 =\beta_k(1+O(D^{1/4}\ell^N)).
\]
Evaluating $T_k$ at $z=0$ and using $a_k^2F_2(q)=4$ gives
\begin{equation}
 a_k^2=\frac{\beta_k\delta_k^4(1+|B_k|^2)}4
              (1+O(D^{1/4}\ell^N)).
 \label{7.51}
\end{equation}
Evaluate the first density at $z=\delta_k$, where $m_{B_k}(1)$ is
uniformly positive for bounded $B_k$. Since \eqref{7.51} gives
$a_k^2/\delta_k^2=O(\delta_k^2)$, \eqref{5.60} and the first component
of \eqref{7.47} yield
\begin{equation}
 \delta_k^4=\frac{e^{c_1}}{16}\,h_1(q)e^{S_B(q)}
                    (1+O(D^{1/4}\ell^N)).
 \label{7.52}
\end{equation}
The positive coefficients in \eqref{7.50}--\eqref{7.52} are
uniformly bounded above and away from zero. Thus
\begin{equation}
 D\asymp e^{c_1/4},\qquad
 s_A\asymp a_B\asymp D^2,\qquad\delta_B\asymp D .
 \label{7.53}
\end{equation}

Let $\mathcal{U}^{\rm app}$ be the approximation \eqref{7.10}
in the same moving charts, with the same centers, the modulated
parameters $B_k=B_{k,\nu}$, normalization $e^{c_1}$, and scalar
branch $w_*(P)$. Use the scales given by the leading terms of
\eqref{7.50}--\eqref{7.52}. The ordinary and mixed local profiles
give Green completions differing from \eqref{7.47} by $O(D^c)$ for
some fixed $c>0$. At a mixed center, write
$d_\delta=\delta^4+|z^2-\delta^2B|^2$ and let $F_\delta,Q_\delta$ be
the Gram polynomials in \eqref{7.8}. Then
$$
 \frac{F_\delta}{d_\delta}=1+O(\delta^2),\qquad
 \frac{Q_\delta}{4\delta^4}
   =(r^2+a^2)(1+\beta r^2/4)+O(\delta^2r^2).
$$
These are uniform on the entire fixed disk, including zero,
because $a^2=\beta\delta^4(1+|B|^2)/4$.
Substitute these expressions into
$n_1=4Q_\delta/F_\delta^2$ and
$n_2=16\beta\delta^8F_\delta/Q_\delta^2$, then subtract the
exact contributions of the limiting second density. This gives the two
logarithmic Green expressions in \eqref{7.47} with uniform error $O(\delta^2)$.
The cutoff and mean-subtraction errors are $O(\delta^2|\log\delta|)$.
At an ordinary center, use the Green identity preceding \eqref{7.44}
with the smooth background. The derivatives with respect to the
logarithmic scales are uniformly bounded, so the changes in
\eqref{7.50}--\eqref{7.52} preserve the estimate. Reducing the
exponent to absorb the logarithmic factors gives
\begin{equation}
 \|u-\mathcal{U}^{\rm app}\|_\infty\le CD^{1/16}.
 \label{7.54}
\end{equation}

\subsubsection{The projected linear estimate}

\begin{proposition}
\label{pr7.6}
Fix a shadow neighborhood where the centers are separated, the scalar
linearization at fixed centers is invertible, and the $B_k$ are bounded.
For the approximation \eqref{7.10}, uniformly for small $t$ and
$|\mu|\leq C_0t^{1/8}$,
\begin{equation}\label{7.55}
 \|\mathcal L_{t,\mu}^{-1}\|_{H^{-1}\to H^1}\leq C\ell_t,
 \qquad \operatorname{ind}\mathcal L_{t,\mu}=m_*+a+4b-1,
 \qquad \ell_t=\log(1/t).
\end{equation}
The Gram matrices of the source representatives are uniformly invertible.
The multipliers satisfy the same $C\ell_t$ bound. These estimates
hold also after the parameter adjustment in \eqref{7.50}--\eqref{7.54}.
\end{proposition}
\begin{proof}
On the compact parameter set, the centers are uniformly separated,
the scale factors are positive, and the scalar inverse is bounded.
At $\mu=0$, Proposition~\ref{pr7.1} with $p=6/5$ gives the
$O(t^{2/3})$ comparison with the model. Corollary~\ref{co6.7}
and \eqref{6.17} give the inverse and source estimates for
\eqref{7.10}, and Corollary~\ref{co6.5} gives the multiplier
estimate. The same $L^{6/5}$ comparison transfers
\eqref{eq:completed-weight-capacity} from the model to the
completed first weight by H\"older's inequality.
The approximate profile \eqref{7.10} is independent of $\mu$.
Since $t^{1/8}\ell_t^2\to0$, Lemma~\ref{le6.8} now gives the
same spectral, inverse and source bounds for $|\mu|\leq C_0t^{1/8}$.
Applying Corollary~\ref{co6.5} again gives the multiplier bound.
For the actual parameters, \eqref{7.53} gives
$t\asymp D^2$ and \eqref{7.42} gives
$|\mu|=O(D^{1/4})=O(t^{1/8})$.
A small parameter adjustment changes $t$ by a factor $1+o(1)$
and keeps the centers and parameters $B_k$ in the same compact set.
\end{proof}

For $|\mu|=O(D^{1/4})$, the residual of \eqref{7.10} satisfies
\begin{equation}
 \|R^{\rm app}\|_{H^{-1}}
       \le C(D^{4/3}+|\mu|\sqrt{\ell}).
 \label{7.56}
\end{equation}
At $\mu=0$, use Proposition~\ref{pr7.1} with $t\asymp D^2$.
In the fixed center neighborhood, the first derivative term at an
ordinary center contributes $\|rm_A^s\|_{6/5}=O(D^{4/3})$. In the
chosen mixed chart, $\psi^B_{j,z}(q_j)=0$.
The first weight at $\mu=0$ satisfies
\eqref{eq:completed-weight-capacity} by the $L^{6/5}$ comparison
in the preceding proof. Since \eqref{7.10} is independent of
$\mu$, the residual part of Lemma~\ref{le6.8} gives the term
$|\mu|\sqrt\ell$.

By the $L^\infty$ estimate in Lemma~\ref{lem:parameter-modulation},
a change of $O(D^{1/16})$ in the scaled parameters enforces the
orthogonality conditions. The estimate \eqref{7.54} still holds.

Let $\eta=u-\mathcal{U}^{\rm app}$ after modulation, and
$E=\|\eta\|_{H^1}$. For normalized exponential densities,
the exact covariance identity gives
$$
 \|N(\eta)\|_{H^{-1}}
      \le C\|\eta\|_\infty\,\ell E .
$$
Differentiate the normalized density along $t\eta$. The ratio of
its weight to the original weight is $1+O(\|\eta\|_\infty)$.
Cauchy--Schwarz for the covariance and the weighted embedding
with constant $C\sqrt{\ell}$ give the stated $H^{-1}$ bound.
By \eqref{7.54}, applying the inverse of norm $C\ell$ gives an
error $CD^{1/16}\ell^2E$ that we can absorb. Thus
\begin{equation}
 E\le CD^{4/3}\ell+C|\mu|\ell^{3/2}.
 \label{7.57}
\end{equation}

\subsubsection{The mass estimate from Green's identities}

We sum the local identities for the first-component mass.
At each mixed center, testing against the dilation kernel function
gives \eqref{5.10} in Section~\ref{sec:inv-chapter}. Its left-hand
side is twice the difference between the actual and model first masses
on the testing disk.

At an ordinary point, put
$$
 U_1^A=\log(4s^2)-2\log(s^2+r^2),\qquad
 U_2^A=\log\beta_A+\log(s^2+r^2),\quad
 z_s=\frac{s^2}{s^2+r^2}.
$$
The ordinary profile satisfies $-\Delta U^A=A(m_A^s,0)^T$.
Set
\begin{equation}
 E_A=s\partial_sU^A=(2-4z_s,2z_s),\quad
 E_{A,\infty}=(2,0),\quad
 K_A=C(E_A-E_{A,\infty})=(-2z_s,0).
 \label{7.58}
\end{equation}
For the correction $\zeta$, test the equation against $CE_A$ and
subtract the pairing with $CE_{A,\infty}$. Green's identity gives
\begin{align}
 2\int_{D_{r_0}}m_A^s(e^{\zeta_1}-1)\dif A
   ={}&\int_{D_{r_0}}K_A\cdot f^\rho\dif A
      +\int_{D_{r_0}}E_{A,1}m_A^s
                    (e^{\zeta_1}-1-\zeta_1)\dif A \notag\\
     &+\int_{D_{r_0}}E_{A,2}F_2\dif A
      -\int_{\partial D_{r_0}}
       (\partial_{\boldsymbol n} K_A\cdot\zeta
        -K_A\cdot\partial_{\boldsymbol n}\zeta)\dif s.
 \label{7.59}
\end{align}
The second-component comparison at the ordinary center and
\eqref{7.54} give $F_2\le C(s^2+r^2)$. Hence
$$
 \begin{aligned}
 0\leq\int_{D_{r_0}}E_{A,2}F_2\dif A
 &\leq Cs^2\int_0^{r_0}r\dif r\leq Cs^2,\\
 \int_{D_{r_0}}|K_A\cdot f^\rho|\dif A
 &\leq C\int_0^{r_0}\frac{s^2r}{s^2+r^2}\dif r
 \leq Cs^2|\log s|.
 \end{aligned}
$$
On a fixed testing circle, $K_A,\nabla K_A=O(s^2)$.
Since $J_{A,1}=O(r)$,
\begin{equation}
 \int_{D_{r_0}}m_A^s|J_{A,1}|^2\dif A
 \leq C\int_0^{r_0}\frac{s^2r^3}{(s^2+r^2)^2}\dif r
 \leq Cs^2|\log s|.
 \label{7.60}
\end{equation}
At a mixed center, use the dilation test $E_B=\epsilon\partial_\epsilon U$
from \eqref{5.9}, with $\epsilon=\delta_B^2$ and $B,\beta$
held fixed. Write $a=a_B$, so $a\asymp\epsilon$, and let
$q_\epsilon=Q/\epsilon^2$. Here $E_B$ denotes a vector test and
$E$ without the subscript is the norm of the correction. On the
fixed testing disk, \eqref{7.8} gives
\[
 \begin{gathered}
 F\asymp(r^2+\epsilon)^2,\qquad q_\epsilon\asymp r^2+\epsilon^2,\qquad
 |\epsilon\partial_\epsilon F|\le C\epsilon(r^2+\epsilon),\\
 |\epsilon\partial_\epsilon q_\epsilon|\le
                  C(\epsilon r^2+\epsilon^2),\qquad
 n_1\le \frac{C\epsilon^2(r^2+a^2)}{(r^2+\epsilon)^4},
 \qquad n_2\le C\left(1+\frac{a^2}{(r^2+a^2)^2}\right).
 \end{gathered}
\]
The constants are uniform for bounded $B$ and for $\beta$ in
a compact subset of $(0,\infty)$. The second component of $E_B$
tends to zero off the center. Quantitatively,
\[
 \begin{aligned}
 |E_{B,1}|&=
 |2+\epsilon\partial_\epsilon\log q_\epsilon
                       -2\epsilon\partial_\epsilon\log F|\le C,\\
 |E_{B,2}|&=
 |\epsilon\partial_\epsilon\log F
                       -2\epsilon\partial_\epsilon\log q_\epsilon|
 \le C\left(\frac{\epsilon}{r^2+\epsilon}+\epsilon+
                  \frac{\epsilon^2}{r^2+\epsilon^2}\right)
 \le \frac{C\epsilon}{r^2+\epsilon}.
 \end{aligned}
\]
The final inequality holds on the fixed disk for $0<\epsilon<1$.
Use $J_1=O(r^2)$ and $J_2=O(r)$ from \eqref{7.15} and consider the
first component. The total mass of $n_1$ is bounded, and the bound
$C\delta_B^4r^{-6}$ on $r\ge\delta_B$ gives
\[
 \int_{D_{r_0}} |E_{B,1}|n_1J_1^2\dif A
 \le C\left(\delta_B^4+\delta_B^4\int_{\delta_B}^{r_0}\frac{\dif r}{r}\right)
 \le C\delta_B^4(1+|\log\delta_B|).
\]
For the second component, estimate the inner peak and the diffuse
density separately.
\[
 \int_{D_{r_0}} |E_{B,2}|n_2J_2^2\dif A
 \le C a^2\int_0^{r_0}\frac{r^3}{(r^2+a^2)^2}\dif r
       +C\epsilon\int_0^{r_0}\frac{r^3}{r^2+\epsilon}\dif r
 \le C\bigl(a^2(1+|\log a|)+\epsilon r_0^2\bigr).
\]
It follows that the moment used in the nonlinear integral in
\eqref{5.10} satisfies
\begin{equation}\tag{7.60a}\label{eq:mixed-dilation-weighted-moment}
 \sum_{i=1}^2\int_{D_{r_0}} |E_{B,i}|n_iJ_i^2\dif A
          \le C\delta_B^2 .
\end{equation}
In the nonlinear term of \eqref{5.10} the correction is
$\zeta_i=J_i+\eta_i+c_i^{\rm nor}$, where
$c_i^{\rm nor}=\log(\rho_i/M_i)-\log\int_M\nu_i e^{\eta_i}\dV$
and $\nu_i=V_i^{\rm app}/\rho_i$. On the testing disk
$V_i^{\rm app}\dV=(\rho_i/M_i)n_i e^{J_i}\dif A$ by
\eqref{7.15} and the normalization. The prenormalized profile is
independent of $\mu$. Proposition~\ref{pr7.1} at $\mu=0$ gives
$M_i=\rho_i^0+O(t)$ with $\rho_i^0>0$. Thus
$|\log(\rho_i/M_i)|\le C(|\mu|+t)$ in the entire rough window,
and $n_i\dif A\le C V_i^{\rm app}\dV$ there.

By \eqref{7.54}, $\|\eta\|_\infty\le CD^{1/16}=o(1)$ after
modulation. The model capacity and the $L^{6/5}$ comparison
used in Proposition~\ref{pr7.6} give the uncentered quadratic bound
for both completed weights. Indeed, the two weight integrals
against $\xi^2$ differ by at most $C t^{2/3}\|\xi\|_{12}^2$
at $\mu=0$, and changing $\mu$ multiplies the first weight by
$\rho_1/\rho_1^0$. In particular,
\[
 \sum_i\int_M V_i^{\rm app}|\eta_i|^2 e^{|\eta_i|}\dV
 \le e^{\|\eta\|_\infty}
                \sum_i\int_M V_i^{\rm app}|\eta_i|^2\dV\le C\ell E^2 .
\]
The same capacity bound gives
\[
 \left|\int_M\nu_i\eta_i\dV\right|
 +\left|\log\int_M\nu_i e^{\eta_i}\dV\right|
 \le C\left(\int_M\nu_i|\eta_i|^2\dV\right)^{1/2}
 \le C\sqrt\ell E .
\]
Here we used $|e^{\eta_i}-1|\le C|\eta_i|$ and the fact that
$\int_M\nu_i e^{\eta_i}\dV=1+o(1)$. Consequently
$|c_i^{\rm nor}|\le C(|\mu|+t+\sqrt\ell E)$. The uniform bounds for
$c_i^{\rm nor}$ also follow directly from
$|\log\int\nu_i e^{\eta_i}\dV|\le\|\eta_i\|_\infty$.

On the fixed disk, $\zeta_i$ is bounded. For the Taylor remainder
$e^{\zeta_i}-1-\zeta_i$, use the integral form of Taylor's formula
to separate the squares and products of $J_i,\eta_i,c_i^{\rm nor}$.
Apply Cauchy--Schwarz with respect to the positive measure
$|E_{B,i}|n_i\dif A$. The $J_i^2$ moment satisfies
\eqref{eq:mixed-dilation-weighted-moment}. The $\eta_i^2$ moment is
bounded by $C\ell E^2$ since $|E_{B,i}|\le C$ and
$n_i\dif A\le C V_i^{\rm app}\dV$. Thus the cross term involving
$J_i\eta_i$ is $O(\delta_B\sqrt\ell E)$.
Young's inequality bounds the terms involving $c_i^{\rm nor}$ by
the same square terms and $C\sum_i|c_i^{\rm nor}|^2$. The latter is
$O(\mu^2+\ell E^2+D^4)$ since $t\asymp D^2$.
At an ordinary center, use the measure $m_A^s\dif A$
and \eqref{7.60} for the first component. The second-density contribution
is already controlled in the test in \eqref{7.59}. In particular,
\begin{equation}\tag{7.60b}\label{eq:mixed-dilation-cross-terms}
 |\mathcal E_{A,\rm cross}|\le Cs_A|\log s_A| E,\qquad
 |\mathcal E_{B,\rm cross}|\le C\delta_B\sqrt\ell E,\qquad
 |\mathcal E_{\rm norm}|\le C(\mu^2+\ell E^2+D^4).
\end{equation}
The bound $|\log s_A|\asymp\ell$ used in the ordinary cross term
follows from \eqref{7.53}.
By \eqref{7.50}--\eqref{7.54} and elliptic estimates
on the fixed matching annuli,
$$
 \begin{gathered}
 \sum_{\rm cores}2\int_{B_{r_0}}(F_1-n_1)\dif A=2\mu+O(D^4),\qquad
 \|\zeta\|_{C^1(\partial B_{r_0})}\leq C,\\
 |\mathcal E_{\partial,A}|\leq Cs_A^2,\qquad
 |\mathcal E_{\partial,B}|\leq C\delta_B^2,\qquad
 s_A\asymp D^2,\quad\delta_B\asymp D.
 \end{gathered}
$$
In the mixed identity \eqref{5.10}, pair the smooth forcing with $K$
from \eqref{5.9}, with $\eps=\delta_B^2$. Pair the constant
first-component forcing $f_1^\rho(0)$ with the signed integral of
$K_1$. Differentiating \eqref{5.11} gives
\[
 \int_{B_{r_0}}K_1\dif A
 =-2\eps\mathcal M(B)+O_{r_0}(\eps^2|\log\eps|).
\]
The bounds on the coefficients and on $B$ are uniform. The
explicit denominators in \eqref{4.1}--\eqref{4.2} give
\[
 |K_1|\leq \frac{C\eps}{\eps+r^2},\qquad
 |K_2|\leq C\left(\eps+\frac{\eps^2}{r^2+\eps^2}\right).
\]
For the first inequality, use $F\geq c(\eps+r^2)^2$.
For the second, use $K_2=-\eps\partial_\eps\log(Q/\eps^2)$
and $Q/\eps^2\geq c(r^2+\eps^2)$. The difference
$f_1^\rho(z)-f_1^\rho(0)=O(r)$ therefore has a pairing with
$K_1$ of order $O_{r_0}(\eps)$. The full second-component pairing
is also $O_{r_0}(\eps)$. This gives
$\int_{B_{r_0}}K\cdot f^\rho\dif A=O_{r_0}(\delta_B^2)$
uniformly even in the present coarse parameter range.
Summing \eqref{7.59} and the mixed identities using
\eqref{eq:mixed-dilation-weighted-moment} and
\eqref{eq:mixed-dilation-cross-terms} yields
$$
 \begin{aligned}
 |\mu|&\leq C\sum_A(s_A^2|\log s_A|+s_A|\log s_A|E)\\
 &\quad+C\sum_B(\delta_B^2+\delta_B\sqrt\ell E)
           +C(\ell E^2+\mu^2+D^4).
 \end{aligned}
$$
Since $b\geq1$, $D^4\ell=o(D^2)$ and $D^2\ell=o(D\sqrt\ell)$,
\begin{equation}\label{7.61}
 |\mu|\leq CD^2+CD\sqrt\ell E+C\ell E^2+C\mu^2.
\end{equation}
Equations \eqref{7.42}, \eqref{7.57} and \eqref{7.61} give
$$
 |\mu|=O(D^{1/4}),\qquad
 E\leq C(D^{4/3}\ell+|\mu|\ell^{3/2}).
$$
Lemma~\ref{le5.18}, with
$(r,M_r,E_r,\alpha)=(D,|\mu|,E,1/4)$, implies
\begin{equation}\label{7.62}
 |\mu|=O(D^2),\qquad \|\eta\|_{H^1}=O(D^{4/3}\ell),
\end{equation}
uniformly on the fixed separated neighborhood with $|B_j|\leq C$.

\begin{theorem}
\label{th7.7}
Let a sequence of Toda solutions satisfy \eqref{2.4}--\eqref{2.6}
with $b\geq1$ and fixed $\rho_2\notin4\pi\mathbb N$. Suppose that
\eqref{7.35} holds and the residual scalar branch is nondegenerate.

In a fixed neighborhood with separated centers, the sequence satisfies
\eqref{7.50}--\eqref{7.53} and \eqref{7.62}.

After modulation, the sequence belongs to the projected family of
Proposition~\ref{pr7.2} for small $D$.
\end{theorem}
\begin{proof}
The profile classification, Green comparison and modulation give
\eqref{7.50}--\eqref{7.54}. Proposition~\ref{pr7.6}
and Lemma~\ref{le5.18} yield
$$
 |\mu|\leq CD^2,\qquad
 \|\eta\|_{H^1}\leq CD^{4/3}\ell<C_1t^{2/3}|\log t|
$$
for sufficiently large fixed $C_1$. Take $C_0$ in
Proposition~\ref{pr7.2} larger than the bound for $|\mu|/t$.
Then $\eta$ solves the projected equation in the contraction ball,
so uniqueness gives
$$
 \eta=\eta^{\rm constructed},\qquad
 F_\rho(u)=0\ \Longrightarrow\ c_\alpha=0\quad\text{for every }\alpha.
$$
\end{proof}

\begin{proof}[Proof of Corollary~\ref{co7.8}]
By Theorem~\ref{th7.7}, the sequence satisfies the reduced equations,
with parameters converging to the prescribed centers and $B_{j,*}$.
Divide \eqref{7.25} by $t$ to obtain $\mu/t\to S_*>0$, since
$\kappa_j\mathcal H(B_{j,*})>0$ for each $j$.

For $\mu>0$, take the center and angular domains of
Theorem~\ref{th7.4} with each $B_{j,*}$ in the interior.
On the closure of these domains, use the bounds $S_-\leq\mathcal S\leq S_+$
from that proof and its scale interval $(\tau_-,\tau_+)$ for $t/\mu$,
with $\tau_-<S_+^{-1}\leq S_-^{-1}<\tau_+$. The resulting inequalities
hold on the whole part of the parameter boundary where the scale is
at an endpoint. Take $C_0$ large enough for this interval and
for the actual sequence. Choose the correction radius
$Kt^{2/3}\ell$ with $K$ larger than the uniform constants for
both the constructed correction and the actual correction.
The projected contraction holds on a still larger fixed
multiple of $t^{2/3}\ell$. The reduced map has no zeros on the
parameter boundary, and uniqueness in the larger ball excludes
zeros on the correction boundary. Lemma~\ref{lem:parameter-modulation}
gives the open neighborhood with degree \eqref{7.64}. The center,
angular, scale and correction estimates put the actual sequence
strictly inside it. For finitely many shadow solutions, take
common bounds in these estimates.
\end{proof}

\begin{lemma}\label{lem:parameter-modulation}
Let $p=(P,Q,B,\log t)$ and fix a base parameter $p_0$. In the
coordinates $\xi$ centered at $p_0$, measure the displacements
of $P$ and $Q$ in units of their respective inner widths.

If $\|u-u^{\rm app}(p_0)\|_\infty=o(1)$, the equations
\[
 \mathcal C_\alpha(p;u)
   =\mathcal L_{\alpha,p}(u^{\rm app}(p)-u)=0
\]
have a unique solution near $p_0$, with
$|\xi|\leq C\|u-u^{\rm app}(p_0)\|_\infty$.
The corresponding difference $u-u^{\rm app}(p)$ has the same order
in $L^\infty$.

For fixed $\mu>0$, take the scale domain $t/\mu\in[\tau_-,\tau_+]$
and compact center and angular domains as above. On the bundle
of corrections $\eta\in E_t$ with
$\|\eta\|_{H^1}\leq Kt^{2/3}\ell$, the map
$(p,\eta)\mapsto u^{\rm app}(p)+\eta$ is one-to-one for small $\mu$.

In the interior, it gives local coordinates, and the matrix of
pairings between source functionals and parameter derivatives has
positive determinant.
The same assertions hold for the family with $b=0$ in
Section~\ref{rev:sec:ordinary}, on domains with comparable scales
and corrections bounded by $Kt^{7/6}\ell$.
\end{lemma}

\begin{proof}
In the coordinates $\xi$, the first and second derivatives of
$u^{\rm app}$ are bounded in $L^\infty$. On a mixed disk, the
denominators of $\log F_j,\log Q_j$ in \eqref{7.8} give the bounds
for derivatives taken in the directions $a_j\partial_{q_j}$,
$\partial_{B_j}$, and $\partial_{\log t}$. On ordinary disks, use
$s_i^2+r^2\geq s_i^2$. The cutoffs and smooth background fields have
bounded derivatives on the fixed annuli.

The same differentiation of the explicit tests and normalized
densities gives
\[
 \|D_\xi\zeta_\alpha\|_\infty
 +\|D_\xi\log V_i^{\rm app}\|_\infty\leq C.
\]
Consequently, for $q_\alpha=
 V_{k(\alpha)}^{\rm app}\widehat\zeta_\alpha$,
\[
 |D_\xi q_\alpha|\leq C V_{k(\alpha)}^{\rm app},\qquad
 \|D_\xi\mathcal L_{\alpha,p}\|_{(L^\infty)^*}\leq C,qquad
 \|D_\xi\mathcal L_{\alpha,p}\|_{(H^1)^*}\leq C\sqrt\ell .
\]
The last estimate is the weighted embedding of
Section~\ref{rev:sec:multiple}. The normalization derivatives
are bounded since $\partial_\xi\log\int h_ie^{u_i^{\rm app}}\dV$
is a weighted average of $\partial_\xi u_i^{\rm app}$.

For a fixed $\mu$, rescale each first-component constraint by the
positive factor $\rho_1^0/\rho_1$. In fact
$V_1^{\rm app}(\mu)=(\rho_1/\rho_1^0)V_1^{\rm app}(0)$ and each
of its weighted averages in the centering formula is unchanged.
The rescaled first-component row therefore equals its value at
$\mu=0$. Second-component rows are unchanged. These row operations
preserve the zero set, uniform local inverse estimates and the
determinant sign. We use this rescaling in the Jacobian below.
For the scaled translations, the leading pairings are those with
$\Psi^A,\Psi$ in Theorem~\ref{th7.3}. The angular derivatives give
$\Phi$, and the common scale derivative gives $E-E_\infty$ on
the inner disks. By centering, the pairing of this scale derivative
with the common source agrees to leading order with the pairing
for $E$. The matrix in real coordinates therefore tends to
$J_0$ in \eqref{7.34}.
For the family with $b=0$, the same calculation gives
$\operatorname{diag}(\frac{4\pi}{3}I_{2a},\frac{16\pi a}{3})$.

For fixed $u$, differentiating $\mathcal C$ yields
\[
 D_\xi\mathcal C
 =\mathcal L_pD_\xi u^{\rm app}
    +(D_\xi\mathcal L_p)(u^{\rm app}(p)-u).
\]
The first term is $J_0+O(|\xi|)+o(1)$ with $J_0$ evaluated
at $p_0$. For the second, use
$\|u^{\rm app}(p)-u^{\rm app}(p_0)\|_\infty=O(|\xi|)$
and the dual bounds above. Thus
\[
 \|D_\xi\mathcal C-J_0\|
 \leq C|\xi|+C\|u-u^{\rm app}(p_0)\|_\infty+o(1).
\]
Since $|\mathcal C(p_0;u)|\leq C\|u-u^{\rm app}(p_0)\|_\infty$,
the quantitative implicit function theorem gives the first assertion.
The same derivative estimate holds with
$C\sqrt\ell\|u-u^{\rm app}(p_0)\|_{H^1}$ in place of the
$L^\infty$ term. It gives local coordinates throughout the
indicated bundle.

To prove injectivity, suppose $u$ has two representations with
parameters $p,p'$. By H\"older's inequality, the exponential Sobolev
inequality, and Lemma~\ref{le6.9}, each normalized density of $u$
differs from its corresponding $V_i^{\rm app}$ by at most
$Ct^{-1/3}\|\eta\|_{H^1}=o(1)$ in $L^1$.

The scales of the two representations are comparable. At an ordinary
center, compare the first densities. At a mixed center, compare the
second densities. If the center displacement divided by the inner
width were unbounded, a fixed large rescaled disk would carry positive
limiting mass for only one representation. This would contradict the
$L^1$ comparison. The displacements are therefore bounded in
inner units. Passing to the limit gives equal translated and dilated
Liouville densities. Their centers and widths agree, so the scaled
displacements tend to zero and the width ratios tend to one.

At the outer scale of a mixed center, the limits of the first
densities satisfy $m_B(y)=\lambda^{-2}m_{B'}(y/\lambda)$,
where $\lambda$ is a limit of outer width ratios. Since
$|y|^6m_B(y)\to16$, it follows that $\lambda=1$.
Comparing the denominators then gives $B=B'$. The scale formulas
\eqref{7.5} and \eqref{7.7} give $t'/t\to1$.
Thus $|\xi|=o(1)$ between any two representations.

Both parameters solve the same constraint equation. Integrating its
derivative along the segment between the parameters gives
$(J_0+o(1))\xi=0$. Thus the parameters agree. Permutations within each
type preserve orientation because they exchange even-dimensional blocks.

For fixed $\mu\ne0$ with a nonempty local domain, the closed parameter
domain is compact. If $u_n=u^{\rm app}(p_n)+\eta_n\to u$ in $H^1$
with $(p_n,\eta_n)$ in the open bundle, pass to a subsequence with
$p_n\to p$. Then $\eta_n\to u-u^{\rm app}(p)$ in $H^1$, and the
limit satisfies the constraints and the closed correction bound.
If $u$ is on the boundary of the image, then $p$ is on the parameter
boundary or the correction attains its radius. Indeed, the local
coordinate assertion gives an open image for the interior.
\end{proof}

\section{Blow-up at ordinary points}\label{rev:sec:ordinary}

We use the wall and scale setting of Section~\ref{sec:scale-conventions}
with $b=0$.

Let $(P^0,w^0)$ be a nondegenerate shadow solution with $m$
ordinary points and an invertible scalar linearization at fixed
centers. We derive the following mass expansion, with $L_0$
defined in \eqref{8.23}.
$$
 \mu=\rho_1-4\pi m
     =L_0t^2\log(1/t)+O(t^2),\qquad
 L_0=\frac\pi2l(P^0,w^0).
$$
For $L_0\ne0$, write $d_T$ for the local degree on the contributing
side and use the two-sided degrees $d_T^\pm$ from
Section~\ref{sec:two-local-degrees}. The mass expansion and the
degree computation give
$$
 \mu L_0>0,\qquad
 d_T=(-1)^m\operatorname{sgn}(\mu)d_S,
 \qquad d_T^+-d_T^-=(-1)^m d_S.
$$

\subsection{The projected equation}
\label{ordinary-sec:shrinking}

We consider type $(m,0)$ with $m\geq1$ and set
$$
 \rho_1^0=4\pi m,\qquad \rho_1=4\pi m+\mu,\qquad
 \rho_2\notin4\pi\mathbb N,\qquad
 A=\begin{pmatrix}2&-1\\-1&2\end{pmatrix},\qquad C=A^{-1}.
$$
Write $E_t$ for the projected space defined below and
$\mathcal L_t$ for the corresponding projected linearization. The orthogonality conditions in Theorem~\ref{th6.1} with $a=m$ and
$b=0$ remove $2m$ translation modes and one common dilation mode.
With the scalar index $m_*$ of Section~\ref{sec:degree-conventions},
this gives
$$
 \operatorname{codim}E_t=2m+1,\qquad
 \operatorname{ind}\mathcal L_t=m_*+m-1.
$$

Use the scalar branch $w(P)$ and the density $n_*^g$ in
\eqref{7.1}, with $\sigma=\rho_2$ and $w(P^0)=w^0$. Write $m_*$ for the Morse index of the scalar linearization and
$F(P)$ for the center field obtained by substituting $w(P)$ into
\eqref{2.8}. Then $F(P^0)=0$ and $DF(P^0)$ is invertible.
With the orientation of Section~\ref{sec:degree-conventions},
\begin{equation}\label{8.2}
 d_S=(-1)^{m_*}\operatorname{ind}_{P^0}F .
\end{equation}
All configurations considered below remain separated.

\subsubsection{Approximate solutions and the residual estimate}

Use the outer pair in \eqref{7.2} and the normalization in
\eqref{7.4} with $Q=\varnothing$ and $\sigma=\rho_2$.
Use the ordinary charts and scales in \eqref{7.5}, and write
\begin{equation}\label{8.4}
 m_i=m_{A,i}.
\end{equation}
By \eqref{2.8}, $\nabla\psi_i(p_i)=F_i(P)$, and $s_i\asymp t$.
Use $D_i$ from \eqref{7.9} and write $\chi_i=\chi_i^A$.
For $b=0$, the ansatz \eqref{eq:ansatz} reduces to
\begin{equation}\label{8.5}
 \begin{aligned}
 v^{\rm app}
 &=\binom{u_1^0}{u_2^0}+\sum_i\chi_i\binom{D_i}{-D_i/2},\\
 u^{\rm app}
 &=v^{\rm app}-\int_Mv^{\rm app}\dV.
 \end{aligned}
\end{equation}
The logarithmic singularities cancel in \eqref{8.5}.
The constant $c_1$ determines all the scales $s_i$.

For fixed constants $C_0,C_1$, restrict the parameters by
\begin{equation}\label{8.6}
 |P-P^0|\leq C_0t^{1/2},\qquad
 |\mu|\leq C_1t^2\ell .
\end{equation}
Then $|F(P)|=O(t^{1/2})$. On the disk where $\chi_i=1$, the
densities $N_j^g=e^{c_j}h_je^{v_j^{\rm app}}$ before normalization satisfy
\begin{equation}\label{8.7}
 \begin{split}
 N_1^{\rm loc}&=m_ie^{J_{1,i}},\qquad
 J_{1,i}=\psi_i-\psi_i(p_i)=F_i(P)\cdot z+O(|z|^2),\\
 N_2^{\rm loc}&=\gamma_i(z)(s_i^2+|z|^2),\qquad
 \gamma_i(z)=n_*^{\rm loc}(z)/|z|^2>0 .
 \end{split}
\end{equation}
The quotient defining $\gamma_i$ extends smoothly across the origin.
The smooth remainders and their derivatives are uniform in the range \eqref{8.6}.

\begin{lemma}
\label{le8.1}
Under \eqref{8.6}, set
$M_j=\int_MN_j^g\dV$ and
$V_j^{\rm app}=(\rho_j/M_j)N_j^g$. Then
\begin{equation}\label{8.8}
 M_1=4\pi m+O(t^2\ell),\qquad
 M_2=\rho_2+O(t^2),\qquad
 \log(\rho_j/M_j)=O(t^2\ell).
\end{equation}
The residual in the normalized equation satisfies
\begin{equation}\label{8.9}
 \left\|-\Delta_gu^{\rm app}
       -A(V^{\rm app}-\rho)\right\|_{H^{-1}}
       \leq Ct^{7/6}.
\end{equation}
The associated normalized variance forms and projected spaces
differ from the corresponding objects in the model with $m$ ordinary
cores by $o(\ell^{-1})$ in the energy norm.
\end{lemma}

\begin{proof}
The local rank-one pair
$U_i=(\log m_i,\log[\gamma_i(0)(s_i^2+r^2)])$
satisfies $-\Delta U_i=A(m_i,0)^T$. At $(\rho_1,\rho_2)=(4\pi m,\rho_2)$,
the prenormalized residual in Euclidean coordinates is
\begin{equation}\label{8.10}
 -A\binom{m_i(e^{J_{1,i}}-1)}{s_i^2\gamma_i}.
\end{equation}
The estimates
$$
 \|m_i\|_{6/5}=O(t^{-1/3}),\qquad
 \|rm_i\|_{6/5}=O(t^{2/3}),\qquad
 \|r^2m_i\|_{6/5}=O(t^{5/3})
$$
together with $|F(P)|=O(t^{1/2})$ bound the residual by $O(t^{7/6})$,
including the $O(t^2)$ cutoff errors. The linear term in $J_{1,i}$
integrates to zero, and $\int_{B_{r_0}} r^2m_i\dif A=O(t^2\ell)$.
This gives the mass estimates in \eqref{8.8}. The correction from
normalization is $O(t^{5/3}\ell)$, and the background correction is
$O(|\mu|)=O(t^2\ell)$.
For the model densities, use $m_i^g=e^{-\phi}m_i$
with the fixed cutoffs and the mass corrections in
Section~\ref{rev:sec:multiple}. The same moments give
\[
 \left\|V_1^{\rm app}-\sum_i m_i^g\right\|_{6/5}
 \leq C(t^{2/3}|F(P)|+t^{5/3}\ell),
 \qquad
 \|V_2^{\rm app}-n_*^g\|_{6/5}\leq Ct^2.
\]
Corollary~\ref{co6.7} now compares the variance forms and the projected
spaces, using test functions centered by their global means.
\end{proof}

In the first component, use the translation and dilation test functions
$$
 \zeta_{i,l}=\chi_i\frac{s_iz_l}{s_i^2+|z|^2},
 \quad l=1,2,\qquad
 \zeta_0=\sum_i\chi_i\left(2-\frac{4s_i^2}{s_i^2+|z|^2}\right).
$$
Subtract the mean with respect to $V_1^{\rm app}\dV$ and set
$$
 \widehat\zeta_\alpha
   =\zeta_\alpha-\rho_1^{-1}\int_MV_1^{\rm app}\zeta_\alpha\dV,
 \qquad
 q_\alpha=V_1^{\rm app}\widehat\zeta_\alpha,\qquad
 \mathcal L_\alpha(\xi)=\int_Mq_\alpha\xi_1\dV .
$$
Let $E_t$ be the common kernel of the $2m+1$ functionals
$\mathcal L_\alpha$ in $\mathcal H$, with energy
$\int_M\nabla\xi^TC\nabla\xi\dV$. The common radial test function
is $-2$ times the normalized radial profile of Section~\ref{rev:sec:multiple}.

\begin{theorem}
\label{th8.2}
For small $t>0$, uniformly under \eqref{8.6}, there is a unique
$\eta\in E_t\cap B_{Ct^{7/6}\ell}^{H^1}$ satisfying
\begin{equation}\label{8.11}
 \begin{aligned}
 -\Delta_g(u^{\rm app}+\eta)
 &=A\begin{pmatrix}
 \displaystyle\rho_1\frac{h_1e^{u_1^{\rm app}+\eta_1}}
      {\int_Mh_1e^{u_1^{\rm app}+\eta_1}\dV}-\rho_1\\[1.2ex]
 \displaystyle\rho_2\frac{h_2e^{u_2^{\rm app}+\eta_2}}
      {\int_Mh_2e^{u_2^{\rm app}+\eta_2}\dV}-\rho_2
 \end{pmatrix}\\
 &\quad+Ae_1\sum_\alpha c_\alpha q_\alpha.
 \end{aligned}
\end{equation}
The resulting family is continuous in the parameters and smooth
for each fixed $t>0$. Moreover,
\begin{equation}\label{8.12}
 E:=\|\eta\|_{H^1}\leq Ct^{7/6}\ell,\qquad
 |c|\leq Ct^{7/6}\ell^2 .
\end{equation}
The projected linearization at the corrected solution has Morse index $m_*+m-1$.
\end{theorem}

\begin{proof}
Apply the contraction argument of Proposition~\ref{pr7.2} with
$b=0$ and the residual estimate \eqref{8.9} of Lemma~\ref{le8.1}.
On the ball of radius $Kt^{7/6}\ell$, apply the nonlinear estimate with
factor $t^{-1/6}$ and then the inverse. This gives a Lipschitz constant
$O(t\ell^2)=o(1)$.
Corollary~\ref{co6.5} yields \eqref{8.12}. The derivative of the
nonlinear remainder has norm $O(t\ell)=o(\ell^{-1})$, so the
projected Morse index remains $m_*+m-1$.
\end{proof}

\subsubsection{The projected center equations}

Let $f^\rho$ and $f^{\rho^0}$ be the canonical forcings in \eqref{2.2}
at the actual parameter $\rho$ and at $\rho^0=(4\pi m,\rho_2)^T$,
respectively. In each coordinate disk, set
$$
 \theta_j=\log(\rho_j/M_j)+\eta_j
                -\log\int_M(V_j^{\rm app}/\rho_j)e^{\eta_j}\dV .
$$
The canonical fields of the corrected profile are $U+J+\theta$.
The term $\theta$ includes the normalization constants. We have
\begin{equation}\label{8.13}
 \begin{split}
 \sum_i\int_{B_{r_0}(p_i)}
       (m_i|\theta_1|^2e^{|\theta_1|}
             +(r^2+s_i^2)|\theta_2|^2e^{|\theta_2|})\dif A
       &\leq C\ell E^2+Ct^4\ell^2,\\
 \|\theta\|_{C^1(\partial B_{r_0}(p_i))}
       &\leq C_{r_0}(\sqrt\ell E+t^2\ell).
 \end{split}
\end{equation}
For the first bound, the weights in the integral have bounded mass
and are bounded above by $Ct^{-2}$. Their $L^{1+1/\ell}$ norms are
uniformly bounded. If $q=1+\ell$, the bounds
$\|\eta_i\|_{4q}\leq C\sqrt q E$ and
$\|e^{|\eta_i|}\|_{2q}\leq Ce^{CqE^2}$ follow from the two-dimensional
Sobolev and Moser--Trudinger inequalities. H\"older's inequality
gives a weighted quadratic remainder of order
$\ell E^2$. Applying the same argument to
$|e^{\eta_i}-1|\leq|\eta_i|e^{|\eta_i|}$ bounds the logarithmic
normalizations by $C\sqrt\ell E$. Equation \eqref{8.8} supplies
the mass errors $O(t^2\ell)$. These give the first bound of
\eqref{8.13}. On a fixed annulus about $\partial B_{r_0}$,
the first density $V_1^{\rm app}$ is $O(t^2)$ and the second is bounded.
For fixed $p>2$, the exponential difference in \eqref{8.11}
has $L^p$ norm at most $C_p(\sqrt\ell E+t^2\ell)$ by
H\"older's inequality and the same exponential estimate at fixed exponents.
The remaining terms on the annulus are $O(t^2\ell+t^2|c|)$.
Interior $W^{2,p}$ estimates and \eqref{8.12} give the
$C^1$ bound on the testing circle. We keep $r_0>0$ fixed
throughout these estimates.

For two vector fields, define the boundary form
$$
 \mathcal B(V,\xi)=
 \int_{\partial B_{r_0}}
       (\xi\cdot\partial_n(CV)-(CV)\cdot\partial_n\xi)\dif s .
$$
The translation kernel functions of the linearized rank-one equation are
$$
 \Psi_l=-s_i\partial_{x_l}U,\quad
 \Psi_{l,1}=\frac{4s_ix_l}{s_i^2+r^2},\quad
 \Psi_{l,2}=-\frac{2s_ix_l}{s_i^2+r^2},\quad
 C\Psi_l=\left(\frac{2s_ix_l}{s_i^2+r^2},0\right).
$$
The smooth field $J$ satisfies
\begin{equation}\label{8.14}
 -\Delta_0J=A(0,n_*^{\rm loc})^T+f^{\rho^0} .
\end{equation}
Let $L=-\Delta_0-A\operatorname{diag}(m_i,0)$ and
$R_1=e^{J_1+\theta_1}-1-J_1-\theta_1$. Subtract the equations
for $U$ and $J$ from the equation for $U+J+\theta$. With
$s_{\rm proj}\dif A=\sum_\alpha c_\alpha q_\alpha\dV$ on this disk,
\[
 L\theta
 =A\binom{m_iJ_1+m_iR_1}{n_{2,\mathrm{actual}}^{\rm loc}-n_*^{\rm loc}}
       +(f^\rho-f^{\rho^0})+Ae_1s_{\rm proj}.
\]
For this operator, $C L$ has a symmetric potential and
$L\Psi_l=0$. Integration by parts in
$\int (C\Psi_l)\cdot L\theta$ gives $\mathcal B(\Psi_l,\theta)$.
Solving for the source pairing gives
\begin{equation}\label{8.15}
 \begin{split}
 \sum_\alpha c_\alpha\int_{B_{r_0}}\Psi_{l,1}q_\alpha\dV
  ={}&-\int_{B_{r_0}} m_i\Psi_{l,1}J_1\dif A+\mathcal B(\Psi_l,\theta)\\
     &-\int_{B_{r_0}} m_i\Psi_{l,1}
                   (e^{J_1+\theta_1}-1-J_1-\theta_1)\dif A\\
     &-\int_{B_{r_0}}\Psi_{l,2}(n_{2,\mathrm{actual}}^{\rm loc}-n_*^{\rm loc})\dif A
       -\int_{B_{r_0}}(C\Psi_l)\cdot(f^\rho-f^{\rho^0})\dif A.
 \end{split}
\end{equation}

The subtraction cancels $\mathcal B(\Psi_l,J)$ and the contribution
of $f^{\rho^0}$.

Expand to second order. The quadratic contribution is odd, so
\begin{equation}\label{8.16}
 \int_{B_{r_0}} m_i\Psi_{l,1}J_1\dif A=4\pi s_iF_{i,l}(P)+O(t^3\ell).
\end{equation}
The cubic remainder has integral $O(s_i^3|\log s_i|)$. For the
nonlinear term we have
$$
 \int_{B_{r_0}} m_i|\Psi_{l,1}|\,|J_1|^2\dif A
       \leq C(s_i^2|F_i|^2+s_i^3)=O(t^3).
$$
The terms involving $\theta$ are controlled by \eqref{8.13} and
Cauchy--Schwarz. The integrand involving the second density in
\eqref{8.15} is
$$
 \Psi_{l,2}(n_{2,\mathrm{actual}}^{\rm loc}-n_*^{\rm loc})
 =-\frac{2s_i^3\gamma_i(x)x_l}{s_i^2+r^2}
                -2s_i\gamma_i(x)x_l(e^{\theta_2}-1).
$$
The constant part of $\gamma_i$ gives an odd integrand in the
first summand, whose integral is $O(t^3)$. The second summand and
$\mathcal B(\Psi_l,\theta)$ each contribute $O(t\sqrt\ell E+t^3\ell)$.
By \eqref{2.2}, the term with $f^\rho-f^{\rho^0}=-A(\mu,0)^Te^\phi$
is $O(t|\mu|)$.
Thus the right-hand sides of the projected center equations are
\begin{equation}\label{8.17}
 -4\pi s_iF_i(P)+
 O(t^3\ell+t\sqrt\ell E+\ell E^2)
       =-4\pi s_iF_i(P)+o(t^{3/2}).
\end{equation}

\subsection{Mass expansion and local degree}\label{sec:ordinary-mass}

\subsubsection{The reduced mass equation}

Since $s_i=td_i(P)$, the change of scale at fixed $P$ is
$t\partial_t=s_i\partial_{s_i}$ in the $i$th disk. The associated
kernel function and its limit at infinity are
$$
 \mathcal E=s_i\partial_{s_i}U
   =\left(2-\frac{4s_i^2}{s_i^2+r^2},
                     \frac{2s_i^2}{s_i^2+r^2}\right),\qquad
 \mathcal E_\infty=(2,0),\qquad
 K=C(\mathcal E-\mathcal E_\infty)
     =\left(-\frac{2s_i^2}{s_i^2+r^2},0\right).
$$
Let $\Delta M_{1,i}=\int_{B_{r_0}}(n_{1,\mathrm{actual}}^{\rm loc}-m_i)\dif A$
and put $X=J+\theta$. The equation is
$LX=A(m_iR_1,n_{2,\mathrm{actual}}^{\rm loc})^T
       +f^\rho+Ae_1s_{\rm proj}$.
Since $L\mathcal E=0$ and $\mathcal E_\infty=(2,0)$, we have
$L(\mathcal E-\mathcal E_\infty)=A(2m_i,0)^T$. The Green
identity thus gives $\mathcal B(\mathcal E-\mathcal E_\infty,X)
+2\int m_iX_1\dif A$ for the integral of $K\cdot LX$.
In solving for the source term use
$\Delta M_{1,i}=\int m_i(X_1+R_1)\dif A$, so that
$2\int m_i X_1\dif A-\int m_i(\mathcal E_1-2)R_1\dif A
=2\Delta M_{1,i}-\int m_i\mathcal E_1R_1\dif A$.
This gives
\begin{equation}\label{8.18}
 \begin{split}
 \sum_\alpha c_\alpha
        \int_{B_{r_0}}(\mathcal E_1-2)q_\alpha\dV
 ={}&2\Delta M_{1,i}
       -\int_{B_{r_0}} K\cdot f^\rho\dif A
       +\mathcal B(\mathcal E-\mathcal E_\infty,J+\theta)\\
    &-\int_{B_{r_0}} m_i\mathcal E_1
                 (e^{J_1+\theta_1}-1-J_1-\theta_1)\dif A
       -\int_{B_{r_0}}\mathcal E_2n_{2,\mathrm{actual}}^{\rm loc}\dif A.
 \end{split}
\end{equation}
The last integral in \eqref{8.18} comes from the second density,
which is absent from the rank-one equation for $\mathcal E$.

Use $\kappa_i(P)$ from \eqref{8.19}. Since $e^{\phi(p_i)}=1$ and the linear part of $f^\rho$ integrates
to zero against $K$, radial integration gives
\begin{equation}\label{8.20}
 \int_{B_{r_0}} K\cdot f^\rho\dif A
     =4\pi\kappa_i(P)s_i^2\log(1/s_i)
             +O(t^2)+O(|\mu|t^2\ell).
\end{equation}
The contribution of $J$ to the boundary integral is $O(t^2)$ because
its linear term has zero circular average. The contribution
of $\theta$ is $O(t^2\sqrt\ell E+t^4\ell)$.
Moreover
$$
 \int_{B_{r_0}} m_i|J_1|^2\dif A
     \leq C(s_i^2+s_i^2|F_i|^2\ell)=O(t^2),
 \qquad
 \int_{B_{r_0}}\mathcal E_2n_{2,\mathrm{actual}}^{\rm loc}\dif A=O(t^2).
$$
Thus the following estimate controls both the nonlinear remainder and
the integral against $n_{2,\mathrm{actual}}^{\rm loc}$ in \eqref{8.18}.
$$
 O(t^2+t\sqrt\ell E+\ell E^2)=o(t^2\ell).
$$
The first-component mass outside the disks and the integrals of
$m_i$ over the complements of their respective disks are $O(t^2)$.
Hence $\sum_i\Delta M_{1,i}=\mu+O(t^2)$. Summing \eqref{8.18} gives
the common mass equation with leading mass difference $2\mu$.

\begin{theorem}
\label{th8.3}
The multipliers in \eqref{8.11} satisfy, uniformly under
\eqref{8.6},
\begin{align}
 c_i&=-3s_iF_i(P)+o(t^{3/2}),
                       \label{8.21}\\
 c_0&=\frac3{8\pi m}
       \left(\mu-2\pi\sum_{i=1}^m
                 \kappa_i(P)s_i^2\log(1/s_i)\right)
          +o(t^2\ell).
                       \label{8.22}
\end{align}
\end{theorem}

\begin{proof}
Pair the translation sources with $\Psi_l$ and sum the pairings of
the common radial source with $\mathcal E_1-2$ over the disks.
The limiting source matrix is diagonal, with center blocks
$(4\pi/3)I_2$ and radial entry
$$
 \Gamma_0=\frac{16\pi m}{3}.
$$
On a core, $\Psi_{l,1}=4\zeta_{i,l}$. The mode
$Z_{i,l}=s_i z_l/(s_i^2+|z|^2)$ before cutoff satisfies
$\int_{\C}m_i Z_{i,l}^2\dif A=\pi/3$. Replacing $Z_{i,l}$ by
$\chi_i Z_{i,l}$ changes this integral by $O(s_i^4)$.
Indeed, the difference is supported at fixed positive radius,
where $m_i Z_{i,l}^2\leq C s_i^4 r^{-6}$. The radial integrals are
$\int_{\C}m_i\mathcal E_1^2\dif A=16\pi/3$ and
$\int_{\C}m_i\mathcal E_1\dif A=0$. The matrix also includes the
source means and cutoff tails.

The matrix error is $O(t^{3/2}+t^2\ell)$. The linear term is bounded
by $|F_i|\int_{B_{r_0}} r m_i\dif A=O(t^{3/2})$. The quadratic term
is bounded by $\int_{B_{r_0}} r^2m_i\dif A=O(t^2\ell)$.
The normalization terms, including the source means, and the cutoff
tails satisfy the same bound.

Apply the uniformly bounded inverse to the right-hand sides of
\eqref{8.15} and \eqref{8.18}. This gives $|c|=O(t^{3/2})$, so the
matrix error contributes $O(t^3\ell)=o(t^2\ell)$ even in the radial
equation.
With the limiting matrix, \eqref{8.17} gives \eqref{8.21}, and
summing \eqref{8.18} and using \eqref{8.20} gives \eqref{8.22}.
\end{proof}

\subsubsection{Nonvanishing of the mass coefficient and local degree}

Use the logarithmic mass coefficient \eqref{8.23} from
Section~\ref{sec:scale-conventions}.
At $P=P^0$, the leading balance in \eqref{8.22} is
$\mu=L_0t^2\ell$.

\begin{proof}[Proof of Lemma~\ref{le8.4}]
Keep $h_2$ fixed and let $x_1,\ldots,x_J$ be the distinct centers
of the shadow solutions. Choose cutoffs $\chi_j=1$ near $x_j$
in disjoint geodesic disks. The perturbations
$$
 b_c(x)=\frac14\sum_{j=1}^J c_j\chi_j(x)d_g(x,x_j)^2,
 \qquad h_{1,c}=h_1e^{b_c},
$$
satisfy
$$
 b_c(x_j)=0,\qquad \nabla b_c(x_j)=0,\qquad
 \Delta_gb_c(x_j)=c_j.
$$
Thus each original shadow solution still solves the equations,
and its scalar linearization at fixed centers is unchanged.
At $x_j$, the Hessian of $b_c$ is $(c_j/2)g$. At any mixed center,
the perturbation preserves the first and pure complex second jets
in the isothermal chart. Thus $\mathcal T$ in \eqref{3.41} is unchanged,
while $\kappa_m$ changes by $c_j$. Sufficiently small coefficients
preserve all strict inequalities.
At a solution with ordinary configuration $P$, the coefficient is
$$
 l_c(P)=l_0(P)+\sum_{j=1}^J a_{Pj}c_j,
 \qquad
 a_{Pj}=\sum_{i:\,p_i=x_j}e^{\psi_i(p_i)}\geq0,
 \qquad \sum_j a_{Pj}>0.
$$
The choices with $l_c(P)=0$ lie on finitely many proper affine hyperplanes
$$
 \mathcal A_P:\quad l_0(P)+a_P\cdot c=0.
$$
For every $\epsilon>0$ choose
$$
 c\in B_\epsilon(0)\setminus\bigcup_P\mathcal A_P.
$$
By the compactness argument of Theorem~\ref{th3.14}\textup{(i)},
the solution sets are uniformly compact for nearby weights. A
collision would give a lower-type limit satisfying \eqref{3.41},
in contradiction with \eqref{3.7}. The implicit function theorem
then excludes additional solutions. Taking $\epsilon$ small
preserves nondegeneracy and \eqref{3.41}. If the shadow system of
type $(m,0)$ has no solutions, take the original weights.

\end{proof}

\begin{proof}[Proof of Theorem~\ref{th8.5}]
\emph{Step 1. Degree of the reduced equations.}
For $\mu L_0>0$, let $t_\mu>0$ be the unique sufficiently
small solution of
$|\mu|=|L_0|t_\mu^2\log(1/t_\mu)$.  In local center coordinates
write
$$
 P=P^0+t_\mu^{1/2}\xi,\qquad t=t_\mu\tau,\qquad
 |\xi|<R,\quad \tfrac12<\tau<2 ,
$$
with $R>0$ fixed and small. This domain satisfies \eqref{8.6}
for suitable $C_0,C_1$. Divide each center block by $-3s_it_\mu^{1/2}$
and multiply the last equation by $8\pi m/(3|\mu|)>0$. The $2m$
negative row factors preserve orientation. The resulting map in
\eqref{8.21}--\eqref{8.22} converges uniformly on the closed domain to
\begin{equation}\label{8.26}
        (\xi,\tau)\longmapsto
        \bigl(DF(P^0)\xi,\
              \operatorname{sgn}(\mu)(1-\tau^2)\bigr).
\end{equation}
The limiting map has a unique nondegenerate zero at $(0,1)$.
Its first block is bounded away from zero on $|\xi|=R$, and its
last component is nonzero on $\tau\in\{1/2,2\}$. By uniform convergence,
the reduced map has no zeros on the boundary. Its Brouwer degree is
$$
              -\operatorname{sgn}(\mu)\operatorname{ind}_{P^0}F .
$$
Equation~\eqref{8.22} excludes zeros on the opposite side and gives
$t/t_\mu\to1$ at every zero in the shrinking neighborhood with $t\to0$.
Equation~\eqref{8.21} gives $P-P^0=o(t^{1/2})$. Hence all such
zeros eventually lie in the chosen parameter domain.

\emph{Step 2. Leray--Schauder degree.}
The parametrization has positive orientation by \eqref{7.34} with
$b=0$. Theorem~\ref{th8.2} gives the contraction estimate required
in \eqref{7.30} and the projected Morse index $m_*+m-1$. Hence
the local degree is
$$
 (-1)^{m_*+m-1}
       [-\operatorname{sgn}(\mu)\operatorname{ind}F]
       =(-1)^m\operatorname{sgn}(\mu)d_S.
$$
Since zeros occur only for $\operatorname{sgn}\mu=\operatorname{sgn}L_0$,
the positive-side degree minus the negative-side degree is $(-1)^m d_S$.

\end{proof}

\subsection{Localization of ordinary blow-up sequences}\label{sec:ordinaryactual-chapter}

\subsubsection{Local profile and pointwise estimates}

With the normalization of Section~\ref{sec:scale-conventions}, write
\begin{equation}
 \begin{gathered}
 A=\begin{pmatrix}2&-1\\-1&2\end{pmatrix},\qquad C=A^{-1},\qquad
 f_i=\rho_i\frac{h_ie^{u_i}}{\int_Mh_ie^{u_i}\dV},\\
 v_i=Gf_i,\qquad u=(2v_1-v_2,\,2v_2-v_1).
 \end{gathered}
 \label{8.27}
\end{equation}
Suppose $\rho_1\to4\pi m$ with $\rho_2\notin4\pi\mathbb N$ fixed,
and that the sequence concentrates with local masses $(4\pi,0)$
at $m$ distinct points $P^0=(p_1^0,\ldots,p_m^0)$. By Lemma~\ref{le2.1},
\begin{equation}
 \begin{gathered}
 v_1\longrightarrow4\pi\sum_jG(\cdot,p_j^0),\qquad
 W:=2v_2\longrightarrow w^0,\\
 -\Delta_gw^0=
 2\rho_2\left(\frac{H_{P^0}e^{w^0}}{\int_MH_{P^0}e^{w^0}\dV}-1\right),
 \qquad H_P=h_2e^{-4\pi\sum_jG(\cdot,p_j)}.
 \end{gathered}
 \label{8.28}
\end{equation}
The translation Pohozaev identity in Lemma~\ref{le2.1} gives the
center equations in \eqref{2.8}. Assume that the resulting shadow
solution is nondegenerate and that its scalar linearization at
fixed centers is invertible. On a fixed separated neighborhood of $P^0$, use the scalar branch
$w(P)$ in \eqref{7.1} and let $F(P)$ be the center field in
\eqref{2.8} evaluated at $w=w(P)$. Then $F(P^0)=0$, $DF(P^0)$ is invertible, and $F(P)=o(1)$ as $P\to P^0$.

Use the same chart normalization as in the construction. In
an initial isothermal chart $z$, with metric $e^{\phi(z)}|dz|^2$,
use the change of coordinates
$$
 f_p(z)=a(p)(z-p)+b(p)(z-p)^2,\qquad
 a(p)=e^{\phi(p)/2},\qquad b(p)=\tfrac12a(p)\phi_z(p).
$$
The coordinate changes are uniformly invertible on a sufficiently
small fixed disk. The new metric factor satisfies $\phi(p)=0$ and
$\nabla\phi(p)=0$. In each fixed disk, choose $p_j$ at a maximum of
$f_1$. For large indices, these points lie in the interior. For the
coordinate density $F_1=e^\phi f_1$, it follows that
$\nabla F_1(p_j)=0$. Set $s_j^2F_1(p_j)=4$ and use radial
cutoffs in these charts.

\begin{lemma}
\label{le8.6}
Consider a sequence as in Lemma~\ref{le2.1} and an ordinary limit
point $p^0$ in a fixed disk separated from the other limiting
centers. Choose $p$ at a maximum of $f_1$ in that disk. In the
normalized chart above, set $s^2F_1(p)=4$ and $r=|z-p|$. Then
$$
 s^2F_1(p+sy)\longrightarrow \frac4{(1+|y|^2)^2}
$$
with smooth local convergence, including after taking logarithms.
The profile accounts for the full local mass $4\pi$. For every
sufficiently small fixed $\varepsilon>0$, on a smaller fixed disk we have
\begin{equation}
 F_1\le C_\varepsilon s^{-2}(1+r/s)^{-4+\varepsilon},\qquad
 F_2\le C_\varepsilon(s+r)^{2-\varepsilon},\qquad
 |\nabla v_1|\le C/(s+r).
 \label{8.29}
\end{equation}
On this smaller disk, $F_2$ is uniformly Lipschitz and converges
in $C^{0,\alpha}$ to the coordinate density of $n_*\dV$ for
$0<\alpha<1$. Moreover, with
$m^1(y)=4(1+|y|^2)^{-2}$,
\begin{equation}
 \left|\log\frac{s^2F_1(p+sy)}{m^1(y)}\right|
       \le Cs^{1/4}\log(2+|y|),\qquad |y|\le s^{-1/2}.
 \label{8.30}
\end{equation}
The estimates are uniform over a fixed finite collection of limiting
shadow solutions.
\end{lemma}
\begin{proof}
\textit{Step 1.} We obtain the profile and pointwise bounds.
Take a fixed disk $D$ about the ordinary point with boundary away
from all limiting centers. By Lemma~\ref{le2.1}, $U_2=\log F_2$ is
bounded above on $\partial D$. By \eqref{2.2} and positivity of the
Dirichlet Green function, $U_2\leq C+2G_DF_2$ on $D$. Fix $q>2$.
The second concentration mass at $p^0$ is zero by \eqref{2.4}, so
$D$ can be chosen with $\limsup_\nu q\int_DF_2\dif A<2\pi$.
Jensen's inequality and the logarithmic Green bound give
$$
 \int_D e^{2qG_DF_2}\dif A\leq C,\qquad \|F_2\|_{L^q}\leq C.
$$
Apply H\"older's inequality to the Green representation to bound
$U_2$ above. On smaller disks, the potentials of $f_i\dV$ outside
$D$ have uniformly bounded derivatives because the kernels are smooth
there and the masses are bounded. The bound for $F_2$ and local
elliptic estimates give a uniform $C^{1,\alpha}$ bound for $v_2$
on smaller disks.

Near the first-component peak, $-\Delta U_1=2e^{U_1}+O(1)$.
The rescaled function $U_1(sy)+2\log s$ has value $\log4$ at
zero. The peak of $f_1$ and the bounded metric factor give an
upper bound on each fixed rescaled disk. The rescaled right-hand
side is bounded there. On a slightly larger fixed disk, subtract the
zero-boundary Poisson potential of this right-hand side. Harnack's
inequality for a constant minus the harmonic remainder gives a lower
bound on the smaller disk from the value at zero. Interior estimates and
a diagonal subsequence give the entire limit. A change of variables on
fixed rescaled disks bounds its finite mass by the local mass.
Rescaling by $s$ gives a regular finite-mass scalar limit, with
maximum density $4$ at zero. By \eqref{eq:scalar-liouville},
its density is $4(1+|y|^2)^{-2}$, of mass $4\pi$.
If the bound $r^2F_1\le C$ failed, the failing radii would tend to
zero by the off-center bounds, and their ratios to $s$ would tend
to infinity along a subsequence by the local scaled convergence.
Apply the selection lemma \eqref{eq:prescribed-center-selection}
in the original chart, retaining the ordinary center. The
bounded forcing, total mass bound and pointwise failure are its
hypotheses. Its new ball of radius $o(r)$ avoids every fixed core
$B_{Rs}$ for large $\nu$. In a fixed disk of radius $\eta$ about
the limiting point, the limiting first mass is $4\pi$ plus the
integral of any diffuse part, and the core mass tends to
$\int_{B_R}4(1+|y|^2)^{-2}\dif A$. Subtract the core mass and
take the sequence limit, then let $R\to\infty$, and finally let
$\eta\downarrow0$. The mass in the new ball from the first density
tends to zero. The second density is bounded throughout
the working disk by the preceding Green estimate, so its mass
in that ball tends to zero as well. This contradicts
\eqref{eq:prescribed-center-selection}. For bounded $r/s$, the
pointwise bound follows directly from the local scaled convergence.
Since the second density is locally bounded, split the Green gradient
integral at distance $r/4$ from a point on the circle. This gives a
$C/r$ bound for $\nabla U_1$ there, as in the proof of
Proposition~\ref{pr4.19}. Integrating around the circle then
bounds its oscillation.

Integration over $D_r$ gives
$$
 r\overline U_1'(r)
   =-\frac1{2\pi}\int_{D_r}(2F_1-F_2+f_1^\rho)\dif A.
$$
The convergence of the local mass bounds the right-hand side by
$-4+\varepsilon$ for $Rs\leq r\leq r_0$, with $R$ large and $r_0$
small. Integration and the oscillation bound give the first estimate
in \eqref{8.29}. Splitting the Green integral for $\nabla v_1$ at
a radius comparable to $s+r$ gives $|\nabla v_1|\leq C/(s+r)$.

The mass in $D_{Rs}$ also gives
$v_1(x)\geq-(2-\varepsilon)\log(s+|x-p|)-C_\varepsilon$.
In $f_2=e^{c_2}h_2e^{2v_2}e^{-v_1}$, the factor
$e^{c_2}h_2e^{2v_2}$ has bounded first derivatives. Thus
$$
 F_2\leq C_\varepsilon(s+r)^{2-\varepsilon},\qquad
 |\nabla F_2|\leq C_\varepsilon(s+r)^{1-\varepsilon}.
$$
Taking $\varepsilon<1$ gives the uniform Lipschitz bound for $F_2$.
Its convergence away from $p^0$ then extends in $C^{0,\alpha}$
across $p^0$ for each $0<\alpha<1$.

\textit{Step 2.} We prove \eqref{8.30}.
Put $U_1=\log F_1$, $g_\rho=f_1^\rho-F_2$, and
$\alpha(r)=(4\pi)^{-1}\int_{D_r}F_1\dif A$.
On $\partial D_r$ with $r=\sqrt s$, Green's representation and
\eqref{8.29} give
\[
 \nabla U_1(z)=-4\alpha(r)\frac{z}{r^2}
     +O\left(1+\frac{s}{r^2}+s^{2-\varepsilon}r^{-3+\varepsilon}\right).
\]
On $D_{r/2}$, the first moment $\int_{D_{r_0}}|y|F_1\dif A=O(s)$
bounds the kernel expansion error. On $D_{2r}\setminus D_{r/2}$,
the density bound $Cs^{2-\varepsilon}r^{-4+\varepsilon}$ and
$\int_{D_{2r}}|z-y|^{-1}\dif A_y=O(r)$ give the last error term.
On $D_{r_0}\setminus D_{2r}$, the bound $O(|y|^{-1})$ on the kernel
gives the same order. Replacing the mass in $D_{r/2}$ by that in
$D_r$ also gives an error of this order.

The radial Pohozaev identity for $-\Delta U_1=2F_1+g_\rho$ is
\[
 r\int_{\partial D_r}
 \left[(\partial_n U_1)^2-\tfrac12|\nabla U_1|^2+2F_1\right]\dif s
 =4\int_{D_r}F_1\dif A-\int_{D_r}g_\rho\,z\cdot\nabla U_1\dif A.
\]
At $r=\sqrt s$, the gradient terms contribute
$16\pi\alpha(r)^2+O(s^{1/4})$, the boundary integral of $2rF_1$
is $O(s^{1-\varepsilon/2})$, and the integral against $g_\rho$
is $O(r^2)$.

It follows that $\alpha(r)-\alpha(r)^2=O(s^{1/4})$. Since the
rescaled profile carries the full local mass, $\alpha(\sqrt s)\to1$.
Thus the mass in $D_{\sqrt s}$ is $4\pi+O(s^{1/4})$. The gradient estimate gives
$\osc_{\partial D_{\sqrt s}}U_1=O(s^{1/4})$.

With $R=s^{-1/2}$, the rescaled error
$\eta(y)=\log(s^2F_1(p+sy)/m^1(y))$ satisfies
$$
 (-\Delta-2m^1)\eta
    =2m^1(e^\eta-1-\eta)+s^2(f_1^\rho-F_2)(sy).
$$
The bounded forcing has weighted norm
$O(s^2R^3)=O(s^{1/2})$ with weight $(1+|y|)^3$.
By the choice of peak and scale, $\eta(0)=0$ and $\nabla\eta(0)=0$.
We first show that $N=\|\eta\|_{\log,R}$ tends to zero. Given
$\epsilon>0$, choose a fixed $R_*$ large enough that the flux
identity gives $|r\overline U_1'(r)+4|\leq C\epsilon$ for
$R_*s\leq r\leq\sqrt s$. The corresponding slope of the model
differs from $-4$ by $O(R_*^{-2})$.

Integrate the slope difference from $R_*s$ and use the bounded
circular oscillation to obtain
\[
 \limsup_{s\to0}N
 \leq C\epsilon+C R_*^{-2}+C/\log(2+R_*).
\]
Local profile convergence makes the error at $R_*s$ and on the
rescaled disk $D_{R_*}$ tend to zero.

Applying Lemma~\ref{le5.12} with \eqref{5.52} and the boundary
oscillation bound gives $N\leq C(s^{1/4}+N^2)$. Since $N=o(1)$,
this proves \eqref{8.30}. Uniformity over a finite set of limiting
shadow solutions follows by subsequence compactness.
\end{proof}

For a sequence satisfying \eqref{8.28}, the local Lipschitz bounds
and convergence away from the ordinary points give global
$C^{0,\alpha}$ convergence of $f_2$ for $0<\alpha<1$. Elliptic
regularity then yields $v_2\to w^0/2$ in $C^{2,\alpha}(M)$.

\subsubsection{Green estimates and the common scale}

Let $S=\max_j s_j$ and $L_S=|\log S|$. Cut off the profiles
$$
 m_j(z)=\frac{4s_j^2}{(s_j^2+|z|^2)^2}
$$
in the fixed disjoint disks when applying $G$. The notation
$Gm_j$ denotes the Green potential of the measure after cutoff,
$m_j\dif A=m_j^g\,dV_g$, where $m_j^g=e^{-\phi}m_j$.
By \eqref{8.29} and \eqref{8.30},
\begin{equation}
 \left\|v_1-\sum_jGm_j\right\|_\infty
                  \le CS^{1/4}L_S^N
 \label{8.31}
\end{equation}
for a fixed integer $N$. On $D_{\sqrt{s_j}}$, \eqref{8.30} bounds
the density error by $Cs_j^{1/4}m_j\log(2+|z|/s_j)$.
The finite logarithmic moments of $m^1$ give a bound on the Green
potential that is uniform in the evaluation point. The tail outside
this disk has logarithmic potential $O(s_j^{1-\eps/2}|\log s_j|)$.
Outside the fixed disks $f_1=O(e^{c_1})$, where
$$
 c_1=\log\frac{\rho_1}{\int_Mh_1e^{u_1}\dV}.
$$
Convergence away from the centers and \eqref{8.29} on a fixed circle give
$$
 e^{c_1}\le C_\eps S^{2-\eps}.
$$
For small $x>0$, the function $x^{1/4}|\log x|^N$ is increasing.
Since $s_j\leq S$, summing these estimates gives \eqref{8.31}.

Set $H_{P,S}=h_2e^{-\sum_jGm_j}$. The expansion
$$
 Gm_j=-\log(s_j^2+r^2)+4\pi R(x,p_j)
                      +O(s_j^2|\log s_j|)
$$
gives the bound
$$
 \|H_{P,S}-H_P\|_\infty\le CS^2L_S .
$$
Near $p_j$, use $H_P=r^2k_j(x)$ and
$r^2k_j(x)(1+s_j^2/r^2)=(r^2+s_j^2)k_j(x)$.
The estimate therefore holds through the center.

Combining this bound with \eqref{8.31} and the normalized equation
for $f_2$ gives, for each fixed $p>2$,
$$
 \left\|f_2-\rho_2\frac{H_Pe^W}{\int_MH_Pe^W\dV}\right\|_{L^p}
       \le CS^{1/4}L_S^N .
$$
The normalization integrals are bounded away from zero.
As a map $W^{2,p}_0\to L^p_0$, the fixed-center differential
operator has index zero because it is a compact perturbation of the
Poisson isomorphism. Elliptic bootstrapping shows that the kernel
elements are smooth. Fixed-center nondegeneracy in $X_s$ makes the
kernel trivial. The inverse is uniformly bounded near the limiting
branch. Apply the inverse to the scalar residual, use $W\to w^0$,
and absorb the quadratic remainder to obtain
\begin{equation}
 \|W-w(P)\|_{W^{2,p}}\le CS^{1/4}L_S^N .
 \label{8.32}
\end{equation}

In the normalized chart at $p_j$, write $\psi_j(\cdot;P)$ for
the regular field in \eqref{7.5} with $Q=\varnothing$. Evaluating \eqref{8.31}--\eqref{8.32} at $p_j$ and using
$s_j^2F_1(p_j)=4$ gives
\begin{equation}
 s_j^2=\frac{e^{c_1}}4e^{\psi_j(p_j;P)}
                  (1+O(S^{1/4}L_S^N)).
 \label{8.34}
\end{equation}
The factors $e^{\psi_j(p_j;P)}$ are bounded above and away from
zero. Set
$$
 t=e^{c_1/2},\qquad L=\log(1/t),\qquad
 \widehat s_j=\tfrac12t\,e^{\psi_j(p_j;P)/2}.
$$
Then
\begin{equation}
 S\asymp t,\qquad s_j\asymp t,\qquad
 \log(s_j/\widehat s_j)=O(t^{1/4}L^N).
 \label{8.35}
\end{equation}
With this choice of $t$, use $u^0(P)$ from \eqref{7.2} with
$Q=\varnothing$ and define
\begin{equation}\label{8.36}
 u^{\rm app}_{t,P}
 =\left[u^0(P)+\sum_j\chi_j
       \binom{-2\log(1+\widehat s_j^2/|z_j|^2)}
             {\log(1+\widehat s_j^2/|z_j|^2)}\right]_0 .
\end{equation}
In the $v$ variables, the second component is $w(P)/2$.
For the first component, the logarithmic scale derivatives are
uniformly bounded, and the cutoff and mean-zero corrections are
$O(t^2L)$. Combining \eqref{8.31}, \eqref{8.32}, and \eqref{8.35}
and decreasing the exponent gives
\begin{equation}
 \quad \|u-u^{\rm app}_{t,P}\|_\infty\le Ct^{1/8},
                    \qquad P\longrightarrow P^0 .\quad
 \label{8.37}
\end{equation}

\begin{lemma}
\label{le8.7}
Suppose the centers vary in a fixed separated neighborhood where
the scalar linearization at fixed centers is invertible. For
$|\mu|\leq C_0t^{1/4}$, the projected linearization at
\eqref{8.36} is invertible, with uniform inverse and multiplier
bounds $CL$ and Morse index $m_*+m-1$.

With $f=\lvert F(P)\rvert$, the residual of \eqref{8.36} satisfies
\begin{equation}\label{8.38}
 \|R^{\rm app}\|_{H^{-1}}
 \leq C(t^{2/3}f+t^{5/3}L+\lvert\mu\rvert\sqrt L).
\end{equation}
\end{lemma}
\begin{proof}
At $\mu=0$, use \eqref{8.10} with $J_j=F_j(P)\cdot z+O(r^2)$
on the fixed neighborhood. The moments in Lemma~\ref{le8.1}
give a residual of order $O(t^{2/3}f+t^{5/3}L)$ and an $O(t^{2/3})$ difference
in $L^{6/5}$ between the model and approximate densities and sources.
Corollary~\ref{co6.7} then compares the two projected operators.
Theorem~\ref{th6.1} and Corollary~\ref{co6.5}, with $b=0$, give
the inverse and multiplier bounds and the Morse index at $\mu=0$.

The $O(t^{2/3})$ comparison in $L^{6/5}$ also transfers
\eqref{eq:completed-weight-capacity} from the ordinary model
to the completed first weight at $\mu=0$. Using the bounded
background in the model, the same argument bounds the uncentered
second quadratic form by $CL\|\xi\|_{H^1}^2$.
The profile \eqref{8.36} is independent of $\mu$, and
$t^{1/4}L^2=o(1)$. Lemma~\ref{le6.8} thus preserves the
spectral and source bounds and adds
$C\lvert\mu\rvert\sqrt L$ to the residual for
$\lvert\mu\rvert\leq C_0t^{1/4}$.
Corollary~\ref{co6.5} gives the corresponding multiplier bound.
\end{proof}

\subsubsection{Mass and center estimates}

Let $(P_0,t_0)$ denote the parameters in \eqref{8.37}. To impose
orthogonality, apply Lemma~\ref{lem:parameter-modulation} with
$b=0$. Then $(t_0^{-1}(P-P_0),\log(t/t_0))=O(t_0^{1/8})$.
The scaled derivatives of the approximate solution are bounded
in $L^\infty$, so \eqref{8.37} still holds, with
$t/t_0=1+O(t_0^{1/8})$.

After modulation, set $f=\lvert F(P)\rvert=o(1)$ and
$E=\|u-u^{\rm app}\|_{H^1}$. We use the inverse bound of
Lemma~\ref{le8.7} to derive the sharper mass estimate \eqref{8.39}
together with bounds on $f$ and $E$.

\begin{lemma}
\label{le8.8}
Every ordinary blow-up sequence with the above nondegenerate
shadow limit satisfies
\begin{equation}\label{8.39}
 |\mu|=O(t^2L),
\end{equation}
and
\begin{equation}\label{8.40}
 |F(P)|+|P-P^0|\leq Ct^{5/3}L^3,\qquad
 \|\eta\|_{H^1}\leq Ct^{5/3}L^2.
\end{equation}
\end{lemma}
\begin{proof}
\emph{Step 1. Energy estimate.}
Integrating \eqref{8.30} over $D_{\sqrt{s_j}}$ gives first-component
mass $4\pi+O(s_j^{1/4})$. By \eqref{8.29}, the rest of each fixed
disk contributes $O(s_j^{1-\varepsilon/2})$. The mass outside
these disks is $O(t_0^2)$ by the definition of $t_0$ and the
bounds away from the centers.

Summing these masses and using $s_j\asymp t_0\asymp t$ gives
$|\mu|=O(t^{1/4})$. By Lemma~\ref{le8.7} and
$\|N(\eta)\|_{H^{-1}}\leq C\|\eta\|_\infty LE$, inversion bounds
the nonlinear term by $Ct^{1/8}L^2E$. Absorbing this term gives
\begin{equation}\label{8.41}
 E\leq C(t^{2/3}Lf+t^{5/3}L^2+\lvert\mu\rvert L^{3/2}).
\end{equation}
In particular, $\sqrt L E=o(1)$.

\emph{Step 2. Estimate of the mass difference.}
First estimate the normalization and the fixed-circle terms in the
present coarse parameter range. The completed unnormalized masses
satisfy $M_1=4\pi m+O(t^2L)$ and $M_2=\rho_2+O(t^2)$ on the
fixed separated neighborhood. Indeed, the linear term of
$e^{J_1}-1$ integrates to zero against $m_i$. The second moment
is $O(s_i^2L)$, and the first-density tails and cutoff errors
are $O(t^2)$. The local second-density difference is $s_i^2\gamma_i$
by \eqref{8.7}, and the cutoff annuli contribute $O(t^2)$.
Write $\nu_i=V_i^{\rm app}/\rho_i$. The weighted quadratic
bounds used in Lemma~\ref{le8.7} and $\|\eta\|_\infty=o(1)$ give
\[
 \left|\log\int_M\nu_i e^{\eta_i}\dV\right|
 \leq C\int_M\nu_i|\eta_i|\dV\leq C\sqrt L E.
\]
The constants in $\theta-\eta$ therefore satisfy
\[
 |\theta_1-\eta_1|\leq C(\sqrt L E+|\mu|+t^2L),
 \qquad |\theta_2-\eta_2|\leq C(\sqrt L E+t^2).
\]
In particular $\|\theta\|_\infty=o(1)$ by \eqref{8.37},
the coarse bound on $\mu$ and \eqref{8.41}.

On a fixed annulus containing the testing circle,
subtract the equations for $U+J$ and $U+J+\theta$ to get
\[
 -\Delta_0\theta
 =A\binom{m_i(e^{J_1+\theta_1}-1)}
          {(r^2+s_i^2)\gamma_i(e^{\theta_2}-1)+s_i^2\gamma_i}
          -A\binom{\mu}{0}e^\phi .
\]
The coefficients on this annulus are bounded, with $m_i=O_{r_0}(t^2)$.
For fixed $p>2$, the $L^p$ norm of the right-hand side is at most
$C_{r_0,p}(t^2+|\mu|+\|\theta\|_{L^p})$.
Interior $W^{2,p}$ estimates on a smaller annulus and Sobolev
embedding give
\[
 \|\theta\|_{C^1(\partial B_{r_0}(p_i))}
 \leq C_{r_0}(\sqrt L E+|\mu|+t^2L).
\]
These bounds use only the coarse parameter range.

Sum \eqref{8.18} with all multipliers zero. Since
$K=(-2s_i^2/(s_i^2+r^2),0)$, the smooth forcing contributes
$O_{r_0}(s_i^2\log(1/s_i))$. The test
$\mathcal E-\mathcal E_\infty$ and its first derivatives are
$O_{r_0}(s_i^2)$ on the testing circle, so the boundary term
against $J+\theta$ is $O_{r_0}(s_i^2)$. The second-density
integrand is $2s_i^2\gamma_i e^{\theta_2}$, also giving
$O(s_i^2)$. The exponential remainder is bounded by
\[
 C\int_{B_{r_0}}m_i
           \bigl(J_1^2+\theta_1^2e^{|\theta_1|}\bigr)\dif A
 \leq C(t^2L+LE^2+\mu^2).
\]
The last inequality follows from the second moment, the uncentered
weighted bound for $\eta_1$, and the constant normalization error.
The actual first-density tail
is $O(t^2)$ by the fixed-scale approximation and its uniform
$L^\infty$ error. The model tail is also $O(t^2)$, so
$\sum_i\Delta M_{1,i}=\mu+O(t^2)$. Therefore
\[
 \lvert\mu\rvert\leq C(t^2L+LE^2+\mu^2).
\]
Substituting \eqref{8.41} and absorbing
$C\mu^2L^4=o(\lvert\mu\rvert)$ gives
\begin{equation}\label{8.42}
 \lvert\mu\rvert\leq C(t^2L+t^{4/3}L^3f^2),\qquad
 E\leq C(t^{2/3}Lf+t^{5/3}L^2).
\end{equation}
The second inequality follows by substituting the first into
\eqref{8.41}.

\emph{Step 3. Estimate of the center displacement.}
Use the translation identities \eqref{8.15} with all multipliers
zero. Taylor expansion of $J_1$ and the oddness of its quadratic
term after multiplication by $\Psi_{l,1}$ give the leading
pairing $4\pi s_i F_{i,l}(P)+O(t^3L)$ as in \eqref{8.16}.
The remainder is bounded using
\[
\begin{aligned}
 \int_{B_{r_0}}m_i|\Psi_{l,1}|J_1^2\dif A
       &\leq C(t^2f^2+t^3),\\
 \int_{B_{r_0}}m_i\theta_1^2e^{|\theta_1|}\dif A
       &\leq C(LE^2+\mu^2+t^4L^2).
\end{aligned}
\]
For the second density, we have the exact integrand
\[
 \Psi_{l,2}(n_{2,\mathrm{actual}}^{\rm loc}-n_*^{\rm loc})
 =-\frac{2s_i^3\gamma_i(x)x_l}{s_i^2+r^2}
   -2s_i\gamma_i(x)x_l(e^{\theta_2}-1).
\]
In the first summand, the contribution of the constant part of
$\gamma_i$ integrates to zero. The rest contributes $O(t^3)$ after
integration. The second summand has integral
$O(t\sqrt L E+t^3)$ by the $L^1$ bound and the uniform smallness
of $\theta_2$. On the fixed testing circle, $\Psi_l$ and its
first derivatives are $O_{r_0}(t)$. The preceding $C^1$
estimate gives a boundary contribution
$O(t\sqrt L E+t|\mu|+t^3L)$. Finally,
$f^\rho-f^{\rho^0}=-A(\mu,0)^T e^\phi$ and the $L^1$
norm of $C\Psi_l$ on the disk is $O(t)$, so the forcing
difference contributes $O(t|\mu|)$. Since $s_i\asymp t$,
these bounds imply
\[
 tf\leq C(t^3L+t^2f^2+t^2Lf+t\sqrt L E+LE^2+t|\mu|+\mu^2).
\]
Substitute \eqref{8.42} and divide by $t$. We can then absorb the
linear and quadratic terms in $f$. For the mass terms,
$$
\lvert\mu\rvert+\mu^2/t\leq C(t^2L+t^3L^2+
                      t^{4/3}L^3f^2+t^{5/3}L^6f^4)
             \leq Ct^{5/3}L^3+o(1)f.
$$
Absorbing the $o(1)f$ terms gives $f\leq Ct^{5/3}L^3$.
The invertibility of $DF(P^0)$ and \eqref{8.42} give
\eqref{8.40} and then \eqref{8.39}.
\end{proof}

\begin{theorem}
\label{th8.9}
Let a solution sequence with $m\geq1$ ordinary points satisfy
\eqref{8.28}. Assume that its separated shadow limit $(P^0,w^0)$
is nondegenerate and that the scalar linearization at fixed
centers is invertible.

Then \eqref{8.39} and \eqref{8.40} hold. After modulation, the
sequence belongs to the family of Theorem~\ref{th8.2}. Both the estimates
and membership in the family hold for $\mu=0$. If $l(P^0,w^0)\ne0$, the sequence cannot
concentrate at this shadow limit on the wall $\rho_1=4\pi m$.
\end{theorem}
\begin{proof}
After modulation the correction belongs to $E_t$.
Choose the constant $C_1$ in \eqref{8.6} strictly larger
than the uniform bound for $|\mu|/(t^2L)$ in \eqref{8.39}.
Moreover, Lemma~\ref{le8.8} gives
\[
 \frac{|P-P^0|}{t^{1/2}}=O(t^{7/6}L^3)=o(1),\qquad
 \frac{\|\eta\|_{H^1}}{t^{7/6}L}=O(t^{1/2}L)=o(1).
\]
For sufficiently small $t$, these estimates put the center and the
correction strictly inside the center domain and contraction ball of
Theorem~\ref{th8.2}, respectively. The actual solution satisfies the
projected equation, so uniqueness makes its correction equal to the
constructed correction. This argument also applies for $\mu=0$.

Since all multipliers vanish, the sequence satisfies the reduced
equations, also at $\mu=0$. If $L_0\ne0$, dividing the radial
equation by $t^2L$ at $\mu=0$ would give $0=-L_0+o(1)$.
\end{proof}

\begin{proposition}
\label{pr8.10}
Under the hypotheses of Theorem~\ref{th8.9},
\begin{equation}\label{8.43}
 \mu=L_0t^2\log(1/t)+O(t^2),\qquad L_0=\frac\pi2l(P^0,w^0).
\end{equation}
The estimate is uniform over a fixed finite set of shadow limits.
If $L_0\ne0$, concentration requires $\mu L_0>0$.
When $\mu L_0>0$ and $|\mu|$ is sufficiently small, let $t_\mu$
be the unique root in $(0,e^{-1/2})$ of
$\mu=L_0t_\mu^2\log(1/t_\mu)$. Then
\begin{equation}\label{8.44}
 \frac t{t_\mu}=1+O\!\left(\frac1{\log(1/t_\mu)}\right).
\end{equation}
\end{proposition}
\begin{proof}
Apply \eqref{8.18} with zero multipliers at fixed radius and use
\eqref{8.40}. The model tails and the first-component mass outside
the disks are $O(t^2)$. Equation~\eqref{8.20} gives the forcing
term. At the actual parameters, the bound for $\mu$ in \eqref{8.39}
gives the required size for the parameter correction
$O(|\mu|t^2L)$, including at $\mu=0$. Also, the
normalization bounds preceding \eqref{8.42}, together with
\eqref{8.39} and \eqref{8.40}, give
$|\theta_i-\eta_i|=O(\sqrt L E+|\mu|+t^2L)=o(t)$.
Their squared contribution is therefore $O(t^2)$.
The remaining terms are bounded by
$$
 O(t^2+t\sqrt L E+LE^2+t^2L|F(P)|^2)=O(t^2).
$$
For the boundary term, use $\mathcal E-\mathcal E_\infty$.
For the second density, use
$\mathcal E_2n_{2,\rm actual}=2s_i^2\gamma_i e^{\theta_2}$.
The estimates leading to \eqref{8.22}, with \eqref{8.40},
then give
$$
 \mu=2\pi\sum_i\kappa_i(P)s_i^2\log(1/s_i)+O(t^2)
     =L_0t^2L+O(t^2),
$$
since $s_i=td_i(P)$ and $L|P-P^0|=o(1)$.
Equation~\eqref{8.43} gives the sign. Since $t^2\log(1/t)$ is strictly
increasing on $(0,e^{-1/2})$ and tends to zero as $t\downarrow0$,
the asserted root exists and is unique. When $\mu L_0>0$,
apply the mean value theorem using
$t\partial_t\log(t^2\log(1/t))=2-1/L$ to obtain \eqref{8.44}.
All estimates are uniform on a fixed finite set of shadow limits.
\end{proof}

\begin{proof}[Proof of Corollary~\ref{co8.11}]
By Theorem~\ref{th8.9} and Proposition~\ref{pr8.10},
$\mu L_0>0$, $(P-P^0)/t_\mu^{1/2}\to0$, $t/t_\mu\to1$, and
$\|\eta\|_{H^1}=o(t^{7/6}L)$. Thus every concentrating sequence
eventually lies strictly inside the neighborhood of
Theorem~\ref{th8.5}. The reduced map has no zeros on the boundary of the parameter domain,
and uniqueness keeps the projected correction inside its ball.

Apply the $b=0$ assertion of Lemma~\ref{lem:parameter-modulation}
to this closed parameter domain. Here $\mu\ne0$ is fixed,
$t/t_\mu\in[1/2,2]$, and the correction radius is a fixed
multiple of $t^{7/6}L$. The lemma gives the open image and describes
its boundary. The correction radius may be
taken strictly between the bound for the constructed correction and
the radius of a larger contraction ball. Uniqueness then excludes
zeros on the correction boundary. The $b=0$ matrix in
\eqref{7.34} provides the orientation. The degree is given by
\eqref{8.24}, which yields \eqref{8.45}. Since there are finitely
many shadow solutions, the constants can be chosen uniformly.
\end{proof}

\section*{Acknowledgement}
Y. Lee was supported by a National Research Foundation of Korea (NRF) grant funded by the Korea government (MSIT), No. RS-2026-25489601. 
The research of J. Wei is partially supported by the Hong Kong General Research Fund of the Research Grants Council of Hong Kong under the project ``New frontiers in singularity formations in nonlinear partial differential equations.'' W. Yang was supported by the National Key Research and Development Program of China (Grant No. 2022YFA1006800), the National Natural Science
Foundation of China (Grant No. 12531010), the Science and Technology Development Fund of the Macao  SAR(FDCT,GrantNo.0070/2024/RIA1), the Multi-Year Research Grants of the University of Macau (Grant Nos. MYRG-GRG2025-00051-FST and MYRG-GRG2026-00034-FSC), and the University of Macau Development Foundation (Grant No. TISF/2025/006/FST).
L. Zhang was partially supported by Simons Foundation Grant SFI-MPS-TSM-00013752.

We declare that AI assistance (OpenAI ChatGPT 6.1 Sol and Codex) was used in the preparation of this manuscript. Specifically, AI assistance was used to address the decoupling process in Subsection~\ref{typed-sec:rotation} and Subsection~\ref{subsection:ProofofTheorem2.6}.

\section*{Competing interests}
The authors declare that they have no competing interests.

\bigskip

\medskip

\noindent
\textsc{Youngae Lee}\\
Department of Mathematical Sciences,
Ulsan National Institute of Science and Technology (UNIST),
Ulsan, Republic of Korea\\
\textit{Email address: youngaelee@unist.ac.kr}

\medskip

\noindent
\textsc{Haoyu Li}\\
Department of Mathematics, Faculty of Science and Technology, University of Macau, Taipa, Macau, China\\
\textit{Email address: hyli1994@hotmail.com}

\medskip

\noindent
\textsc{Juncheng Wei}\\
Department Of Mathematics, The Chinese University Of Hong Kong, Shatin, Nt, Hong Kong\\
\textit{Email address:  wei@math.cuhk.edu.hk} 

\medskip

\noindent
\textsc{Wen Yang}\\
Department of Mathematics, Faculty of Science and Technology, University of Macau, Taipa, Macau, China\\
\textit{Email address: wenyang@um.edu.mo} 

\medskip

\noindent
\textsc{Lei Zhang}\\
Department of Mathematics, University of Florida, Gainesville, FL 32611-8105, USA\\
\textit{Email address: leizhang@ufl.edu}


\begingroup\small
\begin{thebibliography}{99}
\setlength{\itemsep}{0pt}\setlength{\parskip}{0pt}
\bibitem{Aro57} N. Aronszajn,
A unique continuation theorem for solutions of elliptic partial
differential equations or inequalities of second order,
J. Math. Pures Appl. 36 (1957), 235--249.

\bibitem{BJLY} D. Bartolucci, A. Jevnikar, Y. Lee, W. Yang,
Uniqueness of bubbling solutions of mean field equations,
J. Math. Pures Appl. 123 (2019), 78--126.

\bibitem{BJMR15} L. Battaglia, A. Jevnikar, A. Malchiodi, D. Ruiz,
A general existence result for the Toda system on compact surfaces,
Adv. Math. 285 (2015), 937--979.

\bibitem{BM91} H. Brezis, F. Merle,
Uniform estimates and blow-up behavior for solutions of
$-\Delta u=V(x)e^u$ in two dimensions,
Comm. Partial Differential Equations 16 (1991), 1223--1253.

\bibitem{CL15} C.-C. Chen, C.-S. Lin,
Mean field equation of Liouville type with singular data: topological degree,
Comm. Pure Appl. Math. 68 (2015), 887--947.

\bibitem{CL02} C.-C. Chen, C.-S. Lin,
Sharp estimates for solutions of multi-bubbles in compact Riemann surfaces,
Comm. Pure Appl. Math. 55 (2002), 728--771.

\bibitem{CL03} C.-C. Chen, C.-S. Lin,
Topological degree for a mean field equation on Riemann surfaces,
Comm. Pure Appl. Math. 56 (2003), 1667--1727.

\bibitem{ChenLi91} W. Chen, C. Li,
Classification of solutions of some nonlinear elliptic equations,
Duke Math. J. 63 (1991), 615--622.

\bibitem{DEM} M. del Pino, P. Esposito, M. Musso,
Nondegeneracy of entire solutions of a singular Liouville equation,
Proc. Amer. Math. Soc. 140 (2012), 581--588.


\bibitem{Deimling85} K. Deimling,
Nonlinear Functional Analysis,
Springer-Verlag, Berlin, 1985.

\bibitem{FK92}
H.~M. Farkas and I.~Kra,
\emph{Riemann Surfaces}, 2nd ed.,
Graduate Texts in Mathematics, vol.~71,
Springer-Verlag, New York, 1992.


\bibitem{GT01} D. Gilbarg, N. S. Trudinger,
Elliptic Partial Differential Equations of Second Order, 2nd ed.,
Classics in Mathematics, Springer, Berlin, 2001.

\bibitem{JKM15} A. Jevnikar, S. Kallel, A. Malchiodi,
A topological join construction and the Toda system on compact
surfaces of arbitrary genus,
Analysis \& PDE 8 (2015), 1963--2027.

\bibitem{JLW06} J. Jost, C.-S. Lin, G. Wang,
Analytic aspects of the Toda system. II. Bubbling behavior and
existence of solutions,
Comm. Pure Appl. Math. 59 (2006), 526--558.

\bibitem{Kato} T. Kato,
Perturbation Theory for Linear Operators,
Die Grundlehren der mathematischen Wissenschaften 132,
Springer-Verlag, Berlin--Heidelberg--New York, 1966.

\bibitem{Kiel12} H. Kielh\"ofer,
Bifurcation Theory: An Introduction with Applications to Partial
Differential Equations, 2nd ed., Applied Mathematical Sciences 156,
Springer, New York, 2012.

\bibitem{LLWY18} Y. Lee, C.-S. Lin, J. Wei, W. Yang,
Degree counting and shadow system for Toda system of rank two: One bubbling,
J. Differential Equations 264 (2018), 4343--4401.

\bibitem{LLY22} Y. Lee, C.-S. Lin, S. Yan,
Degree counting formula for non-symmetric Toda systems of rank two,
J. Funct. Anal. 283 (2022), 109642.

\bibitem{LLY19} Y. Lee, C.-S. Lin, W. Yang,
Existence of bubbling solutions without mass concentration,
Ann. Inst. Fourier 69 (2019), 895--940.

\bibitem{LLYZ20} Y. Lee, C.-S. Lin, W. Yang, L. Zhang,
Degree counting for Toda system with simple singularity: One point blow up,
J. Differential Equations 268 (2020), 2163--2209.



\bibitem{LLWYZ}H. Li, Y. Lee, J. Wei, W. Yang, L. Zhang,
Degree counting for singular \(SU(3)\) Toda systems.
Preprint.


\bibitem{Li99} Y. Y. Li,
Harnack type inequality: The method of moving planes,
Comm. Math. Phys. 200 (1999), 421--444.

\bibitem{LS94} Y. Y. Li, I. Shafrir,
Blow-up analysis for solutions of $-\Delta u=V(x)e^u$ in dimension two,
Indiana Univ. Math. J. 43 (1994), 1255--1270.

\bibitem{LWY} C.-S. Lin, J. Wei, D. Ye,
Classification and nondegeneracy of $SU(n+1)$ Toda system with singular sources,
Invent. Math. 190 (2012), 169--207.

\bibitem{LWYZ18} C.-S. Lin, J. Wei, W. Yang, L. Zhang,
On rank-2 Toda systems with arbitrary singularities: local mass and new estimates,
Analysis \& PDE 11 (2018), 873--898.

\bibitem{LWZ} C.-S. Lin, J. Wei, L. Zhang,
Classification of blowup limits for $SU(3)$ singular Toda systems,
Analysis \& PDE 8 (2015), 807--837.

\bibitem{LWZhao12} C.-S. Lin, J. Wei, C. Zhao,
Sharp estimates for fully bubbling solutions of a $SU(3)$ Toda system,
Geom. Funct. Anal. 22 (2012), 1591--1635.

\bibitem{MR13} A. Malchiodi, D. Ruiz,
A variational analysis of the Toda system on compact surfaces,
Comm. Pure Appl. Math. 66 (2013), 332--371.

\bibitem{MR15} A. Malchiodi, D. Ruiz,
On Leray--Schauder degree of the Toda system on compact surfaces,
Proc. Amer. Math. Soc. 143 (2015), 2985--2990.

\bibitem{Moser71} J. Moser,
A sharp form of an inequality by N. Trudinger,
Indiana Univ. Math. J. 20 (1971), 1077--1092.


\bibitem{Sard42} A. Sard,
The measure of the critical values of differentiable maps,
Bull. Amer. Math. Soc. 48 (1942), 883--890.

\bibitem{Smale65} S. Smale,
An infinite dimensional version of Sard's theorem,
Amer. J. Math. 87 (1965), 861--866.

\bibitem{Tar23} G. Tarantello,
On the blow-up analysis at collapsing poles for solutions of singular
Liouville-type equations,
Comm. Partial Differential Equations 48 (2023), 150--181.

\bibitem{WWZ} J. Wei, L. Wu, L. Zhang,
Estimates of bubbling sequences of $SU(3)$ Toda systems at critical parameters: Part 2,
J. London Math. Soc. 107 (2023), 2150--2196.
\end{thebibliography}
\endgroup
\end{document}